\documentclass[phd,tocprelim]{cornell}
\usepackage{amsmath,amssymb,amsthm}
\usepackage{graphicx}
\usepackage{graphics}
\usepackage{epsfig}
\usepackage{hangcaption}
\usepackage{txfonts}
\usepackage{palatino}

\usepackage{pictexwd}
\usepackage{dcpic}
\usepackage[all]{xy}
\usepackage{amstext}
\usepackage{amsfonts}
\usepackage{latexsym}

\newtheorem{defin}{Definition}[section]
\newtheorem{thm}[defin]{Theorem}
\newtheorem{cor}[defin]{Corollary}
\newtheorem{rem}[defin]{Remark}
\newtheorem{lem}[defin]{Lemma}
\newtheorem{prop}[defin]{Proposition}
\newtheorem{quest}[defin]{Question}
\newtheorem{prob}[defin]{Problem}
\newtheorem{axi}[defin]{Axiom}
\newtheorem{ex}[defin]{Example}
\newtheorem{notice}[defin]{Note}
\newtheorem{conven}[defin]{Convention}
\newtheorem{conj}[defin]{Conjecture}

\newcommand{\theorem}[1]{\begin{thm} \sl{#1} \end{thm}}
\newcommand{\theoremname}[2]{\begin{thm}[#1] \sl{#2} \end{thm}}
\newcommand{\proposition}[1]{\begin{prop} \sl{#1} \end{prop}}
\newcommand{\propositionname}[2]{\begin{prop}[#1] \sl{#2} \end{prop}}
\newcommand{\lemma}[1]{\begin{lem} \sl{#1} \end{lem}}
\newcommand{\definition}[1]{\begin{defin} \emph{#1} \end{defin}}
\newcommand{\question}[1]{\begin{quest} \emph{#1} \end{quest}}
\newcommand{\remark}[1]{\begin{rem} \emph{#1} \end{rem}}
\newcommand{\remarkname}[2]{\begin{rem}[#1] \emph{#2} \end{rem}}
\newcommand{\corollary}[1]{\begin{cor} \sl{#1} \end{cor}}
\newcommand{\corollaryname}[2]{\begin{cor}[#1] \sl{#2} \end{cor}}
\newcommand{\lemmaname}[2]{\begin{lem}[#1] \sl{#2} \end{lem}}

\newcommand{\example}[1]{\begin{ex} \emph{#1} \end{ex}}
\newcommand{\problem}[1]{\begin{prob} \emph{#1} \end{prob}}
\newcommand{\note}[1]{\begin{notice} \emph{#1} \end{notice}}
\newcommand{\convention}[1]{\begin{conven} \emph{#1} \end{conven}}
\newcommand{\conjecture}[1]{\begin{conj} \emph{#1} \end{conj}}
\newcommand{\ov}[1]{\overline{#1}}
\newcommand{\wh}[1]{\widehat{#1}}
\newcommand{\PL}{\mathrm{PL}}
\newcommand{\Homeo}{\mathrm{Homeo}_+}
\newcommand{\mf}[1]{\mathfrak{#1}}
\DeclareMathOperator{\Diff}{\mathrm{Diff}}
\DeclareMathOperator{\rot}{\mathrm{rot}}
\DeclareMathOperator{\PLog}{\mathrm{PLog}}
\DeclareMathSymbol{\Rb}{\mathbin}{AMSb}{"52}
\DeclareMathSymbol{\Zb}{\mathbin}{AMSb}{"5A}

\renewcommand{\caption}[1]{\singlespacing\hangcaption{#1}\normalspacing}

\title {Algorithms and Classification in Groups of Piecewise-Linear Homeomorphisms}
\author {Francesco Matucci}
\conferraldate {August}{2008}
\degreefield {Ph.D.}
\copyrightholder{Francesco Matucci}
\copyrightyear{2008}

\begin{document}

\maketitle
\makecopyright

\begin{abstract}

The first part (Chapters 2 through 5) studies decision problems in Thompson's groups $F,T,V$
and some generalizations. The \emph{simultaneous conjugacy problem} is determined to be
solvable for Thompson's group $F$ and suitable larger groups of piecewise-linear homeomorphisms
of the unit interval. We describe a conjugacy invariant both from the piecewise-linear point of view
and a combinatorial one using \emph{strand diagrams}.
We determine algorithms to compute roots and centralizers in these groups
and to detect periodic points and their behavior by looking at the \emph{closed strand diagram} associated to an element.
We conclude with a complete cryptanalysis of an encryption protocol based on the \emph{decomposition problem}.

In the second part (Chapters 6 and 7),
we describe the structure of subgroups of the group of all homeomorphisms of the unit circle,
with the additional requirement that they contain no non-abelian free subgroup. It is shown that in this setting
the \emph{rotation number map} is a group homomorphism. We give a classification of such subgroups as
subgroups of certain wreath products and we show that such subgroups can exist by building examples. Similar techniques
are then used to compute centralizers in these groups and to provide the base to generalize the techniques
of the first part and to solve the simultaneous conjugacy problem.

\end{abstract}

\begin{biosketch}
Francesco Matucci was born on May 5th, 1977 in Florence, Italy to Annamaria Savi and Patrizio Matucci.
He attended Liceo Guido Castelnuovo in Florence and
went on to earn a \emph{Laurea} in Mathematics \emph{cum Laude} in April 2002
at the University of Florence. In December 2005 he earned
a \emph{Dottorato di Ricerca} at the University of Florence
while studying abroad at Cornell University in Ithaca, NY, USA. He was admitted to the Mathematics graduate
program at Cornell University in January 2006 and completed his Ph.D. studies in August 2008 under
the supervision of Dr. Kenneth S. Brown and Dr. Martin D. Kassabov.
\end{biosketch}

\begin{dedication}
To whoever offered me a smile.
\end{dedication}

\begin{acknowledgements}
I would like to thank my two advisors Ken Brown and
Martin Kassabov for their guidance, their time
and many laughs.
They have supported me through adverse times and shown me much
deep mathematical thinking. I want to thank Ken who welcomed me
into the graduate program as a father and lifted the many burdens of these years.
I want to thank Martin, who has been a great friend and has always
been available for help. It has been both a pleasure and an honor to work with them. I also
want to thank the other members of my committee, Karen Vogtmann and R. Keith Dennis,
for giving me important advice while deciding my future directions.

I want to thank all of the faculty and staff of the Cornell Department of Mathematics,
but especially Maria Terrell for showing me the beauty of teaching and
encouraging me to develop this passion. An important mention goes also to Donna Smith who has been the ``glue''
throughout these years, providing tremendous support in so many ways.

I would like to thank Indira Chatterji and Collin Bleak who have taken care of me remotely and
shared their own experiences with me. I would also like to thank Matt Brin and Jos\'e Burillo
for giving me support throughout the enduring process of job searching.

I would like to thank Carlo Casolo who introduced me to the world
of Mathematics as an undergraduate and who followed me even when I crossed the ocean.

Many friends also supported me in this process: my fellow graduate students who welcomed me and made me feel at home,
in particular, thanks go to Brad, Jim, Jessica, Treven and Will.
I want to thank Chris, Franco and
Showey for just being my friends, entertaining me and showing me a new world.
And from the old country
I would like to thank Luca and Massimo for reaching out at all times.

Finally I want to thank Hiromi for listening, looking after me
and chasing me through the world.
My parents come last in this list, but they
know well they are first, because they are the ones who have taught me to smile.
\end{acknowledgements}

\setcounter{tocdepth}{1}

\contentspage
\tablelistpage
\figurelistpage

\normalspacing \setcounter{page}{1} \pagenumbering{arabic}
\pagestyle{cornell} \addtolength{\parskip}{0.5\baselineskip}

\chapter*{Introduction}
\label{chapter0}
\addcontentsline{toc}{chapter}{Introduction}

\section*{Decision Problems in Thompson's groups}

 My research started from studying decision problems for the important Thompson's
groups $F,T$ and $V$. These groups are piecewise-linear
homeomorphism groups of a 1-dimensional space and were introduced in 1965 by R.
Thompson in connection with his work in logic.
They were introduced during the creation of a finitely generated group with unsolvable word problem,
and later rediscovered in many other contexts. Thompson's groups provided the first known examples of
finitely presented infinite simple groups and they are still at the center of many geometric group theory questions.
What makes Thompson's groups an interesting starting point is that they are considered a test case for many conjectures.
Even though the groups have a simple definition, many questions prove to be a challenge.
Thompson's groups have many models: they can be described using generators and relations, or by
their action on a 1-dimensional space or by representing elements as combinatorial diagrams. This characteristic
often allows hard questions in one language to be transformed into easy questions in another one.
They are often used as instruments to measure the understanding of a certain property.
For example, it is an outstanding open question whether or not $F$ is an amenable group.

Richard Thompson's group $F$ can be seen as the group $\PL_2([0,1])$
of piecewise linear orientation-preserving homeomorphisms of the unit interval $[0,1]$,
with finitely many breakpoints such that:
\begin{itemize}
\item
all slopes are integral powers of $2$, and
\item
all breakpoints are in $\mathbb{Z}[\frac{1}{2}]$, the ring of dyadic rational numbers.
\end{itemize}
The product of two elements is given by the composition of functions.
The group $F$ is finitely presented (with two generators and two relations)
and torsion-free.

In addition to $F$, Thompson introduced two other finitely-generated
groups known as $T$ and $V$.  Briefly, $T$ is the set of piecewise-linear
self-homeomorphisms of the circle $\left[0,1\right]/\left\{0,1\right\}$
satisfying the two conditions above, while $V$ is the set of piecewise-linear
\emph{bijections} of the interval (or self-homeomorphisms of the Cantor set)
satisfying the above conditions.
We will recall all relevant definitions and
properties of Thompson's groups in Chapter \ref{chapter1}.

We say that a group $G$ has \emph{solvable conjugacy problem}
if there is an algorithm such that, given any two elements $x,y \in G$, we can determine whether there is,
or not, a $g \in G$ such that $g^{-1}xg=y$. The conjugacy problem for $F$ was addressed by Guba and Sapir~\cite{gusa1}, who
solved it for general diagram groups in 1997, observing that $F$ itself is a diagram group.
In Chapter \ref{chapter2}, we give a version of Guba and Sapir's
solution using \emph{strand diagrams}, and generalize it to $T$ and $V$. To the best of our knowledge,
the solution for $T$ is entirely new. The material of Chapter \ref{chapter2} represents joint work with James Belk.

In Chapter \ref{chapter3} we
we derive an explicit correspondence between strand diagrams and piecewise-linear functions. Specifically,
we show that strand diagrams can be interpreted as \emph{stack machines} acting on binary expansions.  Using this correspondence,
we obtain a complete understanding of the dynamics of elements, describing the behavior of fixed points and their
slopes. As a byproduct of our techniques,
we obtain simple proofs of previously known results.
In addition, we describe a completely dynamical solution to the conjugacy problem for one-bump functions in $F$,
similar to the Brin-Squier \cite{brin2} dynamical criterion for conjugacy in $\PL_+\left(I\right)$,
the group of all piecewise-linear orientation-preserving homeomorphisms of the unit interval with finitely many breakpoints.
The material of Chapter \ref{chapter3} represents joint work with James Belk.

For a fixed $k \in \mathbb{N}$, we say that the group $G$ has \emph{solvable $k$-simultaneous
conjugacy problem} if there is an algorithm such that, given any two $k$-tuples of elements in $G$,
$(x_1,\ldots,x_k),(y_1,\ldots,y_k)$, one can determine whether there is, or not, a $g \in G$ such that
$g^{-1} x_i g = y_i$ for all $i=1, \ldots, k$. We say that there is an
\emph{effective solution} if the algorithm produces such an element $g$, in addition to proving its existence.
In 1999, Guba and Sapir~\cite{gusa2} posed the question
of whether or not the simultaneous conjugacy problem was solvable for diagram groups.
In Chapter \ref{chapter4} we prove

\noindent \textbf{Theorem A.} {\sl Thompson's group $F$ has a solvable
$k$-simultaneous conjugacy problem, for every $k \in \mathbb{N}$.
There is an algorithm which produces an effective solution.}

With similar techniques we can prove that
the same result holds for larger groups of
piecewise linear homeomorphisms, containing $F$ as a subgroup.
The material of Chapter \ref{chapter4} represents joint work with Martin Kassabov.

\section*{Decision Problems and Cryptography}

Recent advances in public key cryptography have underlined the need to find alternatives
to the RSA cryptosystem. It has been proposed to use algorithmic problems
in non-commutative group theory as possible ways to build new protocols.
The \emph{conjugacy search problem} was introduced in several papers
as a generalization of the \emph{discrete logarithm problem} in the research of a new safe encryption scheme.
The former problem asks whether or not, given a group $G$ and two elements $a,b \in G$ that are conjugate,
we can find at least one $x \in G$ with $a^x:=x^{-1}ax=b$. It is thus important to look for a platform
group $G$ where this problem is computationally hard. Seminal works by Anshel-Anshel-Godlfeld \cite{aag}
and Ko-Lee et al. \cite{kolee}
have proposed the braid group $B_n$ on $n$ strands as a possible platform group.

It has been observed that Thompson's group $F$ and the braid groups $B_n$ have some similarities.
Belk proved in his thesis \cite{belkthesis} that $F$ and the braid groups have a similar classifying
space. Strand diagrams for elements of $F$ (introduced in Chapter \ref{chapter2}) are similar to braids,
but with merges and splits instead of twists.
However, for cryptographic purposes, $F$ has still not proved to be a
good platform. Theorem A proves
that the simultaneous conjugacy problem is solvable, making it insecure to apply
protocols based on the simultaneous conjugacy problem.

Shpilrain and Ushakov in \cite{su} have proposed using a particular version of the
\emph{decomposition problem} as a protocol and the group $F$ as a platform. The new problem is: given a
group $G$, a subset $X \subseteq G$ and two elements $w_1,w_2 \in G$ with the information that there exist
$a,b \in X$ such that $a w_1 b =w_2$, find at least one such pair $a,b$.
In Chapter \ref{chapter5} we show how to recover efficiently the shared secret key of this protocol.

\section*{Structure Theorems for $\Homeo(S^1)$ and Centralizers}

Let $\Homeo(S^1)$ denote the full group of orientation-preserving homeomorphisms of the unit circle and $G$
be one of its subgroups. Many papers have studied the structure of subgroups under particular assumptions.
Plante and Thurston have discovered that sufficient smoothness imposes a strong condition on nilpotent
groups of orientation-preserving diffeomorphisms.

\noindent \textbf{Theorem (Plante-Thurston, \cite{plathu}).}
{\sl Any nilpotent subgroup of $\mathrm{Diff}_+^2(S^1)$ must be abelian.}

\noindent On the other hand,
Farb and Franks showed that reducing the smoothness produces a contrasting situation, where every possibility
can occur.

\noindent \textbf{Theorem (Farb-Franks, \cite{farbfra1}).}
{\sl Every finitely-generated, torsion-free nilpotent group is isomorphic to a subgroup of
$\mathrm{Diff}_+^1(S^1)$.}

In Chapter \ref{chapter6} we relax the hypotheses on the regularity of the homeomorphisms and on the group. We explore the
dynamics of Poincar\`e's rotation number map $rot:G\to \mathbb{R}/\mathbb{Z}$: under certain conditions,
it is possible to prove that the $rot$ map becomes a homomorphism of groups. In particular, we obtain a result in the flavor
of the Tits Alternative:

\noindent \textbf{Theorem B.} {\sl Let $G \le \Homeo(S^1)$. Then the following alternative holds:

\noindent (i) $G$ has a non-abelian free subgroup, or

\noindent (ii) the map $rot:G \to (\mathbb{R}/\mathbb{Z},+)$ is a group homomorphism.}

Part (ii) of this first result allows us to write subgroups of $\Homeo(S^1)$ as extensions of
the kernel by a subgroup of $\mathbb{R}/\mathbb{Z}$ and hence to reduce the classification to studying
the kernel of the $rot$ map. As a byproduct,
we obtain Margulis' Theorem on the existence of a $G$-invariant measure on the unit circle (see \cite{marg1}).

\noindent \textbf{Theorem C.} {\sl Let $G \le \Homeo(S^1)$ with no non-abelian free subgroups. Then:

\noindent (i) $G$ is abelian, or

\noindent (ii) $G \hookrightarrow H_0 \wr K$, the standard unrestricted wreath product,
where $K:= G/G_0$ is isomorphic to a countable subgroup
of $\mathbb{R}/\mathbb{Z}$ and $H_0 \le \prod \Homeo(I_i)$ has no non-abelian free subgroups.}

We will show that such wreath products do exist in $\Homeo(S^1)$ by providing embedding theorems.
The material of Chapter \ref{chapter6} represents joint work with Collin Bleak and Martin Kassabov.

The techniques developed in Chapter \ref{chapter6} are then employed in Chapter \ref{chapter7} to obtain some results on centralizers
of elements and subgroups in $\PL_+(S^1)$ the group of orientation-preserving piecewise-linear homeomorphisms of the unit circle
with finitely many breakpoints. Centralizers are used in Chapter \ref{chapter4} as an intermediate step to go from the solution of the
ordinary conjugacy problem to the solution of the simultaneous conjugacy one. It is thus of interest to classify centralizers
in groups of homeomorphisms of the unit circle to generalize our results to this setting.
The material of Chapter \ref{chapter7} represents joint work with Collin Bleak and Martin Kassabov.

\section*{Estimating the size of balls in Thompson's group $F$}

In the last Chapter we describe a recurrence formula relating suitable slices of the $n$-sphere of Thompson's group $F$
with those of spheres of smaller radius, where elements are written with respect to the standard finite generating set of $F$.
The algorithm is based on the length formula for elements developed in \cite{bebr} by
Belk and Brown. We study their formula to detect what is the correct pattern to shorten a word to the identity element.
To the best of our knowledge, no other formula existed before to estimate the size of balls (besides direct counting of elements).

\chapter{Thompson's groups $F,T$ and $V$}
\label{chapter1}

In this Chapter we recall the main definitions and results about Thompson's groups $F,T$ and $V$ and some of their overgroups.
The proofs of all the stated results of this Chapter can be found either in \cite{cfp} or in
\cite{belkthesis}, unless otherwise stated.

\section{Background on $F$}

Let $I$ denote the unit interval $[0,1]$. Thompson's group $F$ is the group
of all piecewise-linear homeomorphisms of the unit interval with
finitely many breakpoints and satisfying the following conditions:
\begin{enumerate}
\item Every slope is a power of two, and
\item Every breakpoint has dyadic rational coordinates.
\end{enumerate}
The group $F$ is finitely presented (with two generators and two relations) and torsion-free.
It can be thought of as a ``lattice'' in the full group $\PL_+(I)$ of orientation-preserving piecewise-linear
homeomorphisms of $[0, 1]$ with finitely many breakpoints, and indeed it shares many properties with
this larger group.

\begin{figure}[0.5\textwidth]
\centering
\includegraphics{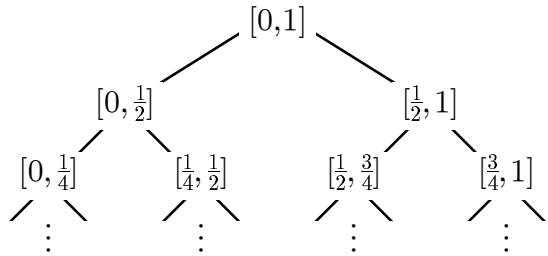}
\caption{The binary tree of dyadic intervals}
\label{fig:tree-dyadic-intervals}
\end{figure}

We are now going to describe how to see the elements of $F$ as diagrams. Consider the
subintervals of $I$ obtained by repeatedly cutting in half (see figure \ref{fig:tree-dyadic-intervals}).

\begin{figure}[0.5\textwidth]
\centering
\includegraphics{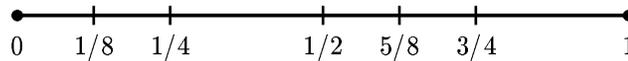}
\caption{An example of a dyadic subdivision}
\label{fig:example-of-dyadic-subdivision}
\end{figure}

\noindent These are the \em{standard dyadic intervals}. A \em{dyadic subdivision} of $I$
is any partition into finitely many standard dyadic intervals (see figure \ref{fig:example-of-dyadic-subdivision}).

\begin{figure}[0.5\textwidth]
\centering
\includegraphics{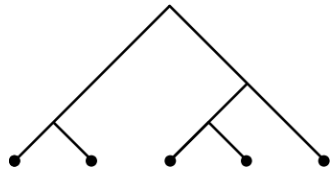}
\caption{The tree corresponding to the subdivision in figure \ref{fig:example-of-dyadic-subdivision}.}
\label{fig:dyadic-subdivision-into-tree}
\end{figure}

\noindent Dyadic subdivisions correspond to finite subtrees of the infinite binary tree
(the dyadic subdivision of figure \ref{fig:tree-dyadic-intervals} becomes the tree represented in
figure \ref{fig:dyadic-subdivision-into-tree}).
A \emph{dyadic rearrangement} is a homeomorphism $f: I \to I$
that maps intervals of one dyadic subdivision linearly to the intervals of another.

\begin{center}
\begin{minipage}{.5\textwidth}
\includegraphics{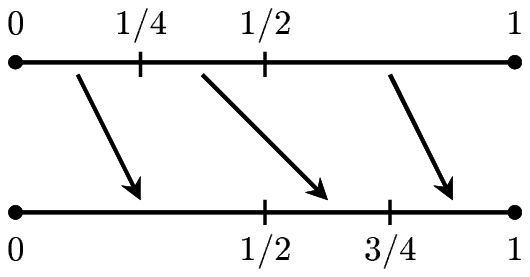}
\end{minipage}
\begin{minipage}{.25\textwidth}
$\begin{aligned}
.00\alpha\,&\mapsto\,.0\alpha\\[6bp]
.01\alpha&\mapsto.10\alpha\\[6bp]
.1\alpha&\mapsto.11\alpha
\end{aligned}$
\end{minipage}
\bigskip
\\ Figure 1.3.1: a dyadic rearrangement
\end{center}

\noindent If we represent elements of $[0, 1]$ in binary,
a dyadic rearrangement acts as a prefix replacement rule on binary sequences,
as illustrated in figure 1.3.1.


\proposition{The elements of Thompson's group $F$ are precisely the dyadic rearrangements of $I$.}

\begin{figure}[0.5\textwidth]
\centering
\includegraphics{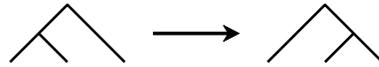}
\caption{An example of a tree diagram}
\label{fig:example-tree-diagram}
\end{figure}

A \emph{Tree diagram} for an element $f\in F$ is a pair of rooted,
binary trees that describe the dyadic subdivisions of the domain and range (see figure \ref{fig:example-tree-diagram}).

\begin{figure}[0.5\textwidth]
\centering
\includegraphics{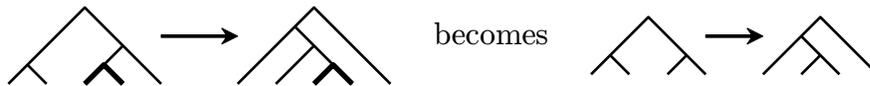}
\caption{The reduction of a tree diagram}
\label{fig:tree-diagram-reduction}
\end{figure}

\noindent The tree diagram for an element of $F$ is not entirely unique.
Specifically, we can reduce a tree diagram by canceling a corresponding pair of bottom carets
(see figure \ref{fig:tree-diagram-reduction}).

\noindent This corresponds to removing an unnecessary ``cut'' from the domain and range
subdivisions. We say that two tree diagrams are equivalent if one can obtained from the other
via a sequence of reductions and inverse reductions.

\theorem{Every element of $F$ has a unique reduced tree diagram.
Moreover, the elements of Thompson's group $F$ are in 1-to-1 correspondence with reduced tree diagrams.}

It is possible to define a product in the set of equivalence classes of
tree diagrams corresponding to the product in $F$.
Given two representatives for equivalence classes of tree diagrams $f:T_1 \to T_2$ and $g:T_3 \to T_4$,
it is possible to unreduce them until the range tree of $f$ is the same as the domain tree for $g$, that is
we can write $f:T_1' \to T_2'$ and $g:T_2' \to T_4'$ and then the product $fg$ is defined as the equivalence class of
the tree diagram $fg:T_1' \to T_4'$. This product agrees with the product of piecewise-linear homeomorphisms, hence
Thompson's group $F$ is isomorphic with the group of equivalence classes of tree diagrams.

\begin{figure}[0.5\textwidth]
\centering
  \includegraphics[height=6cm]{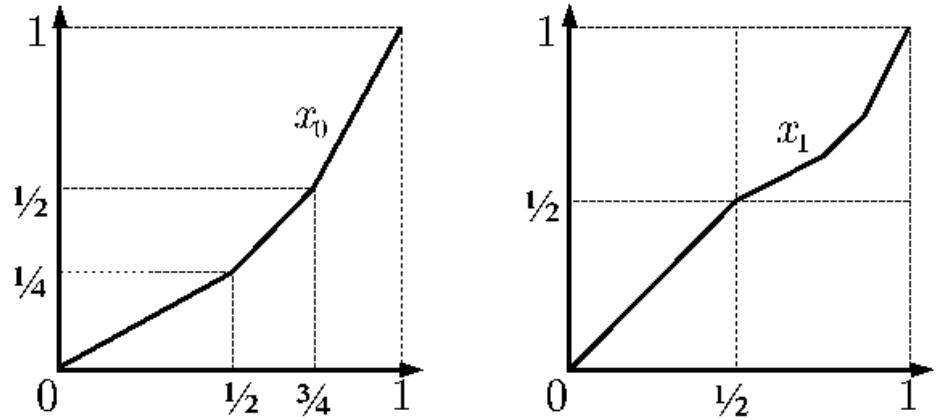}
 \caption{Two of the elements of the generating set of $F$.}
 \label{fig:generators-F}
\end{figure}

As we stated above Thompson's group admits a finite presentation, however we will not need it.
Instead, in Chapter \ref{chapter5}
we will make use of the infinite presentation described in the next result.
\theorem{Thompson's group $F$ is described by the following presentation
\[
F = \langle x_0, x_1, x_2, \ldots \mid x_n x_k = x_k x_{n+1}, \forall \, k<n\rangle.
\] \label{thm:thompson-presentation}}
This presentation has the advantage that the elements of $F$ can be uniquely written in the following \emph{normal form}
\[
x_{i_1}\ldots x_{i_u} x_{j_v}^{-1} \ldots x_{j_1}^{-1}
\]
such that $i_1 \le \ldots \le i_u$, $j_1 \le \ldots \le j_v$ and if both $x_i$ and
$x_i^{-1}$ occur, then either $x_{i+1}$ or $x_{i+1}^{-1}$ occurs, too. Since $x_k=x_0^{1-k}x_1x_0^{k-1}$
for $k \ge 2$, the group $F$ is generated by the elements $x_0$ and $x_1$.
The generators $x_k$ of the infinite presentation can be represented as piecewise-linear homeomorphisms
by shrinking the function $x_0$ shown in figure \ref{fig:generators-F} onto the interval
$[1-\frac{1}{2^k},1]$ and extending it as the identity on $[0,1-\frac{1}{2^k}]$.

\lemma{If $0=a_0< a_1 < a_2 < \ldots < a_n=1$
and $0=b_0 < b_1 < b_2 < \ldots < b_n=1$ are two partitions
of $[0,1]$ consisting of dyadic rational numbers, then we can build an $f \in F$, such that $f(a_i)=b_i$.
In particular, $F$ acts transitively on $k$-tuples of dyadic points in $(0,1)$ for any $k$.
\label{thm:cfp-extension}}

For an interval $J \subseteq \mathbb{R}$ with dyadic endpoints,
it is possible to define Thompson-like analogues of Thompson's group $F$. We can consider
the groups $\PL_2(J)$ of piecewise-linear homeomorphism on the interval $J$ onto itself
with finitely many breakpoints and the same requirements on slopes and breakpoints as $F$.

\theorem{For any two dyadic rationals $\alpha<\beta$ in the real line $\mathbb{R}$, there
exists a Thompson-like homeomorphism $\varphi : [\alpha, \beta] \to [0,1]$, i.e.
a piecewise-linear homeomorphisms with finitely many breakpoints occurring at dyadic rationals and whose
slopes are integral powers of $2$. \label{thm:thompson-like}}

\noindent \emph{Proof.} This is a well known fact, but we provide a proof for the sake of completeness.
Let $m<n$ be a pair of integers such that $[0,1] \cup [\alpha,\beta] \subseteq (m,n)$ and such that $m-n=2^v$,
for some integer $v$. If we consider the straight line map $\psi:[m,n] \to [0,1]$, it is straightforward
to verify that the choice of $m,n$ implies that the map $\rho(G):=\psi G \psi^{-1}$ yields an
isomorphism $\rho:\PL_2([m,n]) \to \PL_2(I)$.
To construct the required Thompson-like map it is sufficient to
consider $m<\alpha<\beta<n$ and $m<0<1<n$ as partitions of the interval $[m,n]$ and bring them to the interval $[0,1]$
through $\psi$. Hence, we can apply Lemma \ref{thm:cfp-extension} to find an element $f \in \PL_2(I)$ that sends
the partition $0=\psi(m)<\psi(\alpha)<\psi(\beta)<\psi(n)=1$ into the partition
$0<\psi(0)<\psi(1)<1$ to find an element $f \in \PL_2(I)$. To conclude,
it is sufficient to define $\varphi$ as the restriction of $\rho^{-1}(f)$ to the interval $[\alpha,\beta]$. $\square$

\corollary{For any two dyadic rationals $\alpha<\beta$ in the interval $[0,1]$,
the groups $\PL_2([\alpha,\beta])$ and $\PL_2(I)$ are isomorphic. \label{thm:thompson-like-isomorphic}}

\noindent \emph{Proof.} By Theorem \ref{thm:thompson-like}, there is a Thompson-like map $\varphi : [\alpha, \beta] \to [0,1]$.
Now define
\[
\begin{array}{ccc}
\PL_2([\alpha,\beta]) & \longrightarrow & \PL_2(I) \\
         f            & \longmapsto     & \varphi f \varphi^{-1} \; \; \; \square
\end{array}
\]

The previous two results are used at many points in this thesis: for example we will use them in Chapters \ref{chapter3},
\ref{chapter4} and \ref{chapter7}.

\section{The generalized groups $\PL_{S,G}(J)$ and the group $\PL_+(J)$}

In Chapter \ref{chapter4} we will work with a generalization of Thompson's group $F$ by relaxing the hypotheses
on breakpoints and slopes. Let $S$ be a subring of $\mathbb{R}$, let $U(S)$ denote the group of invertible elements of $S$ and let
$G$ be a subgroup of $U(S) \cap \mathbb{R}_+$. For any interval $J$ with endpoints in $S$,
we define $\PL_{S,G}(J)$ to be the
group of piecewise linear homeomorphism from the interval $J$ into itself,
with only a finite number of breakpoints and such that
\begin{itemize}
\item
all breakpoints are in the subring $S$,
\item
all slopes are in the subgroup $G$,
\end{itemize}
the product of two elements being given by the composition of functions. Unlike the
case of $F$ and Corollary \ref{thm:thompson-like-isomorphic}, it is not true anymore that for any
two intervals $J_1,J_2$ with endpoints in $S$ the groups $\PL_{S,G}(J_1)$ and $\PL_{S,G}(J_2)$ are isomorphic.
For any interval $J$, we denote by $\PL_+(J)$ the group of piecewise-linear orientation-preserving homeomorphisms
of the interval $J$, with finitely many breakpoints. Since there are no requirements for the breakpoints and the slopes
of elements of $\PL_+(J)$ then, for any subring $S$ and
$G \le U(S) \cap \mathbb{R}_+$, we have $\PL_{S,G}(J) \subseteq \PL_+(J)$.

\section{Thompson's group $T$ and $V$}

\begin{figure}[0.5\textwidth]
\centering
  \includegraphics[height=3cm]{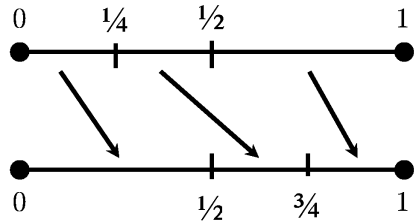}
 \caption{An element of $F$ as a dyadic rearrangement of the unit interval.}
 \label{fig:F-element-partitions}
\end{figure}

Thompson's group $F$ is the group of dyadic rearrangements on the unit interval. The important requirement is that
these rearrangements preserve the order of the intervals of the partition (see figure \ref{fig:F-element-partitions}).

\begin{figure}[0.5\textwidth]
\centering
  \includegraphics[height=6cm]{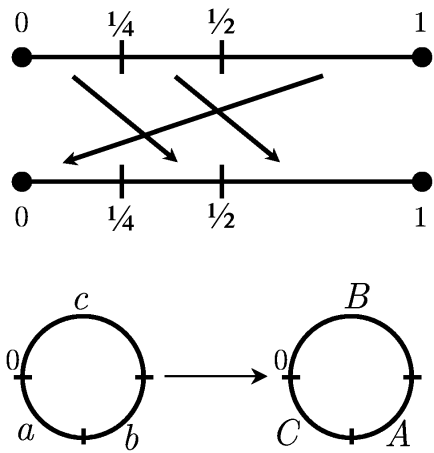}
 \caption{A dyadic rearrangement of the unit circle.}
 \label{fig:T-element-partitions}
\end{figure}

We can relax the last requirement on the order of intervals and require that these rearrangements
preserve the cyclic order of intervals (see figure \ref{fig:T-element-partitions}). \emph{Thompson's group} $T$
is the group of dyadic rearrangements of $[0,1]$ that preserve the cyclic order. Since the elements of
$T$ preserve the cyclic order of subdivisions of $[0,1]$ they can be viewed as homeomorphisms
of the unit circle.

\begin{figure}[0.5\textwidth]
\centering
  \includegraphics[height=2cm]{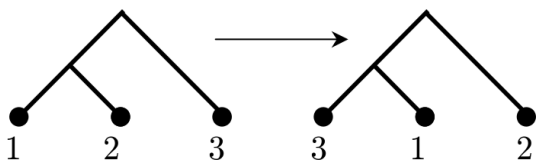}
 \caption{A tree diagram for the element of $T$.}
 \label{fig:T-tree-diagram}
\end{figure}

It is possible to define tree diagrams for elements of $T$. They correspond to dyadic rearrangements of the
unit circle and are thus represented by a pair of rooted, binary trees along with a cyclic permutation of the leaves
(see figure \ref{fig:T-tree-diagram}).

\begin{figure}[0.5\textwidth]
\centering
  \includegraphics[height=3cm]{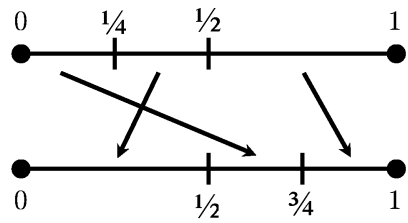}
 \caption{An element of $V$ seen as a dyadic rearrangement.}
 \label{fig:V-element-partitions}
\end{figure}

If we allow a dyadic rearrangement to permute the order of intervals arbitrarily, we obtain
a larger group containing $F$ and $T$. \emph{Thompson's Group $V$} is the group of dyadic rearrangements
of $[0,1]$ that may permute the order of the subdivisions (see figure \ref{fig:V-element-partitions}).

Note that this produces bijections
$[0,1] \rightarrow [0,1]$ that are \emph{not} continuous.  (By convention, all functions
in $V$ are required to be continuous from the right.  Alternatively, one can define $V$
as a group of homeomorphisms of the Cantor set.) The set of all elements of $V$ forms a group under composition.
\begin{figure}[0.5\textwidth]
\centering
  \includegraphics[height=3.5cm]{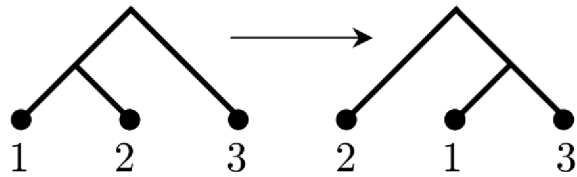}
 \caption{A tree diagram for the element of $V$ in figure \ref{fig:V-element-partitions}.}
 \label{fig:V-tree-diagram}
\end{figure}

Even for $V$ it is possible to define tree diagrams,
which will appear as a pair of rooted, binary trees along with a permutation of the leaves (see figure \ref{fig:V-tree-diagram}).


\chapter{Conjugacy in Thompson's groups}
\label{chapter2}

In this Chapter we give a unified solution for the conjugacy problem in Thompson's groups $F,T$ and $V$.
We introduce \emph{strand diagrams}, a modification of tree diagrams for these groups, and show how
identifying the roots of the two trees defines a conjugacy invariant in all cases. This
reduces the conjugacy problem to the study of the isomorphism problem for certain classes of graphs
and gives us elementary proofs of some known results. Strand diagrams were first introduced by Pride
in his study of the homotopy of relations using the term \emph{pictures} in \cite{Pr1},
\cite{Pr2} and \cite{BoPr} and are dual to the diagrams introduced by Guba and Sapir.

In 1997 Guba and Sapir showed that $F$ can be viewed as a \emph{diagram group} for
the monoid presentation $\langle x\mid x^2=x\rangle$ \cite{gusa1}. They give a solution
for each diagram group, and in particular for $F$. Their solution amounted to an algorithm
which had the same complexity as the isomorphism problem for planar graphs. This last problem has been solved in linear
time in 1974 by Hopcroft and Wong~\cite{hopwon}, thus proving the
Guba and Sapir solution of the conjugacy problem for diagram groups optimal. We mention here relevant
related work: in 2001 Brin and Squier in \cite{brin2} produced a criterion for describing conjugacy classes in $\PL_+(I)$.
In 2007 Gill and Short \cite{ghisho1} extended this criterion to work in $F$,
thus finding another way to characterize conjugacy classes
from a piecewise linear point of view.

The conjugacy problem in $V$ was previously solved by Higman
\cite{Hig} by combinatorial group theory methods
and again by Salazar-Diaz \cite{saldiazthesis} by using the techniques introduced by Brin
in his paper \cite{brin4}. On the other hand, to the best of our knowledge, the solution for $T$ is entirely new.

This Chapter is organized as follows.  In section 1 we give a simplified
solution to the conjugacy problem in $F$.  We extend this solution to $T$ in
section 2, and to $V$ in section~3, and in section 4 we analyze the running
time of the algorithm.  Finally, we have relegated to the appendix a proof
that every closed strand diagram for a conjugacy class in $F$, $T$, or $V$
possesses a cutting path. The material of this Chapter represents
joint work with James Belk. It can also be found in \cite{bema1}.

\section{Conjugacy in Thompson's group $F$}

\subsection{Strand Diagrams}

\begin{figure}[0.5\textwidth]
 \centering
  \includegraphics[height=3cm]{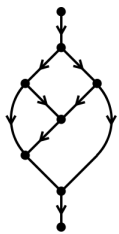}
  \caption{A $(1,1)$-strand diagram}
  \label{fig:example-strandiag}
\end{figure}

In this section, we describe Thompson's group $F$ as a group of \emph{strand
diagrams}. A strand diagram is similar to a braid, except instead of twists,
there are splits and merges (see figure \ref{fig:example-strandiag}).

\begin{figure}[0.5\textwidth]
 \centering
  \includegraphics{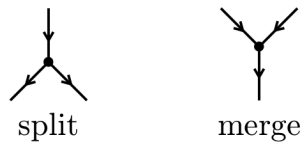}
  \caption{A split and a merge}
  \label{fig:split-merge-def1}
\end{figure}

To be precise, a \emph{strand diagram} (or a $(1,1)$-\emph{strand diagram})
is any directed, acyclic graph in the unit square
satisfying the following conditions:

\begin{enumerate}
\item There exists a unique univalent source along the top of the square,
and a unique univalent sink along the bottom of the square.
\item Every other vertex lies in the interior of the square, and is either
a \emph{split} or a \emph{merge} (see figure \ref{fig:split-merge-def1}).
\end{enumerate}

\begin{figure}[0.5\textwidth]
 \centering
  \includegraphics[height=3cm]{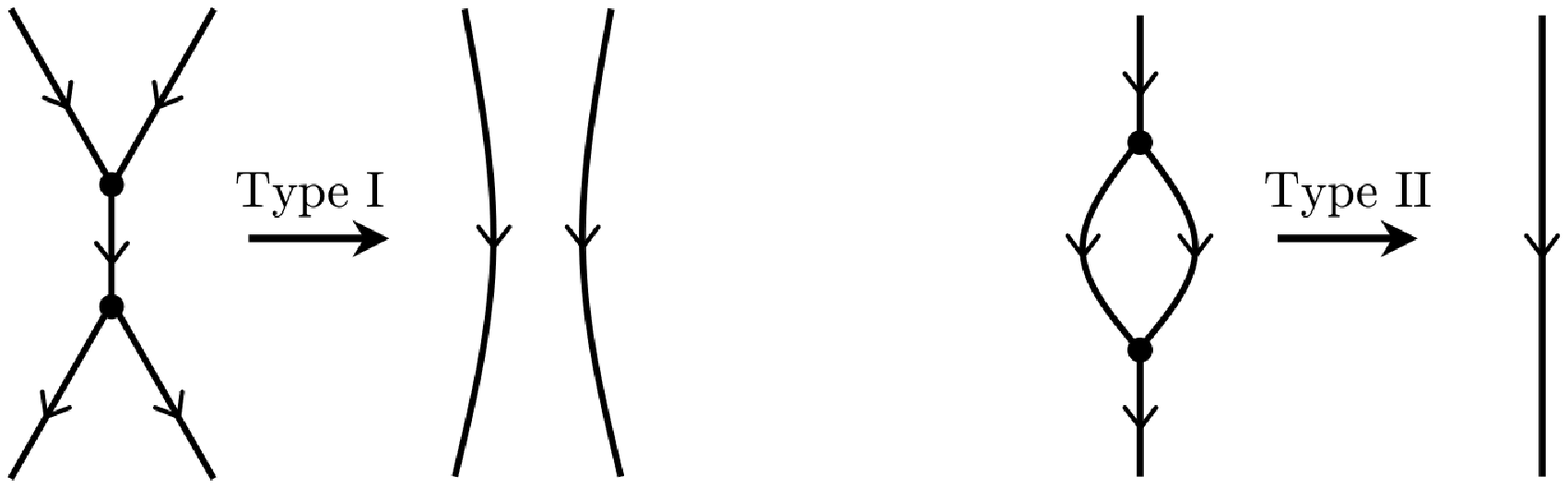}
  \caption{Reductions}
  \label{fig:reductions-strand-F}
\end{figure}

As with braids, isotopic strand diagram are considered equal.
A \emph{reduction} of a strand diagram is either of the moves
shown in figure \ref{fig:reductions-strand-F}.

Two strand diagrams are \emph{equivalent} if one can be obtained from the other via
a sequence of reductions and inverse reductions. A strand diagram is
reduced if it is not subject to any reductions.

\proposition{Every strand diagram is equivalent to a unique reduced strand diagram.
\label{thm:equivalent-unique-diagram}}

\noindent \emph{Proof:} This result was first proved by Kilibarda \cite{kilibardathesis}, and
appears as lemma 3.16 in \cite{gusa1}. We repeat the proof here, for
we must prove several variations of this result later.
Consider the directed graph $\mathcal{G}$ whose vertices are strand
diagrams, and whose edges represent reductions.  We shall use Newman's
Diamond Lemma (see \cite{new1}) to show that each component of $\mathcal{G}$ contains a unique
terminal vertex.

\begin{figure}[0.5\textwidth]
 \centering
  \includegraphics[height=3cm]{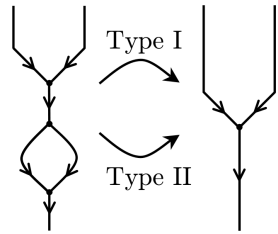}
  \caption{Diamond Lemma}
  \label{fig:diamond-lemma-F}
\end{figure}

Clearly $\mathcal{G}$ is terminating, since each reduction decreases
the number of vertices in a strand diagram.  To show that $\mathcal{G}$ is
locally confluent, suppose that a strand diagram is subject to two different
reductions, each of which affects a certain pair of vertices.  If these two
pairs are disjoint, then the two reductions simply commute.  The only other
possibility is that the two pairs have a vertex in common, in which case the
two reductions have the same effect (figure \ref{fig:diamond-lemma-F}). $\square$

The advantage of strand diagrams over tree diagrams is that \emph{multiplication} is the same as concatenation
(see figure 2.4.1).
\begin{center}
\begin{minipage}{0.75in}
\centering
\includegraphics{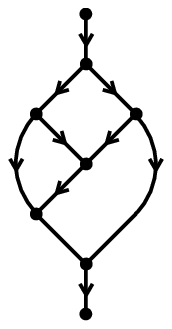}\\
$f$
\end{minipage}%
\hspace{0.6in}
\begin{minipage}{0.75in}
\centering
\includegraphics{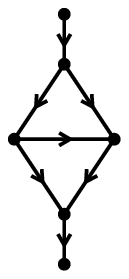}\\
$g$
\end{minipage}%
\hspace{0.6in}
\begin{minipage}{0.75in}
\centering
\includegraphics{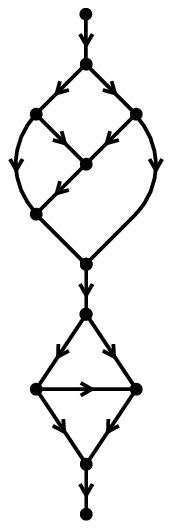}\\
$g \circ f$
\end{minipage}%
\bigskip
\\ Figure 2.4.1: Product is given by concatenating diagrams
\end{center}
This algorithm is considerably simpler than the standard multiplication algorithm for tree diagrams (see
Chapter \ref{chapter1}). The \emph{inverse} of a strand diagram is obtained by reflection
across a horizontal line and by inverting the direction of all the edges.
Note that the product of a strand diagram with
its inverse is always equivalent to the identity.

\theorem{Thompson's group $F$ is isomorphic with the group of all equivalence
classes of strand diagrams, with product induced by concatenation.
\label{thm:correspondence-F-strand-diagrams}}

\noindent \emph{Proof.} There is a close relationship between
strand diagrams and
the well-known tree pair diagrams for elements of $F$. In particular, a
strand diagram for an element $f\in F$ can be constructed by gluing the two
trees of a tree pair diagram together along corresponding leaves,
after turning one tree upside down (see figure 2.4.2).
\begin{center}
\begin{minipage}{2.25in}
\centering
\includegraphics{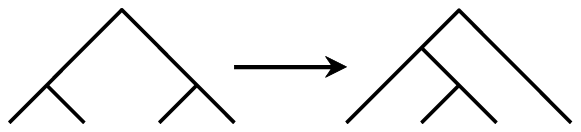}
\end{minipage}%
\hspace{0.4in}becomes\hspace{0.4in}
\begin{minipage}{1.4in}
\centering
\includegraphics{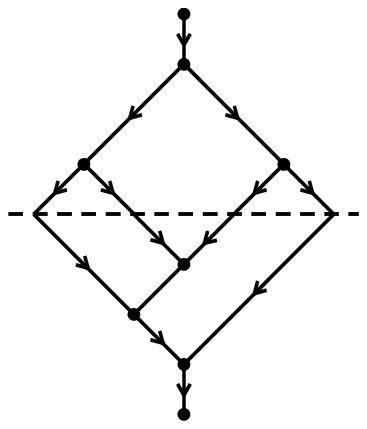}\\[0in]
\end{minipage}
\bigskip
\\ Figure 2.4.2: Gluing a tree diagram
\end{center}
Conversely, any reduced strand diagram can be ``cut'' in a unique way to
obtain a reduced pair of binary trees (see figure 2.4.3).
\begin{center}
\begin{minipage}{2.5in}
\centering
\includegraphics{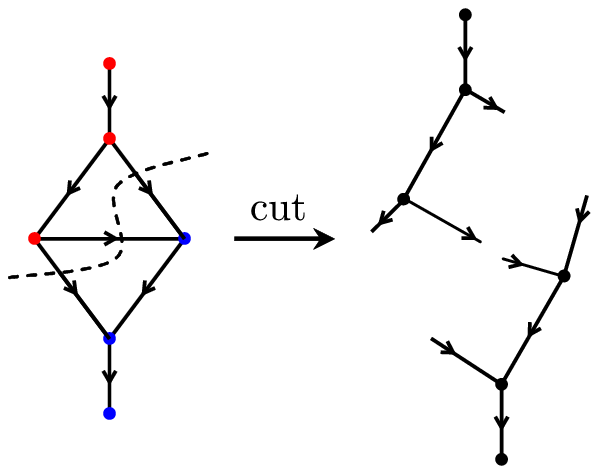}
\end{minipage}%
\hspace{0.2in}becomes\hspace{0.2in}
\begin{minipage}{1.5in}
\centering
\includegraphics{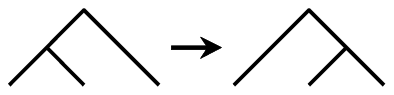}
\end{minipage}
\bigskip
\\ Figure 2.4.3: Cutting a strand diagram
\end{center}
We will now give a more formal proof of the intuition given by the previous two figures.
We define the following three sets: $\mathcal{SD}=\{$reduced strand diagrams$\}$,
$\mathcal{TD}=\{$reduced tree diagrams$\}$ and $PL_2(I)$ is the piecewise linear description of $F$.
We want to prove that the previous three groups are isomorphic. Define the following map
\[
\begin{array}{cccc}
\varphi: & \mathcal{SD} &    \longrightarrow & PL_2(I) \\
         &     D        & \longmapsto & f_D
\end{array}
\]
To build $f_D(t)$ we construct a path through the diagram which carries along a number
that will change at each vertex.
Given a number $t\in [0,1]$ we write it in binary expansion $t=0.b_1b_2 \ldots$.
We  start with $t$ at the source and follow the outgoing edge. If we meet a split, we remove
the first bit $b_1$ from $t$ obtaining a new number $t'=0.b_2 \ldots$. We then exit from the left
if $b_1=0$ and from the right otherwise. We keep this going until we meet a merge.
If we meet a merge with a number $0.a_1\ldots$ we add a new first bit,
obtaining a number $0.a_0a_1\ldots$. The number $a_0$ that we add is a $0$
if we arrived to the merge from the left and a $1$ otherwise. We now iterate this
pattern for all the splits and merges that we encounter, until we exit through the
unique output of the strand diagram (the sink).

The number that remains at the end of this procedure is
called $f_D(t)$. It is not difficult to see that each map $f_D(t)$ is a homeomorphism and
that the map $\varphi$ is a homomorphism (we will describe this construction in more detail in Chapter \ref{chapter3}).

Given a reduced strand diagram $D$ there is a way to cut it into two halves which are both
directed trees. We consider the set of edges that leave all the splits above and all the merges below and cut each
of them. More precisely, to determine this set of edges, we use the following procedure: start with the set containing
only the edge leaving the source, replace it with its two children. Now, if an edge terminates in a split, replace it
with its two children. We repeat until we remain with a set of edges each of which terminates in a merge.
It is not difficult to prove that this procedure is well defined, using the Diamond Lemma. Moreover, it is clear that
what lies above these edges is a tree, and that what lies below must also be a tree (otherwise
the diagram would not be reduced, or the set of edges is not minimal according to the previous procedure)

Thus, there exists at least one cut dividing a strand diagram into a tree diagram. A priori, there might be more than
one way to cut the diagram in two trees, hence we need to prove
that there is only one such cut. We build the following diagram
\begin{displaymath}
\xymatrix{
\mathcal{SD} \ar[r]^{\varphi}  &  PL_2(I) \\
             & \ar[u]_{\mu} \ar[ul]^{\sigma} \mathcal{TD}
}
\end{displaymath}
\noindent where $\sigma$ is the map obtained by taking a
tree diagram and gluing all its corresponding leaves to get
a strand diagram and $\mu$ is the standard map
which associates a piecewise-linear homeomorphism to a tree diagram (see Chapter \ref{chapter1}).
By definition of the maps, the diagram is commutative.

\noindent Let $D \in \mathcal{SD}$ and choose some cut $v$ on the edges of $D$ which divides $D$ into two directed trees
and define $T_v$ to be the tree diagram associated to this cut.

\noindent \emph{Claim:} The tree diagram $T_v$ is independent of the choice of the cut $v$.

\noindent \emph{Proof of the Claim.} Since the map $\sigma$ glues back the points where we have cut $D$,
we have that $\sigma(T_v)=D$. Hence we get
\[
\mu(T_v)=\varphi \sigma (T_v) = \varphi(D).
\]
By applying $\mu^{-1}$ to both sides of the previous equality, we get
$T_v=\mu^{-1}\varphi(D)$ and therefore $T_v$ is the image of a map and is then
defined independently of the chosen cut. $\square$

\noindent The previous Claim allows us to well-define a map
\[
\begin{array}{cccc}
\psi: & \mathcal{SD} & \longrightarrow & \mathcal{TD} \\
      &      D       & \longmapsto     &     T_v
\end{array}
\]
By the proof of the Claim, we have shown that $\mu \psi = \varphi$. Moreover, the relations
$\sigma \psi = id_{\mathcal{SD}}$ and $\psi \sigma = id_{\mathcal{TD}}$ follow easily.
Therefore $\psi$ is bijective, and it is clearly a homomorphism. $\square$

\note{\label{thm:groupoid} We will sometimes need to consider more general strand diagrams, with more than one source and sink
(see figure \ref{fig:m,n-strand-diagram1}).}

\begin{figure}[0.5\textwidth]
 \centering
  \includegraphics[height=3cm]{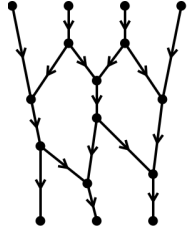}
  \caption{An $(m,n)$-strand diagram}
  \label{fig:m,n-strand-diagram1}
\end{figure}
\begin{figure}[0.5\textwidth]
 \centering
  \includegraphics[height=2cm]{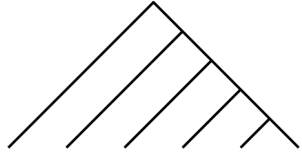}
  \caption{The right vine}
  \label{fig:right-vine1}
\end{figure}

We call an object like this an
$\left(m,n\right)$-strand diagram, where $m$ is the number of sources and $n$
is the number of sinks. This graph is built with the same conditions as (1,1)-strand diagrams,
except that it is allowed to have multiple sources and sinks.
We observe that, for every positive integer $k$ the equivalence classes of $(k,k)$-strand diagrams
equipped with the product given by concatenation returns a group.
It is possible to prove that the group of all $(1,1)$-strand diagrams is isomorphic to the group
of all $(k,k)$-strand diagrams, for every positive integer $k$.
In fact, if we denote by $v_m$ the right vine with $m$ leaves (figure \ref{fig:right-vine1}).

then any $\left(k,k\right)$-strand
diagram can be identified with the corresponding $(1,1)$-strand diagram given by $v_k^{-1}fv_k$
of Thompson's group $F$.In particular, we can compose
two elements of $F$ by concatenating any corresponding pair of strand
diagrams. The previous description can be seen more formally from a categorical point of view.
We consider a category $\mathcal{C}$, where
$Obj(\mathcal{C})=\mathbb{N}$, $Mor(\mathcal{C})=\{$morphisms $i \to j$ are labeled binary forests
with $i$ trees and $j$ total leaves$\}$ (notice
that the labeling on the trees induces a labeling on the leaves). The composition of two morphisms
$f:i\to j$ and $g:j \to k$,
is the morphism $fg:i \to k$ obtained by attaching the roots of the trees of $g$ to the leaves of $f$ by respecting
the labeling
on the roots of $g$ and the leaves of $f$. With this definition,
the equivalence classes of strand diagrams with any number of sources and sinks
is the groupoid of fractions of the category $\mathcal{C}$
(more details on this construction can be found in \cite{belkthesis}). We call this \emph{Thompson's groupoid}
$\mathcal{F}$. From the categorical point of view, we see that
the projection $f\mapsto v_n^{-1}fv_m$ is an epimorphism from the groupoid $\mathcal{F}$ to the group $F$.

\subsection{Annular Strand Diagrams}

\definition{An \emph{annular strand diagram} is a directed graph embedded
in the annulus with the following properties:
\begin{enumerate}
\item Every vertex is either a merge or a split.
\item Every directed cycle has positive winding number around the central hole.
\end{enumerate}}

\begin{figure}[0.5\textwidth]
 \centering
  \includegraphics[height=4cm]{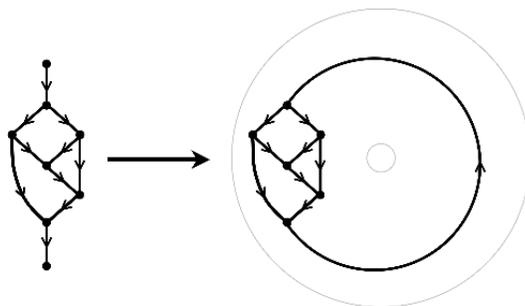}
  \caption{An annular strand diagram}
  \label{fig:annular-strand}
\end{figure}

Our definition of graph allows the existence of \emph{free loops},
i.e. directed cycles with no vertices on them.
Every element of $F$ gives an annular strand diagram:
given a strand diagram in the square, we can identify the top and bottom
and delete the resulting vertex to
get an annular strand diagram (figure \ref{fig:annular-strand}).

\begin{figure}[0.5\textwidth]
 \centering
  \includegraphics[height=4cm]{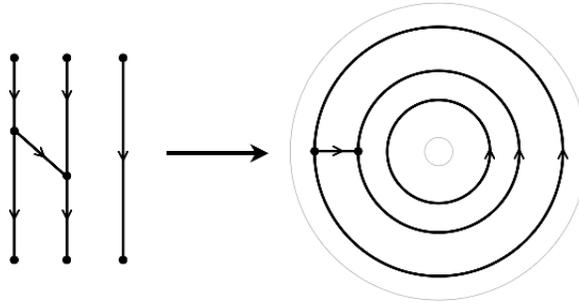}
  \caption{An example of a free loop in an annular strand diagram}
  \label{fig:free-loops-exist}
\end{figure}

More generally, you can obtain an annular strand diagram from any
$\left(k,k\right)$-strand diagram in the square, for any $k \ge 1$. We observe that we may obtain free
loops with no vertices (figure \ref{fig:free-loops-exist})

\definition{A \emph{cutting path} for an annular strand diagram is a
continuous path in the annulus that satisfies the following conditions:}

\begin{figure}[0.5\textwidth]
 \centering
  \includegraphics[height=3cm]{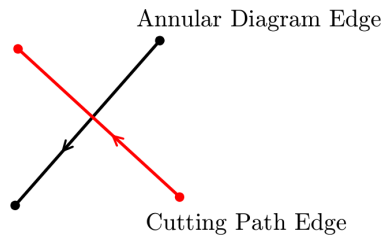}
  \caption{Orientation of the cutting path}
  \label{fig:annular-cutting-orientation}
\end{figure}

\begin{enumerate}
\item The path begins on the inner circle of the annulus, and ends on the outer circle.
\item The path does not pass through any vertices of the strand diagram.
\item The path intersects edges of the strand diagram transversely, with the orientation shown in figure
\ref{fig:annular-cutting-orientation}.
\end{enumerate}

Cutting an annular strand diagram along a cutting path yields a
$\left(k,k\right)$-strand diagram embedded in the unit square (thus an element of Thompson's group $F$,
by Note \ref{thm:groupoid}). Conversely, given an element of $F$ as a
$(k,k)$-strand diagram we can build an \emph{associated annular strand diagram} by gluing the $i$-th source and the $i$-th sink, for
$i=1, \ldots, k$. This gluing also defines a cutting path for the associated annular diagram.
On the other hand, it can be shown that every annular strand diagram has at least
one cutting path, hence any annular strand diagram is the associated
annular strand diagram for some $(k,k)$-strand diagram. We sketch a proof of this fact, even though it will not be used in the characterization
of conjugacy for elements of $F$ seen as strand diagrams.

\theorem{Every annular strand diagram has a cutting path.}

\noindent \emph{Sketch of a Proof:} Let $S$ be an annular strand diagram, and let $c$ be the class in
$H^1\left(S\right)$ induced by winding number on the annulus. By theorem
\ref{thm:appendix} in the appendix, there exists a cochain $\alpha$ on $S$ representing
$c$ which takes a non-negative value on each directed edge.

If we regard the directed graph $S$ as embedded in the plane, we observe that $S$
divides the plane into regions and we can define
$S^\ast$ the directed dual graph to $S$. That is, $S^\ast$
is the graph with one vertex for each region of $S$---including a vertex $i$
for the inner region and a vertex $o$ for the outer region---and with
directed edges that transversely intersect the directed edges of $S$ in the
same manner as a cutting path. The cochain $\alpha$ on $S$ can be viewed as
a chain $\alpha^\ast$ on $S^\ast$, which is a positive linear combination of
directed edges.  In particular, the boundary of $\alpha^\ast$ must be the
difference $o-i$.  Then $\alpha^\ast$ must be the sum of directed cycles and
a single directed path from $i$ to $o$, the latter being the desired cutting
path. $\square$

\begin{figure}[0.5\textwidth]
 \centering
  \includegraphics[height=6cm]{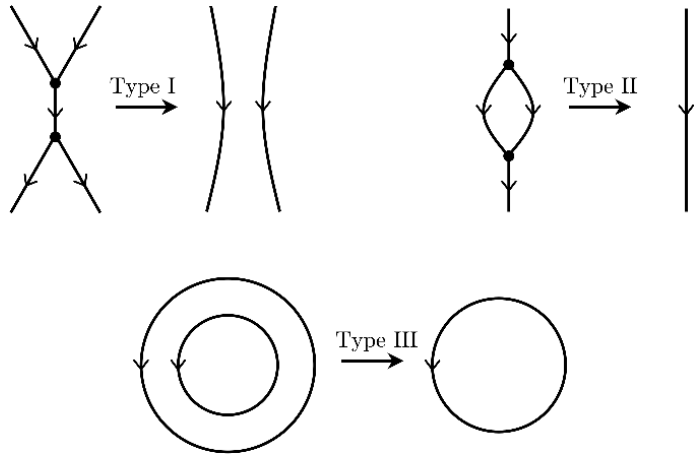}
  \caption{Reductions of an annular strand diagram}
  \label{fig:annular-reductions}
\end{figure}

\definition{A reduction of an annular strand diagram is any of the
three types of moves shown in figure \ref{fig:annular-reductions}.}

In the third move, two concentric free loops with nothing in between are combined into one.
Note that a reduction of an annular strand diagram yields an annular strand
diagram. Note also that any two annular strand diagrams for the same
element of $F$ are equivalent.

\proposition{Every annular strand diagram is equivalent to a
unique reduced annular strand diagram. \label{thm:unique-reduced-diagram}}

\noindent \emph{Proof:} We shall use Newman's Diamond Lemma (see \cite{new1}).
Clearly the process of
reduction terminates, since any reduction reduces the number of edges.  We
must show that reduction is locally confluent.

\begin{figure}[0.5\textwidth]
 \centering
  \includegraphics[height=4cm]{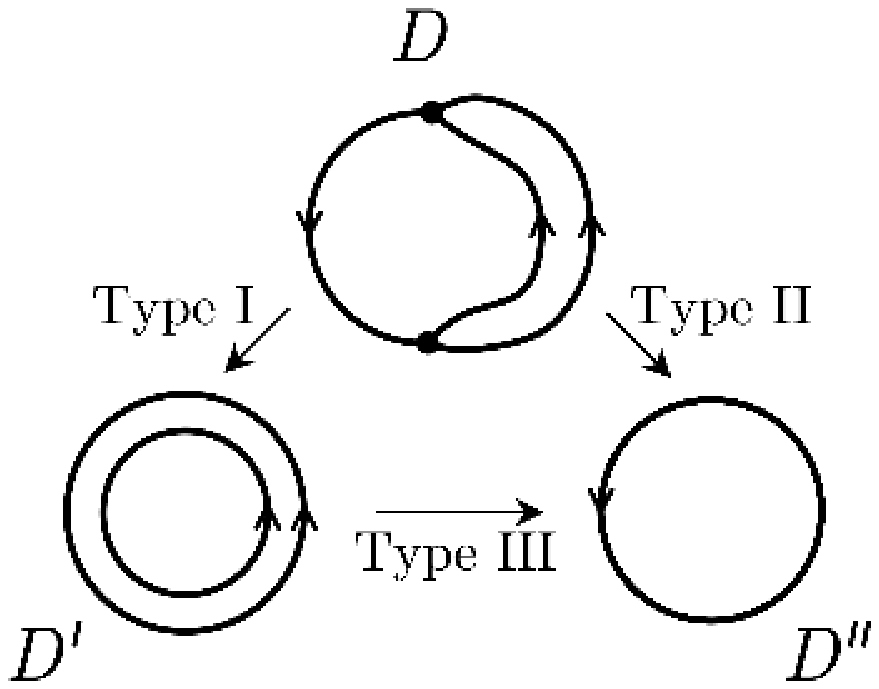}
  \caption{Diamond lemma in the annular case}
  \label{fig:annular-diagramlemma1}
\end{figure}

Suppose that a single annular strand diagram is subject to two
different reductions.  If one of these reductions is of type III, then the
two reductions commute: if the other one is of type I or II, then the reductions
must act on disjoint connected components of the diagram,
while if the other is of type III too, we can collapse all adjacent free loops
in any given order. Otherwise, both of the reductions involve the
removal of exactly two trivalent vertices.  If the reductions remove
disjoint sets of vertices, then they commute.  If the reductions share a
single vertex, then the results of the two reductions are the same (see figure
\ref{fig:diamond-lemma-F}).  Finally, it is possible for the reductions
to involve the same pair of vertices, in which case they can be resolved with
a reduction of type III (see figure \ref{fig:annular-diagramlemma1}).

\subsection{Characterization of Conjugacy in $F$}

The goal of this section is to prove the following theorem:

\theorem{Two elements of $F$ are conjugate if and only if they
have the same reduced annular strand diagram. \label{thm:characterize-conjugacy-F}}

It is not hard to see that conjugate elements of $F$ yield the same reduced annular strand
diagram. The task is to prove that two elements of $F$ with the same
reduced annular strand diagram are conjugate.

We begin with the following proposition, whose proof closely follows
the arguments of Guba and Sapir regarding conjugacy \cite{gusa1}.

\proposition{Any two cutting paths for the same annular strand
diagram yield conjugate elements of $F$. \label{thm:cutting-path-same-annular-diagram}}

\noindent \emph{Proof:} Let $\sigma_1$ and $\sigma_2$ be cutting paths for the same
annular strand diagram, and let $g_1,g_2$ be the resulting strand diagrams.
Consider the universal cover of the annulus, with the iterated preimage of
the annular strand diagram drawn upon it. Any path $\sigma$ in the annulus
lifts to a collection $\bigl\{\sigma^{\left(i\right)}:i\in\mathbb{Z}\bigr\}$
of disjoint paths in the universal cover:

\begin{figure}[0.5\textwidth]
 \centering
  \includegraphics[height=7cm]{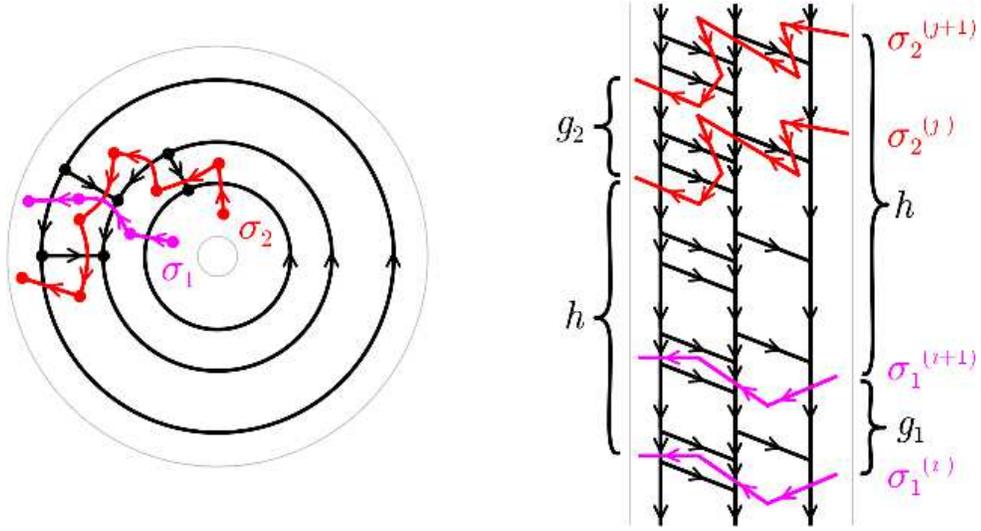}
  \caption{Creating a conjugator}
  \label{fig:universal-conjugator}
\end{figure}
Then $g_1h=hg_2$, where $h$ is the
strand diagram bounded by $\sigma_1^{\left(i\right)}$ and
$\sigma_2^{\left(j\right)}$ for some $j\gg i$, that is we choose $j$ big enough so that
the two paths $\sigma_1^{\left(i\right)}$ and $\sigma_1^{\left(i+1\right)}$ do not intersect
any of the two paths $\sigma_2^{\left(j\right)}$ and $\sigma_2^{\left(j+1\right)}$
(see figure \ref{fig:universal-conjugator}). Assume now that $g_1$ and $g_2$ are, respectively,
$(k,k)$-strand diagram and an $(m,m)$-strand diagram. We have proved that they are conjugate in Thompson's groupoid $\mathcal{F}$
(see Note \ref{thm:groupoid}).
To conclude the proof we can rewrite $g_1,g_2,h$ as $(1,1)$-strand diagrams using the right vine, that is
\[
v_k g_1 v_k^{-1} \left( v_k h v_m^{-1} \right) = \left( v_k h v_m^{-1} \right)   v_m g_2 v_m^{-1}. \; \; \square
\]

\noindent Therefore, any annular strand diagram determines a conjugacy class in $F$.

\proposition{Equivalent strand diagrams determine the same
conjugacy class.}
\noindent \emph{Proof:} Recall that a type III reduction is the composition of a type II
reduction and an inverse reduction of type I.  Therefore, it suffices to
show that the conjugacy class is unaffected by reductions of types I and II.

\begin{figure}[0.5\textwidth]
 \centering
  \includegraphics[height=5cm]{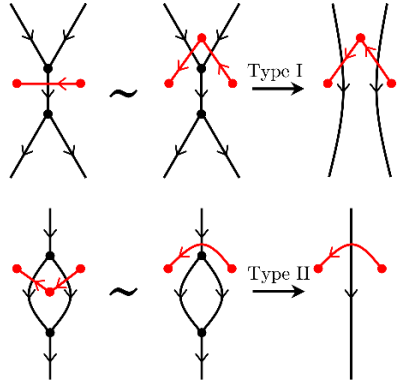}
  \caption{Moving the cutting path past the reduction area}
  \label{fig:moving-cutting-path}
\end{figure}

Given any reduction of type I or type II, it is possible to find a
cutting path that does not pass through the affected area. In particular,
any cutting path that passes through the area of reduction can be moved
(figure \ref{fig:moving-cutting-path}).

If we cut along this path, then we are performing a reduction of the
resulting strand diagram, which does not change the corresponding element
of $F$. $\square$

\noindent This proves Theorem \ref{thm:characterize-conjugacy-F}. The reduced annular strand diagram is a
computable invariant, so this gives a solution to the conjugacy problem in
$F$. We will discuss in Section \ref{sec:running-time} that the complexity of this algorithm can be implemented in
linear time.

\subsection{Structure of Annular Strand Diagrams \label{sec:structure-annular-diagrams}}

\begin{figure}[0.5\textwidth]
 \centering
  \includegraphics[height=5cm]{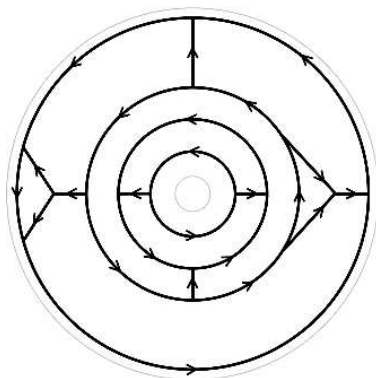}
  \caption{A reduced annular strand diagram}
  \label{fig:typical-reduced-diagram}
\end{figure}

Figure \ref{fig:typical-reduced-diagram} shows an example of a reduced annular strand diagram.

The main feature of this
diagram is the large directed cycles winding counterclockwise around the
central hole. We begin by analyzing the structure of these cycles:

\proposition{Let $L$ be a directed cycle in a reduced annular
strand diagram. Then either:
\begin{enumerate}
\item $L$ is a free loop, or
\item Every vertex on $L$ is a split, or
\item Every vertex on $L$ is a merge.
\end{enumerate}}
\noindent \emph{Proof:} Suppose $L$ has both splits and merges. Then if we trace around
$L$, we must eventually find a merge followed by a split, implying that the
annular strand diagram is not reduced. $\square$

We shall refer to $L$ as a \emph{split loop} if its vertices are all splits, and as
a merge loop if its vertices are all merges.

\proposition{For any reduced annular strand diagram:
\begin{enumerate}
\item Any two directed cycles are disjoint, and no directed cycle can
intersect itself.
\item Every directed cycle winds exactly once around the central hole.
Hence, any cutting path intersects each directed cycle exactly once.
\item Every component of the graph has at least one directed cycle.
\item Any component with only one directed cycle is a free loop.
\item By following the cutting path within a component, it is possible to order
all the directed cycles touched by the path. This order is independent of the choice of the cut.
Moreover, these concentric cycles must alternate between merge loops and split loops.
\end{enumerate}
\label{thm:combinatorial-annular-properties}}

\remark{In the next Chapter we will analyze again in more detail the connection between
strand diagrams and piecewise linear functions.
The order of directed loops defined in part 5 of Proposition \ref{thm:combinatorial-annular-properties}
follows naturally from the order of the unit interval. Compare with Theorem \ref{thm:directed-loops-fixed-point}.}

\noindent \emph{Proof:} For statement (1), observe that intersecting directed cycles
would have to merge together and then subsequently split apart, implying that the diagram
is not reduced.  For (2), recall that the directed cycles are required to wind around at least once,
and, since the graph is embedded in the plane, any closed curve
that wound around more than once would have a self-intersection.

For (3), observe that any vertex in an annular strand diagram has at
least one outgoing edge, and therefore any directed path can be extended
indefinitely. If we start a path at a vertex $p$, then the path must
eventually intersect itself as there are only finitely many vertices in the component,
which proves the existence of a directed loop in the component containing $p$.

For (4), suppose that a component of an annular strand diagram has a
split loop.  Any path that begins at a split can never again intersect the
split loop, and must therefore eventually intersect a merge loop, proving that this component
has at least two directed cycles.
Similarly, any path followed backwards from a merge loop must eventually
intersect a split loop.

For (5), observe that two adjacent concentric directed cycles in the same
component cannot both be split loops: a path starting in the region between them
must eventually cycle, for it cannot end on any of the two split loops.
Similarly, it is not possible to have two concentric
merge loops. To prove that the order does not depend on the cutting path
we start by observing that, by the proof of
Proposition \ref{thm:cutting-path-same-annular-diagram}, any two cutting paths bound a conjugator $h$
in Thompson's groupoid. By Proposition 7.2.1 in \cite{belkthesis}, $h$ must be a product
of merges and splits, so to conclude we must observe that if one cutting path can
be obtained from another by passing through a merge or a split, the order of directed cycles does not change.
This is immediately clear by looking at the moves in figure
\ref{fig:moving-cutting-path}. $\square$

\begin{figure}[0.5\textwidth]
 \centering
  \includegraphics[height=4cm]{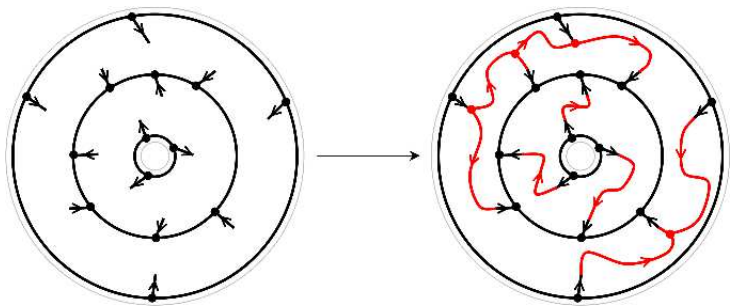}
  \caption{Constructing an annular strand diagram}
  \label{fig:construct-annular}
\end{figure}

In the next section we will define \emph{cylindrical strand diagrams} for elements of Thompson's group $T$.
With this definition and the previous proposition, we can construct a component of a reduced
annular strand diagram by drawing alternating split and merge loops, and then
filling the connections between them with unlabeled reduced cylindrical
strand diagrams (see figure \ref{fig:construct-annular}).

A general reduced annular strand diagram consists of several
concentric rings, each of which is either a free loop or a component of this
form.

\section{Conjugacy in Thompson's group $T$}

\subsection{Strand Diagrams for $T$}

\begin{figure}[0.5\textwidth]
 \centering
  \includegraphics[height=3cm]{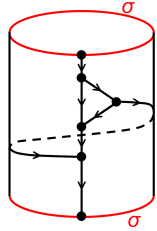}
  \caption{A cylindrical strand diagram}
  \label{fig:cylindrical-strand}
\end{figure}

We are now going to generalize to Thompson's group $T$ the diagrams and the characterization of conjugacy
that we have found for $F$. As many parts of this section are similar to the previous one, we are going to
omit some details to avoid repetition.
A \emph{cylindrical strand diagram} is a strand diagram drawn on the cylinder
$S^1\times\left[0,1\right]$, instead of on the unit square (figure
\ref{fig:cylindrical-strand}).

\begin{figure}[0.5\textwidth]
 \centering
  \includegraphics[height=3cm]{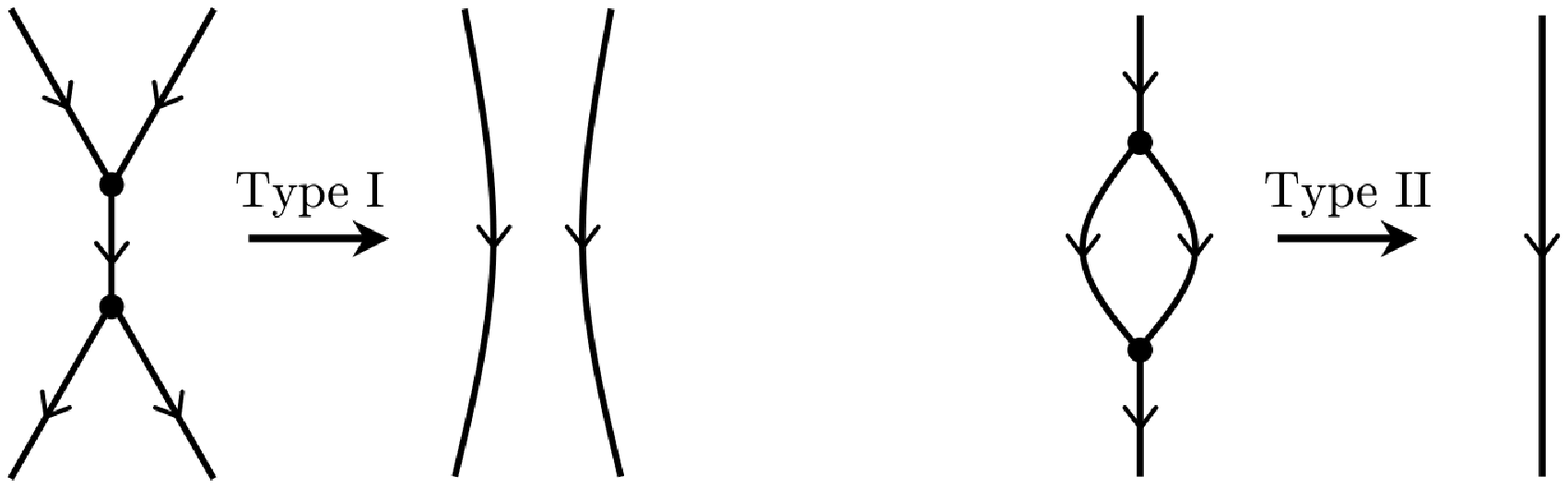}
  \caption{Reductions for a cylindrical strand diagram}
  \label{fig:cylindrical-reductions}
\end{figure}

As with strand diagrams on the square, isotopic cylindrical strand diagrams are considered
equal. We remark that isotopies of the cylinder include Dehn twists.
We recall that a \emph{Dehn twist} of
the cylinder is a homeomorphism obtained by holding the top circle rigid
while rotating the bottom circle through an angle of $2\pi$. Hence
two diagrams are equal if we can get from one to the other through a Dehn twist on the bottom.
A \emph{reduction} of a cylindrical strand diagram is either of the
moves shown in figure \ref{fig:cylindrical-reductions}.

\begin{figure}[0.5\textwidth]
 \centering
  \includegraphics[height=3cm]{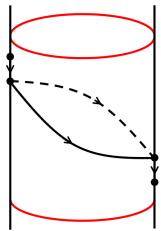}
  \caption{A cylindrical strand diagram that is not reducible}
  \label{fig:not-reducible-strand}
\end{figure}

For the second move, the two parallel edges are required to span a disc on the cylinder. In
particular, the diagram shown in figure \ref{fig:not-reducible-strand} cannot be reduced.

\begin{figure}[0.5\textwidth]
 \centering
  \includegraphics[height=4.5cm]{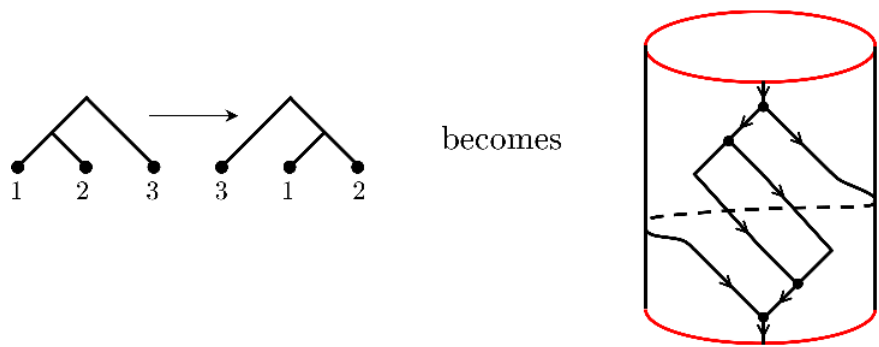}
  \caption{From a tree diagram to a cylindrical strand diagram}
  \label{fig:tree-to-cylinder}
\end{figure}

Any cylindrical strand diagram is equivalent to a unique reduced cylindrical strand
diagram. Cylindrical strand diagrams represent elements of Thompson's group $T$.
Given an element of $T$, we can construct a cylindrical strand diagram by
attaching the two trees of the tree diagram along corresponding leaves
(figure \ref{fig:tree-to-cylinder}).

\begin{figure}[0.5\textwidth]
 \centering
  \includegraphics[height=5cm]{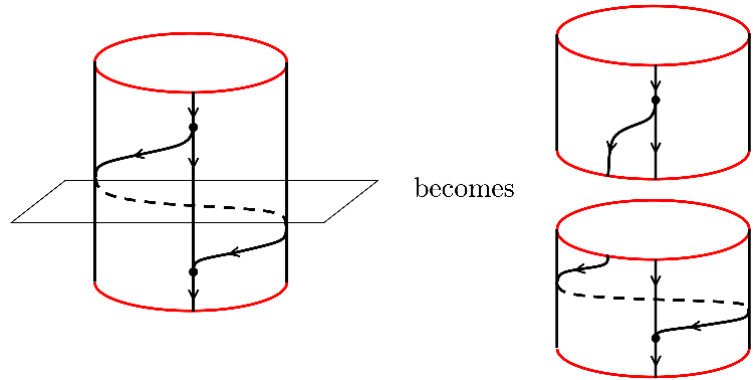}
  \caption{From a cylindrical strand diagram to a tree diagram}
  \label{fig:cylinder-to-tree}
\end{figure}

Conversely, we can cut any reduced cylindrical strand diagram along all the
edges that go from a split to a merge. This cuts the diagram into two
trees, each of which is contained in its own cylinder (figure \ref{fig:cylinder-to-tree}).

The leaves of each tree lie along a circle, and therefore the correspondence
between the leaves must be a cyclic permutation.

\begin{figure}[0.5\textwidth]
 \centering
  \includegraphics[height=4cm]{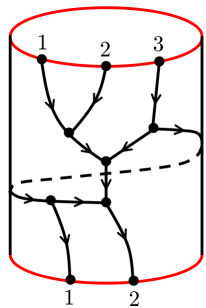}
  \caption{An $(m,n)$-strand diagram}
  \label{fig:T-mn-strand}
\end{figure}

\note{There is a slight difficulty in the definition of
cylindrical $\left(m,n\right)$-strand diagrams (figure \ref{fig:T-mn-strand}).}

If we want concatenation of
cylindrical $\left(m,n\right)$-strand diagrams to be well-defined, we must
insist on a labeling of the sources and sinks (as in the figure above).
Assuming this requirement, the set of cylindrical $\left(m,n\right)$-strand
diagrams forms a groupoid, with the group based at $1$ being Thompson's group
$T$.  Using the canonical embedding of the right vine on a cylinder, we can
then view any cylindrical $\left(m,n\right)$-strand diagram as representing
an element of $T$.

\subsection{Characterization of Conjugacy in $T$}

If we glue together the top and bottom of a cylindrical strand diagram, we
obtain a strand diagram on the torus. The common image of the top and
bottom circles is called a \emph{cutting loop}.

\definition{A \emph{toral strand diagram} is a directed graph embedded on
the torus $S^1 \times S^1$ with the following properties:
\begin{enumerate}
\item Every vertex is either a merge or a split.
\item Every directed cycle has positive index around the central
hole.
\end{enumerate}}

To make the second requirement precise, let $c$ be the cohomology class $(1,0)$ in
$H^1(S^1 \times S^1) = \mathbb{Z} \times \mathbb{Z}$.
Then a toral strand diagram is required to satisfy the condition
$c\left(\ell\right)>0$ for every directed loop $\ell$.  For a toral strand
diagram obtained from a cylinder, $c$ is precisely the cohomology class
determined by counting intersection number with the cutting loop.  For this
reason, we shall refer to $c$ as the \emph{cutting class}.

\begin{figure}[0.5\textwidth]
 \centering
  \includegraphics[height=7cm]{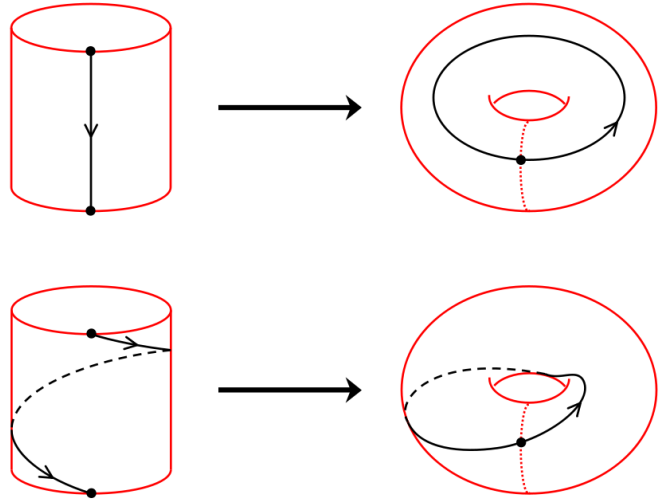}
  \caption{Toral strand diagrams that are not isotopic}
  \label{fig:Dehn-Twist-Picture}
\end{figure}

The cutting class is related to a slight difficulty in defining the
notion of equality for toral strand diagrams. Because a Dehn
twist of the cylinder is isotopic to the identity map, two cylindrical strand
diagrams that differ by a Dehn twist are isotopic and hence considered equal.
However, the resulting toral strand diagrams are not isotopic (for example,
see figure \ref{fig:Dehn-Twist-Picture}).

 This
difficulty arises because the Dehn twist descends to a nontrivial
homeomorphism of the torus (i.e. a homeomorphism that is not isotopic to the
identity). Using the standard basis for the first cohomology group of the torus (since
$c=\left(1,0\right)$), this Dehn twist acts as
$\scriptstyle \begin{pmatrix}
1&1\\
0&1
\end{pmatrix}$.
Therefore, we must consider two toral strand diagrams equal if their
isotopy classes differ by a Dehn twist of the form
$\scriptstyle \begin{pmatrix}
1&n\\
0&1
\end{pmatrix}$. This is equivalent to the following convention:

\convention{Let $S_1$ and $S_2$ be two strand diagrams embedded on
the torus $\mathbb{T}$.  We say that $S_1$ and $S_2$ are \emph{equal} if there
exists an orientation-preserving homeomorphism
$h\colon\mathbb{T}\rightarrow\mathbb{T}$ such $h^\ast\left(c\right)=c$ and
$h\left(S_1\right)=S_2$. \label{thm:convention-toral}}

\begin{figure}[0.5\textwidth]
 \centering
  \includegraphics[height=3cm]{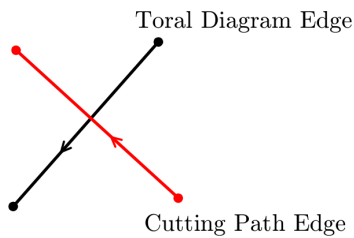}
  \caption{Orientation of the cutting class}
  \label{fig:toral-cutting-orientation}
\end{figure}

\definition{A \emph{cutting loop} for a toral strand diagram is a simple
continuous loop in the torus that satisfies the following conditions:
\begin{enumerate}
\item The loop is dual to the cohomology class $c$.
\item The loop does not pass through any vertices of the strand
diagram.
\item The loop intersects edges of the strand diagram transversely,
with the orientation shown in figure \ref{fig:toral-cutting-orientation}.
\end{enumerate}}

\proposition{Cutting a toral strand diagram along a cutting loop
yields a $\left(k,k\right)$-strand diagram embedded on the
cylinder. $\square$}

\theorem{Every toral strand diagram has a cutting
loop. $\square$}

\begin{figure}[0.5\textwidth]
 \centering
  \includegraphics[height=6cm]{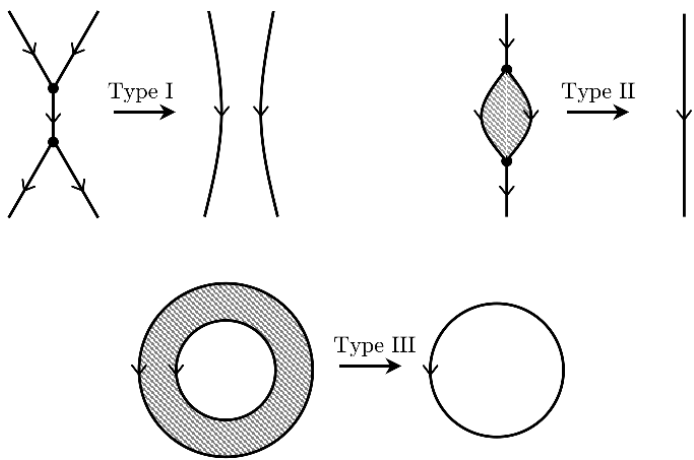}
  \caption{Reductions for a toral strand diagram}
  \label{fig:toral-reductions}
\end{figure}

\definition{A reduction of a toral strand diagram is any of the
three types of moves shown in figure \ref{fig:toral-reductions}}

In the second move, the two edges of the bigon are required to span a disc,
and in the third move the two loops must be the boundary of an annular region.
Two toral strand diagram are \emph{equivalent} if one can obtained from the other
via a sequence of reductions and inverse reductions.

\proposition{Every toral strand diagram is equivalent to a unique
reduced toral strand diagram.}

\noindent \emph{Proof:} The argument on reductions used in the proof of Proposition
\ref{thm:unique-reduced-diagram} can be extended to this case without additional details and so we omit it. $\square$

Gluing the top and the bottom circles of a cylindrical strand diagram
is not well defined in general. However by Convention \ref{thm:convention-toral}, all
resulting toral diagrams are equal.

\theorem{Two elements of $T$ are conjugate if and only if they
have the same reduced toral strand diagram.}

\noindent \emph{Proof:} Our convention for equality of strand diagrams guarantees that
any two conjugate elements of $T$ yield the same reduced toral strand diagram.

We claim that any two cutting loops for the same toral strand diagram
yield conjugate elements of $T$.  Suppose we are given cutting loops
$\ell_1$ and $\ell_2$, and consider the cover of the torus corresponding to
the subgroup ker$\left(c\right)\leq\pi_1\left(\mathbb{T}\right)$. This
cover is an infinite cylinder, with the deck transformations
$\pi_1\left(\mathbb{T}\right)/$ker$\left(c\right)\cong\mathbb{Z}$ acting as
vertical translation. Each of the loops $\ell_i$ lifts to an infinite
sequence $\bigl\{\ell_i^{\left(j\right)}\bigr\}_{j\in\mathbb{Z}}$ of loops in
this cover, and the region between $\ell_i^{\left(j\right)}$ and
$\ell_i^{\left(j+1\right)}$ is the cylindrical strand diagram $f_i$ obtained
by cutting the torus along $\ell_i$. \ It follows that $f_1g=gf_2$, where $g$
is the cylindrical strand diagram between $\ell_1^{\left(j\right)}$ and
$\ell_2^{\left(k\right)}$ for some $k\gg j$.

Clearly reductions do not change the conjugacy class described by a
toral strand diagram, and therefore any two elements of $T$ with the same
reduced toral strand diagram are conjugate. $\square$

\subsection{Structure of Toral Strand Diagrams}

The following section is an analogue for $T$ of Section ref{sec:structure-annular-diagrams} for $F$.
Given an element $f\in T$, the structure of the toral strand diagram for $f$
is closely related to the dynamics of $f$ as a self-homeomorphism of the
circle. In this section we analyze the structure of toral strand diagrams,
and in the next we show how this structure is related to the dynamics of an
element.

We begin by noting some features of annular strand diagrams that remain
true in the toral case:

\proposition{For any reduced toral strand diagram.
\begin{enumerate}
\item Any directed cycle is either a free loop, a split loop, or a merge loop.
\item Any two directed cycles are disjoint, and no directed cycle can intersect itself.
\item Every component of the graph has at least one directed cycle, and any
component with only one directed cycle is a free loop. $\square$
\end{enumerate}}

In an annular strand diagram, each directed cycle winds around the central
hole exactly once, and the components of the diagram form concentric rings.
The structure of a toral strand diagram is more complicated.

\proposition{Let $c \in H^1\left(\mathbb{T}\right)$ denote the cutting class.
Without loss of generality, we can assume that $c=\left(1,0\right)$.
Then any two directed cycles represent the same
element $\left(n,k\right)\in H_1\left(\mathbb{T}\right)$, where $n>0$ and $k$
and $n$ are relatively prime.}

\noindent \emph{Proof:} By the definition of a toral strand diagram, $n>0$ for any
directed cycle. \ Any two disjoint nontrivial loops on a torus are homotopic,
and therefore any two directed cycles must have the same $\left(n,k\right)$.
Furthermore, since a directed cycle cannot intersect itself, $n$ and $k$
must be relatively prime. $\square$

Note that the number $k$ is not uniquely determined. Specifically, recall
that two strand diagrams that differ by the Dehn twist
$\scriptstyle \begin{pmatrix}
1&0\\
1&1
\end{pmatrix}$ are equal. (This matrix is the transpose of the earlier
matrix, since we are now considering the action on homology.) Applying this
Dehn twist to a diagram whose directed cycles are $\left(n,k\right)$ yields a
diagram whose cycles are $\left(n,k+n\right)$, so the number $k$ is only
well-defined modulo $n$.

We will always assume that $0\leq k<n$.  The reduced fraction
$k/n\in[0,1)$ is called the \emph{rotation number} of a toral strand diagram.  It is possible to
show that this corresponds to the dynamical rotation number of a homeomorphism $f\in T$
(see Chapter \ref{chapter6} for the definition).

\propositionname{Ghys-Sergiescu, \cite{ghys2}}{Every element of $T$ has a periodic point. \label{thm:ghys-sergiescu}}

\noindent \emph{Proof.} We delay this proof to Section \ref{sec:mn-strand-diagrams-dynamical}.
See Proposition \ref{thm:ghys-sergiescu-2}. $\square$

We remark that the previous result has been recently proved again by Calegari in \cite{caleg1}.
Bleak and Farley \cite{bleak-farley} also have a proof of this result using ``revealing
tree-pair diagrams'' as introduced by Brin \cite{brin4}.

\propositionname{Burillo-Cleary-Stein-Taback, \cite{burclearysteintaback}}{For any
positive integer $n$, let $c_n$ be the $(n,n)$-strand diagram
\begin{figure}[0.5\textwidth]
 \begin{center}
  \includegraphics[height=4cm]{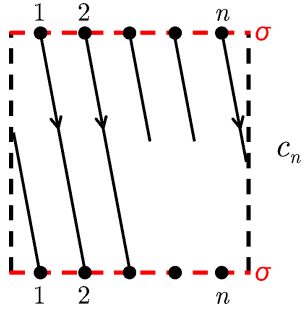}
 \end{center}
 \caption{The torsion element $c_n$}
 \label{fig:torsion-element}
\end{figure}
Then any torsion element of $T$ is conjugate to a power of some $v_n^{-1} c_n v_n$ for some integer $n$. \label{thm:T-torsion}}

\noindent \emph{Proof.} If $f \in T$ is torsion, then $f$ cannot have any merge or split loops. Hence the reduced toral
strand diagram for $f$ has only free loops. By opening the reduced toral strand diagram through the cutting line,
we obtain a strand diagram which is a power $c_n^k$ for some $n$. We can obtain the corresponding $(1,1)$-strand diagram
by writing $v_n^{-1} c_n^k v_n$.  $\square$

It is not too hard to see that, for any $1 \le k < n$, the
element $c_n^k$ has rotation number $k/n$. This proves that, for any rational number $k/n \pmod{1}$ there is
an element of $T$ with rotation number $k/n$ (another result due to Ghys and Sergiescu in \cite{ghys2}).
%
%

\section{Conjugacy in Thompson's group $V$}

\subsection{Strand Diagrams for $V$}

\definition{An \emph{abstract strand diagram} is an acyclic directed
graph, together with a cyclic ordering of the edges incident on each vertex,
and subject to the following conditions:
\begin{enumerate}
\item There exists a unique univalent source and a unique univalent sink.
\item Every other vertex is either a split or a merge.
\end{enumerate}}

\begin{figure}[0.5\textwidth]
 \centering
  \includegraphics[height=4cm]{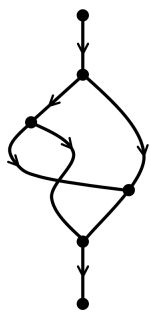}
  \caption{An abstract strand diagram}
  \label{fig:V-diagram}
\end{figure}

The cyclic orderings of the edges allow us to distinguish between the left
and right outputs of a split, and between the left and right inputs of a
merge. We can draw an abstract strand diagram as a directed graph in the plane
with edge crossings (see figure \ref{fig:V-diagram}).

\begin{figure}[0.5\textwidth]
 \centering
  \includegraphics[height=5cm]{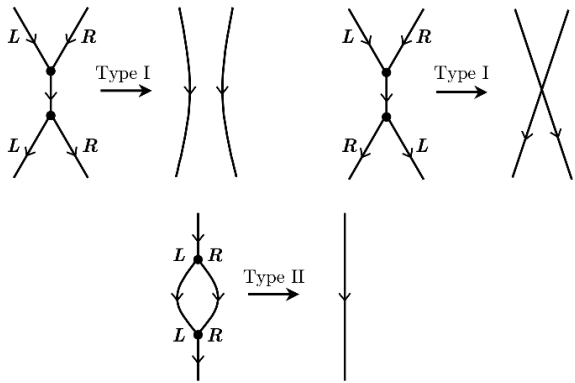}
  \caption{Reductions for abstract strand diagrams}
  \label{fig:V-reductions}
\end{figure}

By convention, the edges incident on a vertex are always
drawn so that the cyclic order is counterclockwise. Reductions in this setting are defined via
the drawing of the graph in the plane, because we need the vertices to be oriented in the same way of the plane.
A \emph{reduction} of an abstract strand diagram (drawn in the plane) is either of the
moves of figure \ref{fig:V-reductions}.

\begin{figure}[0.5\textwidth]
 \centering
  \includegraphics[height=2.5cm]{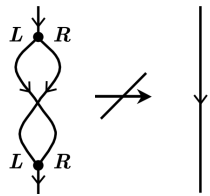}
  \caption{Non valid reduction.}
  \label{fig:newbad-things}
\end{figure}

The first two moves are the same kind of move, drawn differently depending on the embedding in the plane.
The cyclic order of the vertices must be exactly as shown above. The move shown in figure
\ref{fig:newbad-things} is not valid.

\begin{figure}[0.5\textwidth]
 \centering
  \includegraphics[height=4cm]{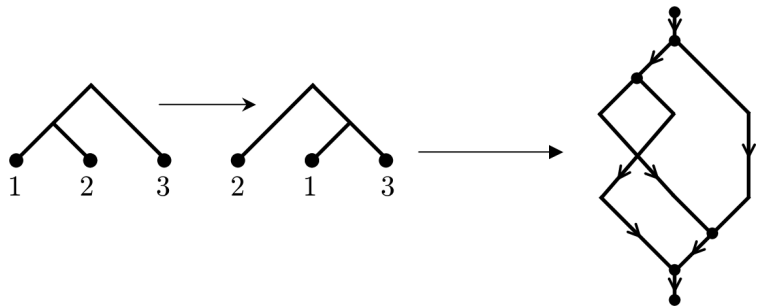}
  \caption{From a tree diagram to an abstract strand diagram}
  \label{fig:tree-to-abstract-strand}
\end{figure}

Every abstract strand diagram is equivalent to a unique reduced abstract
strand diagram. Abstract strand diagrams represent elements of Thompson's group $V$.
Given an element $f\in V$, we can construct an abstract strand diagram for
$V$ by attaching the two trees of a tree diagram for $f$ along corresponding
leaves (figure \ref{fig:tree-to-abstract-strand}).

Conversely, any reduced abstract strand diagram can be cut along all the
edges that go from splits to merges to yield a tree diagram.
Assuming we label the sources and sinks, the set of abstract
$\left(m,n\right)$-strand diagrams forms a groupoid, and elements of this
groupoid can viewed as representing elements of Thompson's group $V$.

\subsection{Characterization of Conjugacy in $V$}

If we glue together the sources and sinks of an abstract strand diagram, we
obtain a directed graph whose vertices are all merges and splits. The
images of the original sources and sinks now fall in the interiors of certain
edges, and are called the \emph{cut points}.  Note that a single edge may contain
more than one cut point. The function that measures the number of cut
points in each edge is a $1$-cocycle, and therefore yields a cohomology class
$c$, which we call the \emph{cutting class}.

\definition{A \emph{closed strand diagram} is a triple
$\left(D,o,c\right)$, where
\begin{enumerate}
\item $D$ is a directed graph composed of splits and merges,
\item $o$ is a cyclic ordering of the edges around each vertex of
$D$, and
\item $c$ is an element of $H^1\left(D\right)$ satisfying
$c\left(\sigma\right)>0$ for every directed cycle $\sigma$.
\end{enumerate}}

\begin{figure}[0.5\textwidth]
 \centering
  \includegraphics[height=5cm]{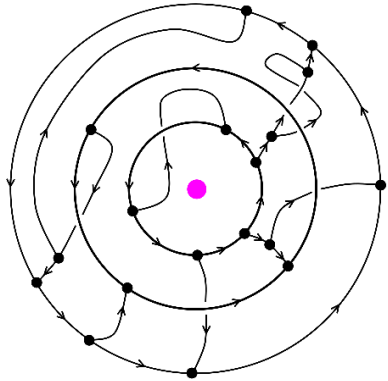}
  \caption{A closed strand drawn on the punctured plane}
  \label{fig:diagram-punctured-plane}
\end{figure}

The cohomology class $c$ is called the \emph{cutting class}.  To make our arguments
as accessible as possible, we will use a very geometric approach to
cohomology.  In particular, we will make heavy use of the following well known result:
for any CW-complex $X$, there is a natural one-to-one
correspondence between elements of $H^1\left(X\right)$ and homotopy classes
of maps from $X$ to the punctured plane.
Using the above theorem, we can represent a closed strand diagram as a graph
with crossings drawn on the punctured plane (figure \ref{fig:diagram-punctured-plane}).

The cohomology class $c$ is given
by winding number around the puncture. By convention, we always draw closed
strand diagrams so that the cyclic order of the edges around each vertex is
counterclockwise.

\begin{figure}[0.5\textwidth]
 \centering
  \includegraphics[height=3cm]{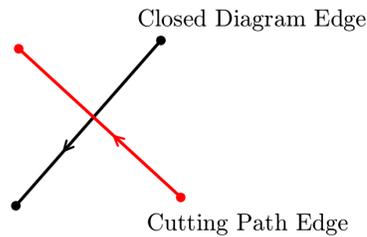}
  \caption{Orientation of the cutting class}
  \label{fig:closed-cutting-orientation}
\end{figure}

\definition{Given a drawing of a closed strand diagram, a \emph{cutting
line} is a continuous path going from the center to the outer region so that it does not
intersect any vertex but it intersects the edges of the diagram transversely,
with the orientation shown in figure \ref{fig:closed-cutting-orientation}}.

The sequence $p_1,\ldots,p_n$ of points
on the graph cut by the line is called a \emph{cutting sequence}.
Note that we can ``cut'' along a cutting sequence to obtain an ordered
abstract $\left(k,k\right)$-strand diagram.
The above definition is very geometric.  Here is a combinatorial
description of cutting sequences:

\proposition{Let $p_1,\ldots,p_n$ be a sequence of points lying in
the interiors of the edges of a closed strand diagram. Then
$p_1,\ldots,p_n$ is a cutting sequence if and only if the function
\begin{equation*}
e\;\longmapsto\;\text{\#}\left\{i:p_i\in e\right\}
\end{equation*}
is a $1$-cochain representing the cutting class $c$. $\square$}

\theorem{Every closed strand diagram has a cutting sequence.}

\noindent \emph{Proof:} From theorem \ref{thm:appendix} in the appendix, there exists a non-negative,
integer-valued cochain $\alpha$ representing $c$. Then the sequence
$p_1,\ldots,p_n$ can be constructed by choosing $\alpha\left(e\right)$ points
from each edge $e$. $\square$

\begin{figure}[0.5\textwidth]
 \centering
  \includegraphics[height=6cm]{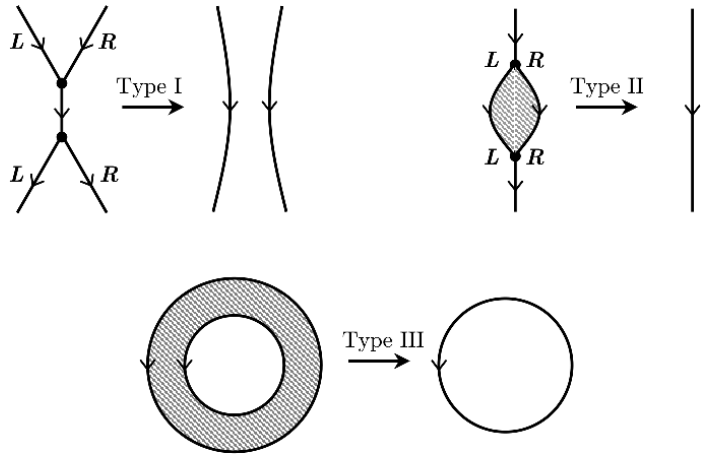}
  \caption{Reductions for closed strand diagrams}
  \label{fig:closed-strand-reductions}
\end{figure}

\definition{A \emph{reduction} of a closed strand diagram is any of the
three moves shown in figure \ref{fig:closed-strand-reductions}.}

In the second move, the loop spanned by the bigon must lie in the kernel of
$c$, i.e. the parallel edges must be homotopic in the punctured plane.  In
the third move, we require that the difference of the two loops lie in the
kernel of $c$, or equivalently that the two loops have the same winding
number around the puncture.

\begin{figure}[0.5\textwidth]
 \centering
  \includegraphics[height=3cm]{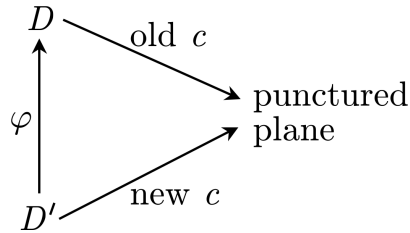}
  \caption{Cutting classes and reductions}
  \label{fig:punctured-plane-commutative-diagram}
\end{figure}

In each of the three cases, the reduced graph $D'$ inherits a cutting
class in the obvious way. For a type I reduction, the new cutting class is
$\varphi^\ast\left(c\right)$, where $\varphi$ is the obvious map
$D'\rightarrow D$ (figure \ref{fig:punctured-plane-commutative-diagram}).

For a reduction of type II, there are two obvious maps
$D'\rightarrow D \to \{$punctured plane$\}$: we send
the reduced edge to any of the two sides of the bigon.
These maps are homotopic, and therefore yield the same
homomorphism $H^1\left(D\right)\rightarrow H^1\left(D'\right)$. The same
holds for reductions of type III.

\proposition{Every closed strand diagram is equivalent to a unique reduced
closed strand diagram.}

\begin{figure}[0.5\textwidth]
 \centering
  \includegraphics[height=5cm]{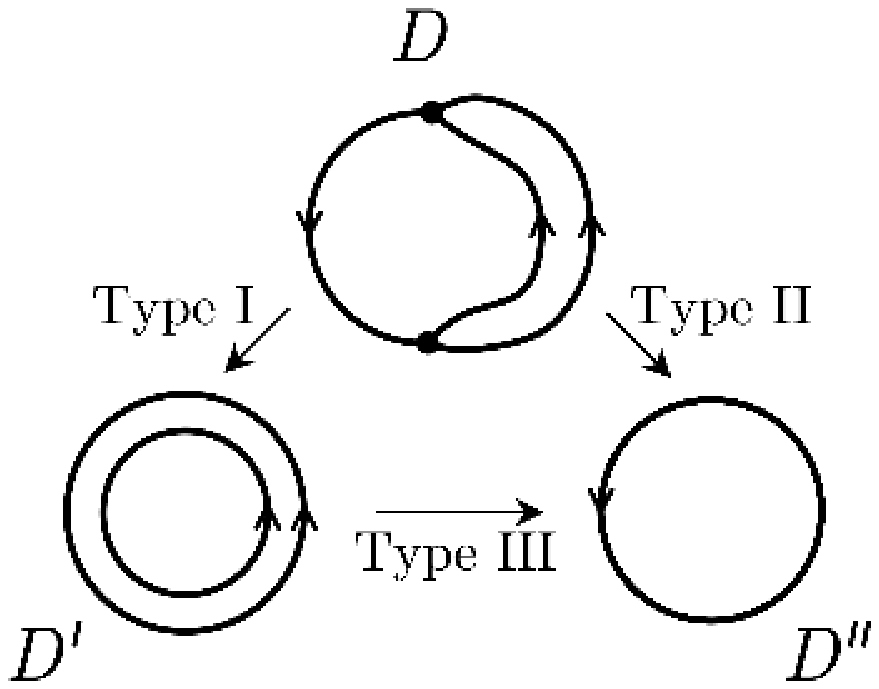}
  \caption{Diamond Lemma}
  \label{fig:annular-diamondlemma2}
\end{figure}

\noindent \emph{Proof:} We must show that reduction is locally confluent, keeping careful
track of the fate of the cohomology class $c$. Suppose that a single closed strand diagram is subject to two different
reductions.  If one of these reductions is of type III or they remove disjoint sets
of vertices, then they commute. If the reductions share a single vertex, then the results of the two
reductions are the same, as seen in previous cases (see figure \ref{fig:annular-diagramlemma1}
in Proposition \ref{thm:unique-reduced-diagram}).
Note in particular that the map $D'\rightarrow D$ obtained from the type I
reduction is homotopic in the punctured plane to the pair of maps
$D'\rightarrow D$ obtained from the type II reduction.
Finally, it is possible for the reductions to involve the same pair of
vertices, in which case they can be resolved with a reduction of type
III (figure \ref{fig:annular-diamondlemma2}).

Again, observe that the two maps $D''\rightarrow D$ obtained from the
type II reduction are homotopic to the two composite maps
$D'' \rightrightarrows D'\rightarrow D$. $\square$

\lemma{Conjugate elements of $V$ yield isomorphic reduced closed strand diagram.}

\begin{figure}[0.5\textwidth]
 \centering
  \includegraphics[height=3cm]{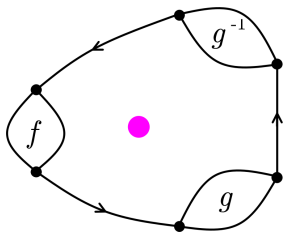}
  \caption{Same reduced closed strand diagram}
  \label{fig:conjugate-same-diagram}
\end{figure}

\noindent \emph{Proof.} Let $f,g\in V$.  Then figure \ref{fig:conjugate-same-diagram}
is a closed strand diagram for both $f$ and $g^{-1}fg$.

Since $f$ and $g^{-1}fg$ share a closed strand
diagram, they must have the same reduced closed strand diagrams. $\square$

\theorem{Two elements of $V$ are conjugate if and only if they
have isomorphic reduced closed strand diagram.}

\noindent \emph{Proof:} We claim that any two cutting sequences
$\left\{p_1,\ldots,p_m\right\},\left\{q_1,\ldots,q_n\right\}$ for isomorphic
closed strand diagram $S$ yield conjugate elements of $V$.  Consider the
infinite-sheeted cover of the strand diagram obtained by lifting to the
universal cover of the punctured plane.  (Abstractly, this is the cover
corresponding to the subgroup ker$\left(c\right)$ of $\pi_1\left(D\right)$.)
If we arrange $S$ on the punctured plane so that the points
$\left\{p_1,\ldots,p_n\right\}$ lie on a single radial line $\ell$, then the
lifts of this line cut the cover into infinitely many copies of the abstract
strand diagram $f$ obtained by cutting $S$ along
$\left\{p_1,\ldots,p_n\right\}$.  Specifically, the points
$\left\{p_1,\ldots,p_m\right\}$ have lifts
$\left\{p_1^{\left(i\right)},\ldots,p_m^{\left(i\right)}\right\}_{i\in%
\mathbb{Z}}$, with the $i$th copy of $f$ having
$\left\{p_1^{\left(i\right)},\ldots,p_m^{\left(i\right)}\right\}$ as its
sources and
$\left\{p_1^{\left(i+1\right)},\ldots,p_m^{\left(i+1\right)}\right\}$ as its
sinks. \ Similarly, if we homotope $S$ so that
$\left\{q_1,\ldots,q_n\right\}$ lie on a single radial line, we obtain a
decomposition of the cover into pieces isomorphic to the abstract strand
diagram $g$ obtained by cutting $S$ along $\left\{q_1,\ldots,q_n\right\}$.
It follows that $fh=hg$, where $h$ is the abstract strand diagram lying
between $\left\{p_1^{\left(i+1\right)},\ldots,p_m^{\left(i+1\right)}\right\}$
and $\left\{q_1^{\left(j\right)},\ldots,q_n^{\left(j\right)}\right\}$ for
some $i\ll j$. $\square$

\subsection{Structure of Abstract Closed Strand Diagrams}

Most of the results seen before, generalize to this setting. For example, reduced abstract
closed strand diagram have the same combinatorial structure as toral strand diagrams
(i.e. They must contain a directed cycle, all cycles must be disjoint, etc.).

\theoremname{Brin, \cite{brin4}}{Let $f \in V$, then:
\begin{enumerate}
\item $f$ has a periodic point.
\item If $f$ is torsion, it is conjugate to a permutation.
\item There is an integer $n(f)$ so that every finite orbit of $f$ has no more than $n(f)$ elements.
\end{enumerate}}

\noindent \emph{Proof.} (i) follows along the same lines as Proposition \ref{thm:ghys-sergiescu}.
(ii) follows along the same lines as Proposition \ref{thm:T-torsion}. For (iii) we
recall that the reduced tree diagram of $f$ gives a partition of the interval $[0,1]$ and
let $p$ be a point in a finite orbit. Then $p$ is in some interval $[a,b]$ of the partition and
there is a power $f^k$ such that $f^k$ sends $[a,b]$ into itself, and so the abstract closed strand
diagram has a directed cycle passing through the vertex corresponding to
$[a,b]$. Thus the orbit of $p$ cannot have more than $k$ points,
where $k$ is the length of the cycle containing $[a,b]$. Thus the length of any finite orbit
is bounded by $n(f)$ maximum length of a directed cycle of the abstract closed strand diagram. $\square$

\subsection{Generalizations and Conjectures}

It seems possible to take this unified point of view and generalize it
to a more abstract setting. We start by conveying an intuition. We observe that,
in the case of $F$, diagrams were embedded
in the unit square $I \times I$ while in the case of $T$ they were embedded in the cylinder
$S^1 \times I$. Similarly, diagrams in $V$ could be regarded as
embedded in the space $\mathbb{R}^3 \times I$. We can try to define strand diagrams
embedded in a space $M \times I$ where $M$ is a suitable space (for example, a differentiable manifold
or a CW-complex). However, issues arise when we try to define reductions since isotopy
can move the cyclic ordering of a vertex and switch left with right.
We remember that all of our definitions of reductions took place on a particular surface.

\begin{figure}[0.5\textwidth]
 \centering
  \includegraphics[height=4cm]{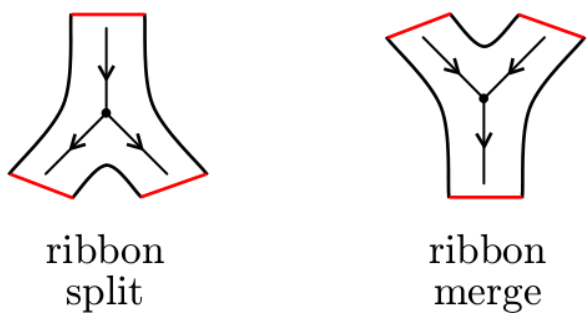}
  \caption{Ribbon splits and merges}
  \label{fig:ribbon-split-merge}
\end{figure}

One can define \emph{ribbon surfaces} instead of strand diagrams. They are diagrams obtained
by ``fattening'' a strand diagram to become an orientable surface with boundary. More precisely,
we can take an abstract strand diagram as defined in $V$ and embed it in $M \times I$ so that
the source is the unique point contained in $M \times \{0\}$ and the sink is the unique point contained
in $M \times \{1\}$, then we attach the following
\emph{ribbon splits} and \emph{ribbon merges} at each vertex (see figure \ref{fig:ribbon-split-merge}).

\begin{figure}[0.5\textwidth]
 \centering
  \includegraphics[height=4cm]{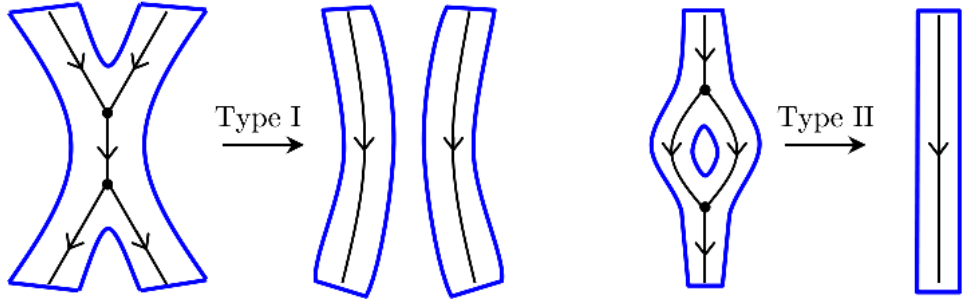}
  \caption{Ribbon reductions}
  \label{fig:ribbon-reductions}
\end{figure}

We can complete the surface by attaching rectangles along the edges and gluing them on the inputs and the outputs
of the ribbon splits and merges. Now we have a well defined surface inside $M \times I$ and we can orient it, starting
from the top. We can still talk about ``trivalent parts'' of the surface, to mean the corresponding ``vertices'' of
this surface. We can now define isotopy and equivalence of the surface inside $M \times I$ and define a product of surfaces.
Finally we can define reductions of a surface: they will be defined by the motions shown in figure \ref{fig:ribbon-reductions}.

\noindent where the depicted pieces of the ribbon surface are assumed to have been embedded in the plane, according to the orientation
of these pieces. We observe that the second reduction assumes that there exists a $2$-cell contained in $M \times I$
that can be attached to the ribbon surface and whose interior does not intersect the ribbon surface at any point.
We can continue this line of thinking by defining \emph{closed ribbon surfaces} inside the space $M \times S^1$,
by gluing the ends of a ribbon surface in $M \times I$. We can thus define a third kind of reduction if we can glue two circular
strips by another circular strip that is entirely contained in $M \times S^1$ and whose interior does not intersect the ribbon
surface at any point.

It is still possible to prove that each of the closed diagrams is equivalent to a unique reduced diagram and
that any two conjugate ribbon surfaces in $M \times I$ give rise to the same reduced closed ribbon surface.
Proving the converse seems also possible, although it will depend on the properties of the space $M$,
however, one can still theoretically follow the ``infinite sheeted cover'' argument and try to prove it in this setting.

\begin{figure}[0.5\textwidth]
 \centering
  \includegraphics[height=4cm]{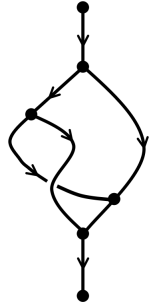}
 \caption{An element of Thompson's group $BV$}
\label{fig:element-BV}
\end{figure}

Using this new generalization, if we choose the space $M$ to be $\mathbb{R}^2$ we get a group
of ribbon surfaces which resembles the braided Thompson's group $BV$.
Brin has studied some properties and presentations of Thompson's group $BV$ in the paper \cite{brin3}
while Burillo and Cleary have described some of its metric properties in \cite{burcleary1}. The descriptions
given previously by Brin and Burillo were describing elements as
as pairs of trees with a braid connecting the leaves (see figure \ref{fig:element-BV}).

It would be interesting to study the particular case of ribbon surfaces in $\mathbb{R}^2 \times I$.
Our procedure is a generalization of
the point of view of binary tree diagrams for Thompson's groups. However, one could extend
it to define strand diagrams in $M \times I$ where all the vertices are $n$-ary splits and merges
to recover the same conjugacy results for higher dimensional Thompson's groups as defined
by Brin in \cite{brin4}.

We want to conclude this section, by describing another tempting generalization.
Guba and Sapir have introduced a class of groups called \emph{diagram groups}
in their monograph \cite{gusa1}. These groups are defined by elements that can be represented
as diagrams in the plane and Thompson's group $F$ can be described as one such group.
Our solution of the conjugacy problem for $F$ is heavily inspired by and similar to
the one they give for diagrams groups. Guba and Sapir describe other families of
diagram groups that live in other spaces and they recover $T$ and $V$ inside these classes.
They suggest that their work could be extended to diagram groups in these settings.

It is possible to define general \emph{ribbon-diagram groups} in $M \times I$ and it would be interesting to see
if their solution for the conjugacy problem could be extended to this general setting.

\conjecture{The strategy to solve the conjugacy problem for strand diagram groups
can also be used to solve the conjugacy problem for ribbon-diagram groups defined in $M \times I$.}

\section{Running Time \label{sec:running-time}}

In this section, we study the complexity of our solution of the conjugacy problem
for Thompson's groups $F, T$ and $V$. We start by sketching a proof of the following result:

\theorem{There exists a linear-time algorithm to determine whether
two elements of $F$ are conjugate.}

We assume that the two elements of $F$ are given as words in the generating
set $\left\{x_0,x_1\right\}$.  ``Linear time'' means that the algorithm
requires $O\left(N\right)$ operations, where $N$ is the sum of the lengths of
these words. We shall use the algorithm of Hopcroft and Wong (see \cite{hopwon}):

\theoremname{Hopcroft and Wong}{There exists a linear-time algorithm to
determine whether two planar graphs are isomorphic. $\square$ \label{thm:hopwon}}

We remark that Guba and Sapir had already proven that their solution to the
conjugacy problem for diagram groups had the same complexity of the isomorphism problem
for planar graphs (private communication). Thus their solution along with Theorem \ref{thm:hopwon}
give again a linear time algorithm.

\proposition{There exists a linear-time algorithm to determine
whether two (reduced) annular strand diagrams are isotopic.}

\begin{figure}[0.5\textwidth]
 \centering
  \includegraphics[height=3cm]{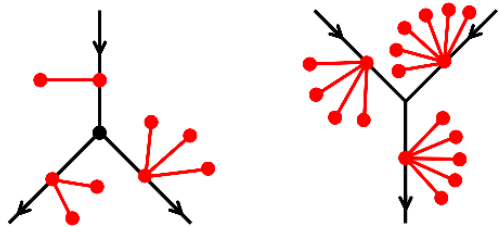}
 \caption{Decorating the annular strand diagram.}
\label{fig:running-time-decoration}
\end{figure}

\noindent \emph{Proof:} We must show that isotopy of connected annular strand diagrams
reduces to isomorphism of planar graphs in linear time.  If the given strand
diagrams are disconnected, then we may check isotopy of the components
separately.  It therefore suffices to prove the proposition in the connected
case. Given a strand diagram, subdivide each edge into three parts, and
attach new edges around each merge and split as drawn in figure \ref{fig:running-time-decoration}.

This new graph can be
constructed in linear time, and its isomorphism type completely determines
the isotopy class of the original reduced annular strand diagram.  In
particular, the decorations determine both the directions of the original
edges and the cyclic order of the original edges around each merge or
split. $\square$

All that remains is to show that the reduced annular strand diagram for an
element of $F$ can be constructed in linear time. This requires two steps:

\begin{enumerate}
\item Construct a strand diagram for the element.
\item Reduce the resulting annular strand diagram.
\end{enumerate}

The first step is easy to carry out in linear time: given a word in
$\left\{x_0,x_1\right\}$, simply concatenate the corresponding strand
diagrams for the generators and their inverses. No reduction is necessary
in this phase. For the second step, observe that any reduction of a strand diagram
reduces the number of vertices, and therefore only linearly many reductions
are required.  However, it is not entirely obvious how to search for these
reductions efficiently.

\proposition{Suppose that any one reduction can be performed in
constant time.  Then a given annular strand diagram $G$ can be reduced in linear
time.}

\noindent \emph{Proof:} We give a linear-time algorithm for performing all the necessary
type I and type II reductions.  Any required type III reductions can be
performed afterwards.

We can write $G$ as a set of vertices $V=\{v_1, \ldots, v_k \}:=V_1$.
We build inductively new sets of vertices $V_i$. To build the sequence $V_{i+1}$,
we read and classify every vertex of $V_i$. We let $R_i$ be the set of vertices
of $V_i$ which are reducible. By definition $R_i$ will have an even number of vertices,
since a vertex is reducible if there is an adjacent vertex which forms a reduction in the diagram.
Then we define $V_{i+1}$ to be all the vertices of $V \setminus \left(\bigcup_{j=1}^i R_j \right)$ which are adjacent to a vertex in $R_i$.

This algorithm goes to look for vertices which were not reducible at the $i$-th step, but might have become
reducible at the $i+1$-th step. We repeat this process until we find an $m$ such that $R_m= \emptyset$.
By construction,
\[
|R_{i+1}| \le |V_{i+1}| \le 4 |R_i|
\]
since every point involved in a reduction is adjacent to at most 3 vertices, one of which will be reduced.
In other words, for each pair of vertices that we reduce, we might have to reinsert up to 4 vertices which were
previously not reducible. On the other hand, it is obvious that
\[
\sum_{i=1}^m |R_i| \le |V|.
\]
The final cost of the computation is thus given by

\[
\sum_{i=1}^{m} |V_i|=  |V|+ \sum_{i=1}^m |V_{i+1}| \le |V|+ 4\sum_{i=1}^m |R_i| \le |V| + 4 |V| = 5 |V|. \; \; \square
\]

Though it may seem that we are done, we have not yet specified the time needed to
perform a reduction.  To do this we must choose a specific
data structure to represent an annular strand diagram, and this choice is
fraught with difficulty.  We have worked out the details, and it suffices to keep track of either the
dual graph (i.e. the cell structure) or of the sequence of edges crossed by
some cutting path.  In neither case can reductions actually be performed in
constant time, but one can show that the amount of time required for linearly
many reductions is indeed linear.

Unfortunately, the algorithm may not be as fast for the groups $T$ and
$V$.  Checking whether two closed strand diagrams are the same involves a
comparison of the cutting cohomology classes.  This requires a Gaussian
elimination, for it must be determined whether the difference of the two
classes lies in the subspace spanned by the coboundaries of the vertices.
Gaussian elimination has cubic running time, thus:

\theorem{Let $X$ be Thompson's group $T$ or $V$ described through their
standard generating sets.
There exists a cubic time algorithm to determine whether
two elements of $X$ are conjugate.}

\chapter{Dynamics in Thompson's Group $F$}
\label{chapter3}

In Chapter \ref{chapter2}, we used \emph{strand diagrams} to give a unified
solution to the conjugacy problems in Thompson's groups $F$, $T$, and $V$.
In the present Chapter, we derive an explicit correspondence between strand diagrams for $F$
\footnote{It is expected that many of the results of this Chapter can also be extended
to Thompson's groups $T$ and $V$ and to \emph{Generalized Thompson's groups} (see Example
\ref{thm:example-rationals}).}
and piecewise-linear functions and
we obtain a complete understanding of the dynamics of elements.
In particular we are able to give simple proofs of several previously known results.
In addition, we describe a completely dynamical solution
to the conjugacy problem for one-bump
functions in $F$, similar to the dynamical criterion for conjugacy in $\PL_+(I)$
derived by Brin and Squier \cite{brin2}.
The material of this Chapter represents
joint work with James Belk. It can also be found in \cite{bema2}.

\section{Strand Diagrams}

We present here a new interpretation of strand diagrams as stack machines.
This provides a direct link between strand diagrams and piecewise-linear
functions, and paves the way for a dynamical understanding of conjugacy.
This description was inspired by a similar description of $F$ in \cite{GNS} as an ``asynchronous automata group''.
We have already
introduced this point of view in the proof of Theorem \ref{thm:correspondence-F-strand-diagrams}.

\subsection{Representation of Elements \label{sec:representation-elements}}

\begin{figure}[0.5\textwidth]
\centering
\includegraphics{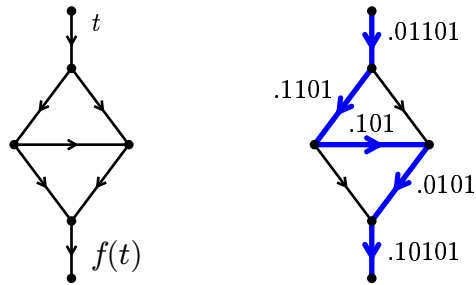}
\caption{A strand diagram as a circuit}
\label{fig:strand-diagram-circuit}
\end{figure}

Each strand diagram represents a certain piecewise-linear homeomorphism
$f\colon I\rightarrow I$. The strand diagram is like a computer circuit:
whenever a binary number $t\in[0, 1]$ is entered into the top, the
signal winds its way through the circuit and emerges from the bottom as
$f\left(t\right)$ (see figure \ref{fig:strand-diagram-circuit}).

\begin{figure}[0.5\textwidth]
\centering
\includegraphics{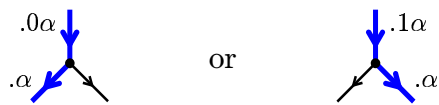}
\caption{Split rule}
\label{fig:split-rule}
\end{figure}

\begin{figure}[0.5\textwidth]
\centering
\includegraphics{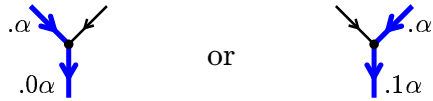}
\caption{Merge rule}
\label{fig:merge-rule}
\end{figure}

\noindent During the computation, the binary number changes each time that the signal
passes through a vertex. For a split, the signal travels either left or
right based on the first digit of the number (figure \ref{fig:split-rule}).
The first digit is lost after the signal passes through the split. For a
merge, the number gains an initial $0$ or a $1$, depending on whether it
enters from the left or from the right (figure \ref{fig:merge-rule}).
This describes the action of a strand diagram on the unit interval. We will
show in the next section that every strand diagram acts as an element of $F$.

\example{The following figure shows the three different paths that numbers might take through a certain strand diagram:
\begin{center}
\begin{minipage}{0.3\columnwidth}
\centering
\includegraphics{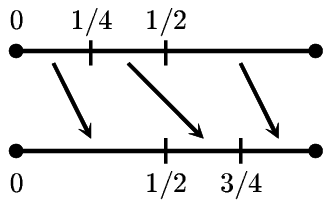}\\[6bp]
$\begin{aligned}
.00\alpha\,&\mapsto\,.0\alpha\\[0bp]
.01\alpha&\mapsto.10\alpha\\[0bp]
.1\alpha&\mapsto.11\alpha
\end{aligned}$
\end{minipage}%
\begin{minipage}{0.6\columnwidth}
\centering
\includegraphics{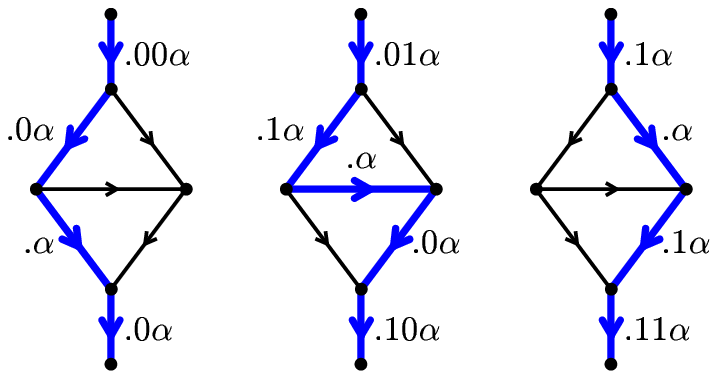}
\end{minipage}
\bigskip
\\ Figure: Three paths through a strand diagram
\end{center}
As you can see, this strand diagram acts as the element of $F$ shown on the left.}

\note{The scheme above is really the description of a \emph{stack
machine} represented by a strand diagram. A stack machine is similar to a
finite-state automaton, except that the input and output are replaced by one
or more stacks of symbols. Each state of a stack machine is either a \emph{read
state}, \emph{write state}, or a \emph{halt state}.  A read state pops a symbol from a stack, and then
moves to another state determined by which symbol was read. A write state pushes
a symbol onto a stack and then moves to a specified other state.
The process ends when the machine moves to a halt state.
A strand diagram can be interpreted as a stack machine with one stack.
Each edge represents a state of the stack machine. Edges that end
with a split are read states, edges that end with a merge are write
states, and the edge that ends with the sink is a halt state.}

\subsection{Reductions}

Recall that a \emph{reduction} of a strand diagram is either of the following moves:
\begin{center}
\begin{minipage}{2in}
\centering
\includegraphics{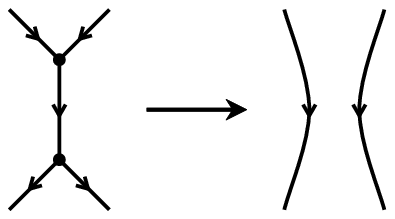}\\[0in]
Type I
\end{minipage}%
\hspace{0.5in}
\begin{minipage}{2in}
\centering
\includegraphics{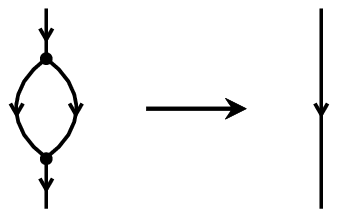}\\[0in]
Type II
\end{minipage}
\bigskip
\\ Figure: Reductions for strand diagrams
\end{center}

\begin{figure}[0.5\textwidth]
\centering
\includegraphics{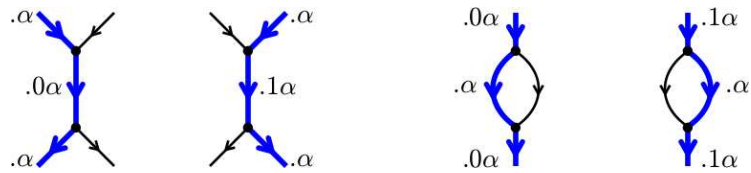}
\caption{Reductions do not change the underlying map}
\label{fig:reductions-equivalent}
\end{figure}

Neither of these simplifications changes the action of the strand diagram on
binary sequences (see figure \ref{fig:reductions-equivalent}).

\subsection{$(m,n)$-Strand Diagrams \label{sec:mn-strand-diagrams-dynamical}}

\begin{figure}[0.5\textwidth]
\centering
\includegraphics{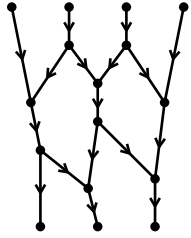}
\caption{An element of Thompson's groupoid}
\label{fig:mn-groupoid}
\end{figure}

We look at the groupoid of $(m,n)$-Strand Diagrams from the dynamical point of view (see figure \ref{fig:mn-groupoid}).
Recall that a strand diagram with $m$ sources and $n$ sinks is called an \emph{$(m, n)$-strand diagram}.
Such a strand diagram can receive input along any of its sources;  the signal then travels through
the diagram according to the rules in section \ref{sec:representation-elements}, eventually emerging from one of the sinks.

We can interpret an $(m,n)$-strand diagram as a piecewise-linear homeomorphism $[0, m]\rightarrow[0, n]$.
Specifically, a number of the form $k + 0.\alpha$ corresponds to an input of .$\alpha$ entered into the $k$th source,
or an output of .$\alpha$ emerging from the $k$th sink. The set of piecewise-linear functions determined in this way is
precisely the set of dyadic rearrangements from $[0, m]$ to $[0, n]$, i.e. the orientation-preserving homeomorphisms
$[0, m]\rightarrow[0, n]$ whose slopes are powers of two, and whose breakpoints have dyadic rational coordinates.

The set of homeomorphisms described above is closed under compositions and inverses,
and therefore forms a \emph{groupoid} with objects $\{ [0, 1], [0, 2], [0, 3], \ldots \}$.
Indeed two homeomorphisms $f:[0,m] \to [0,m]$ and $g:[0,n] \to [0,n]$ from Thompson's groupoid are conjugate
if and only if they have the same reduced annular strand diagram, by the results of Chapter \ref{chapter2}.

It is immediate to generalize this description to the case of cylindrical $(m,n)$-Strand Diagrams
and hence to give a proof to Proposition \ref{thm:ghys-sergiescu}.

\propositionname{Ghys-Sergiescu, \cite{ghys2}}{Every element of $T$ has a periodic point. \label{thm:ghys-sergiescu-2}}

\noindent \emph{Proof.} Let $f$ be a $(1,1)$-cylindrical strand diagram. Up to conjugacy,
we can select a reduced $(k,k)$-cylindrical strand diagram representing an element $g$ conjugate to $f$,
with $k$ as above.
We consider the toral strand diagram associated to $g$.
Every toral strand diagram has a merge cycle $\lambda$, which thus intersects the cutting line $\sigma$ in at least
one point $p$. This point must thus correspond to one of the $k$-sources of $g$, say the $i$-th source.

Since the point $p$ lies on a merge cycle, there must be a power $g^r$ such that the strand leaving the $i$-th source of $g^r$
ends into the $i$-th sink. By using
the right vine, we know that the $i$-th source corresponds to some dyadic subinterval $[a,b]$ of $[0,1]$.
Thus $g^r$ is a continuous map that sends $[a,b]$ into itself, thus $g^r$ must have a fixed point
and so $g$ has a periodic point. $\square$

\section{Dynamics of Annular Strand Diagrams}

\subsection{Fixed points and ``chaos'' }
\begin{figure}[0.5\textwidth]
\centering
\includegraphics{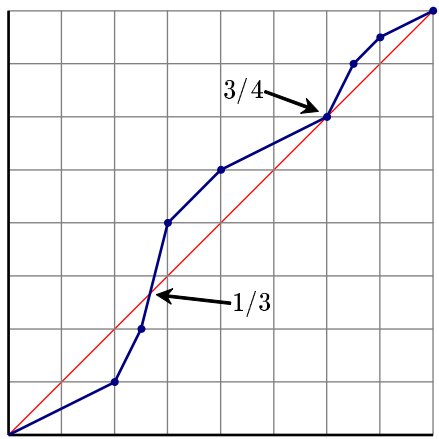}
\caption{An example of an element of $F$}
\label{fig:example-graph-F}
\end{figure}

In this section we survey some known results on dynamics in $F$.
Figure \ref{fig:example-graph-F} is the graph for an element of $F$
The main dynamical features of this element are the four fixed points
at $0$, $\frac{1}{3}$, $\frac{3}{4}$, and $1$.  Every element of $F$ fixes $0$ and $1$, but
not every element has \emph{interior fixed points} like $\frac{1}{3}$ and $\frac{3}{4}$.
We are going to observe the properties of the fixed points of this element by studying the
\emph{local replacement rule}: we look at a one-sided neighborhood $U_p$ of a fixed point $p$
that is small enough so that the map $x \to f(x)$ is linear for any $x \in U_p$ and hence,
if $x$ is written in binary expansion, then $f(x)$ is obtained by adding some digits in front of $x$
or subtracting some of the first digits of $x$, with the tail of the binary expansion of $x$ and $f(x)$ remaining the same.
\begin{enumerate}
\item The fixed point at $0$ is attracting, since the slope is $\frac{1}{2}$.
The local replacement rule is $.\alpha\mapsto.0\alpha$, which causes points near zero to converge to zero:
\begin{equation*}
.\alpha \;\;\mapsto\;\; .0\alpha \;\;\mapsto\;\; .00\alpha \;\;\mapsto\;\; .000\alpha \;\;\mapsto\;\; \cdots
\end{equation*}
\item Fixed points do not have to be dyadic. In fact, the fixed point at $\frac{1}{3}$ is not a dyadic fraction.  In binary, the local replacement
rule is $.10\alpha \mapsto .\alpha$, with a fixed point at $.101010\ldots = \frac{1}{3}$.
The slope here is $4$, so the fixed point is repelling:
\begin{equation*}
.101010\alpha \;\;\mapsto\;\; .1010\alpha \;\;\mapsto\;\; .10\alpha \;\;\mapsto\;\; .\alpha \;\;\mapsto\;\; \cdots
\end{equation*}
\item The fixed point at $\frac{3}{4}$ is dyadic, and has two local replacement
rules: $.10\alpha\mapsto.101\alpha$ on the left, and $.1100\alpha\mapsto.110\alpha$
on the right.  This makes $\frac{3}{4} = .101111\ldots = .110000\ldots$ attracting from the left:
\begin{equation*}
.10\alpha \;\;\mapsto\;\; .101\alpha \;\;\mapsto\;\; .1011\alpha \;\;\mapsto\;\; .10111\alpha \;\;\mapsto\;\; \cdots
\end{equation*}
and repelling from the right:
\begin{equation*}
.110000\alpha \;\;\mapsto\;\; .11000\alpha \;\;\mapsto\;\; .1100\alpha \;\;\mapsto\;\; .110\alpha \;\;\mapsto\;\; \cdots\text{.}
\end{equation*}
Only an interior dyadic fixed point can have different behavior from the
left and from the right, because only a dyadic rational can be a breakpoint for an element of $F$.

If we think of $F$ as acting on the Cantor set, then $\frac{3}{4}$ corresponds to \emph{two} fixed points of $f$:
one at $.101111$ and the other at $.110000$.  Each of these fixed points has a well-defined slope.
\item The fixed point at $1$ is attracting, with local replacement rule $.\alpha\mapsto.1\alpha$.
\end{enumerate}

If we think of $F$ as acting on the Cantor set, then each fixed point of an
element of $F$ has a well-defined slope, because dyadic rational fixed points are counted twice
(as they can have different slopes on the right and on the left).
The possible values of this slope depend on the tail of the fixed point:

\proposition{Suppose that $f\in F$ has a fixed point at $t$,
and let $n$ be the eventual period of the binary expansion for $t$.
Then the slope of $f$ on each side of $t$ is an integer power of $2^n$.
If $t$ is non-dyadic, the slopes at the two sides must be equal.}
\noindent \emph{Proof.} By hypothesis, $t = .\mu \overline{\rho}$, where $\rho$ is a binary sequence of length $n$.
If $\mu$ is as short as possible, then any element of $F$ with a fixed point at $t$ must have the local replacement rule
\begin{equation*}
.\mu\rho^k\alpha \;\longmapsto\; .\mu\alpha \;\;\;\;\text{or}\;\;\;\; .\mu\alpha \;\longmapsto\; .\mu\rho^k\alpha
\end{equation*}
near $t$, for some $k \geq 0$.  The first case gives a slope of $(2^n)^k$, and the second a slope of $(2^n)^{-k}$. $\square$

For example, any element of $F$ that fixes $1/3$ must have slope $4^n$ at the fixed point.
Because a dyadic rational has eventual period $1$, the left and right slopes at a dyadic fixed point can be any powers of $2$.

Most of the properties of the fixed points are preserved under conjugation:

\proposition{Let $f,g\in F$, and suppose that $f$ has fixed points at
$$ 0 = t_0  < t_1 < \cdots < t_n = 1. $$
Then $gfg^{-1}$ has fixed points at
$$0 = g(t_0) < g(t_1) < \cdots < g(t_n) = 1\text{.} $$
Moreover, the slopes of $gfg^{-1}$ on the left and on the right of $g(t_i)$ are the same as the slopes of $f$
on the left and on the right of $t_i$.}
\noindent \emph{Proof.} This is very elementary.  The statement about slopes follows from the chain rule. $\square$

\begin{figure}[0.5\textwidth]
\centering
\includegraphics{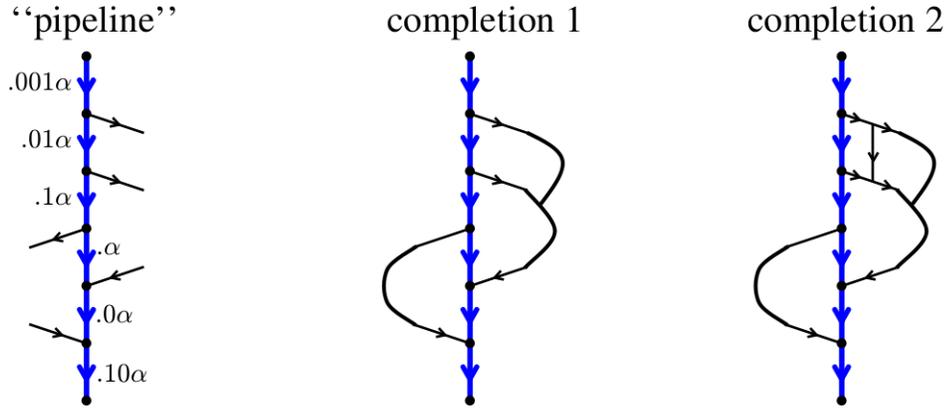}
\caption{Completing an element out of a ``pipeline''}
\label{fig:pipeline}
\end{figure}

Thus it makes sense to talk about the ``number of fixed points''
for a conjugacy class of $F$, as well as the ``slope at the 5th fixed point''.
The following proposition lets us talk about the ``tail of a fixed point'':

\proposition{\label{F-orbits}Let $t,u\in(0, 1)$. Then $t$ and $u$ are in the same orbit of $F$
if and only if $t$ and $u$ have binary expansions with the same tail---that is, if and only if
\begin{equation*}
t=.\mu\omega\hspace{0.5in}\text{and}\hspace{0.5in}u=.\nu\omega
\end{equation*}
for some finite binary sequences $\mu,\nu$ and some infinite binary sequence $\omega$.
\footnote{This result cannot be extended to generalized Thompson's groups (see Example \ref{thm:example-rationals}).
In fact, while Thompson's group $F$ is transitive on all dyadic rational points, this is not
true anymore for generalized Thompson's groups and $n$-adic rational points:
we will see in Chapter \ref{chapter4}, Remark \ref{thm:multiple-orbits-of-points}.}
\label{thm:pipeline}}
\noindent \emph{Proof.}
For the forward direction, observe
that any replacement rule preserves the tail of a binary sequence.  For the backwards direction,
it is easy to draw a ``pipeline'' that implements the rule $.\mu\alpha\mapsto .\nu\alpha$
(see figure \ref{fig:pipeline}).

Up to taking the common tail to start from a further digit, we can assume
$\mu$ and $\nu$ each have both $0$'s and $1$'s (i.e.~both left
and right connections), otherwise
$.\mu \alpha$ or $.\nu \alpha$ is $0$ or $1$.
This drawing can easily be extended to a complete
strand diagram by adding strands
on the left and on the right so that all the outgoing
strands can be suitably arranged to get into the ingoing ones.
Figure \ref{fig:pipeline} shows two possible ways to complete the pipeline, leading
to two distinct elements of $F$. $\square$

For example, the image of $\frac{3}{4}$ under an element $g\in F$ can be any dyadic fraction,
and the image of $\frac{1}{3}$ can be any rational number whose binary
expansion ends in $010101\ldots$ (i.e. any number whose difference from $\frac{1}{3}$ is dyadic).
The previous result can be obtained using the language of piecewise-linear homeomorphisms and we will
do so in Chapter \ref{chapter4} to get similar results (see Lemma \ref{thm:extension-partial-maps} and
Proposition \ref{rationals}).


The following proposition shows that there are no further constraints on the positions of the fixed points within a conjugacy class:

\proposition{\label{F-multiorbits} Let $0 = t_0 < \cdots < t_n = 1$ and $0 = u_0 < \cdots < u_n = 1$,
and suppose that each $t_i$ is in the same $F$-orbit as the corresponding $u_i$. Then
there exists an element of $F$ that maps $(t_0, \ldots, t_n)$ to $(u_0, \ldots, u_n)$.}
\noindent \emph{Sketch of the Proof.} A strand diagram for the required element can be constructed using a method similar
to the proof of the previous proposition.  See Corollary \ref{thm:overlap-intersection-diagonal}
for another proof using different techniques. $\square$


\subsection{Cut Paths and Thompson's Groupoid}
\begin{figure}[0.5\textwidth]
\centering
\includegraphics{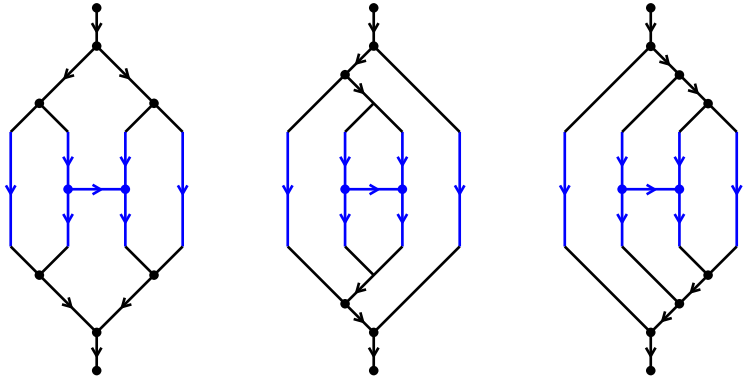}
\caption{Three conjugate elements}
\label{fig:three-conjugate}
\end{figure}

\begin{figure}[0.5\textwidth]
\centering
\includegraphics{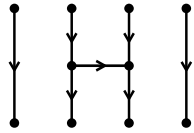}
\caption{A minimal representative}
\label{fig:conjugate-groupoid-element:}
\end{figure}

\begin{figure}[0.5\textwidth]
\centering
\includegraphics{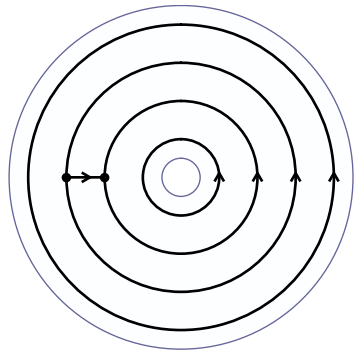}
\caption{The corresponding reduced annular strand diagram}
\label{fig:groupod-example-annulus}
\end{figure}

Thompson's groupoid is fundamental to the study of conjugacy in $F$.  For example,
figure \ref{fig:three-conjugate} shows three strand diagrams that represent conjugate elements of $F$.
Each of these elements begins by partitioning $[0,1]$ into four subintervals, and ends by recombining
these four subintervals into $[0,1]$.  They differ only in the choice of the partition.
These elements are all conjugate to the element of Thompson's groupoid shown in figure \ref{fig:conjugate-groupoid-element:}.
As you can see, this homeomorphism $[0,4]\rightarrow [0,4]$ is simpler than any of
the elements of $F$ above.  Indeed, this element is a \emph{minimal} representative for its
conjugacy class, in the sense that it is reduced (it has the fewest possible splits and merges).
The reason is that any element of this conjugacy class must have at least as many splits and
merges as the reduced annular strand diagram of figure \ref{fig:groupod-example-annulus}.

By Propositions \ref{thm:unique-reduced-diagram} and \ref{thm:cutting-path-same-annular-diagram}
we know that any two cut paths of an annular strand diagrams yield conjugate elements. Hence,
minimal representatives of a conjugacy class are precisely those obtained by cutting the reduced annular strand diagram
along some cut path.

\subsection{Directed Loops and Fixed Points}

It is possible for an element of $F$ to have infinitely many fixed points.
For example, the identity element fixes the entire interval $[0,1]$, and any element of $F$ can have a
linear segment that coincides with the identity on some interval $[d,e]$ ($d$ and $e$ dyadic).
If $f \in F$, a \emph{fixed interval} of $f$ is either
\begin{enumerate}
\item An isolated fixed point $\{ t \}$ of $f$, or
\item A maximal open interval of fixed points,
\item An endpoint of a maximal open interval of fixed points.
\end{enumerate}
\convention{Each isolated interior dyadic fixed point of $f$ corresponds to two fixed intervals.}

\theorem{Let $f\in F$, and let $S$ be the reduced annular strand diagram for $f$.  Then the directed loops
$L_0, \ldots, L_n$ of $S$ (ordered from outside to inside) are in one-to-one correspondence
with the fixed intervals $I_0 < \cdots < I_n$ of $f$.  This correspondence has the following properties:
\begin{enumerate}
\item Every free loop corresponds to a maximal interval of fixed points.
\item Every split loop corresponds to an isolated repeller.  In particular, a split
loop with $n$ splits corresponds to a fixed point with slope $2^n$.
\item Every merge loop corresponds to an isolated attractor.  In particular, a merge loop
with $n$ merges corresponds to a fixed point with slope $2^{-n}$.
\end{enumerate}
In the latter two cases, the pattern of outward and inward connections around the loop determines
the tail of the binary expansion of the fixed point.
Specifically, each outward connection corresponds to a $1$, and each inward connection corresponds to a $0$.
\label{thm:directed-loops-fixed-point}}

\remark{The previous result induces
an order on the components and the directed cycles for annular strand diagrams, going from the inside to the outside.
Compare this result with part 5 of Proposition \ref{thm:combinatorial-annular-properties}.}

\noindent \emph{Proof.} We have already shown that all of the
information outlined in the statement of the theorem is conjugacy invariant.  Therefore, we may replace $f$ by
any element whose reduced annular strand diagram is $S$.  Specifically, we may assume that $f$ is the
dyadic rearrangement $[0,k] \rightarrow [0,k]$ obtained by cutting $S$ along a cutting path $c$.

\begin{figure}[0.5\textwidth]
\centering
\includegraphics{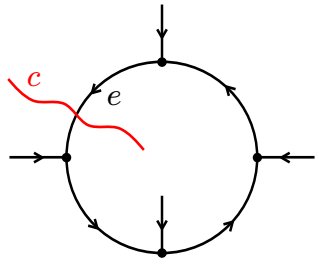}
\caption{An example of a merge loop}
\label{fig:example-merge-loop}
\end{figure}

$S$ contains a merge loop: some of the vertices on this loop are coming from the inner part of
the loop, while some are coming from the outer part of the loop.
We work out an example in detail. The general procedure follows closely from it, as it will become apparent that
the general case does not depend on the number of vertices on the loops.
Suppose that $S$ contains the merge loop in figure \ref{fig:example-merge-loop}.
The cutting path $c$ cuts through this loop exactly once, along some edge $e$. If we place a binary number .$\beta$
along $e$, the number will trace a directed path through the annular strand diagram, changing in value every time it passes
through a vertex. Assuming that $c$ crosses $i$ edges before crossing $e$, this corresponds to feeding $i+.\beta$ into
the strand diagram for $f$.

\begin{figure}[0.5\textwidth]
\centering
\includegraphics{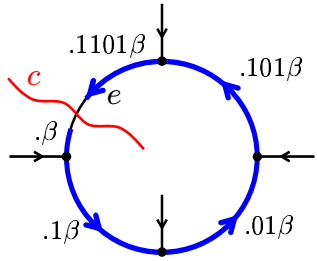}
\caption{Traveling through the merge loop}
\label{fig:example-merge-loop2}
\end{figure}

In the case we are considering, the number will simply travel around the merge loop (see figure
\ref{fig:example-merge-loop2}).
By the time it returns to $e$, its value will be the fractional part of $f\bigl(i+.\beta\bigr)$. If we continue following the
number along the merge loop, the values it has when it passes through $e$ will be the fractional parts of the iterates
$f^n (i + .\beta)$. In the case that we are considering, it follows that:
\begin{equation*}
f(i + .\beta) = i + .1101\beta  \quad  f^2(i + .\beta) = i + .1101\,1101\beta  \quad  \text{etc.}
\end{equation*}
In particular, the number $\alpha=i+.\overline{1101}$ is a fixed point of $f$.

Note that the sequence $1101$ is determined by the counterclockwise pattern of inward and outward edges,
exactly as stated in the theorem. In addition, we have shown that $f$ is linear on $[i, i+1]$, with formula:
\begin{equation*}
f(i + .\beta) \,=\, i + .1101\, \beta
\end{equation*}
This linear function has slope $2^{-4}$. This implies that $\alpha$ is an attracting fixed point---indeed, for any
$i + .\beta \in [i, i+1]$, the first $4n$ digits of $f^n(i + .\beta)$ are the same as the first $4n$ digits of $\alpha$.

\begin{figure}[0.5\textwidth]
\centering
\includegraphics{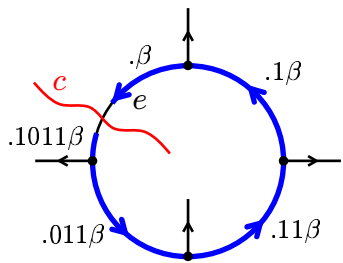}
\caption{An example of a split loop}
\label{fig:example-split-loop}
\end{figure}

A split loop works in roughly the same way, except that a split loop is repelling (see figure \ref{fig:example-split-loop}).
Note that every fixed point of $f$ arises from either a split loop or merge loop.
In particular, suppose that $i+.\beta$ is a fixed point of $f$, and let $e$ be the $(i + 1)$'st edge crossed by $c$.
If we place the binary number .$\beta$ along $e$, then the resulting path of motion must wind once around the
central hole and then return to $e$ with value .$\beta$. It follows that .$\beta$ must have traveled around a directed
loop, and $i+.\beta$ is the unique fixed point determined by the loop. $\square$

Note that the outermost loop of an annular strand diagram for $f\in F$ corresponds to the fixed point
$0=.0000\cdots$, while the innermost loop corresponds to the fixed point $1=.1111\cdots$.
Within each connected component of $S$, the outermost and innermost loops correspond to dyadic fixed points,
while the interior loops correspond to non-dyadic fixed points.

\corollary{Let $S$ be the reduced annular strand diagram for an element $f\in F$.
Then every component of $S$ corresponds to exactly one of the following:
\begin{enumerate}
\item A maximal open interval of fixed points of $f$ (for a free loop), or
\item A maximal interval with no dyadic fixed points of $f$ in its interior.
\end{enumerate}}

\begin{figure}[0.5\textwidth]
\centering
\includegraphics{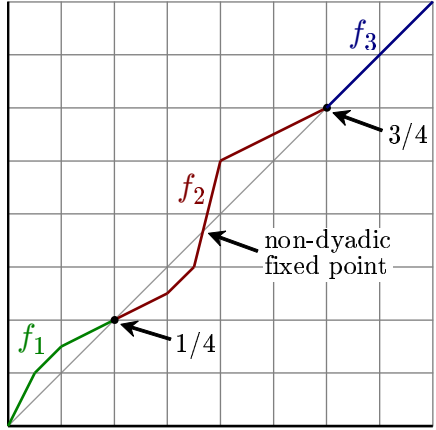}
\caption{Components of a function}
\label{fig:graph-of-components}
\end{figure}

\begin{figure}[0.5\textwidth]
\centering
\includegraphics{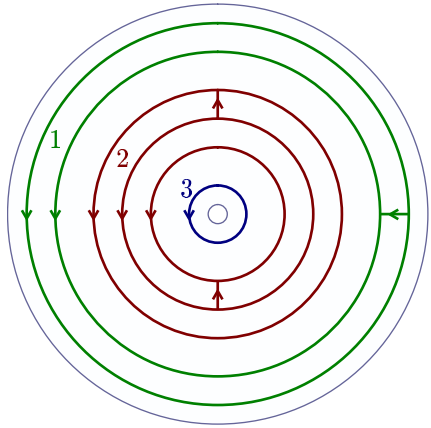}
\caption{Annular strand diagram for a component}
\label{fig:annular-strand-diagram-for-a-component}
\end{figure}

If $f\in F$, a \emph{cut point} of $f$ is either an isolated dyadic fixed point of $f$,
or an endpoint of a maximal interval of fixed points.
If $0 = \alpha_0 < \alpha_1 < \cdots < \alpha_n = 1$ are the cut points of $f$, then the restrictions
$f_i\colon [\alpha_{i-1}, \alpha_i] \rightarrow [\alpha_{i-1}, \alpha_i]$ are called the \emph{components}
of $f$ (see figure \ref{fig:graph-of-components}).
Each component of $f$ corresponds to one connected component of the reduced annular strand diagram
(figure \ref{fig:annular-strand-diagram-for-a-component}).
We recall if $\alpha < \beta$ are any dyadic rationals, there exists a Thompson-like homeomorphism
$\varphi\colon [\alpha, \beta]\rightarrow [0,1]$ such that any map in $\PL_2([\alpha, \beta])$
can be conjugated by $\varphi$ to become an element of $F$ by Corollary \ref{thm:thompson-like-isomorphic}.

\proposition{Let $f\in F$ have components $f_i\colon [\alpha_{i-1}, \alpha_i] \rightarrow [\alpha_{i-1}, \alpha_i]$,
and let $S$ be the reduced annular strand diagram for $f$.  Then for each $i$, the component of $S$ corresponding to $f_i$
is the reduced annular strand diagram for any element of $F$ conjugate to $f_i$.}

\noindent \emph{Proof.} Suppose $f$ has $n + 1$ cut points $0 = \alpha_0 < \alpha_1 < \cdots < \alpha_n = 1$.
Then we can conjugate $f$ to an
element of Thompson's groupoid whose cut points are at $0, 1, 2, \ldots, n$.  The resulting $(n, n)$-strand diagram
has $n$ connected components which, when reduced, yield the $n$ components of $S$. $\square$

\corollary{Let $f, g\in F$ have components $f_1, \ldots, f_n$ and $g_1, \ldots, g_n$.  Then $f$ is conjugate to $g$
in $F$ if and only if each $f_i$ is conjugate to $g_i$ through some Thompson-like homeomorphism.}

\section{Mather Invariants \label{sec:mather-invariant}}

Conjugacy in $F$ was first investigated by Brin and Squier \cite{brin2}, who successfully
found a criterion for conjugacy in the full group of piecewise-linear
homeomorphisms of the interval. This solution was based on some ideas of
Mather \cite{Math} for determining whether two given diffeomorphisms
of the unit interval are conjugate.

In this section we show that solution we have proved in Chapter \ref{chapter2}
can be described in a way similar to the solutions given by Mather for
Diff$_{+}\left(I\right)$ and by Brin and Squier for $\PL_+\left(I\right)$. Specifically, we define a Mather-type invariant
for elements of $F$, and show that two one-bump functions in $F$ are conjugate if and only if they have the same
Mather invariant.

A somewhat different dynamical description of conjugacy in $F$ has been obtained independently by Gill and Short \cite{ghisho1}.

\subsection{Background on Mather Invariants}

Consider the group $\Diff_+ (I)$ of all orientation-preserving diffeomorphisms of~$[0, 1]$.

\definition{A \emph{one-bump function} is an element $f \in \Homeo(I)$ such that $f(x) > x$ for all $x \in (0, 1)$.}

\begin{figure}[0.5\textwidth]
\centering
\includegraphics{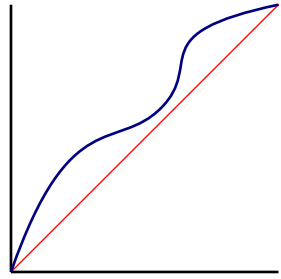}
\caption{A one-bump function}
\label{fig:one-bump-function}
\end{figure}

Figure \ref{fig:one-bump-function} shows an example of a one-bump function.
By the chain rule, two one-bump functions $f,g \in \Diff_+(I)$
can only be conjugate if $f'(0) = g'(0)$ and $f'(1) = g'(1)$, but this condition
is not sufficient.  In 1973, Mather constructed a more subtle conjugacy invariant of one-bump functions $f$
such that $f'(0)>1$ and $f'(1)<1$,
and proved that two such one-bump functions in $\Diff_+ (I)$ are conjugate if and only if they have the same
slopes at $0$ and $1$ and the same Mather invariant.  In 1995, Yocozz extended this to a complete criterion
for conjugacy in $\Diff_+ (I)$ \cite{Yoc}.  Similar invariants are used for conjugacy of diffeomorphisms
in \cite{Belt}, \cite{You}, and \cite{AfYo}, the last of which introduces the term ``Mather invariant''.

In 2001 \cite{brin2}, Brin and Squier
\footnote{Brin and Squier originally developed this theory in 1987, but it was published in 2001.}
extended Mather's analysis to the group $\PL_+(I)$ of all
orientation-preserving piecewise-linear homeomorphisms of $[0, 1$].
Specifically, they defined a Mather invariant for one-bump functions in $\PL_+(I)$, and
showed that two one-bump functions are conjugate if and only if they have the same slopes at $0$ and $1$ and
the same Mather invariant.  Using this result, they went on to describe a complete criterion for conjugacy in $\PL_+(I)$.

The Mather invariant is simpler to describe in the piecewise-linear case.
The following description is based on the geometric viewpoint introduced in \cite{You} and \cite{AfYo}, so the
language differs considerably from that used in \cite{brin2} or \cite{Math}.

\begin{figure}[0.5\textwidth]
\includegraphics{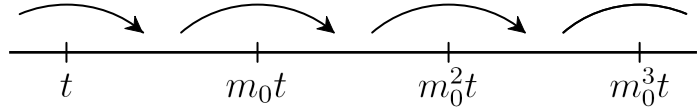}
\centering
\caption{Action of $f$ in a neighborhood of $0$}
\label{fig:mather-circle-identification}
\end{figure}

Consider a one-bump function $f\in \PL_+(I)$, with slope $m_0$ at $0$ and
slope $m_1$ at $1$.  In a neighborhood of zero, $f$ acts as multiplication by $m_0$;
in particular, for any sufficiently small $t > 0$, the interval $[t, m_0t]$ is a fundamental domain for the action of $f$
(see figure \ref{fig:mather-circle-identification}).
If we make the identification $t \sim m_0t$ in the interval $(0, \epsilon)$, we obtain a circle $C_0$, with partial covering map
$p_0 \colon (0, \epsilon)\rightarrow C_0$. Note that the restriction of $f$ is a deck transformation of this cover:
$$
\begindc{\commdiag}[12]
\obj(0,4){$(0,\epsilon)$}
\obj(6,4){$(0, \epsilon)$}
\obj(3,0){$C_0$}
\mor{$(0,\epsilon)$}{$(0, \epsilon)$}{$f$}
\mor{$(0,\epsilon)$}{$C_0$}{$p_0$}[\atright, \solidarrow]
\mor{$(0, \epsilon)$}{$C_0$}{$p_0$}[\atleft, \solidarrow]
\enddc
$$
Similarly, if we identify $(1 - t) \sim (1 - m_1 t)$ on the interval $(1 - \delta, 1)$, we obtain a circle $C_1$, with
partial covering map $p_1\colon (1 - \delta, 1)\rightarrow C_1$.

If $N$ is sufficiently large, then $f^N$ will take some lift of $C_0$ to $(0, \epsilon)$
and map it to the interval $(1 - \delta, 1)$.
This induces a map $f^\infty\colon C_0\rightarrow C_1$, making the following diagram commute:
$$
\begindc{\commdiag}[1]
\obj(-5,55){$(0, \epsilon)$}
\obj(70,55){$(1-\delta, 1)$}
\obj(0,0){$C_0$}
\obj(55,0){$C_1$}
\mor(0,55)(55,55){$f^N$}[\atleft, \solidarrow]
\mor(0,55)(0,0){$p_0$}[\atright, \solidarrow]
\mor(55,55)(55,0){$p_1$}[\atleft, \solidarrow]
\mor{$C_0$}{$C_1$}{$f^\infty$}[\atleft, \dasharrow]
\enddc
$$
\definition{The map $f^\infty$ defined above is the Mather invariant for $f$.}

We note that $f^\infty$ does not depend on the specific value of $N$ chosen. Any map $f^m$,
for $m \ge N$, induces the same map $f^\infty$. This is because $f$ acts as the identity on $C_1$ by construction
and $f^m$ can be written as $f^{m-N}(f^N(t))$, with $f^N(t) \in (1-\delta,1)$.
If $k > 0$, then the map $t \mapsto kt$ on $(0, \epsilon)$ induces a ``rotation'' $\rot_k$ of $C_0$. In particular, if we use the
coordinate $\theta = \log t$ on $C_0$, then
\begin{equation*}
\rot_k (\theta) \,=\, \theta + \log k
\end{equation*}
so rot$_k$ is an actual rotation.

\theoremname{Brin and Squier\label{MatherTheorem}}{Let $f, g \in \PL_+ (I)$ be one-bump functions with
$f'(0) = g'(0) = m_0$ and $f'(1) = g'(1) = m_1$, and let $f^\infty, g^\infty\colon C_0\rightarrow C_1$ be the corresponding Mather
invariants. Then $f$ and $g$ are conjugate if and only
if $f^\infty$ and $g^\infty$ differ by rotations of the domain and range circles:
\begin{equation*}
\begindc{\commdiag}[50]
\obj(0,1){$C_0$}
\obj(1,1){$C_1$}
\obj(0,0){$C _0$}
\obj(1,0){$C _1$}
\mor{$C_0$}{$C_1$}{$f^\infty$}[\atleft, \solidarrow]
\mor{$C _0$}{$C _1$}{$g^\infty$}[\atright, \solidarrow]
\mor{$C_0$}{$C _0$}{$\rot_k$}[\atright, \solidarrow]
\mor{$C_1$}{$C _1$}{$\rot_\ell$}[\atleft, \solidarrow]
\enddc
\end{equation*}}

\noindent \emph{Proof.} We will show here that conjugate elements have similar
Mather invariants.  See \cite{brin2} for the converse.

Suppose that $f = h^{-1} g h$ for some $h \in \PL_+ (I)$. Then the following diagram commutes,
where $k = h'(0)$ and $\ell = h'(1)$:
\begin{equation*}
\begindc{\commdiag}[20]
\obj(0,0){$C_0$}
\obj(5,0){$C_1$}
\obj(2,2){$C _0$}
\obj(7,2){$C _1$}
\obj(0,5){$(0,\epsilon)$}
\obj(5,5){$(1-\delta,1)$}
\obj(2,7){$(0, \epsilon)$}
\obj(7,7){$(1-\delta, 1)$}
\mor(0,0)(5,0){$f^\infty$}[\atright, \solidarrow]
\mor(0,0)(2,2){$\rot_k$}[\atright, \solidarrow]
\mor(5,0)(7,2){$\rot_\ell$}[\atright, \solidarrow]
\mor(2,2)(5,2){$\quad g^\infty$}[\atleft, \solidline]
\mor(5,2)(7,2){}[\atleft, \solidarrow]

\mor(0,5)(0,0){$p_0$}[\atright, \solidarrow]
\mor(2,7)(2,5){}[\atleft, \solidline]
\mor(2,5)(2,2){$p_0$}[\atleft, \solidarrow]
\mor(5,5)(5,0){$\begin{matrix} p_1 \\ \, \end{matrix}$}[\atleft, \solidarrow]
\mor(7,7)(7,2){$p_1$}[\atleft, \solidarrow]

\mor{$(0,\epsilon)$}{$(1-\delta,1)$}{$\qquad f^N$}[\atleft, \solidarrow]
\mor{$(0,\epsilon)$}{$(0, \epsilon)$}{$h$}[\atleft, \solidarrow]
\mor{$(1-\delta,1)$}{$(1-\delta, 1)$}{$h$}[\atleft, \solidarrow]
\mor{$(0, \epsilon)$}{$(1-\delta, 1)$}{$g^N$}[\atleft, \solidarrow]
\enddc
\tag*{\qedhere} \; \; \; \square\end{equation*}

For diffeomorphisms, one-bump functions are not linear in neighborhoods of $0$ and $1$,
but it is still possible to define the Mather invariant by taking a limit as $t \to 0$ and $t \to 1$.
(Essentially, a one-bump function in $\Diff_+(I)$ acts linearly on infinitesimal neighborhoods of $0$ and $1$.)
In this case, the Mather invariant is a  $C^\infty$ function $C_0\rightarrow C_1$.

\theoremname{Mather, Young}{Two one-bump functions $f, g \in \Diff_+ (It)$ with the same slopes at $0$ and $1$ are
conjugate if and only if $f^\infty$ and $g^\infty$ differ by rotations of the domain and range.}

We conclude this section by remarking that Mather invariants have been defined for diffeomorphisms
acting on higher dimensional manifolds.
In \cite{AfYo}, Afraimovich and Young extend this result to a certain class of diffeomorphisms of the sphere $S^2$.
Specifically, they consider diffeomorphisms $f$ of the sphere with two fixed points, one a hyperbolic
attractor and the other a hyperbolic repeller, with the property that all of the orbits are heteroclinic from the
repeller to the attractor. By choosing fundamental annuli for $Df$ in the tangent spaces of the two fixed points,
one can construct a Mather invariant for such diffeomorphisms which is a smooth map between two tori.

\subsection{Mather Invariants for $F$ \label{ssec:mather-invariant-F}}
In this section, we show that the reduced annular strand diagram for a one-bump function in $F$ can be interpreted
as a Mather invariant.  Therefore, two one-bump functions in $F$ are conjugate in $F$ if and only if they have the same
Mather invariant.  We also briefly describe the dynamical meaning of reduced annular strand diagrams for more
complicated elements, thereby giving a completely dynamical description for conjugacy in $F$.

\definition{The \emph{piecewise-linear logarithm} $\PLog\colon (0, \infty)\rightarrow (-\infty, \infty)$
is the piecewise-linear function that maps the interval $\left[ 2^k, 2^{k+1} \right]$ linearly onto
$[k, k+1]$ for every $k \in \Zb$ (see figure \ref{fig:PLog-map}).
\begin{figure}[0.5\textwidth]
\centering
\includegraphics{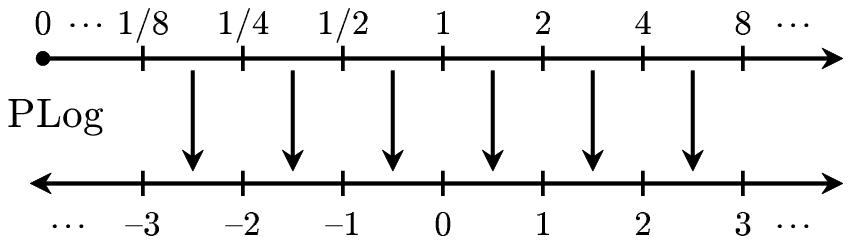}
\caption{The $\PLog$ map}
\label{fig:PLog-map}
\end{figure}}

Suppose that $f\in F$ is a one-bump function with slope $2^m$ at $0$ and slope $2^{-n}$ at $1$,
and let $f^\infty\colon C_0 \rightarrow C_1$ be the corresponding Mather invariant.
In a neighborhood of $0$, the function $f$ acts as multiplication by $2^m$.  In particular,
$\PLog f(t) = m + \PLog t$ for all $t \in (0,\epsilon)$, so we can identify $C_0$ with the
circle $\Rb / m \Zb$.  Figure \ref{fig:dyadic-circle-identification} shows the case $m = 3$:
\begin{figure}[0.5\textwidth]
\centering
\includegraphics{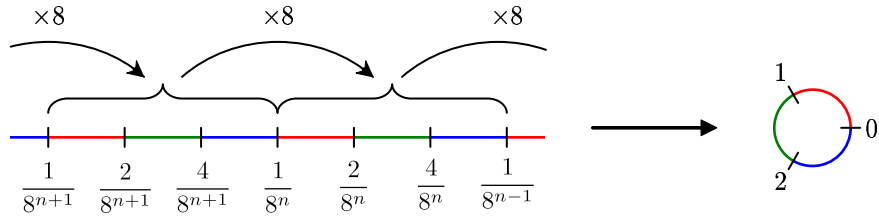}
\caption{Construction of the circle $C_0$}
\label{fig:dyadic-circle-identification}
\end{figure}
In a similar way, we can use the function $t \mapsto -\PLog(1- t)$ to identify $C_1$ with the circle $\Rb / n \Zb$.
This lets us regard the Mather invariant for $f$ as a function $f^\infty\colon \Rb /m\Zb \rightarrow \Rb /n \Zb$.
Because $f^N$ and $\PLog$ are piecewise-linear, the Mather invariant $f^\infty$ is a piecewise-linear function.
Moreover, $f^\infty$ is Thompson-like: all the slopes are powers of $2$, and the breakpoints are dyadic rational
numbers of $\Rb/m \Zb = [0,m]/\{0,m\}$.

Now, if $k \in \Zb$, then the map $t \mapsto 2^kt$ on $(0, \epsilon)$ induces a "rotation" of $C_0$.  Using
our new scheme, this is precisely an integer rotation of $\Rb / m \Zb$:
\begin{equation*}
\rot_k (\theta) = \theta + k \mod m
\end{equation*}
We are now ready to state the main theorem:

\theorem{\label{MainMatherTheorem}
Let $f, g \in F$ be one-bump functions with $f'(0) = g'(0) = 2^m$ and $f'(1) = g'(1) = 2^{-n}$,
and let $f^\infty, g^\infty\colon \Rb / m \Zb \rightarrow \Rb / n \Zb$ be the corresponding Mather invariants.
Then $f$ and $g$ are conjugate if and only if $f^\infty$ and $g^\infty$ differ by integer rotations of the domain and range circles:
\begin{equation*}
\begindc{\commdiag}[17]
\obj(0,3) {$\Rb / m\Zb$}
\obj(4,3){$\Rb / n\Zb$}
\obj(0,0){$\Rb / m   \Zb$}
\obj(4,0){$\Rb / n   \Zb$}
\mor{$\Rb / m\Zb$}{$\Rb / n\Zb$}{$f^\infty$}[\atleft, \solidarrow]
\mor{$\Rb / m   \Zb$}{$\Rb / n   \Zb$}{$g^\infty$}[\atleft, \solidarrow]
\mor{$\Rb / m\Zb$}{$\Rb / m   \Zb$}{$\rot_k$}[\atright, \solidarrow]
\mor{$\Rb / n\Zb$}{$\Rb / n   \Zb$}{$\rot_\ell$}[\atleft, \solidarrow]
\enddc
\end{equation*}}

The forward direction follows from the same argument given for proposition \ref{MatherTheorem}.
The converse is more difficult: we must show that any two one-bump functions whose Mather invariants
differ by integer rotation are conjugate in $F$.  To prove this, we describe an explicit correspondence between
Mather invariants and reduced annular strand diagrams.

\begin{figure}[0.5\textwidth]
\centering
\includegraphics{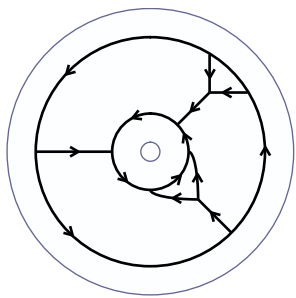}
\caption{Annular strand diagram for a one-bump function}
\label{fig:mather-annular-strand}
\end{figure}

\begin{figure}[0.5\textwidth]
\centering
\includegraphics{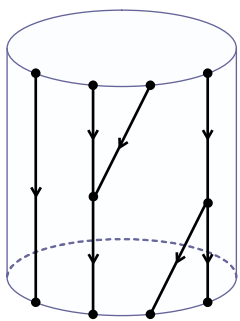}
\caption{From an annular strand diagram to a cylindrical one}
\label{fig:annulus-to-cylinder}
\end{figure}

If $f\in F$ is a one-bump function, then the only fixed points of $f$ are at $0$ and $1$.  Therefore, the reduced annular
strand diagram for $f$ has only two directed cycles (see figure \ref{fig:mather-annular-strand}).
Since $f'(0) > 1$, the outer cycle (corresponding to $0$) must be a split loop, and the inner cycle (corresponding to $1$)
must be a merge loop.  If we remove these two cycles, we get an  $(m,n)$-strand diagram drawn on a cylinder
(see figure \ref{fig:annulus-to-cylinder}).
In Chapter \ref{chapter2}, this is referred to as a \emph{cylindrical strand diagram}.
Such a diagram can be used to describe a Thompson-like map between two circles.

\proposition{There is a one-to-one correspondence between
\begin{enumerate}
\item Reduced cylindrical $(m,n)$-strand diagrams, and
\item Thompson-like functions $\Rb / m \Zb \rightarrow \Rb / n \Zb$, with two functions considered equivalent if they differ
by integer rotation of the domain and range circles.
\end{enumerate}}

\begin{figure}[0.5\textwidth]
\centering
\includegraphics{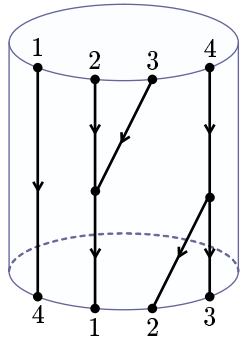}
\caption{Labeling of a cylindrical strand diagram}
\label{fig:labeled-cylinder}
\end{figure}

\noindent \emph{Proof.} A \emph{labeling} of a cylindrical $(m,n)$ strand diagram is a counterclockwise assignment of the numbers $1, 2, \ldots m$
to the sources, and a counterclockwise assignment of the numbers $1, 2, \ldots n$ to the sinks (see figure
\ref{fig:labeled-cylinder}).
Given a labeling, we can
interpret the cylindrical strand diagram as a function $\Rb / m \Zb \rightarrow \Rb / n \Zb$, with the source labeled $k$
corresponding to the interval $[k-1,k] \subset \Rb / \Zb$, and so forth.  We claim that labeled reduced cylindrical
$(m,n)$-strand diagrams are in one-to-one correspondence with Thompson-like functions $\Rb / m \Zb \rightarrow \Rb / n \Zb$.

\begin{figure}[0.5\textwidth]
\centering
\includegraphics{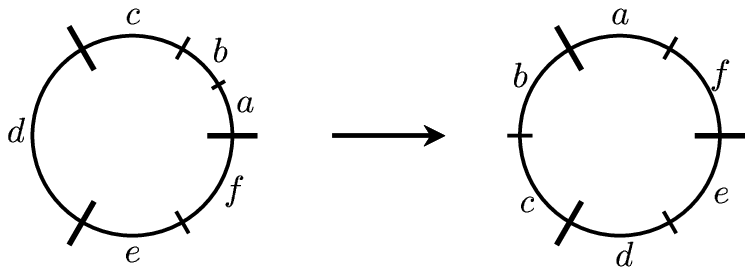}
\caption{A circle map}
\label{fig:circular-function}
\end{figure}

\begin{figure}[0.5\textwidth]
\centering
\includegraphics{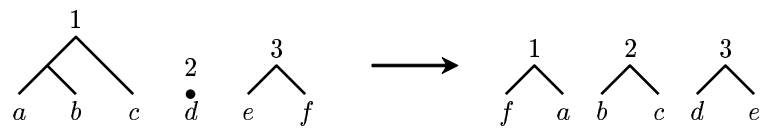}
\caption{A forest diagram for the circle map}
\label{fig:cylindrical-tree-diagram}
\end{figure}

\begin{figure}[0.5\textwidth]
\centering
\includegraphics{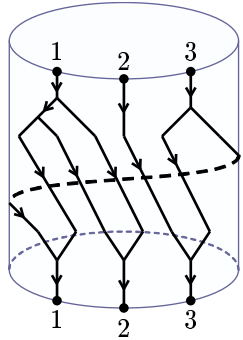}
\caption{The constructed labeled cylindrical strand diagram}
\label{fig:constructed-cylinder}
\end{figure}

The argument is similar to the proof of Theorem \ref{thm:correspondence-F-strand-diagrams}.
Suppose we are given a Thompson-like homeomorphism $f \colon \Rb / m \Zb \rightarrow \Rb / n \Zb$
(see figure \ref{fig:circular-function}).
Then we can construct a pair of binary forests representing the dyadic subdivisions of
the domain and range circles (see figure \ref{fig:cylindrical-tree-diagram}).
The forest for the domain has $m$ trees (corresponding to the subdivisions of the intervals $[0, 1], [1,2], \ldots [m-1,m]$ in
$\Rb / m \Zb$), and the forest for the range has $n$ trees.  Since the function $f$ is continuous, it must preserve the cyclic
order of the intervals.  Therefore, we can construct a strand diagram for $f$ by attaching the leaves of the top forest to the leaves
of the bottom forest via some cyclic permutation (see figure \ref{fig:constructed-cylinder}).
This gives a labeled cylindrical strand diagram for $f$.  Conversely, given any reduced labeled cylindrical $(m,n)$-strand diagram,
we can cut along every edge that goes from a split to a merge.  This decomposes the cylindrical strand diagram into two forests,
and therefore specifies a Thompson-like homeomorphism $f$.

Finally, note that changing the labeling of the sources of a cylindrical $(m,n)$-strand diagram has the effect of
performing an integer rotation on the domain of the corresponding function.  Similarly, changing the labeling of
the sinks performs an integer rotation on the range. $\square$

All that remains is the following:

\proposition{Let $\mathcal{A}$ be the reduced annular strand diagram for a one-bump function $f \in F$, and let
$\mathcal{C}$ be the cylindrical $(m,n)$-strand diagram obtained by removing the merge and split loops from $\mathcal{A}$.
Then $\mathcal{C}$ is the cylindrical strand diagram for the Mather invariant $f^\infty\colon \Rb /m\Zb \rightarrow \Rb /n\Zb$.}

\noindent \emph{Proof.} Let $f\colon [0,k]\rightarrow[0,k]$ be the one-bump function obtained by cutting a reduced annular
strand diagram $\mathcal{A}$ along a cutting path $c$. Let $e_0$ and $e_1$ be the edges on the inner and
outer loops crossed by $c$.

\begin{figure}[0.5\textwidth]
\centering
\includegraphics{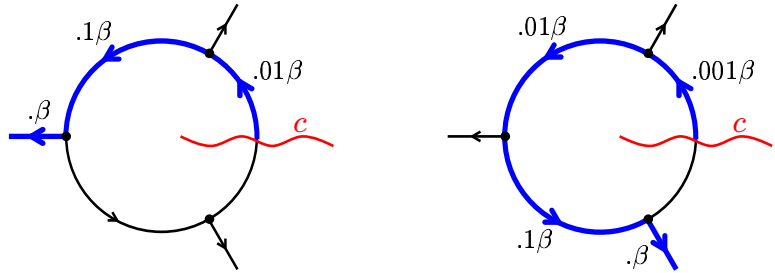}
\caption{Traveling through a split loop}
\label{fig:mather-split-loop}
\end{figure}

If we place a binary number along $e_0$, it will circle the split loop for a while, eventually exiting along some edge.
This edge depends on the length of the initial string of zeroes in the binary expansion of the number
(figure \ref{fig:mather-split-loop}).
In particular, a number leaves along the $i$th edge with value .$\beta$ if and only if the image of the number in $\Rb / m\Zb$ is
$(i - 1) + .\beta$.

After leaving the split loop, the number travels through the cylindrical strand diagram for the circle map,
eventually entering the merge loop. If we stop the number when it reaches the edge $e_1$, it will have the form
.$11\cdots 10\gamma$, where $\gamma$ is the fractional part of the image of $(i - 1) + .\beta$
under the circle map, and the length of the string of $1$'s determines the integer part. $\square$

\noindent This completes the proof of theorem \ref{MainMatherTheorem}.

 \chapter{The Simultaneous Conjugacy Problem in Groups of Piecewise Linear Functions}
 \label{chapter4}

In this chapter we look at Thompson's group $F$ as a group of piecewise-linear homeomorphisms and solve
the simultaneous conjugacy problem for $F$ and suitable $F$-like groups of piecewise linear homeomorphisms containing $F$.

For a fixed $k \in \mathbb{N}$, we say that the group $G$ has \emph{solvable $k$-simultaneous
conjugacy problem} if there is an algorithm such that, given any two $k$-tuples of elements in $G$,
$(x_1,\ldots,x_k),(y_1,\ldots,y_k)$, one can determine whether there is, or not, a $g \in G$ such that
$g^{-1} x_i g = y_i$ for all $i=1, \ldots, k$. We say that there is an
\emph{effective solution} if the algorithm produces such an element $g$, in addition to proving its existence.

This problem was studied before for various classes of groups.
The $k$-simultaneous conjugacy problem has been proved to be solvable for the matrix
groups $\mathrm{GL}_n(\mathbb{Z})$ and $\mathrm{SL}_n(\mathbb{Z})$
by Sarkisyan in 1979 in \cite{sarkisyan} and independently by
Grunewald and Segal in 1980 in \cite {grunsegal}.
In 1984 Scott constructed examples of finitely presented infinite simple groups that
have an unsolvable conjugacy problem in her paper \cite{scottelizabeth}. In their 2005 paper \cite{bridhowie}
Bridson and Howie constructed examples of finitely presented groups
where the ordinary conjugacy problem is solvable, but the $k$-simultaneous
conjugacy problem is unsolvable for every $k \ge 2$.

We will give a solution of the $k$-simultaneous conjugacy problem for
Thompson's group $F$ and then generalize it to the groups $\PL_{S,G}(J)$ (defined in Chapter \ref{chapter1}),
for an interval $J$ with endpoints in $S$. We observe that in order to make some calculations possible
inside the ring $S$ and its quotients, we need to
impose some computability requirements in $S$.
These will be clearly stated in Remark \ref{thm:ring-requirements} and will be assumed from then on.
The material of this Chapter represents
joint work with Martin Kassabov. It can also be found in \cite{kama}.

%

\section{The Ordinary Conjugacy Problem for $\PL_2(I)$}
\label{sec:ordcon}

We begin our investigation with the special case of Thompson's group $F$, seen as the group $\PL_2(I)$.
Most of the techniques that we develop for this case will extend to the general case of $\PL_{S,G}(I)$.

We prove a sequence of lemmas which will yield the solution to the ordinary conjugacy problem, that is, the
$k$-simultaneous conjugacy problem with $k=1$. To attack the ordinary conjugacy problem,
we will split the study into that of some families of functions inside $\PL_2(I)$.
The reduction to these subfamilies will come from the study of the fixed point subset of the interval $I$ for
a function $f$. For an interval $J=[\eta,\zeta] \subseteq I$, a function $f \in \PL_2(J)$ can be extended to the interval $I$ by
$f(t)=t$ for $t \in I \setminus J$, which allows us to consider $\PL_2(J)$ as a subgroup of $\PL_2(I)$.
Throughout the chapter we will assume the interval $J$ to have dyadic endpoints, so that $\PL_2(J) \cong \PL_2(I)$.
If one of the two endpoints is not dyadic, we define $\PL_2(J)$ to be the
group of restrictions of functions in $\PL_2(I)$ fixing the endpoints of $J$, that
is
\[
\PL_2(J)=\Big\{f|_J \, \, \big| \, \, f \in \PL_2(I), f(\eta)=\eta, f(\zeta)=\zeta \Big\}.
\]
We state the following interesting question:

\question{Let $J$ be an interval such
that at least one of its endpoints is non-dyadic. Is the group
$\PL_2(J)$ finitely generated?}

\noindent For a function $f \in \PL_2(J)$ we define the following closed set:
$$
D_J(f):=\{t \in J \;\big\vert\; f(t)=t \},
$$
where, to simplify the notation, we will often drop the subscript $J$.
The motivation for introducing this subset is easily explained --- If $y,z \in \PL_+(J)$ are conjugate through
$g \in \PL_+(I)$ and $s \in (\eta,\zeta)$ is such that $y(s)=s$ then
$z(g^{-1}(s))=(g^{-1}yg)(g^{-1}(s))=g^{-1}(s)$, that is,
if $y$ has a fixed point then $z$ must have a fixed point. For a subset $S \subseteq J$, we denote by
$\partial S$ the usual boundary of $S$ in $J$.

\definition{We define $\PL_+^<(J)$ (respectively, $\PL_+^>(J)$) to be the set of all functions
in $\PL_+(J)$ with graph strictly below the diagonal (respectively, above the diagonal).
Similarly, we can define $\PL_2^<(J)$ (respectively
$\PL_2^>(J)$) as the set of all functions of $\PL_2(J)$ with graph strictly below
the diagonal (respectively, above the diagonal).}

\noindent Since $x \in \PL_2(I)$ has only finitely many breakpoints, $D(x)$
consists of a disjoint union of a finite number of closed intervals and
isolated points. It is easy to see that $\partial D(x) \subseteq \mathbb{Q}$.
As mentioned before, if $g^{-1}yg=z$, then $D(y)=g(D(z))$. Thus, as a first step
we need to know if, given $y$ and $z$, there exists a $g \in \PL_2(I)$ such that $D(y)=g(D(z))$
and, in particular, $\partial D(y)=g(\partial D(z))$.

\begin{figure}[0.5\textwidth]
 \begin{center}
  \includegraphics[height=6.5cm]{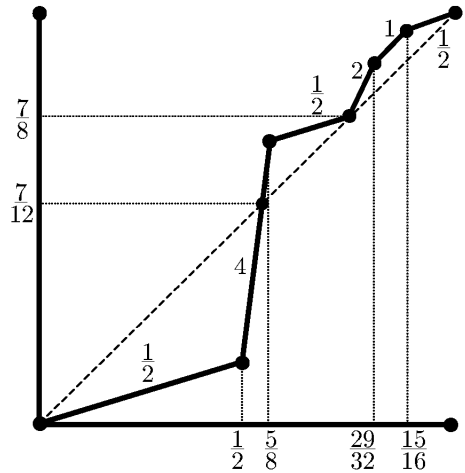}
 \end{center}
 \caption{A function with a non-dyadic fixed point.}
 \label{fig:rational-intersection}
\end{figure}

\noindent Our strategy will be the following: first we will find a way to verify if we can make $\partial D(y)$
coincide with $\partial D(z)$ through conjugation.
Then we reduce the problem to
$\partial D(y)=\partial D(z) = \{\alpha_1,\ldots,\alpha_n\}$ and so we can focus on solving the conjugacy
problem on each group $\PL_2([\alpha_i,\alpha_{i+1}])$. If $y=z=id$ on the interval $[\alpha_i,\alpha_{i+1}]$
there is nothing to prove, otherwise we can suppose that both $y,z$ are below/above the diagonal on
$[\alpha_i,\alpha_{i+1}]$. This case will be dealt with through a procedure called the
``stair algorithm'' that we provide in section \ref{ssec:stair}.
However, we observe that the $\alpha_i$'s described above need not be dyadic. The example in figure
\ref{fig:rational-intersection} shows a function with a non-dyadic rational fixed point.
In order to avoid working in intervals $J$ where the endpoints may not be dyadic,
we introduce a new definition of boundary which deals with this situation: for a subset $S$, we define

$$
\textstyle
\partial_2 S:=\partial S \cap \mathbb{Z}\left[\frac{1}{2}\right]
$$

With this definition, the set $\partial S \setminus \partial_2 S$ becomes the set of isolated non-dyadic points of $S$.

\definition{We define $\PL_2^0(J) \subseteq \PL_2(J)$ to be the set of functions $f \in \PL_2(J)$ such that
the set $D(f)$ does not contain dyadic rational points other than the endpoints of $J$, i.e., $D(f)$ is discrete and
$\partial_2 D(f)=\partial_2 J=\partial J$.}

\subsection{Making $D(y)$ and $D(z)$ coincide}
\label{ssec:diagonls}

\theorem{Given $y,z \in \PL_2(I)$, we can determine if there is (or not) a $g \in \PL_2(I)$ such
that $g(D(y)) = D(g^{-1}yg)=D(z)$. If such an element exists, it can be constructed.\label{congdiagonals}}

\noindent To start off, we need a tool to decide if this
can be proved for the boundaries of the fixed point sets. In other
words, we need to decide if it is possible to make $\partial D(y)$ coincide with $\partial D(z)$
(see figure \ref{fig:dif-fixed}).
The first step is to see how, given two rational
numbers $\alpha$ and $\beta$, we can find a $g \in \PL_2(I)$ with $g(\alpha)=\beta$.
The next two results are well known:

\begin{figure}[0.5\textwidth]
 \begin{center}
  \includegraphics{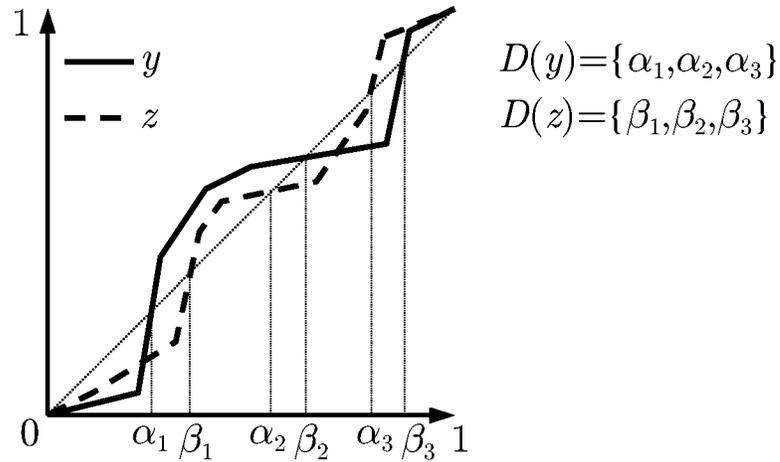}
 \end{center}
 \caption{An example with $\partial D(y) \ne \partial D(z)$.}
 \label{fig:dif-fixed}
\end{figure}

\lemmaname{Extension of Partial Maps}{Suppose $I_1,\ldots,I_k \subseteq[0,1]$ is a family of disjoint compact
intervals  $I_i=[a_i,b_i]$, with $b_i < a_{i+1}$
for all $i=1,\ldots,k$ and $a_i,b_i \in \mathbb{Z}[\frac{1}{2}]$. Let
$J_1, \ldots, J_k \subseteq [0,1]$, with $J_i=[c_i,d_i]$, be another family of intervals with the same property.
Suppose that $g_i:I_i \to J_i$ is a piecewise-linear function with a finite number
of breakpoints, occurring at dyadic rational points, and such that all slopes are integral powers of $2$.
Then there exists an $\widetilde{g} \in \PL_2(I)$ such that $\widetilde{g}|_{I_i}=g_i$.\label{thm:extension-partial-maps}}

\noindent \emph{Proof.} By our hypotheses we have that
$0<a_1<b_1<\ldots<a_k<b_k<1$ and $0<c_1<d_1<\ldots<c_k<d_k<1$
are two partitions of $[0,1]$ with the same number of points.
By Lemma \ref{thm:cfp-extension}, there exists an $h \in \PL_2(I)$ with
$h(a_i)=c_i$ and $h(b_i)=d_i$. Define
\[
\widetilde{g}(t):=
\begin{cases}
h(t) & t \not \in I_1 \cup \ldots \cup I_k \\
g_i(t) & t \in I_i
\end{cases}
\]
This function satisfies the extension condition. $\square$

\noindent We observe that this proof is constructive and produces easily an element of $F$ seen as a piecewise-linear function.
The previous result is an analogue of the proof of Proposition \ref{thm:pipeline}. In fact, another way to build an
extension of a partial map would be to write down the strand diagrams for the various given pieces and then fill
them in between by adding strands.

\proposition{Let $\alpha, \beta \in \mathbb{Q}\cap(0,1)$. Then there is a $g \in F$ such that $g(\alpha)=\beta$ iff
$$
\alpha=\frac{2^t m}{n}, \quad \beta=\frac{2^k u}{n},
$$
with $t,k \in \mathbb{Z}$, $m,n,u$ odd integers, $(m,n)=(u,n)=1$,
and the following holds
\begin{equation} \label{eq:prob}
u \equiv 2^Rm \pmod{n}
\end{equation}
for some $R \in \mathbb{Z}$.\label{rationals}}

\noindent \emph{Proof.} Suppose that there is $g \in F$ such that $g(\alpha)=\beta$.
If $\alpha$ is a dyadic rational then $\beta$ is also a dyadic rational and the conclusion of the lemma holds.
Otherwise $g(t)= 2^r t + 2^s w$ inside a small open neighborhood of $\alpha$, for some $r,s,w \in \mathbb{Z}$.
Let $\alpha=\frac{2^t m}{n}, \beta=\frac{2^k u}{v}$, for some $t,k \in \mathbb{Z}$, $(m,n)=(u,n)=1$, $m,n,u,v$ odd. Then
\[
\frac{2^k u}{v}=\beta=g(\alpha)= 2^r \frac{2^t m}{n} + 2^s w= \frac{2^{r+t}m + 2^s w n}{n}.
\]
Now the numerator of $\frac{2^{r+t}m + 2^s w}{n}$ and $n$ may not be coprime any more,
in which case we may cancel the common part
and get a new odd part $n'$ of the denominator of the right hand side.
Moreover we have $v|n$. Applying the same argument for $g^{-1}$ we have that
$n|v$, i.e., $v=n$. Thus, if there is a $g$ carrying $\alpha$ to $\beta$, then
\[
u = 2^{r+t-k}m + 2^{s-k}w n
\]
Now we can rename $R:=r+t-k$ so that the equation becomes
\[
u \equiv 2^R m \pmod{n}
\]

\begin{figure}[0.5\textwidth]
 \begin{center}
  \includegraphics[width=14cm]{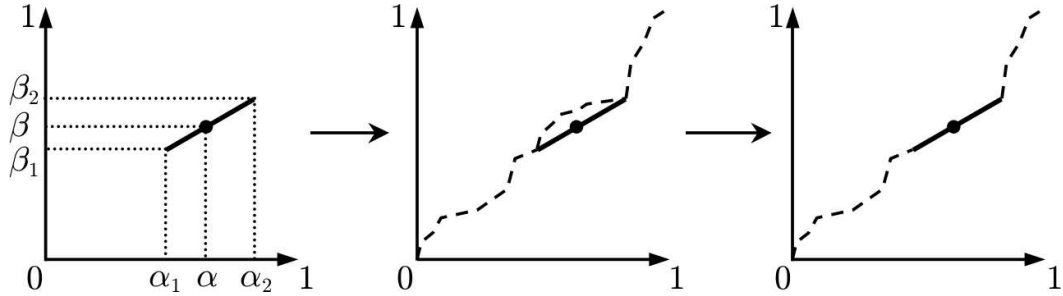}
 \end{center}
 \caption{How to build a $g \in \PL_2(I)$, with $g(\alpha)=\beta$.}
 \label{fig:rationals}
\end{figure}

\noindent Conversely, suppose $u$ satisfies (\ref{eq:prob}).
Then we can find $r,s,w$ such that, by going backwards in the "only if" argument, there is a small
open interval $(\gamma,\delta) \subset [0,1]$ containing $\alpha$ and a function $g(t)= 2^r t + 2^s w$,
with $g(\alpha)=\beta$ and we can choose $\gamma, \delta$ so that they are dyadic rationals and
$g(\gamma),g(\delta) \in I$.
Now we just apply the extension Lemma \ref{thm:extension-partial-maps}
and extend $g$ to the whole interval $[0,1]$ (see figure \ref{fig:rationals}).
$\square$

\medskip

\noindent \textbf{Example} Let $\alpha=\frac{1}{17}$, $\beta=\frac{13}{17}$ and $\gamma=\frac{3}{17}$.
It is easy to see that we can find a $g \in \PL_2(I)$ with $g(\alpha)=\beta$, but there is no
$h \in \PL_2(I)$ with $h(\alpha)=\gamma$. The same can be determined applying Proposition
\ref{thm:pipeline}, since the binary expansions of the previous numbers are
$\alpha=0.\ov{00001111}, \beta=0.11\ov{00001111}, \gamma=001011\ov{01001011}$. $\square$

\corollary{Given $\alpha,\beta \in \mathbb{Q}\cap(0,1)$ there is an algorithm to determine
whether or not there is a $g \in \PL_2(I)$ such
that $g(\alpha)=\beta$. \label{thm:how-to-identify-rationsls}}

\noindent \emph{Proof.} In order to apply the previous proposition we need
to check whether the odd parts of the denominator of $\alpha$ and $\beta$
are the same and whether they satisfy condition (\ref{eq:prob}). Equation (\ref{eq:prob}) is solvable
if and only the equation $2^X u = 2^Y m + 2^Z w n$, for some $X,Y,Z \in \mathbb{N}$, $w \in \mathbb{Z}$, is solvable.
This last equation in turn is solvable if and only if we can solve
\begin{equation} \label{eq:mod}
2^{X-Y} u \equiv m \pmod{n},
\end{equation}
because $2$ and $n$ are coprime integers, and so $2$ is invertible in $\mathbb{Z}/n\mathbb{Z}$.
If $\phi$ denotes Euler's 
function then we have that $2^{\phi(n)} \equiv 1 \pmod{n}$.
Thus, to see if (\ref{eq:mod}) is solvable, we just need to plug in all the possible
$X-Y \in \{0,1,\ldots, \phi(n)\}$. $\square$

\noindent The previous Corollary is another interpretation of Proposition \ref{thm:pipeline}:
in fact it tells us how to determine if the tail of the binary expansions of two rational numbers are the same.
We now  state the same results for a
finite number of points. Its proof uses the extension Lemma \ref{thm:extension-partial-maps}
on a number of disjoint intervals, one around each point.

\begin{cor}
\textsl{Let $0<\alpha_1 < \ldots < \alpha_r<1$ and $0<\beta_1 < \ldots < \beta_r<1$ be two
rational partitions of $[0,1]$. There exists a
$g \in \PL_2(I)$ with $g(\alpha_i)=\beta_i$ if and only if there are $g_1, \ldots, g_r \in \PL_2(I)$ such
that  $g_i(\alpha_i)=\beta_i$. Moreover, if such element $g$ exists it can be constructed.
\label{thm:overlap-intersection-diagonal}}
\end{cor}

\noindent \emph{Proof of Theorem~\ref{congdiagonals}.}
Using the previous Lemma we can determine whether or not we can make $\partial D(y)$ and
$\partial D(z)$ coincide. First we have to check if $\# \partial D(y) = \# \partial D(z)$. Then we use the previous
Corollary to find a $g \in \PL_2(I)$, with $g(\partial D(y))=\partial D(z)$, if it exists.
Let $\widehat{y}:=g^{-1}yg$. Now we just have to check if the sets where the graphs of the two functions $\widehat{y}$
and $z$ intersect the diagonal are the same. In fact, we know that the boundary points of these sets are the same, so
it is enough to check whether $D(\widehat{y})$ contains the same intervals as $D(z)$. $\square$


\subsection{The Linearity Boxes}
\label{ssec:linearitybox}

\medskip

\noindent The very first thing to check, if $y$ and $z$ are to be conjugate through a $g \in \PL_2(J)$, is
whether they can be made to coincide in neighborhoods of the endpoints of $J=[\eta,\zeta]$. This subsection
and the following one will deal with functions in $\PL_+(J)$: we will reuse them in the discussion
on $\PL_{S,G}(I)$. We start by making the following observation: the map $\PL_+(J) \to \mathbb{R}_+$
which sends a function $f$ to $f'(\eta^+)$ is a group homomorphism.

\begin{figure}[0.5\textwidth]
 \begin{center}
  \includegraphics{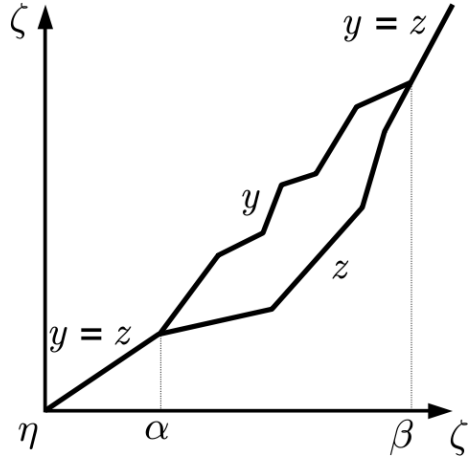}
 \end{center}
 \caption{$y$ and $z$ coincide around the endpoints.}
 \label{fig:startendequal}
\end{figure}

\lemma{Given three functions $y,z,g \in \PL_+(J)$ such that $g^{-1} y g= z$,
there exist $\alpha, \beta \in (\eta,\zeta)$ such that $y(t)=z(t)$, for all
$t \in [\eta,\alpha] \cup [\beta,\zeta]$ (refer to figure \ref{fig:startendequal}).\label{thm:startendequal}}

\noindent \emph{Proof.} We prove the Lemma for the first interval. Let $\varepsilon>0$
be a number small enough that
$$
\begin{array}{lll}
g(t) - \eta = a (t - \eta), &\,\,& \mbox{for} \,\, t \in [\eta,\eta + \varepsilon],\\
y(t) - \eta = b (t - \eta), &    & \mbox{for} \,\, t \in [\eta,g(\eta+\varepsilon)], \\
g^{-1}(t)- \eta =a^{-1} (t - \eta), & &  \mbox{for} \,\, t\in [\eta,y g(\eta+\varepsilon)].
\end{array}
$$
for some $a,b>0$. Let $\alpha=\min \{\eta+\varepsilon,g(\eta+\varepsilon),y g(\eta+\varepsilon)\}$. Then,
for $t \in [\eta,\alpha]$, we have
$$
z(t) = g^{-1} y g (t)  - \eta = a^{-1} b a (t- \eta) = b(t-\eta) =y(t).
$$
\noindent The second interval is found in the same way, after recentering the axis at the point $(\zeta,\zeta)$. $\square$

\medskip

\noindent If two functions coincide at the beginning and at the end,
then a candidate conjugator $g$ will have to be linear in certain particular ``boxes'',
which depend only on $y$ and $z$.

\begin{figure}[0.5\textwidth]
 \begin{center}
  \includegraphics{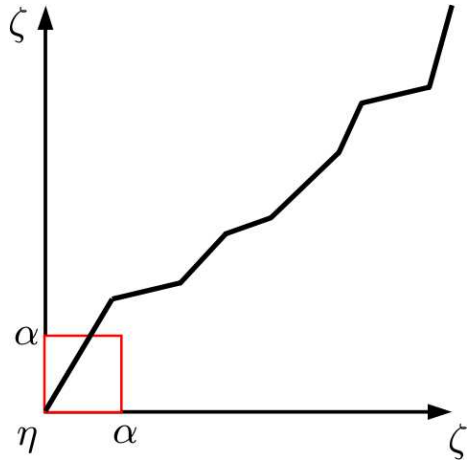}
 \end{center}
 \caption{Initial linearity box.}
 \label{fig:linearityboxes}
\end{figure}

\lemmaname{Initial Box}{Suppose $y,z,g \in \PL_+(J)$ and $g^{-1}yg=z$.
Let $\alpha>0$ and $y'(\eta^+)=z'(\eta^+)=c > 1$ satisfy
$$
y(t) - \eta = z(t)-\eta = c (t - \eta) \,\, \mbox{for } t\in[\eta,\eta+\alpha].
$$
Then the graph of $g$ is linear inside the square $[\eta,\eta+\alpha]\times[\eta,\eta+\alpha]$, i.e.,
the graph of $g$ is linear in some neighborhood of the point
$(\eta,\eta)$ in $J \times J$ depending only on $y$ and $z$ (see figure \ref{fig:linearityboxes}).
\label{thm:initial-box}}

\noindent \emph{Proof.} We can rewrite the conclusion of this lemma, by saying that, if we define
$$
\varepsilon = \sup \{r \mid g \mbox{ is linear on } [\eta,\eta+r]\},
$$
then $\eta +\varepsilon\ge \min \{g^{-1}(\eta+\alpha),\eta+\alpha\}$.
Assume the contrary, let
$\varepsilon < \alpha$ and $\eta + \varepsilon < g^{-1}(\eta + \alpha)$ and write $g(t)- \eta =\gamma (t-\eta)$
for  $t\in [\eta,\eta+\varepsilon]$, for some constant $\gamma > 0$. Let $0 \le \sigma <1$ be any number.
Since $\sigma < 1$ and $\varepsilon < \alpha$, we have $\eta+\sigma \varepsilon<\eta + \alpha$ and
so $y$ is linear around $\eta+\sigma \varepsilon$:
$$
g(y(\eta+\sigma \varepsilon))=g(\eta+c\sigma \varepsilon).
$$
On the other hand, since $\eta + \varepsilon < g^{-1}(\eta + \alpha)$,
it follows that $g(\eta+\sigma\varepsilon) < g(\eta+\varepsilon)< \eta+\alpha$ and so $z$
is linear around the point $g(\eta+ \sigma \varepsilon)=\eta +\gamma \sigma \varepsilon$:
$$
z(g(\eta+ \sigma \varepsilon))= z(\eta +\gamma \sigma \varepsilon)=\eta+ c \gamma \sigma \varepsilon.
$$
Since $gy=zg$, we can equate the previous two equations and write
$g(\eta+c \sigma \varepsilon)=\eta+\gamma c \sigma \varepsilon$,
for any number $0 \le \sigma <1$.
If we choose $1/c<\sigma<1$, we see that $g$ must be linear on the interval $[0, c \sigma \varepsilon]$,
where  $c \sigma \varepsilon> \varepsilon$. This is a contradiction to the definition of $\varepsilon$. $\square$

\noindent Notice that the square neighborhood depends only on $y$ and $z$. We observe that
the Lemma also holds when $z'(\eta^+)=y'(\eta^+)=c<1$ and the proof is given by applying
the previous proof to the homeomorphisms $y^{-1},z^{-1}$. Thus we can remove any requirement on the
initial slopes of $y$ and $z$. Note that the Initial Box Lemma has an analogue for the points close to $\zeta$:

\remarkname{Final Box}{Let $y,z,g \in \PL_+(J)$. Suppose $(g^{-1}yg)(t)=y(t)$, for all $t \in J$.
If there exist $\beta,c \in (0,1)$ such that $y(t)=z(t)=c \cdot (t-\zeta) + \zeta$
on $[\beta,\zeta]$, then the graph of $g$ is linear inside the square $[\beta,\zeta]\times[\beta,\zeta]$.
\label{thm:final-box}}


\subsection{The Stair Algorithm for $\PL_+^<(J)$}
\label{ssec:stair}

\medskip

\noindent This subsection will deal with the main construction of this Chapter. We show how, under certain hypotheses,
if there is a conjugator, then it is unique. On the other hand, we give a construction of such a conjugator, if
it exists. Given two elements $y,z$ the set of their conjugators is a coset of the centralizer of one of them,
thus it makes sense to start by deriving properties of centralizers.

\lemma{Let $z \in \PL_+(J)$. Suppose there exist $\eta \le \lambda \le \mu \le \zeta$
such that $z(t) \le \lambda$, for every $t \in [\eta,\mu]$. Suppose further that
$g \in \PL_+(I)$ is such that

\noindent (i) $g(t)=t$, for all $t \in [\eta,\lambda]$ and

\noindent (ii) $g^{-1}z g(t)=z(t)$, for all $t\in [\eta,\mu]$.

\noindent Then $g(t)=t$, for all $t \in [\eta,\mu]$.}

\noindent \emph{Proof.}
Suppose, by contradiction, that there exist points $\lambda \le \theta_1 < \theta_2 \le \mu$
such that $g(t)=t$, for all $t \in [\eta,\theta_1]$ and $g(t) \ne t$ and $g$ is linear, for
$t \in (\theta_1, \theta_2]$. Recenter the axes in the point $(\theta_1,\theta_1)$ through $T=t-\theta_1$ and $Z=z-\theta_1$.
Then $g(t)=\alpha t$,  for $t \in [0, \theta_2 - \theta_1]$, for some positive $\alpha \ne 1$ and
$z(t)=\beta t - \gamma$, for $t \in [0, \varepsilon]$, for $\beta,\gamma \in \mathbb{R},\varepsilon>0$
suitable numbers. Observe that now $-\theta_1 \le z(t) \le z(\theta_2-\theta_1)\le \lambda -\theta_1 \le 0$
and that due to the recentering $g(t)=t$ on $[-\theta_1,0]$.
For any $0<t<\min\{\theta_2-\theta_1,\varepsilon,\varepsilon/\alpha\}$ the following equalities hold:
$$
\beta t - \gamma = z(t)= gz(t)=z g(t)=z(\alpha t)=\alpha\beta t - \gamma,
$$
and so this implies $\beta t = \alpha \beta t$, hence $\alpha = 1$. Contradiction. $\square$

\corollary{Let $z \in \PL_+^<(J)$ and $g \in \PL_+(J)$ be such that

\noindent (i) $g'(\eta^+)=1$,

\noindent (ii) $g^{-1}z g(t)=z(t)$, for all $t\in J$.

\noindent Then $g(t)=t$, for all $t \in J$.\label{thm:unique-partial}}

\noindent \emph{Proof.} Since $g'(\eta^+)=1$, we have $g(t)=t$ in an open neighborhood of $\eta$.
Suppose, to set a contradiction, that $g(t_0) \ne t_0$,
for some $t_0 \in (\eta,\zeta)$. Let $\lambda$ be the first point after which $g(t) \ne t$. It is obvious
that $\eta<\lambda<\zeta$. Thus $z(\lambda) < \lambda$ and we let $\mu = z^{-1}(\lambda)>\lambda$. So we have that
$z(t) \le \lambda$ on $[0,\mu]$, $g(t)=t$ on $[\eta,\lambda]$ and $g^{-1}zg=z$ on $I$. By the previous Lemma,
$g(t)=t$ on $[\eta,\mu]$, with $\mu>\lambda$. Contradiction. $\square$

\lemma{Let $z \in \PL^<_0(J)$. Let $C_{\PL_+(J)}(z)$ be the centralizer of $z$ in $\PL_+(J)$. Define the map
$$
\begin{array}{lrcl}
\varphi_z: & C_{\PL_+(J)}(z) & \longrightarrow & \mathbb{R_+} \\
           & g              & \longmapsto     & g'(\eta^+).
\end{array}
$$
Then $\varphi_z$ is an injective group homomorphism.\label{sec:unque}}

\noindent \emph{Proof}
Let $y \in \PL^<_0(J)$ and suppose that there exists two elements $g_1,g_2 \in C_{\PL_+(J)}(y)$ such that
$\varphi_y(g_1)=\varphi_y(g_2)$, then $g^{-1}_1 g_2$ has a slope $1$ near $\eta$ and by the previous Lemma
is equal to the identity. Therefore $g_1=g_2$, which proves the injectivity. Clearly this is a group
homomorphism. $\square$

\medskip

\noindent The Lemma implies the following:

\lemma{Let $y,z \in \PL^<_0(J)$, let $C_{\PL_+(J)}(y,z)=\{g \in \PL_+(J) \; | \; y^g=z\}$ be the set
of all conjugators and let $\lambda$ be in the interior of $J$. We define the following two maps
$$
\begin{array}{lrcl}
\varphi_{y,z}: & C_{\PL_+(J)}(y,z) & \longrightarrow & \mathbb{R_+} \\
               & g                & \longmapsto     & g'(\eta^+)
\\
\psi_{y,z,\lambda}:      & C_{\PL_+(J)}(y,z)   & \longrightarrow & J \\
             & g                & \longmapsto     & g(\lambda).
\end{array}
$$
Then

\noindent (i) $\varphi_{y,z}$ is an injective map.

\noindent (ii) There is a map $\rho_\lambda:J \to \mathbb{R}_+$ such that the following diagram commutes:

\begin{displaymath}
\xymatrix{
C_{\PL_+(J)}(y,z) \ar[rr]^{\varphi_{y,z}} \ar[drr]_{\psi_{y,z,\lambda}} & & \mathbb{R}_+ \\
                                                                    & & J \ar[u]_{\rho_\lambda}
}
\end{displaymath}

\noindent (iii) $\psi_{y,z,\lambda}$ is injective. \label{cong:unique}}

\noindent \emph{Proof.} (i) is an immediate corollary of Lemma \ref{sec:unque}.
(ii) Without loss of generality we can assume that the initial slopes
of $y,z$ are the same (otherwise the set $C_{\PL_+(J)}(y,z)$ is obviously empty and any map will do).
We define the map $\rho_\lambda:J \to \mathbb{R}_+$ as
\[
\rho_\lambda(\mu)= \lim_{n \to \infty} \frac{y^n(\mu)-\eta}{z^n(\lambda)-\eta}
\]
We observe that the limit exists, i.e. the sequence stabilizes under these assumptions.

To prove that the diagram commutes we define $\mu=g(\lambda)$ and observe that
$y^{n}(\mu)\underset{n \to \infty}{\longrightarrow}\eta$ and
$z^{n}(\lambda)\underset{n \to \infty}{\longrightarrow}\eta$.
By hypothesis $y(\mu)=g(z(\lambda))$ so that $g(z^n(\lambda))=y^n(\mu)$, for every $n \in \mathbb{Z}$.
Since $g$ fixes $\eta$ we have
$$
g(t)= g'(\eta^+)(t-\eta)+\eta
\mbox{ on a small interval } [\eta,\eta + \varepsilon],
$$
where $\varepsilon$ depends on $g$. Let $N=N(g) \in \mathbb{N}$
be large enough, so that the numbers $y^N(\lambda), z^N(\lambda) \in (\eta,\eta+\varepsilon)$.
This implies that, for any $n\ge N$
$$
y^n(\mu)=g(z^n(\lambda))=g'(\eta^+)(z^n(\lambda)-\eta)+\eta
$$
and so then
\[
\varphi_{y,z}(g)=g'(\eta^+)=
\frac{y^n(\mu) -\eta}{z^n(\lambda)-\eta}=\rho_{\lambda}(\psi_{y,z,\lambda}(g)).
\]
\noindent (iii) Since $\varphi_{y,z}=\rho_{\lambda}\psi_{y,z,\lambda}$ is injective by part (i),
then $\psi_{y,z,\lambda}$ is also injective. $\square$

\medskip

\noindent Our strategy will be to construct a ``section'' of the map $\varphi_{y,z}$, if it exists.
Then as a consequence we will build a ``section'' of the map $\psi_{y,z,\lambda}$ too.
The main tool of this subsection is the {\bf Stair Algorithm}. This procedure builds a conjugator (if it exists) with
a given fixed initial slope. The idea of the algorithm is the following. In order for $y$ and
$z$ to be conjugate, they must have the same initial slope; by the initial linearity box Lemma this determines
uniquely the first piece of a possible conjugator. Then we ``walk up the first step of the stair'', with the
Identification Trick, that is basically identifying $y$ and $z$ inside a rectangle next to the linearity box,
by taking a suitable product of $y$ and $z$ as a conjugator. Then we repeat and walk up more rectangles,
until we ``reach the door'' (represented by the final linearity box) and this happens when a rectangle
that we are building crosses the final linearity box.

\begin{figure}[0.5\textwidth]
 \begin{center}
  \includegraphics[height=5.5cm]{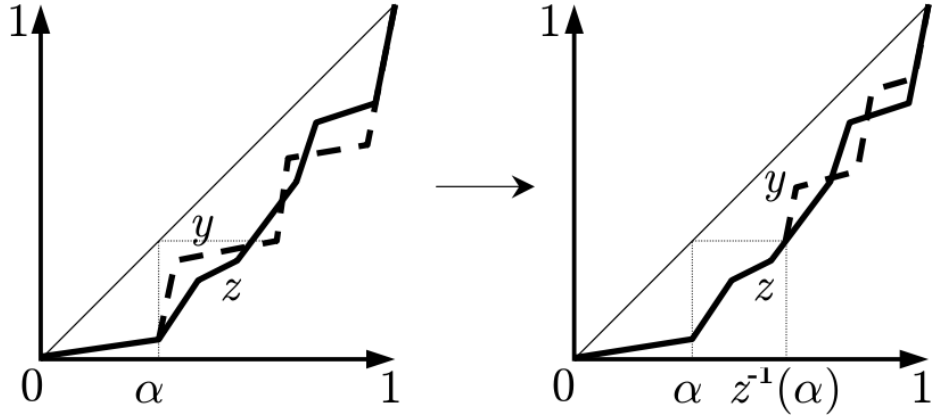}
 \end{center}
 \caption{The identification trick}
 \label{fig:identificationtrick}
\end{figure}

\lemmaname{Identification Trick}{Let $y,z \in \PL_+^<(J)$
and let $\alpha \in (\eta,\zeta)$ be such that $y(t)=z(t)$ for $t \in [\eta,\alpha]$.
Then there exists a $g \in \PL_+(I)$ such that $z(t)=y^g(t)$ for $t \in [\eta,z^{-1}(\alpha)]$
and $g(t)=t$ in $[\eta,\alpha]$. The element $g$ is uniquely defined up to the point
$z^{-1}(\alpha)$. If $y,z \in \PL_2^<(J)$ then $g$ can be chosen in $\PL_2(J)$ (see figure \ref{fig:identificationtrick}).}

\noindent \emph{Proof.} If such $g$ exists then we have that, for
$t \in [\eta,z^{-1}(\alpha)]$
$$
y(g(t))= g(z(t))=z(t)
$$
since $z(t) \le \alpha$ in $[\eta,z^{-1}(\alpha)]$. Thus, for $t \in [\alpha,z^{-1}(\alpha)]$, we have that
$g(t)=y^{-1} z(t)$. Now, that we have derived this necessary condition, we are ready to prove that
such a $g$ exists. Now define
$$
g(t):=
\begin{cases}
t & t \in [\eta,\alpha] \\
y^{-1} z(t) & t \in [\alpha,z^{-1}(\alpha)].
\end{cases}
$$
and extend it to $J$ as a line from the point $(z^{-1}(\alpha),y^{-1}(\alpha))$ to $(\zeta,\zeta)$. If
$y,z \in \PL_2(J)$ then we extend $g$ to $J$, through the extension Lemma.
A direct computation verifies that $y^g(t)=z(t)$ for $t \in z^{-1}(\alpha)$.
$\square$

\propositionname{Stair Algorithm for $\PL_+^<(J)$}{Let $y,z \in \PL_+^<(J)$
and let $q$ be a fixed positive real number.
We can decide whether or not there is a $g \in \PL_+(I)$ with initial slope $g'(\eta^+)=q$
such that $y^g=z$. If $g$ exists, it is unique and can be constructed. \label{thm:stair-algorithm}}

\noindent \emph{Proof.} Assume $y \ne z$ and, up to taking inverses, suppose
$0<g'(\eta^+)=q <1$. Let $[\eta,\alpha]^2$ the initial
linearity box and $[\beta,\zeta]^2$ the final one. Then, for $y$ and $z$ to be conjugate we must have that
$g$ has is linear in $[\eta,\alpha]^2$ and in $[\beta,\zeta]^2$.  Since $q<1$ we must have
$g$ linear on the interval $[\eta,\alpha]$ and so we define:
$$
g_0(t):= q(t-\eta)+ \eta \qquad t \in [\eta,\alpha].
$$
and extend it to the whole $J$. Now take the function
$y_1=g_0^{-1} y g_0$, which is still below the diagonal. Our goal now is to see if
$y_1$ and $z$ are conjugate. What is different now is that the new conjugator we will try to build is the identity
on $[\eta,\alpha]$, where we already know that the functions $y_1$ and $z$ coincide.
We use the Identification Trick under the diagonal and build
$$
g_1(t):=
\begin{cases}
t                & \; \; t \in [\eta,\alpha] \\
y_1^{-1} z(t)       & \; \; t \in [\alpha,z^{-1}(\alpha)]
\end{cases}
$$
then extending it to $J$. Again, we want to see we can find a conjugator
of $y_2$ and $z$ such that it is the identity on $[\eta,z^{-1}(\alpha)]$. Thus if we iterate this process and we
build a sequence $g_2, y_3, g_3 \ldots, y_r, g_r, \ldots$. By construction, we always have that
$g_r$ is the identity on $[\eta,z^{-r}(\alpha)]$ and that $y_r(t)=z(t)$  for all $t \in [\eta.z^{-r}(\alpha)]$.
We apply Lemma \ref{thm:same-lin-boxes} and choose the smallest integer $r$ so that
\[
\min\{z^{-r}(\alpha),y^{-r}(\eta + q(\alpha-\eta))\} > \beta
\]
and notice that this $r$ depends \emph{only} on $y,z$ and $q$.
Observe now that the Identification Trick tells us that, if the function $g$ of the statement exists,
it must coincide with the function $h(t):=g_0 \ldots g_r(t)$,
for $t \in [\eta,z^{-r}(\alpha)]$. If we prove that the part of the graph of $h$ on the interval
$[z^{-r}(\alpha),1]$ is inside the final box,
then we can build $g$ by extending it linearly up to the point $(\zeta,\zeta)$. Recall that,
by construction $g_{i-1}y_i^{-1}=y_{i-1}^{-1}g_{i-1}$ and $g_i(z^{-i}(\alpha))=y_i^{-1}(z^{-i+1}(\alpha))$,
for all $i=1,\ldots,r$. Then

\begin{eqnarray*}
h(z^{-r}(\alpha))=g_0 \ldots  g_{r-2} y_{r-1}^{-1} g_{r-1}(z^{-r+1}(\alpha))=\\
=g_0 \ldots g_{r-2} y_{r-1}^{-2}(z^{-r+2}(\alpha))= \ldots = \\
= y^{-r} g_0(\alpha)=y^{-r}(\eta+q(\alpha-\eta))>\beta.
\end{eqnarray*}

\noindent Since $z^{-r}(\alpha)>\beta$ by our choice of $r$ then
$(z^{-r}(\alpha),h(z^{-r}(\alpha))) \in [\beta,\zeta]^2$ and therefore we can define
$g$ by extending it linearly in the last segment, i.e. joining $(z^{-r}(\alpha),h(z^{-r}(\alpha)))$ with $(1,1)$.


\noindent If the function $h$ is not linear on
$[\beta,z^{-r}(\alpha)]$, then there is no conjugator for $y$ and $z$.
Otherwise, we have to check whether $g^{-1} y  g =z$
and we are done. To prove the uniqueness of $g$, we just apply Lemma~\ref{cong:unique}. $\square$

\lemma{Let $y,z \in \PL_+^<(J)$, $g \in \PL_+(J)$ and $n \in \mathbb{N}$. Then $g^{-1}y g=z$ if and only if
$g^{-1}y^n g=z^n$. \label{thm:unique-roots}}

\noindent \emph{Proof.} The ``only if'' part is obvious. The ``if'' part follows from the injectivity
of $\varphi_x$ of Lemma \ref{sec:unque} since $g^{-1}yg$ and $z$ both centralize the element
$g^{-1}y^n g=z^n$. $\square$

\lemma{Let $J=[\eta,\zeta]$ be a compact interval, let
$y,z \in \PL_+^<(J)$ and $g \in \PL_+(J)$ be such that $g^{-1}yg=z$. Suppose moreover that
$[\eta,\alpha] \times [\eta,\alpha]$ is the initial
linearity box and $[\beta,\zeta] \times [\beta,\zeta]$ is the final one for $y$ and $z$.
For every positive real number $q$ there
is a $k \in \mathbb{N}$ such that $y^k(\beta)<\eta + q(\alpha-\eta), z^k(\beta)<\alpha$. Moreover $y^k$ and
$z^k$ are still conjugate through $g$, so $g$ must still be linear in the same linearity boxes of $y$ and $z$.
\label{thm:same-lin-boxes}}

\noindent \emph{Proof.} Since $y(\beta)<\beta$ and $y \in \PL_+^<(J)$ then
$y^n(\beta) \underset{n \to \infty}{\longrightarrow} \eta$.
Similarly this is true for $\{z^n(\beta)\}$ and so we can pick any number $r \in \mathbb{N}$ big
enough to satisfy the statement. Moreover,
we have $g^{-1}y^kg=(g^{-1}yg)^k=z^k$. Finally we observe that the linearity box of $y^r$ and $z^r$
is smaller than that of $y$ and $z$, but that we already know that $g$ has to be linear
on $[\eta,\alpha]$ and on $[\beta,\zeta]$. $\square$

\medskip

\noindent The stair algorithm can also be proved in a slightly different way. We can apply Lemma
\ref{thm:same-lin-boxes} at the beginning and work with $y^r$ and $z^r$ instead of $y$ and $z$.
This gives a proof which concludes in two steps, although it yields the same complexity for a machine
which has to compute immediately the powers $y^r$ and $z^r$.

\noindent \emph{``Short'' Proof of Proposition \ref{thm:stair-algorithm}.}
Assume the same setting of the Proposition \ref{thm:stair-algorithm}.
We choose $r$ to be the smallest
number satisfying Lemma \ref{thm:same-lin-boxes}, so that
\[
\min\{z^{-r}(\alpha),y^{-r}(\eta + q(\alpha-\eta))\} > \beta
\]
If we call $\widehat{z}=z^r$ and $\widehat{y}=y^r$ then we have:
\[
\min\{\widehat{z}\,{}^{-1}(\alpha),\widehat{y}\,{}^{-1}(\eta + q(\alpha-\eta))\} > \beta.
\]
With this assumption, the algorithm we are going to define will need only two steps to end.
We define $g_0$ as before. Then we define
$\widehat{y}_1=g_0^{-1} \widehat{y} g_0$ and we define a map $g_1$ as in the previous proof out of $\widehat{y}_1$.
Now we observe that $g_0g_1$ is a conjugator up to the point $\widehat{z}\,{}^{-1}(\alpha)$
and that it enters the final linearity box, as in the previous proof.  Now we extend it by linearity and we check whether
this is a conjugator. If it is, it is the unique one. $\square$

\remark{By the uniqueness of the conjugator (Lemma \ref{cong:unique}) we remark that both the proofs
of the stair algorithm do not depend on the choice of $g_0$. The only requirements on $g_0$ are that it must
be linear in the initial box and $g_0'(\eta^+)=q$.}

\corollaryname{Explicit Conjugator}{Let $y,z \in \PL_+^<(J)$, let $[\eta,\alpha]$ be the initial linearity box
and let $q$ be a positive real number. There is an $r \in \mathbb{N}$ such that the unique candidate conjugator
with initial slope $q<1$ is given by
\[
g(t)=y^{-r}g_0z^r(t) \qquad \forall t \in [\eta,z^{-r}(\alpha)]
\]
and linear otherwise, where $g_0$ is any map in $\PL_+(J)$ which is linear in the initial box and
such that $g_0'(\eta^+)=q$. \label{thm:explicit-conjugator}}

\noindent \emph{Proof.} We run the short stair algorithm and let $g=g_0g_1$ be defined as above. By the short proof
of the stair algorithm and the previous Remark,
we have $g=g_0g_1=y^{-1}g_0g_1z$ on $[\eta,z^{-r}(\alpha)]$ for some $r$. Therefore
\[
g(t)=y^{-r}g_0g_1z^{-r}(t)=y^{-r}g_0 z^r(t) \qquad \forall t \in [\eta,z^{-r}(\alpha)]
\]
and it is linear on $[z^{-r}(\alpha),\zeta]$. $\square$

\corollary{Let $y,z \in \PL_+^<(J)$, and let $\lambda$ be in the interior of $J$. The map
$$
\begin{array}{lrcl}
\psi_{y,z,\lambda}:      & C_{\PL_+(J)}(y,z)   & \longrightarrow & J \\
                 & g                  & \longmapsto     & g(\lambda).
\end{array}
$$
admits a \emph{section}, i.e. if $\psi_{y,z,\lambda}(g)=\mu \in J$ , then $g$ is unique and can be constructed.
\label{thm:psi-surjective}}

\remark{Suppose $y,z \in \PL_+^<(J) \cup \PL_+^>(J)$, then in order to
be conjugate, they will have to be both in $\PL_+^<(J)$ or both in $\PL_+^>(J)$,
because by Lemma \ref{thm:startendequal} they will have to coincide
in a small interval $[\eta,\alpha]$. Moreover, $g^{-1}yg=z$ if and only
if $g^{-1}y^{-1}g=z^{-1}$, and so, up to working
with $y^{-1},z^{-1}$, we may reduce to studying the case where they are both in $\PL_+^<(J)$.}

\remarkname{Backwards Stair Algorithm}{The stair algorithm
for $\PL_+^<(J)$ can be reversed. This is to say that, given $q$
a positive real number, we can determine whether or not there is a conjugator $g$ with final
slope $g'(\zeta^-)=q$. The proof is the same: we simply start building $g$ from the final box.}

\remark{All the results of subsections \ref{ssec:linearitybox} and \ref{ssec:stair} can be
stated and proved by subsituting $\PL_2(J)$ and
$\PL_2^<(J)$ for every appearance of $\PL_+(J)$ and $\PL_+^<(J)$.
Only a few more remarks must be made in order to prove it.
In the Identification Trick we need to observe that $\alpha$ and $z^{-1}(\alpha)$ are dyadic and to take
all the extensions in $\PL_2(J)$ through the extension Lemma.}

\noindent The stair algorithm gives a practical way to find conjugators if they exist and we have chosen
a possible initial slope. By modifying the algorithm we can see that, if two elements are in $\PL_2^<(J)$
and they are conjugate through an element with initial slope a power of $2$ then the conjugator
is an element of $\PL_2(J)$.

\corollary{Let $y,z \in \PL_2^<(J)$, $g \in \PL_+(J)$ such that $y^g=z$ and $g'(\eta^+)$ is a power of $2$.
Then $g \in \PL_2(J)$.}

\subsection{The Stair Algorithm and the Mather Invariant}

\begin{figure}[0.5\textwidth]
 \begin{center}
  \includegraphics[height=4cm]{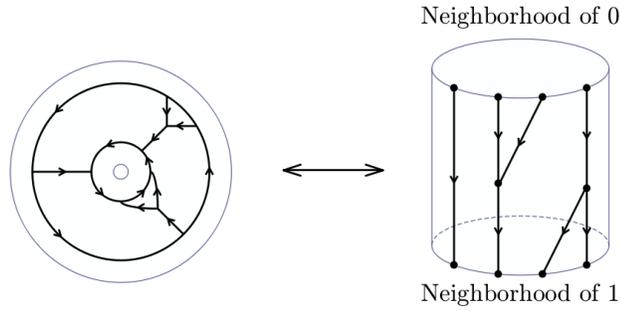}
 \end{center}
 \caption{Mather invariant as an unlabeled cylindrical strand diagram}
 \label{fig:mather-stair-algorithm-1}
\end{figure}

In Subsection \ref{ssec:mather-invariant-F} we have defined the Mather invariant for elements of $\PL_2^>(I)$.
For an element $f \in \PL_2^>(I)$, the invariant $f^\infty$ is defined by taking large powers of $f$ and then taking
a quotient so that $f^\infty$ is a map from the quotient space of a neighborhood of $0$ to the quotient space
of a neighborhood of 1 (they both become circles). The Mather invariant can be represented as an
annular strand diagram or an unlabeled cylindrical strand diagram (see figure \ref{fig:mather-stair-algorithm-1}).

\begin{figure}[0.5\textwidth]
 \begin{center}
  \includegraphics[height=6cm]{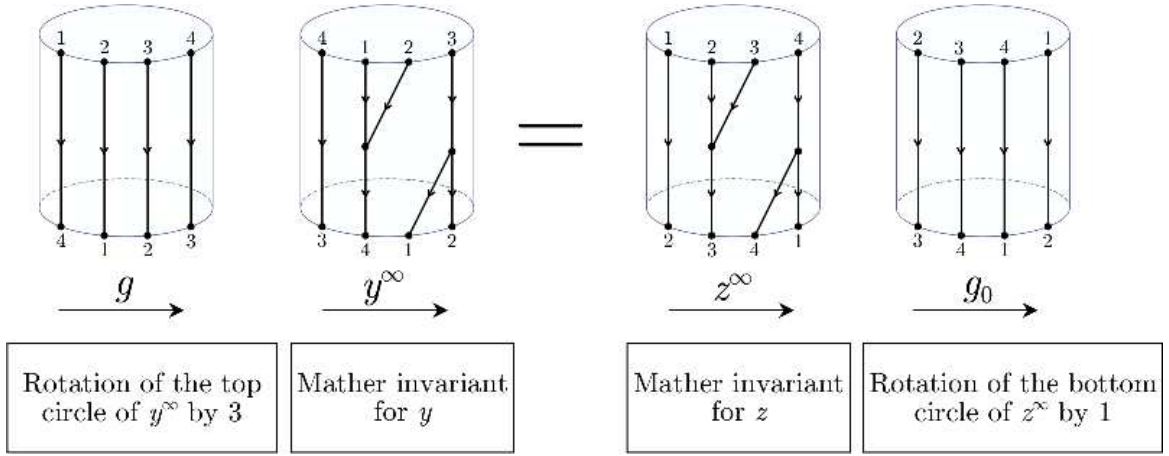}
 \end{center}
 \caption{Cylindrical strand diagrams ``differ'' by a rotation on the top or on the bottom}
 \label{fig:mather-stair-algorithm-2}
\end{figure}

In Corollary \ref{thm:explicit-conjugator} the Stair Algorithm yields that two elements $y,z \in \PL_2^>(I)$ are conjugate if
and only if the map $y^{-r}g_0z^r$ is a conjugator, for an integer $r$ large enough and for any element $g_0 \in \PL_2^>(I)$
with a given initial slope. We observe that if there is a conjugator $g$ it is thus equal to $y^{-N}g_0z^N$ for any
integer $N \ge r$.
The parallel between the two descriptions is now more apparent. If we take $N$ very large, the two maps $y^N$ and $z^N$
can be seen as the Mather invariants for $y$ and $z$. We can rewrite the equation as
\[
y^N g = g_0 z^N.
\]
If we pass to quotients, what we see on the left hand side is the composition of the Mather invariant $y^\infty$ by
the map $g$ which acts as a rotation on the domain circle of $y^\infty$. On the right hand side, we see
the composition of a rotation of the range circle of $z^\infty$ by the Mather invariant $z^\infty$.
This can also be visualized in figure \ref{fig:mather-stair-algorithm-2}.

We recall that, by Theorem \ref{MainMatherTheorem}, $y$ and $z$ are conjugate
if and only if their Mather invariants differ by a rotation in the domain and the range,
and this is precisely the same result that we obtain from Corollary \ref{thm:explicit-conjugator},
and the two points of view agree.

\remark{The previous discussion does not depend from the point of view of strand
diagrams (they only provide a different way to visualize it). The parallel between the formula for
the explicit conjugator of Corollary \ref{thm:explicit-conjugator} and the definition of Mather Invariants
for functions of $\PL_+^>(I)$ is shown by the cube diagram contained in the proof of
the Brin-Squier Theorem \ref{MatherTheorem}. Hence the parallel between the two points
of view can be generalized to $\PL_+(I)$ and later on to the groups $\PL_{S,G}(I)$
(in Section \ref{sec:stair-algorithm} we will generalize the Stair Algorithm to the groups $\PL_{S,G}(I)$).}

\subsection{The Stair Algorithm for $\PL^0_2(J)$}
\label{ssec:stairfull}

\medskip

\noindent Subsection~\ref{ssec:diagonls} proves that we can reduce our study to $y$ and $z$ such that
$\partial D(y) = \partial D(z)$. It is now important to notice that an intersection point $\alpha $
of the graph of $z$ with the diagonal
may not be a dyadic rational. If this is the case then $\alpha$ cannot be a breakpoint for $y,z,g$.
This means that, for these
$\alpha$'s, we have that $y'(\alpha)$, $z'(\alpha)$ and $g'(\alpha)$ are defined, i.e.,
the left and right derivatives coincide. Recall that a function $z$ belongs to the set $\PL_2^0(J)$ if its graph does
not have dyadic intersection points with the diagonal.

\propositionname{Stair Algorithm for $\PL^0_2(J)$}
{Let $y,z \in \PL^0_{2}(J)$ and suppose that $\partial D(y)= \partial D(z)$.
Let $q$ be a fixed power of $2$. We can decide whether or not there is a
$g \in \PL_2(J)$ with initial slope $g'(\eta^+)=q$ such that $y$ is conjugate to $z$ through $g$. If $g$ exists
it is unique.\label{thm:stairfull}}

\noindent \emph{Proof.} This proof will be essentially the same as the previous stair algorithm with a few more remarks.
We assume therefore that such a conjugator exists and build it. Let
$\partial D(y)= \partial D(z)=\{\eta=\alpha_0  < \alpha_1 < \ldots < \alpha_s < \alpha_{s+1}=\zeta\}$.
We restrict our attention to $\PL_2([\alpha_i,\alpha_{i+1}])$, for each $i=0,\ldots,s$. If $y$ and
$z$ are conjugate on $[\alpha_i,\alpha_{i+1}]$ then we can speak of linearity
boxes: let $\Gamma_i:=[\alpha_i,\gamma_i] \times [\alpha_i,\gamma_i]$
be the initial linearity box and $\Delta_i:=[\delta_i,\alpha_{i+1}] \times
[\delta_i,\alpha_{i+1}]$ the final one for $\PL_2([\alpha_i,\alpha_{i+1}])$.
Now what is left to do is to repeat the procedure of the stair algorithm for elements in $\PL_2^<(U)$,
for some interval $U$. We build a conjugator $g$ on $[\alpha_0,\alpha_1]$
by means of the stair  algorithm.
We observe that $\alpha_1$ is not a breakpoint, hence $g'(\alpha_1^+)=g'(\alpha_1^-)$. Thus we are given an initial
slope for $g$ in $[\alpha_1,\alpha_2]$, then we can repeat the same procedure and repeat the stair algorithm
on $[\alpha_1,\alpha_2]$. We keep repeating the same procedure until we
reach $\alpha_{s+1}=\zeta$. Then we check whether the $g$ we have found conjugates
$y$ to $z$. Finally, we observe that in each square $[\alpha_i,\alpha_{i+1}]\times[\alpha_i,\alpha_{i+1}]$
the determined function is unique, since we can apply Lemma \ref{cong:unique} on it. $\square$

\noindent An immediate consequence of the previous result is the following Lemma:

\lemma{Suppose $z \in \PL^0_2(J)$ and $g \in \PL_2(J)$ are such that

\noindent (i) $g'(\eta^+)=1$,

\noindent (ii) $(g^{-1}zg)(t)=z(t)$, for all $t\in J$.

\noindent Then $g(t)=t$, for all $t \in J$.}

\remarkname{Backwards and Midpoint Stair Algorithm}{It is possible to run
a backwards version of the stair algorithm also for $\PL_2^0(J)$.
Moreover, in this case it also possible to run a midpoint version of it:
if we are given a point $\lambda$ in the interior of $J$ fixed by $y$ and $z$
and $q$ a fixed power of $2$, then, by running
the stair algorithm at the left and the right of $\lambda$ we determine whether there is or not
a conjugator $g$ such that $g'(\lambda)=q$. \label{backwardstairalgorithm}}

\noindent From the previous Lemma and Remark it is immediate to derive:

\corollary{Let $y,z \in \PL^0_2(J)$ such that $D(y)=D(z)$
and let $C_{\PL_2(J)}(y,z)=\{g \in \PL_2(J) \, | \, y^g=z\}$ be the set of all conjugators. For any $\tau \in D(y)$
define the map
$$
\begin{array}{lrcl}
\varphi_{y,z,\tau}: & C_{\PL_2(I)}(y,z) & \longrightarrow & \mathbb{R_+}  \\
               &  g               & \longmapsto     & g'(\tau),
\end{array}
$$
where if $\tau$ is an endpoint of $J$ we take only a one-sided derivative. Then

\noindent (i) $\varphi_{y,z,\tau}$ is an injective map.

\noindent (ii) If $\varphi_{y,z,\tau}$ admits a section, i.e. if there is a map $\mathbb{R_+} \to C_{\PL_2(I)}(y,z)$,
$\mu \to g_\mu$ such that $\varphi_{y,z,\tau}(g_\mu)=\mu$
then $g_\mu$ is unique and can be constructed. \label{con:unique}}

\proposition{Let $y,z \in \PL_2^0(J)$ such that $D(y)=D(z)$
and let $\lambda$ be in the interior of $J$ such that
$y(\lambda) \ne \lambda$. Define
$$
\begin{array}{lrcl}
\psi_{y,z,\lambda}:      & C_{\PL_+(J)}(y,z)   & \longrightarrow & J \\
                         & g                  & \longmapsto     & g(\lambda).
\end{array}
$$
Suppose $y^{n}(\lambda)\underset{n \to \infty}{\longrightarrow} \tau$. Then

\noindent (i) There is a map $\rho_\lambda:J \to \mathbb{R}_+$ such that the following diagram commutes:

\begin{displaymath}
\xymatrix{
C_{\PL_+(J)}(y,z) \ar[rr]^{\varphi_{y,z,\tau}} \ar[drr]_{\psi_{y,z,\lambda}} & & \mathbb{R}_+ \\
                                                                    & & J \ar[u]_{\rho_\lambda}
}
\end{displaymath}
\noindent (ii) $\psi_{y,z,\lambda}$ is injective.

\noindent (iii) If $\psi_{y,z,\lambda}$ admits a section, i.e. if there is a map $J \to C_{\PL_2(I)}(y,z)$,
$\mu \to g_\mu$ such that $\psi_{y,z,\lambda}(g_\mu)=\mu$
then $g_\mu$ is unique and can be constructed. \label{thm:psi-surjective2}}

\noindent \emph{Proof.} Let $D(y)=D(z)=\{\eta=\mu_0 < \mu_1 < \ldots < \mu_k < \mu_{k+1} = \zeta\}$
and suppose $\mu_i < \lambda < \mu_{i+1}$ for some $i$. We define the partial map $\rho_\lambda:J \to \mathbb{R}_+$ as
\[
\rho_\lambda(\mu)=\begin{cases}
\lim_{n \to \infty} \frac{y^n(\mu)-\tau}{z^n(\lambda)-\tau} & \mu \in [\mu_i,\mu_{i+1}] \\
1 & \mbox{otherwise}
\end{cases}
\]
Since $D(y)=D(z)$, $z^{n}(\lambda)\underset{n \to \infty}{\longrightarrow} \tau$ and
$\tau$ is fixed by $g$. Thus if $\mu=g(\lambda)$,
then $y^n(\mu)=g(z^n(\lambda))\underset{n \to \infty}{\longrightarrow} \tau$.
With this definition, the proof follows closely that of Lemma \ref{cong:unique}(ii),
Proposition \ref{thm:psi-surjective} and by applying Corollary \ref{con:unique} and
the previous Remark. $\square$

\noindent We conclude this subsection with a technical lemma which we will need later on:

\lemma{Let $\tau,\mu \in J$, $h \in \PL_+(J)$. Then:

\noindent (i) The limit $\varphi_\pm=\underset{n \to \infty}{\lim}
h^{\pm n}(\tau)$ exists and $h(\varphi_\pm)=\varphi_\pm$,

\noindent (ii) We can determine whether there is or
not an $n \in \mathbb{Z}$, such that $h^n(\tau) = \mu$.\label{conginZ}}

\noindent \emph{Proof.} If $h(\tau)=\tau$ then it is clear. Otherwise, without loss of generality,
we can assume $h(\tau)>\tau$. The two sequences
$\{h^{\pm n}(\tau)\}_{n \in \mathbb{N}}$ are strictly monotone, and they have a limit
$\underset{n \to \infty}{\lim} h^{\pm n}(\tau)= \varphi_\pm \in [0,1]$. Thus, by continuity of $h$
\[
\varphi_\pm = \underset{n \to \infty}{\lim}h^{n+1}(\tau)=
\underset{n \to \infty}{\lim}h(h^n(\tau))= h(\varphi_\pm).
\]
Thus we have that $\{h^{n}(\tau)\}_{n \in \mathbb{Z}} \subseteq (\varphi_-,\varphi_+)$ and we have that
$\varphi_+$ is the closest intersection of $h$ with the diagonal on the right of $\tau$ (similarly
for $\varphi_-$), so we can compute $\varphi_+,\varphi_-$ directly, without using the limit.
As a first check, we must see if $\mu \in (\varphi_-,\varphi_+)$.
Then since the two sequences $\{h^{\pm n}(\tau)\}_{n \in \mathbb{N}}$
are monotone, then after a finite number of steps we find $n_1,n_2 \in \mathbb{Z}$ such that $h^{-n_1}(\tau) < \mu
< h^{n_2}(\tau)$ and so this means that either there is an integer $-n_1 \le n \le n_2$ with $h^n(\tau)=\mu$ or
not, but this is a finite check. $\square$

\subsection{The conjugacy problem for $\PL_2(I)$}
\label{ssec:con}

\medskip

\noindent We are now ready to give an alternative proof of the solvability
of the ordinary conjugacy problem (compare it with Theorem \ref{thm:characterize-conjugacy-F}).

\lemma{For any $y,z \in \PL^0_2(I)$ we can decide whether there is (or not) a $g \in \PL_2(I)$ with $y^g=z$.}

\noindent \emph{Proof.} Let $y,z \in \PL_2(I)$, $y \ne z$.
We use Theorem \ref{congdiagonals} to make $\partial D(y) = \partial D(z)$, if possible.
In order to be conjugate, we must have $y'(0^+)=z'(0^+)$ and $y'(1^-)=z'(1^-)$.
Up to taking inverses of $y$ and $z$, we can assume that $2^{u} = y'(0^+)=z'(0^+) < 1$.
Now observe that $g^{-1}yg=z$ is satisfied if and only if $(y^v g)^{-1}y (y^v g)=z$ is satisfied for every
$v \in \mathbb{Z}$. If $2^{\rho(g)}$ is the initial slope of $g$, then $2^{v u + \rho(g)}$ is the initial
slope of $y^v g$. Thus, up to taking powers of $y$, we can assume
that the initial slope of $g$ is between $2^u$ and $2^{-1}$.
Now we choose all $q \in U:=\{2^u, 2^{u+1}, \ldots, 2^{-1}\}$ as possible initial slopes for $g$, therefore
we apply the stair algorithm for $\PL^0_2(I)$ for all the elements of $U$ and check if we find a solution or not.
There is only a finite number of ``possible'' initial slopes, so the algorithm will terminate. Moreover, by Lemma 2.22
we can derive the uniqueness of each solution, for a given initial slope. $\square$

\noindent The previous Lemma provides a way to find all
possible conjugators, however it is not an efficient way to do it because
we are taking all possible slopes into consideration.

\theorem{The group $\PL_2(I)$ has solvable conjugacy problem.\label{solofcongproblem}}

\noindent \emph{Proof.} We use Theorem \ref{congdiagonals} again and suppose that
$\partial_2 D(y) = \partial_2 D(z) = \{ 0 = \alpha_0  < \alpha_1 < \ldots < \alpha_r < \alpha_{r+1} = 1\}$.
Now we restrict to an interval $[\alpha_i,\alpha_{i+1}]$ and consider $y,z \in \PL_2^0([\alpha_i,\alpha_{i+1}])$.
If $D(y)$ contains the a subinterval of $[\alpha_i,\alpha_{i+1}]$, then we must have $y=z=id$ on
the whole interval $[\alpha_i,\alpha_{i+1}]$
and so any function $g \in \PL_2([\alpha_i,\alpha_{i+1}])$ will be a conjugator. Otherwise,
$D(y)$ does not contain any subinterval of $[\alpha_i,\alpha_{i+1}]$ and so we can apply the previous Lemma on it.
If we find a solution on each such interval, then the conjugacy problem is solvable. Otherwise, it is not.
$\square$

\noindent The argument given to solve the conjugacy problem in $F$ also works,
in much the same way, to solve the power conjugacy problem. We say
that a group $G$ has \emph{solvable power conjugacy problem} if there is an algorithm
such that, given any two elements $y,z \in G$, we can determine whether there is,
or not, a $g \in G$ and two non-zero integers $m,n$ such that $g^{-1}y^m g=z^n$,
that is, there are some powers of $y$ and $z$ that are conjugate.

\theorem{The group $\PL_2(I)$ has solvable power conjugacy problem. \label{thm:power-conjugacy}}

\noindent \emph{Proof.} Again, we can use Theorem \ref{congdiagonals}, suppose that
$\partial_2 D(y) = \partial_2 D(z)$ and restrict to a smaller interval $J=[\eta,\zeta]$ with
dyadic endpoints and such that $y,z \in \PL_2^0(J)$. If $g \in \PL_2(J)$ and $m,n$ exist then we must
have that the initial slope of $y^m$ and $z^n$ must coincide. A simple argument on
the exponent of these slopes, implies that this can happen if and only if $y^m$ and $z^n$
are both powers of a common minimal power $(y^\alpha)'(\eta) = (z^\beta)'(\eta)$. Hence
the problem can be reduced to finding whether there is a $g \in \PL_2(J)$ and an integer $k$
such that $g^{-1}y^{k \alpha} g=z^{k \beta}$. By Lemma \ref{thm:unique-roots} (that
can be naturally generalized to $\PL_2^0(J)$), we have that this is equivalent
to finding a $g \in \PL_2(J)$ such that $g^{-1}y^{\alpha} g=z^{\beta}$. Hence
solving the power conjugacy problem is equivalent to solving the conjugacy problem
for $y^\alpha$ and $z^\beta$. $\square$

\section{Roots and Centralizers in $\PL_2(I)$}
\label{sec:applications}

In this section we show how the techniques developed so far allow us to obtain two previously known results.
The first of these results was first proved by Brin and Squier in \cite{brin1} in 1985 and
later proved again by Guba and Sapir
in \cite{gusa1} in 1997. Most of the results of this section are proved in \cite{gusa1} using different methods.

\propositionname{Computing Roots}{Let $id \ne x \in \PL_2(I)$, then the function $x$ has only a finite number of
roots and every root is constructible, i.e., there is an algorithm to compute it.\label{thm:compute-roots}}

\noindent \emph{Proof.} We suppose that
$\partial_2 D(z) = \{ 0 = \alpha_0  < \alpha_1 < \ldots < \alpha_r < \alpha_{r+1} = 1\}$
and we restrict again to an interval $[\alpha_i,\alpha_{i+1}]$. So we can suppose
$\partial_2 D(z)=\{0,1\}$. Let $m = x'(0^+)$ and let $n \in \mathbb{N}$ such that $\sqrt[n]{m}$
is still an integral power of $2$ (otherwise it does
not make sense to look for a $n$-th root). We want to determine whether there is, or not,
a $g \in \PL_2(I)$ such that $g^{-1}xg=x$ and such that $g'(0^+)=\sqrt[n]{m}$.
Suppose that there is such a $g$, then $g^{-k} x g^k = x$ and $(g^k)'(0^+)=m$.
Then, by the uniqueness of the solution of
the conjugacy problem with initial slope $m$ (Corollary \ref{con:unique}), we have that $g^n=x$. Conversely,
if we have $h$ such that $h^n=x$, then $h'(0^+)=\sqrt[n]{m}$. But $h^{-1}xh= h^{-1} h^n h = h^n =x$.
Thus an element $h$ is a $n$-th root of $x$ if and only if it
is the solution the ``differential type'' equation with a given initial condition
$$
\begin{cases}
h^{-1}xh=x \\
h'(0^+)=\sqrt[n]{m}.
\end{cases}
$$

\noindent So we can decide whether or not there is a $n$-th root, by solving the equivalent conjugacy problem.
Moreover, if the $n$-th root of $g$ exists,  it is computable by Theorem~\ref{solofcongproblem} and unique
by Corollary~\ref{con:unique}. $\square$

\begin{figure}[0.5\textwidth]
 \begin{center}
  \includegraphics[height=6cm]{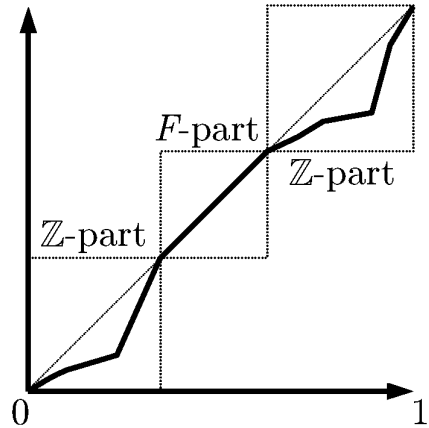}
 \end{center}
 \caption{The structure of centralizers in $F$}
 \label{fig:structure-centralizers-F}
\end{figure}

\propositionname{Centralizers}{Suppose $x \in F$, then its centralizer is
$C_F(x) \cong F^m \times \mathbb{Z}^n$, for some positive integers $m,n$ such that
$0\le m \le n+1$ (see figure \ref{fig:structure-centralizers-F}). \label{thm:centralizers-F}}

\noindent \emph{Proof.} Consider the conjugacy problem with $y=z=x$ and
let $\partial_2 D(x)=\{\eta=\alpha_0  < \alpha_1 < \ldots < \alpha_s < \alpha_{s+1}=\zeta\}$. Since
all the points of $\partial_2 D(x)$ are fixed by $x$ then $g \in C_{\PL_2(I)}(x)$
must fix the set $\partial_2 D(x)$ and thus each of the $\alpha_i$'s.
This implies that we can restrict to solve the conjugacy problem
in each of the subgroups $\PL_2([\alpha_i,\alpha_{i+1}])=\PL_2^0([\alpha_i,\alpha_{i+1}])$ and so we
can assume that $x \in \PL_2^0(I)$. If $x=id$, then it is immediate $C_{\PL_2(I)}(x)=\PL_2(I)$.
Suppose $x \ne id$ on $[0,1]$, then the map
$\varphi_{x,x}$ of Corollary \ref{con:unique} is a non-trivial injective group homomorphism. Thus
$C_{\PL_2(I)}(x) \cong \log_2(\textrm{Im} \, \varphi_{x,x}) \le \mathbb{Z}$, and
so $C_{\PL_2(I)}(x)$ is isomorphic to a subgroup
of $\mathbb{Z}$. Therefore $C_{\PL_2(I)}(x) \cong \mathbb{Z}$.
Let $[\alpha_{i_1},\alpha_{i_1+1}], \ldots,[\alpha_{i_n},\alpha_{i_n+1}]$ be
the family of intervals such that $x|_{[\alpha_{i_j},\alpha_{i_j+1}]} \ne id$, then the number of intervals
where there restriction of $x$ is trivial cannot be more than $n+1$: $x$ can be trivial only on
the intervals $[\eta,\alpha_{i_1}], [\alpha_{i_1+1},\alpha_{i_2}],
\ldots,[\alpha_{i_{n-1}+1},\alpha_{i_n}],[\alpha_{i_n+1},\zeta]$. $\square$

\corollary{Suppose $x \in \PL_2(I)$ is such that $C_{\PL_2(I)}(x) \cong \mathbb{Z}$,
then
$C_{\PL_2(I)}(x) = \langle \sqrt[k]{x} \rangle$,
for some $k \in \mathbb{Z}$.}

\noindent \emph{Proof.} Let $\varphi_{x,x}$ be as in
Corollary~\ref{con:unique}, then $\log_2(\textrm{Im} \, \varphi_{x,x}) = M\mathbb{Z}$,
for some $M \in \mathbb{Z}$. Let $2^n=\varphi_{x,x}(x)$ and let $\widehat{x}=\varphi_{x,x}^{-1}(2^M)$.
Thus there is a $k \in \mathbb{Z}$ with
$2^n=\varphi_{x,x}(x)=\varphi_{x,x}(\widehat{x}^k)=2^{kM}$. This implies
that $k=n/M$ and that $\widehat{x}=\sqrt[k]{x}$, since $\varphi_{x,x}$ is injective.
Thus $C_{\PL_2(I)}(x)=\langle \sqrt[k]{x} \rangle$. $\square$

\propositionname{Intersection of Centralizers}
{Let $x_1, \ldots, x_k \in \PL_2(I)$ and define
$C:=C_{\PL_2(I)}(x_1) \cap \ldots \cap C_{\PL_2(I)}(x_k)$.
If the interval $I$ is divided by the points in the union
$\partial_2 D(x_1) \cup \dots \cup \partial_2 D(x_k)$
into the intervals $J_i$ then
$$
C = C_{J_1} \cdot C_{J_2} \cdot \ldots \cdot C_{J_r},
$$
where $C_{J_i}:= \{f \in C \mid f(t)=t, \forall t \not \in J_i\} = C \cap \PL_2(J_i)$.
Moreover, each
$C_{J_i}$ is isomorphic to either $\mathbb{Z}$, or $\PL_2(J_i)$ or
the trivial group.\label{intersectionofcetralizers}}

\noindent \emph{Proof.}
The set $\partial_2 D(x_i)$ is fixed by all elements in $C_{\PL_2(I)}(x_i)$,
therefore all elements in $C$ fix the end points of the intervals $J_i$.
The decomposition of $C$ as $C_{J_1}\cdot \ldots \cdot C_{J_r}$ follows from the observation:

\noindent \emph{Claim:} Let $J$ and $J'$ be intervals such that $J' \subseteq J$. Then for any
$x\in \PL_2(J)$, such that $\partial_2 D(x)$ does not contain any points in the interior of $J'$ we have
the restriction of
$$
C_{\PL_2(J)}(x) \cap \{g \in \PL_2(J) \mid g(J')=J' \}
$$
to the interval $J'$ is either trivial in the case that $x$ does not preserves the interval $J'$
or $C_{\PL_2(J')}(x)$ otherwise.

\medskip

\noindent \emph{Proof of the Claim.}
Let $g \in C_{\PL_2(J)}(x) \cap \{g \in \PL_2(J) \mid g(J')=J' \}$.
If $x(J')=J'$ then it is immediate that $g|_{J'} \in C_{\PL_2(J')}(x)$.
Suppose now that $x(J') \ne J'$ and $g|_{J'} \ne id$ and say that
$J'=[\gamma_1,\gamma_2]$.
Thus $x(\gamma_1) \ne \gamma_1$ or $x(\gamma_2) \ne \gamma_2$.
Without loss of generality we can assume that $x(\gamma_1) \ne \gamma_1$.
Let $[c,d]$ be the largest interval containing $\gamma_1$ such that
$x(t) \ne t$ for any $t \in \mathbb{Z}[1/2] \cap (c,d)$.
The proof of the previous Proposition implies that $g$ coincides with
$(\sqrt[M]{x})^k$
for some root of $x$ and some integer $k$ on the interval $[c,d]$.
Since $\sqrt[M]{x}(\gamma_1) \ne \gamma_1$, $k$ must be $0$ and so $g$
is the identity map
on $[c,d]$. The restrictions on $J'$ and $J$ imply that $J' \subseteq [c,d]$,
which completes the proof. $\square$

\noindent By the previous claim we see that, for each $i=1,\ldots,r$ and $j=1, \ldots, k$,
the restriction of the subgroup $C_{\PL_2(I)}(x_j) \cap \{g \in \PL_2(I) \mid g(J_i)=J_i \}$ is either trivial
or equal to $C_{\PL_2(J_i)}(x_j)$. Thus $C_{J_i}=id$ or $C_{J_i}=C_{\PL_2(J_i)}(x)$ for some $x \in \PL_2(I)$ which,
by the previous Proposition, is isomorphic with $\mathbb{Z}$ or $\PL_2(J_i)$  $\square$

\corollary{The intersection of any number $k\ge2$ centralizers of elements of $F$ is equal
to the intersection of 2 centralizers. \label{thm:intersection-two}}

\noindent \emph{Proof.} We use the same notation of Proposition \ref{intersectionofcetralizers},
where the $x_i$'s do not denote the standard generators of $\PL_2(I)$ seen in Chapter \ref{chapter1},
but only some arbitrary elements of the group.
Let $C=C_{\PL_2(I)}(x_1) \cap \ldots \cap C_{\PL_2(I)}(x_k)$ be the intersection of
$k\ge 2$ centralizers of elements of $F$. By the previous Proposition we have $I=J_1 \cup \ldots \cup J_r$
and $C=C_{J_1} \cdot \ldots \cdot C_{J_r}$. We want to define $w_1,w_2 \in \PL_2(I)$ such that
$C=C_{\PL_2(I)}(w_1) \cap C_{\PL_2(I)}(w_2)$. We define them on each interval $J_i:=[\alpha_i,\alpha_{i+1}]$,
depending on $C_{J_i}$. \emph{Case 1:} If $C_{J_i}=id$, then we define $w_1,w_2$
as any two elements in $\PL_2^<(J_i)$ such that are not one a
power of another. \emph{Case 2:} If $C_{J_i}\cong \langle x \rangle$ for some $id \ne x \in \PL_2(J_i)$,
then we define $w_1=w_2=x$.
\emph{Case 3:} If $C_{J_i}=\PL_2(J_i)$, then we define $w_1=w_2=id$. $\square$

\question{Corollary \ref{thm:intersection-two} determines that any intersection
of more than one centralizer of elements in $F$ can be expressed as the intersection
$C_F(w_1) \cap C_F(w_2)$ for two suitable elements $w_1,w_2 \in F$. Is it possible
to build the two elements $w_1,w_2$ inside the subgroup $\langle x_1,\ldots,x_k\rangle$?
The current proof does not give an answer to this question.}

\section{The $k$-Simultaneous Conjugacy Problem in $\PL_2(I)$}
\label{sec:simconj}

We will make a sequence of reductions to solve first a particular case. These reductions will use the fact
that we are able to solve the ordinary conjugacy problem. First we notice that, since we know how to solve
the ordinary conjugacy problem, then solving the $(k+1)$-simultaneous conjugacy problem is equivalent to find a positive
answer to the following problem:
\problem{Is there an algorithm such that given
$(x_1, \ldots, x_k, y)$ and $(x_1, \ldots, x_k,z)$ it can decide whether
there is a function $g \in C_{\PL_2(I)}(x_1) \cap \ldots \cap C_{\PL_2(I)}(x_k)$
such that $g^{-1} y g = z$? \label{thm:equiv-prob}}

\noindent Since we understand the structure of the intersection of centralizers,
we are going to work on solving this last question. Our strategy now is to reduce
the problem to the ordinary conjugacy problem and to isolate a very special case that must be dealt with.

\subsection{General case: any $k$ and any centralizer}
\label{ssec:general}

\medskip

\noindent This subsection deals with the general case. We will first extend Theorem \ref{congdiagonals}
and then we will use our description for the intersection of many centralizers to
solve the general problem. The argument of Proposition~\ref{simconfork=2}
will show us that we can  build possible conjugators by
using the stair algorithm and then check if they are in an intersection of
centralizers. This will be verifiable, since we have given a description of such intersection in
Proposition~\ref{intersectionofcetralizers}.

\lemma{Let $x_1,\ldots,x_k,y,z \in \PL_2(J)$. We can determine whether there is, or not,
a $g \in C=C_{\PL_2(J)}(x_1) \cap \ldots \cap C_{\PL_2(J)}(x_k)$ such that $g(D(y))=D(z)$.
\label{simdiagonals}}

\noindent \emph{Proof.} The proof is essentially the same as that of Corollary \ref{thm:overlap-intersection-diagonal}
on each of the intervals between two dyadic fixed points of $y$ and $z$. The only new tool
required is Lemma~\ref{conginZ} on the intervals where $C$ is isomorphic to $\mathbb{Z}$. We omit
the details of this proof. $\square$

\proposition{Let $x_1,\ldots,x_k,y,z \in \PL_2(J)$. We can determine whether there is, or not,
a $g \in C=C_{\PL_2(J)}(x_1) \cap \ldots \cap C_{\PL_2(J)}(x_k)$ with $g^{-1}y g=z$. \label{simconfork=2}}

\noindent \emph{Proof.} Apply Lemma \ref{simdiagonals} to make $D(y)=D(z)$, if possible. Recall
that a candidate conjugator must centralize $x_1,\ldots,x_k$ too, so it has to fix
$\bigcup^k_{i=1} \partial_2 D(x_i)$ and $\partial_2 D(y)=\partial_2 D(z)$. Let
$\bigcup^k_{i=1} \partial_2 D(x_i)= \{\lambda_m \}_m$ and $\partial_2 D(y)=\{\mu_1 < \ldots < \mu_k\}$
and let $J_i$ denote the interval $[\mu_i,\mu_{i+1}]$, for $i=1,\ldots,k-1$.
We build $g$ on each interval $J_i$, depending on how $y$ is defined on $J_i$. We have the following three cases:
\\ \emph{Case 1:} $y$ is the identity on $J_i$. In this case we define $g$ to be the identity on $J_i$.
\\ \emph{Case 2:} $y$ is not the identity on $J_i$ and there is a point
$\lambda_j \in \bigcup^k_{i=1} \partial_2 D(x_i)$
which is in the interior of $J_i$. Since $\mu_i<\lambda_j<\mu_{i+1}$ and $\lambda_j$ is dyadic, then
$\lambda_j \not \in \partial_2 D(y)$ and in particular $\lambda_j$ is not fixed by $y$ and $z$.
Since $g(\lambda_j)=\lambda_j$, the proof of Lemma \ref{cong:unique}(ii) implies that $g'(\mu_i^+)=\lim_{n \to \infty}
\frac{y^n(\lambda_j)-\mu_i}{z^n(\lambda_j)-\mu_i}$, hence the slope of $g$ on the right of $\mu_i$
is uniquely determined.
Therefore we can apply Proposition \ref{thm:psi-surjective2}(iii) to build the unique candidate conjugator $g$.
\\ \emph{Case 3:} $y$ is not the identity on $J_i$ and $\bigcup^k_{i=1} \partial_2 D(x_i)$
does not contain any point of the interior
of $J_i$. More precisely, each $x_r$ does not fix any point in $J_i$ and so, by the Claim contained
in the proof of Proposition~\ref{intersectionofcetralizers} we have that the restriction group
\[
C_{\PL_2(J)}(x_r) \cap \{g \in \PL_2(J) \mid g(J_i)=J_i \}
\]
is the trivial group or $\PL_2(J_i)$ or isomorphic to a copy of $\mathbb{Z}$.
Since $C_{J_i}$ is the intersection of all the restriction groups for $r=1, \ldots, k$,
then $C_{J_i}$ will also be trivial or $\PL_2(J_i)$ or infinite cyclic.
If $C_{J_i}$ is trivial, we choose $g$ to be trivial on $J_i$. If $C_{J_i} = \PL_2(J_i)$ then
the construction reduces to solving
the ordinary conjugacy problem in $\PL_2(J_i)$. The case $C_{J_i} \cong \mathbb{Z}$
will be covered in Subsection ~\ref{ssec:specialcase}.
\\ Finally we have to verify that the element $g$ constructed by the above
procedure conjugates $y$ to $z$ and commutes with $x$. $\square$

\medskip

\noindent The restatement of the $k$-simultaneous conjugacy problem given in Problem
\ref{thm:equiv-prob} and the previous Proposition imply the result of Theorem A in the introduction.

\subsection{A special case: $k=1$ and $C_{\PL_+(J)}(x) \cong \mathbb{Z}$}\label{ssec:specialcase}

\medskip

\noindent This subsection is technical and it deals with a variant of the ordinary conjugacy problem. We want to see
if we can solve it when we have a restriction on the possible conjugators. Thus, given $x,y,z$
we want to see if $g^{-1}yg=z$ for a $g \in C_{\PL_2(J)}(x) \cong \mathbb{Z}$. In particular, if
$\sqrt[M]{x}$ is the ``smallest possible'' root (in the sense of the proof of centralizers in $\PL_2(J)$),
then we need to find if there is a power of $\sqrt[M]{x}$ which conjugates $y$ to $z$.
Since $C_{\PL_2(J)}(x) = C_{\PL_2(J)}(\sqrt[M]{x})= \langle \sqrt[M]{x} \rangle$ then we can substitute
$x$ with $\widehat{x}:=\sqrt[M]{x}$. For simplicity, we assume still call $\widehat{x}$ with $x$.
The plan for this subsection will be to reduce to solve an equation of the type
\[
f^k= w h^k
\]
where $f,h,,w$ are given, $w'(\eta^+)=1$ and we need to find if there is any $k \in \mathbb{Z}$ satisfying
the previous equation. The second step will be to prove that there is only a finite number of $k$'s to that
may solve the equation and so we try all of them.

\noindent We need first to run the usual conjugacy problem on $[\eta,\zeta]$
between $y$ and $z$ to see if they are conjugate. If they are, we continue. Otherwise we stop.
Let $C_{\PL_2(J)}(y,z)=\{g \in \PL_2(J)\mid g^{-1} y g(t)= z(t),$ for all $t \in J\} =g_0 \cdot
C_{\PL_2(I)}(y)$, for some $g_0 \in T$. Now $C_{\PL_2(J)}(y) \cong \mathbb{Z}^s \times
\PL_2(J)^t$. Notice that $s=t=0$ is impossible.

\medskip

\noindent If $s+t \ge 2$, then there must be some
$\tau \in (\eta,\zeta) \cap \mathbb{Z}[\frac{1}{2}]$ fixed point for every element
in $C_{\PL_2(I)}(y)$. So if $y$ and $z$ are conjugate through a power of $x$ then there is a $k$
such that $x^k(\tau)=g_0(\tau)$. Notice $x(\tau) \ne \tau$, so we apply Lemma \ref{conginZ} with $\mu := g_0(\tau)$
and find, if possible a unique integer $\bar{k}$ such that $x^{\bar{k}}(\tau)=\mu$. Now we take $g:=x^{\bar{k}}$ and
we check if it is a conjugator or not.

\medskip

\noindent If $s=0, t=1$, then this would mean that $C_{\PL_2(J)}(y) \cong \PL_2(J)$ and so that
$y=id $ on $[\eta,\zeta]$ and so do not need to check the powers of $x$, but simply if the function
$z=id$ on $[\eta,\zeta]$.

\medskip

\noindent If $s=1, t=0$, then $C_{\PL_2(J)}(y) = \langle \widehat{y} \rangle \cong \mathbb{Z} $,
for $\widehat{y}$ a generator. Thus, $y$ and $z$ are conjugate through an element of  $C_{\PL_2(J)}(x)$,
if and only if there exist $k,m \in \mathbb{Z}$ such that $x^m=g_0 \widehat{y}^n$ in $[\eta,\zeta]$.

\lemma{Let $x,y,z \in \PL_2(J)$ such that $C_{\PL_2(J)}(x)=\langle x \rangle$ and $C_{\PL_2(J)}(y)=\langle
\widehat{y} \rangle$. Then there exists
$X,Y,G_0 \in \PL_2(J)$ such that $G_0'(\eta^+)=1$ and following two problems are equivalent:

\noindent (i) Find powers $k,m \in \mathbb{Z}$ such that $x^m=g_0 \widehat{y}^n$

\noindent (ii) Find a power $k \in \mathbb{Z}$ such that $X^k=G_0 Y^k$. \label{thm:special-case}}

\noindent \emph{Proof.} Suppose we have $x'(\eta^+)=2^\alpha, \widehat{y}'(\eta^+)=2^\beta, g_0'(\eta^+)=2^\gamma$
for some $\alpha,\beta,\gamma \in \mathbb{Z}$, then we must
have $2^{\alpha m}=(x^m)'(\eta^+)=(g_0 \widehat{y}^n)'(\eta^+)=2^{\gamma + \beta n}$
and so $\alpha m = \gamma + \beta n$.
Thus, in order for $y$ and $z$ to be conjugate we must have
that $\gcd(\alpha,\beta)$ divides $\gamma$. That is,
$\gamma = \alpha m_0 - \beta n_0$, for some $m_0,n_0 \in \mathbb{Z}$ which can be computed and thus
$\alpha (m-m_0)=\beta(n-n_0)$. We can change variables and call $\widetilde{m}=m-m_0$ and $\widetilde{n}=n-n_0$. So we have to find
$\widetilde{m},\widetilde{n}$ such that $\alpha \widetilde{m}=\beta \widetilde{n}$ and so that
\[
\frac{\alpha}{\gcd(\alpha,\beta)}\widetilde{m} = \frac{\beta}{\gcd(\alpha,\beta)}\widetilde{n}
\]
Thus there must exist a $k \in \mathbb{Z}$ such that
\[
\widetilde{m} = \frac{\beta}{\gcd(\alpha,\beta)}k \; \; \textrm{and} \; \;
\widetilde{n} = \frac{\alpha}{\gcd(\alpha,\beta)}k.
\]
Going backwards, we write
\[
m := \frac{\beta}{\gcd(\alpha,\beta)}k + m_0 \; \; \textrm{and} \; \;
n := \frac{\alpha}{\gcd(\alpha,\beta)}k + n_0.
\]
By substituting these two values in the equation $x^m=g_0 y^n$ we get
\[
(x^{\frac{\beta}{\gcd(\alpha,\beta)}})^k = x^{-m_0} g_0 \widehat{y}^{n_0}
(\widehat{y}^{\frac{\alpha}{\gcd(\alpha,\beta)}})^k
\]
We rename $X= x^{\frac{\beta}{\gcd(\alpha,\beta)}}$, $G_0 = x^{-m_0} g_0 \widehat{y}^{n_0}$ and
$Y=\widehat{y}^{\frac{\alpha}{\gcd(\alpha,\beta)}}$ and so we are left to find a
$k \in \mathbb{Z}$, if it exists, such that
\begin{equation}
\label{eq:power}
X^k = G_0 Y^k.
\end{equation}
Notice that, with these adjustments, $G_0(\eta^+)=2^0=1$. $\square$

\medskip

\noindent In the last case we are examining, both $x$ and $y$ cannot have fixed
dyadic points, since their centralizers are cyclic groups.
Thus $D(x) \cap (\eta,\zeta)$ and $D(y)\cap(\eta,\zeta)$ must be empty or finite. The same is also
true for the new functions $X$ and $Y$, i.e. $D(X) \cap (\eta,\zeta)$ and $D(Y) \cap (\eta,\zeta)$
must be empty or finite.
For sake of simplicity, we will still call $X,Y,G_0$ with lowercase letters. We will make distinction
in the following cases, by checking what are $D(x) \cap (\eta,\zeta)$ and $D(y) \cap (\eta,\zeta)$
and see if they coincide or not.

\medskip

\noindent $D(x) \cap (\eta,\zeta) \ne D(y) \cap (\eta,\zeta)$. There exists a $\tau \in
(\eta,\zeta)$ with $y(\tau)=\tau \ne x(\tau)$.
Thus, by applying Lemma \ref{conginZ},
we can determine if there is a $k$ such that $x^k(\tau)=g_0(\tau)$. We act similarly if
there is a $\tau \in (\eta,\zeta)$ with $x(\tau)=\tau \ne y(\tau)$.

\medskip

\noindent $D(x) \cap (\eta,\zeta)=D(y) \cap (\eta,\zeta) \ne \emptyset$. Suppose
$D(x)=D(y)=\{r_1<\ldots<r_v\}$. Observe that if the equation has a
solution then $g_0(r_i)=r_i$ for all $r_i$. If these conditions are satisfied, then we can build all the solutions
by solving the equation in each interval $[r_i,r_{i+1}]$. This reduces the problem to the next case.

\medskip

\noindent $D(x) \cap (\eta,\zeta)=D(y) \cap (\eta,\zeta)=\emptyset$,
that is we have that $x,y \in \PL_2^<(J) \cup \PL_2^>(J)$.
We can now assume that both $x,y \in \PL_2^<(J)$.
Define
$$
K:=\{k \in \mathbb{Z} \, \, \mbox{such that $x^k(t)=g_0(y^k(t))$ for all} \, \, t \in J\}.
$$
\noindent Our goal is to find whether or not $K \ne \emptyset$. The first step will be to prove
that the set $K$ is finite, by computing directly its upper and lower bounds. Therefore, we will have
that $K\subseteq \mathbb{Z}\cap [l_0,k_0]$, for some integers $l_0,k_0$, and so we can check
all these integers and see if any satisfies $x^k(t)=g_0(y^k(t))$.

\lemma{Let $x,y \in \PL_2^<(J)$ and let $K:=\{k \in \mathbb{Z}$ such that $x^k=g_0 y^k\}$, then
$K$ is bounded.}

\noindent \emph{Proof.} The first step is to prove that there exists
a $k_0 \in \mathbb{Z}$, upper bound for $K$.
Suppose that $K$ has no upper bound. Let $\theta < \zeta$
be a point such that $g_0(t)=t$ and $x(t)=y(t)$ on $[\eta,\theta]$.
Let $\psi>\theta$ a number such that $x(\psi)<y(\psi)$ and $x(t)\le y(t)$ for $t \le \psi$.
Since $y \in \PL_2^<(J)$ then $\underset{k \to \infty}{\lim}y^{k}(\psi)=\eta$, and so
we can choose a $k_0 \in K$ be a large enough number such that
$y^{k_0}(\psi)<\theta$. Suppose $k \ge k_0$, by definition of $\theta$ and $k_0 \in K$ we have
\[
x^{k}(\psi)=g_0(y^{k}(\psi))= y^{k}(\psi).
\]
Now recall that $x(\psi) < y(\psi) < \theta +
\varepsilon$ and so, since $x \le y$ on $[\eta,\psi]$
\begin{eqnarray*}
x^{k}(\psi)=x^{k-1}(x(\psi))<x^{k-1}(y(\psi)) \\
= x^{k-2}(x(y(\psi)) \le x^{k-2}(y^2(\theta +\varepsilon)) \\
\le \ldots \le x (y^{k-1}(\psi)) \le y^{k}(\psi).
\end{eqnarray*}
By comparing the last two expressions, we get
$x^{k}(\psi) < y^{k}(\psi) = x^{k}(\psi)$. Contradiction.
Therefore $k_0$ is an upper bound for $K$.

\noindent We now want to bound the $K$ from below, and so we use a similar technique.
If $k \in K$ is negative, then we consider the
equation
\[
y^{-k}=x^{-k}g_0 = g_0 (g_0^{-1}x^{-k}g_0)=g_0(g_0^{-1}xg_0)^{-k}=g_0 \widehat{x}^{-k}
\]
where we have set $\widehat{x}:=g_0^{-1}xg_0$. Since $D(x)=\emptyset$, then $D(\widehat{x})=\emptyset$ and
$\widehat{x} \in \PL_2^<(I)$. So we have reduced to the situation of the previous claim (with $\widehat{x}$ and $y$
switched in their role) and we obtain that the set of possible $(-k)$'s is bounded above,
so that $k$ is bounded below. $\square$

\medskip

\noindent Since $K$ is finite the $k$'s to be checked are finite and we can find its bound in finite time.
Now we can check all possible the elements of $K$ and we conclude this case.

\subsection{The twisted conjugacy problem for $\PL_2(I)$}

We conclude this section by describing an interesting open question for Thompson's group $F$.
It has been shown by Bogopolski, Martino, Maslakova and Ventura \cite{bomave1} and \cite{bomave2} that
the conjugacy problem for certain extensions of groups can be reduced to solving the twisted
conjugacy problem for a subgroup. We say that a group $G$ has \emph{solvable $\varphi$-twisted conjugacy problem},
for a given $\varphi \in \mathrm{Aut}(G)$, if there is an algorithm such that,
given any two elements $y,z \in G$, we can determine whether there is, or not, a $g \in G$ such that
$\varphi(g)^{-1} y g = z$. Brin \cite{brin5} has classified the structure of automorphisms
of Thompson's group $F$. Let $\PL_{dis}(\mathbb{R})$ denote the group of
piecewise-linear orientation-preserving homeomorphisms with finitely many breakpoints occurring
at dyadic rational coordinates, such that every slope is an integral power of $2$
and with a discrete set of breakpoints (infinitely many breakpoints are possible):
then $F$ is isomorphic to the subgroup $H$
of $\PL_{dis}(\mathbb{R})$ of elements $f$ such that there is an interval $[a_f,b_f]$
and $f(t)=t+m_f$ on $(b_f,\infty)$ and $f(t)=t+k_f$ on $(-\infty,a_f)$, for
some integers $m_f,k_f$. Then any orientation-preserving automorphism of $H$
is given by the maps $\varphi:H \to H$, defined by $\varphi(g)= \tau^{-1} g \tau$
where $\tau \in A$, the subgroup of $\PL_{dis}(\mathbb{R})$ such that
there is an interval $[m_\tau,n_\tau]$ with $\tau(t+1)=\tau(t)+1$ for any $t \not \in [m_\tau,n_\tau]$.
Hence, if we rewrite the $\varphi$-twisted conjugacy problem for Thompson's group $F$ seen as the subgroup $H$
we can rewrite it as
\[
z=\varphi(g)^{-1} y g= \tau^{-1}g^{-1}\tau y g
\]
and so it becomes
\[
g^{-1} \left(\tau y \right) g = \tau z
\]
that is, a conjugacy problem for the elements $\tau y, \tau z \in A$ with respect to an element $g \in F$.
It is an interesting problem to work on solving this generalization of the conjugacy problem
and see if any of the presented techniques of closed diagrams, Mather invariants or the Stair Algorithm
can be extended to this setting.

\question{Is the twisted conjugacy problem solvable for the group $F$?}

\section{Stair Algorithm in $\PL_{S,G}(I)$ \label{sec:stair-algorithm}}

We now move on to prove the solvability of the simultaneous conjugacy problem to other subgroups of $\PL_+(I)$
whose structure generalizes that of Thompson's group $F$. We remark that Brin and Squier \cite{brin2}
give a criterion for conjugacy in $\PL_+(I)$.
Let $S$ be a subring of $\mathbb{R}$, let $U(S)$ be the group of invertible elements of $S$ and let
$G$ be a subgroup of $U(S) \cap \mathbb{R}_+$. For any subinterval $J$ of $I$, we define $\PL_{S,G}(J)$ to be the
group of piecewise linear orientation-preserving homeomorphism from the interval $J$ into itself,
with only a finite number of breakpoints and such that
\begin{itemize}
\item
all breakpoints are in the subring $S$,
\item
all slopes are in the subgroup $G$,
\end{itemize}
the product of two elements is given by the composition of functions. If $G=U(S) \cap \mathbb{R}_+$ we
write $\PL_S(J)$ instead of $\PL_{S,G}(I)$. Thompson's group $F$ can thus be recovered as the group
$\PL_{\mathbb{Z}{\left[\frac{1}{2}\right]},\langle 2\rangle}(I)$.
We observe that, in order to make some calculations possible inside the ring $S$ and its quotients,
we need to ask for some requirements to be satisfied by $S$ from the computability standpoint.
These assumptions will be clearly stated in Remark \ref{thm:ring-requirements}
and will be assumed throughout the remainder of the chapter.

We introduce briefly the notation to generalize the results obtained on $F$.
For a subset $J \subseteq [0,1]$, we denote with $\partial J$ the usual boundary of $J$ in $[0,1]$.
For an interval $J=[\eta,\zeta] \subseteq I$ such that $\partial J \subseteq S$,
a function $f \in \PL_{S,G}(J)$ can be extended to the interval $I$ by
$f(t)=t$ for $t \in I \setminus J$, which allows us to consider $\PL_{S,G}(J)$ as a subgroup of $\PL_{S,G}(I)$.
For a function $f \in \PL_{S,G}(J)$ we define
$$
D_J(f):=\{t \in J \;\big\vert\; f(t)=t \},
$$
to simplify the notation will often drop the subscript $J$.

\begin{defin}
\emph{We define $\PL_{S,G}^<(J)$ (and respectively. $\PL_{S,G}^>(J)$) to be the set of all functions
in $\PL_{S,G}(J)$ with graph below the diagonal (respectively, above the diagonal).}
\end{defin}

\noindent Given a function $f \in \PL_{S,G}(I)$ and a number $0<t_0<1$ fixed by $f$, it is not
always true that $t_0 \in S$. For any subset $J \subseteq I$ we define
$$
\textstyle
\partial_S J:=\partial J \cap S
$$

\definition{We define $\PL_{S,G}^0(J) \subseteq \PL_{S,G}(J)$, the set of functions $f \in \PL_{S,G}(J)$ such that
the set $D(f)$ does not contain elements of $S$ other than the endpoints of $J$, i.e., $D(f)$ is discrete and
$\partial_S D(f)=\partial_S J$.}

\question{Given two elements $\alpha,\beta \in S$, is it true that there is a $g \in \PL_{S,G}(J)$ such that
$g(\alpha)=\beta$? We will now discuss conditions to verify this generalization of Proposition \ref{thm:pipeline}
and of Corollary \ref{thm:how-to-identify-rationsls}.}

\definition{We define an ideal in $S$ given by $\mathcal{I}_{S,G}=\langle (g-1) \mid g \in G \rangle$.
We denote with $\pi_{S,G}:S \to S/\mathcal{I}$ the natural quotient map.
Unless otherwise stated,
we will drop the subscript and write $\mathcal{I}$ and $\pi$ instead of $\mathcal{I}_{S,G}$ and $\pi_{S,G}$.}

The following two results are used to detect when two points of $S$ are in the same $\PL_{S,G}$-orbit.

\lemma{Let $J \subseteq [0,1]$ be a closed interval with
at least one of the endpoints $\eta$ in $S$ and let $g \in \PL_{S,G}(J)$.
Then, for every $t \in J \cap S$, we have $\pi(g(t))=\pi(t)$.\label{thm:nec-cond-trans}}

\noindent \emph{Proof.} We can assume that the $\eta$ is the left one
and we apply induction on the number of breakpoints before $t$. In case the endpoint in $S$
is the right one, we apply induction on the breakpoints after $t$. Let $\{\eta_1, \ldots, \eta_r\}$
be the set of all breakpoints of $g$ on the interval $[\eta,t)$. Then $g(t)=c_r (t-\eta_r) +g(\eta_r)$ for some suitable $c_i \in G$.
By induction on $r$ we have that $\pi(g(\eta_r))=\pi(\eta_r)$ and thus
\begin{eqnarray*}
\pi(g(t))=\pi(c_r (t-\eta_r)+ g(\eta_r))=\\
\pi(c_r -1)\pi(t-\eta_r) + \pi(1)\pi(t-\eta_r) + \pi(g(\eta_r))= \\
\pi(t-\eta_r)+\pi(\eta_r)=\pi(t). \square
\end{eqnarray*}

\noindent This result gives us a necessary condition on how homeomorphisms can be built.
We want to know what orbits of elements are under the action of $\PL_{S,G}(J)$.

\proposition{Let $J \subseteq [0,1]$ be a closed interval with both endpoints in $S$
and let $u,v \in J \cap S$. Then $\pi(u)=\pi(v)$
if and only if there is a $g \in \PL_{S,G}(J)$ such that
$g(u)=v$.\label{thm:equiv-cond-trans}}

\noindent The proof of this proposition can be found in the Appendix (see Proposition \ref{thm:equiv-cond-trans-appendix}).

\remarkname{Computational Requirements}{We need to add a few requirements
to the ring $S$ in order to make a machine able to work with the algorithm.
It is reasonable to make the following assumptions to work in the ring $S$:
\begin{itemize}
\item There is solution to the membership problem in $S$ (\emph{i.e.}
an algorithm to determine whether an element $s \in \mathbb{R}$ lies in $S$ or not)
\item There is a solution to the membership problem in $\mathcal{I}$
\item There is an algorithm that, for every $q \in S$, is able to determine
whether two elements in the quotient ring $S/q \mathcal{I}$ are equal or not.
\item There is an algorithm such that, given $a,b,c \in G$, it is able to determine whether or not
there exist $x,z \in \mathbb{Z}$ such that $a^x = b c^z$.
\end{itemize}
All these requirements are reasonable to assume in order to make computations inside $S$ and will
be checkable in the special cases that we take as examples in Section \ref{sec:examples}.
\label{thm:ring-requirements}}

\remark{By taking logarithms, we can rewrite
all of the terms of the last requirement on the algorithm in base $b$, so that it becomes equivalent
to the following: given any $\alpha,\beta,\gamma \in \mathbb{R}$,
determine whether or not they are linearly dependent over $\mathbb{Q}$ and, if they are, we can find
$q_1,q_2 \in \mathbb{Q}$ such that $\gamma=q_1 \alpha + q_2 \beta$. This rewriting transforms the equation
$a^x = b^y c^z$ into a $\mathbb{Q}$-linearity dependence relation, hence if there is a solution, it is unique.
\label{thm:rewrite-requirement}}

\remark{In general, given two intervals $J_1,J_2$ with endpoints in $S$, the groups
$\PL_{S,G}(J_1)$ and $\PL_{S,G}(J_2)$ may not be isomorphic (that is, the analogue of Theorem \ref{thm:thompson-like}
may not hold).
Proposition \ref{thm:equiv-cond-trans} tells us that two elements in $S$ are in the same $\PL_{S,G}$-orbit if their image
under the map $\pi$ is the same.
For example in the cases
$S=\mathbb{R},G=\mathbb{R}_+$ and $S=\mathbb{Q},G=\mathbb{Q}^*$ and $S=\mathbb{Z}\left[\frac{1}{2}\right],G=\langle 2 \rangle$,
it is not difficult to see that every two points in $S$ have the same image under $\pi$
and that any two groups $\PL_{S,G}(J_1)$ and $\PL_{S,G}(J_2)$ are thus isomorphic, for any two intervals
$J_1,J_2$ with endpoints in $S$. On the other hand, if we consider generalized Thompson's groups
(see Section \ref{sec:examples}), it can be shown that the number of orbits is finite but more than one,
so that are only finitely many isomorphism classes for the groups $\PL_{S,G}(J)$,
for $S=\mathbb{Z}\big[\frac{1}{n_1}, \ldots, \frac{1}{n_k}\big]$
and $G=\langle n_1, \ldots, n_k \rangle$ for $n_1,\ldots,n_k \in \mathbb{Z}$.
In general, it seems likely that if two elements $\alpha,\beta \in S$ have different image under $\pi$
then the groups $\PL_{S,G}([0,\alpha])$ and $\PL_{S,G}([0,\beta])$ are not isomorphic, but it is not easy to prove it.
\label{thm:multiple-orbits-of-points}}

\subsection{Making $D(y)$ and $D(z)$ coincide}
\label{ssec:diag-intersec}

We start by generalizing Proposition \ref{thm:equiv-cond-trans} to a finite number of points.

\lemma{Let $J=[\eta,\zeta] \subseteq [0,1]$ be a closed interval with endpoints in $S$ and suppose we have
$u_1,v_1,\ldots,u_k,v_k \in J \cap S$ such that $\pi(u_i)=\pi(v_i)$ for all $i=1,\ldots,k$.
Then there exists a $g \in \PL_{S,G}(J)$ such that $g(u_i)=v_i$ for all $i=1,\ldots,k$.
\label{thm:multiple-trans}}

\noindent \emph{Proof.} We can assume that $J=[\eta,\zeta]$ and that
the $u_i$'s are ordered in an increasing sequence $u_1 < \ldots <u_k$
and therefore $v_1 < \ldots < v_k$. By Proposition \ref{thm:equiv-cond-trans}, there is an $g_1 \in \PL_{S,G}(J)$ such that
$g_1(u_1)=v_1$. Now we notice that $v_1=g_1(u_1) < g_1(u_2) < \ldots < g_1(u_k)$ and so we
restrict to the interval $[v_1,\zeta]$ and, since $\pi(g_1(u_2))=\pi(u_2)=\pi(v_2)$ we can
use again Proposition \ref{thm:equiv-cond-trans} to find an $h_2 \in \PL_{S,G}([v_1,\zeta])$ such that $h_2(g_1(u_2))=v_2$.
Define
\[
g_2(t):=\begin{cases}
t & t \in [\eta,v_1] \\
h_2(t) & t \in [v_1,\zeta]
\end{cases}
\]
so that $g_2g_1(u_1)=v_1,g_2g_1(u_2)=v_2$ and $g_2 \in \PL_{S,G}(J)$.
By iterating this procedure, we build functions $g_i \in \PL_{S,G}(J)$ such that $g_ig_{i-1}\ldots g_1(u_j)=v_j$
for all $j=1,\ldots,i$ and $i=1,\ldots,k$. Thus we define $g:=g_k g_{k-1} \ldots g_1 \in \PL_{S,G}(J)$ and we
get a function such that $g(u_i)=v_i$. $\square$

\noindent The previous Lemma yields the following natural generalization of the Extension Lemma
\ref{thm:extension-partial-maps} which we state without proof.

\lemmaname{Extension of Partial Maps}{Let $J \subseteq [0,1]$ be a closed interval
with endpoints in $S$ and suppose
$I_1,\ldots,I_k \subseteq J$ is a finite family of disjoint closed
intervals in increasing order and of the form $I_i=[a_i,b_i]$,
for all $i=1,\ldots,k$ and $a_i,b_i \in S$. Let
$J_1, \ldots, J_k \subseteq J$, with
$J_i=[c_i,d_i]$, be another family of intervals with the same property and such that
$\pi(a_i)=\pi(c_i)$ and $\pi(b_i)=\pi(d_i)$.
Suppose that $g_i:I_i \to J_i$ is a piecewise-linear function with a finite number
of breakpoints, occurring at points in $S$ and with slopes in $G$.
Then there exists a $\widetilde{g} \in \PL_{S,G}(J)$ such that
$\widetilde{g}|_{I_i}=g_i$. $\square$ \label{thm:general-ext-part-maps}}

\noindent Let $g \in \PL_{S,G}(J)$ be equal to $g(t)=at+b$ around a point $q\in \mathbb{R}$
fixed by $f$, for some $a \in G,b \in S$, then $q=b/(1-a)$
and so the intersection points of $f$ with the diagonal lie in $Q_S$, the field of fractions of $S$.
Now that we have a way to recognize whether we can make two elements of $S$ coincide through an element
of $\PL_{S,G}(J)$, we need to see if it is possible to do the same
for the field of fractions $Q_S$.

\proposition{Let $J=[\eta,\zeta] \subseteq [0,1]$ be a closed interval with endpoints in $S$ and
let $\alpha,\beta \in J \cap Q_S$. There is a $g \in \PL_{S,G}(J)$ with $g(\alpha)=\beta$
if and only if we can find $p,q,r \in S$ such that $\alpha=p/q,\beta=r/q$ and
\[
pG=rG \pmod{q \mathcal{I}}
\]
where $q \mathcal{I}$ denotes the product of the ideal generated by $q$ and $\mathcal{I}$.}

\noindent \emph{Proof.}
Suppose there is a map $g \in \PL_{S,G}(J)$ with $g(\alpha)=\beta$
and let $g(t)=ct+d$ in a small neighborhood $J_\alpha$ of $\alpha$.
We can choose representatives $p,q,r \in S$ such that
$\alpha=p/q,\beta=r/q$ and then, since $g \in \PL_{S,G}(J)$, we use Lemma \ref{thm:nec-cond-trans} to get
\[
\pi(t)=\pi(g(t))=\pi(c-1)\pi(t) + \pi(t) + \pi(d)
\]
for all $t \in J_\alpha \cap S$ and therefore $\pi(d)=0$, which
implies $d \in \mathcal{I}$. Conversely, suppose that we can write
$\alpha=p/q,\beta=r/q$, for some $p,q,r \in S$ and that $pG=rG \pmod{q \mathcal{I}}$.
The second condition implies that there exist $c_1,c_2 \in G, d_2 \in \mathcal{I}$ such that
\[
c_1 r = c_2 p + q d_2
\]
and so if we set $c=c_2/c_1$ and $d=d_2/c_1$, we get $r=cp+qd$. Let $f(t)=ct+d$ be a line
through the point $(\alpha,\beta)$ and let $[\gamma,\delta] \subseteq J$ be a small
interval such that $\gamma,\delta \in S$. Finding the interval $[\gamma,\delta]$ can be accomplished this way:
we can assume $G \ne 1$
and pick any $1 \ne c \in G$ such that $0<c<1$. Then we choose $m,n \in \mathbb{N}$ such that
$\eta + c^{m} < \alpha < \eta + nc^{m} < \zeta$ and we set $\gamma:=\eta + c^{m},
\delta:=\eta + nc^{m}$. Since $\pi(d)=0$ we have that $\pi(f(\gamma))=\pi(\gamma)$
and $\pi(f(\delta))=\pi(\delta)$ and so, by the Extension Lemma
\ref{thm:general-ext-part-maps} there is an $g \in \PL_{S,G}(J)$ with $g|_{[\gamma,\delta]}=f$.
By construction $g(\alpha)=\beta$ as required. $\square$

\noindent In a similar fashion, we can get the same result for a
finite number of points. This amounts to finding small segments
passing through the rational pairs $(\alpha_i,\beta_i)$
and then applying the Extension Lemma to obtain a homeomorphism
of the whole interval $J$. We thus state without proof the following Lemma.

\lemma{Let $J=[\eta,\zeta] \subseteq [0,1]$ be a closed interval with endpoints in $S$ and
let $\alpha_i,\beta_i \in J \cap Q_S$ for $i=1,\ldots,k$. There is a $g \in \PL_{S,G}(J)$
with $g(\alpha_i)=\beta_i$ if and only if there exist $g_1, \ldots, g_k \in \PL_{S,G}(J)$
such that $g_i(\alpha_i)=\beta_i$.}


\noindent By the assumptions made in Remark \ref{thm:ring-requirements}, we can detect
whether or not two elements in $Q_S$ are equal, thus we obtain the following generalizations
of Corollary \ref{thm:overlap-intersection-diagonal} and Lemma \ref{congdiagonals}:

\corollary{Let $J=[\eta,\zeta] \subseteq [0,1]$ be a closed interval with endpoints in $S$ and
let $\alpha_i,\beta_i \in J \cap Q_S$ for $i=1,\ldots,k$. We can determine whether there is or not
an $f \in \PL_{S,G}(J)$ such that $g(\alpha_i)=\beta_i$ for every $i=1,\ldots,k$. $\square$ }

\proposition{Given $y,z \in \PL_{S,G}(I)$, we can determine whether there is or not
a $g \in \PL_{S,G}(I)$ such that $g(D(y)) = D(g^{-1}yg)= D(z)$. If such a
$g$ exists, we can construct it. $\square$ \label{thm:same-intersection}}

\subsection{Linearity Boxes and Stair Algorithm \label{ssec:lin-box-stair}}

In this Subsection we generalize the results of Subsections \ref{ssec:linearitybox}, \ref{ssec:stair}
and \ref{ssec:stairfull}. First we observe that
two conjugate elements $y,z$ in $\PL_{S,G}(I)$ must have the same slopes around the two endpoints,
then we build an algorithm which makes these two elements coincide in a sequence of steps. More
precisely, we build a sequence of functions $g_1,g_2,\ldots,g_k$ and of intervals
$J_1 \subseteq J_2 \subseteq \ldots \subseteq J_i \subseteq \ldots$ such that $0 \in J_1$
and $g_i^{-1}\ldots g_1^{-1} y g_1 \ldots g_i=z$ on $J_i$. We prove that the procedure will stop
because the two elements $y,z$ coincide around the second endpoint of $I$.
When the algorithm stops, we have that $J_k=[0,1]$.

\lemmaname{Linearity Boxes}{Suppose $y,z,g \in \PL_{S,G}(J)$ and $g^{-1}yg=z$.

\noindent (i) If there exist two numbers $\alpha>0$ and $c \ge 1$ such that
$y(t) = z(t) =c (t - \eta) + \eta$ for $t\in[\eta,\eta+\alpha]$, then
the graph of $g$ is linear inside the square $[\eta,\eta+\alpha]\times[\eta,\eta+\alpha]$

\noindent (ii) If there exist $\beta,c \in (0,1)$ such that $y(t)=z(t)=c \cdot (t-\zeta) + \zeta$
on $[\beta,\zeta]$, then the graph of $g$ is linear inside the square $[\beta,\zeta]\times[\beta,\zeta]$.
\label{thm:lin-box}}

\noindent \emph{Proof.} These results follow from the proofs of
Lemma \ref{thm:initial-box} and Remark \ref{thm:final-box}. $\square$

\noindent We recall that $\PL_{S,G}^0(J)$ denotes the set of functions $f \in \PL_{S,G}(J)$ such that
the set $D(f)$ does not contain elements of $S$ other than the endpoints of $J$.

\propositionname{Stair Algorithm for $\PL_{S,G}^0(J)$}{Let $J \subseteq [0,1]$ be a closed interval with endpoints
in $S$, let $y,z \in \PL^0_{S,G}(J)$ such that $D(y)=D(z)$ and define
$C_{\PL_{S,G}(J)}(y,z)=\{g \in \PL_{S,G}(J) | y^g=z\}$ the set of all conjugators. For any $\tau \in D(y)$
we define the map

$$
\begin{array}{lrcl}
\varphi_{y,z,\tau}: & C_{\PL_{S,G}(J)}(y,z) & \longrightarrow & \mathbb{R_+}  \\
               &  g               & \longmapsto     & g'(\tau),
\end{array}
$$
where if $\tau$ is an endpoint of $J$ we take only a one-sided derivative. Then

\noindent (i) $\varphi_{y,z,\tau}$ is an injective map. In particular, if we define $\varphi_{z,\tau}:=\varphi_{z,z,\tau}$,
then $\varphi_{z,\tau}$ is a group homomorphism.\label{thm:conj-unique}

\noindent (ii) If $q \in G$ is a fixed number we can decide whether or not there is a
$g \in \PL_{S,G}(J)$ with initial slope $g'(\eta^+)=q$ such that $y^g=z$. If $g$ exists, it is unique.
In other words, if there is a $g$ such that $\varphi_{y,z,\tau}(g)=\mu \in G$
then $g$ is unique and can be constructed.}

\noindent \emph{Proof.} Immediate generalization of Corollary \ref{con:unique}. $\square$

\corollary{Let $y,z \in \PL^<_{S,G}(J)$ and $g \in \PL_+(J)$ such that $y^g=z$ and $g'(\eta) \in G$. Then
$g \in \PL_{S,G}(J)$.}


\subsection{Centralizers and Roots in $\PL_{S,G}(J)$ \label{sec:centralizers}}

This section proves a generalization of Proposition \ref{thm:centralizers-F}. The centralizers $C_{\PL_{S,G}(J)}(z)$
of elements will be direct products of copies of $\mathbb{Z}$'s and of $\PL_{S,G}(U)$'s,
for some suitable intervals $U$. In order to prove
this, we will use the Stair Algorithm to build a ``section'' of the map $\varphi_x$. As in the proof
of Proposition \ref{thm:centralizers-F}, we will reduce the study to functions in $\PL_{S,G}^0(J)$.
Consider the conjugacy problem with $y=z$ and let
$\partial_S D(z)=\{0=\alpha_0  < \alpha_1 < \ldots < \alpha_s < \alpha_{s+1}=1\}$. Since
all the points of $\partial_S D(z)$ are fixed by $z$, then $g \in C_{\PL_{S,G}(I)}(z)$
must fix the set $\partial_S D(z)$ and thus each of the $\alpha_i$'s.
This implies that we can restrict to solving the conjugacy problem
in each of the subgroups $\PL_{S,G}([\alpha_i,\alpha_{i+1}])=\PL_{S,G}^0([\alpha_i,\alpha_{i+1}])$.
If $z=1$, it is immediate that $C_{\PL_{S,G}(J)}(x)=\PL_{S,G}(J)$, so
now we can focus on $1 \ne z \in \PL_{S,G}^0(J)$. Consider $\mathbb{R}_+$ to
be the multiplicative group of positive real numbers.
Let $A \subset \mathbb{R}_+$ be the set of all possible
initial slopes of functions $g$ such that $g^{-1} z g = z$. The set $A$ is not empty,
since $\langle z \rangle \subseteq C_{\PL_{S,G}(J)}(z)$. For a given closed interval $J$ with endpoints in $S$
we define a map
\[
\psi: A \to C_{\PL_{S,G}(J)}(z)
\]
\[
\alpha \mapsto g_\alpha
\]
which sends an initial slope $\alpha$ to its associated conjugating function $g_\alpha$.
By the uniqueness
of a conjugator with a given initial slope, we notice immediately that
$g_\alpha \circ g_\beta = g_{\alpha \cdot \beta}$ and so $A$ is a subgroup of $\mathbb{R}_+$ and $\psi$ is
a homomorphism of groups. Moreover, the uniqueness of the conjugator implies also that $\psi$ is an isomorphism.
The main result of this section is the following:

\theorem{Let $J \subseteq [0,1]$ be a closed interval with endpoints in $S$ and let
$id \ne z \in \PL_{S,G}^0(J)$. Then $C_{\PL_{S,G}(J)}(z)$ is isomorphic with $\mathbb{Z}$.
\label{thm:new-centralizer}}

\noindent \emph{Proof of Theorem \ref{thm:new-centralizer}.} By the discussion above we have that
the group $A=\{g'(\eta^+) \mid g \in C_{\PL_{S,G}(J)}(z)\}$ is isomorphic with
$C_{\PL_{S,G}(J)}(z)$. We start by assuming that
$z \in \PL_{S,G}^<(J)$ and we want to prove that $A$ is discrete. We assume, by contradiction
that $A$ is not discrete.

\noindent \emph{Step 1:} If $A$ is not discrete, then $A$ is dense in $\mathbb{R}_+$.

\noindent \emph{Proof.} This is a standard well known fact (for example see \cite{robertalain}). $\square$

\noindent \emph{Step 2:} Let $[\eta,\alpha]^2$ be the first initial linearity box and
$[\beta,\tau]^2$ be the first final linearity box, for some $\tau \le \zeta$ fixed point for $z$.
Without loss of generality, we can assume that the restriction $z|_{[\eta,\tau]} \in \PL_+^<([\eta,\tau])$.
Let $r$ be a positive integer big enough so that $z^{-r}(\alpha)>\beta$.
Then $z^r$ is linear on $[\beta,z^{-r}(\alpha)]$, say with slope $b$.

\noindent \emph{Proof.} Since $A$ is dense in $\mathbb{R}_+$, we can pick a $c \in C_{\PL_{S,G}(J)}(z)$ such that
$c'(\eta^+) < 1$ is arbitrarily close to $1$. Now, observe that $c \in \PL_+^<([\eta,\tau])$ and
look at the two hand sides of $cz^r=z^rc$, by restricting this equality
to the interval $[\beta,z^{-r}(\alpha)]$. Suppose $\{\mu_1< \ldots < \mu_k \}$ are the breakpoints of
$z^r$ on $[\beta,z^{-r}(\alpha)]$, hence they are also the breakpoints of $cz^r$ on $[\beta,z^{-r}(\alpha)]$,
since $c$ is linear on $[\eta,\alpha]$. On the interval $[\beta,\tau]$ we can write $c^{-1}(t)=\lambda(t-1)+1$,
where $\lambda=c'(\tau^-)$:
if we have chosen $c'(\eta^+) \ne 1$ to be close enough to $1$, then $\lambda<1$ is also arbitrarily close to $1$.
Since $c^{-1}$ is linear on $[\beta,\tau]$ then, if we choose $\lambda$ close enough to 1,
the set of breakpoints of $z^rc$ on $[\beta,z^{-r}(\alpha)]$ will be $c^{-1}(\{\mu_1,\ldots,\mu_k \})=
\{\lambda(\mu_1-1)+1,\ldots,\lambda(\mu_k-1)+1 \}$. As $cz^r=z^rc$ on $[\beta,z^{-r}(\alpha)]$
we must have that $\{\mu_1,\ldots,\mu_k \}=c^{-1}(\{\mu_1,\ldots,\mu_k \})$ and so $\lambda=1$,
which is a contradiction. $\square$

\noindent \emph{Step 3:} Define $a=\frac{d}{dt}z^{r}(t)\Big\vert_{t=\eta^+}<1$ to be the initial slope of $z^r$.
For every positive integer $i$, the map $z^r$ is linear on $[z^{-ir}(\beta),z^{-(i+1)r}(\alpha)]$
with slope $a$.

\noindent \emph{Proof.} We assume by induction that the result is true for any integer less than $i$.
Consider now the map $z^{(i+1)r}$ and apply the chain rule on two intervals, recalling that $\frac{d}{dt}z^r(t)=a$
on the intervals $[z^{-jr}(\beta),z^{-(j+1)r}(\alpha)]$ for any $j < i$:
\begin{center}
$$
\begin{array}{cc}
\frac{d}{dt} z^{(i+1)r}(t)= a^ib  & \qquad t \in [\beta,z^{-ir}(\alpha)] \\
{} \\
\frac{d}{dt} z^{(i+1)r}(t)= a^{i-1}b \frac{d}{dt} z^{r}(t) & \qquad t \in [z^{-ir}(\beta),z^{-(i+1)r}(\alpha)].
\end{array}
$$
\end{center}
\noindent We apply Step 2 using the positive integer $(i+1)r$, hence we have that $z^{(i+1)r}$
must be linear on $[\beta,z^{-(i+1)r}(\alpha)]$ and we can equate the two derivatives computed above to get
$a^i b=a^{i-1}b \frac{d}{dt} z^{r}(t)$ on the interval $[z^{-ir}(\beta),z^{-(i+1)r}(\alpha)]$. We simplify
both sides and get the thesis of the Claim. $\square$

\noindent By sending $i \to \infty$ in Claim 2 we see that the slope of $z^r$ around $\tau^-$ must be equal to $a<1$.
However, since the restriction $z^r|_{[\eta,\tau]} \in \PL_+^<([\eta,\tau])$, we must have that
$\frac{d}{dt}z^r(t)\Big\vert_{t=\tau^-}>1$, which is a contradiction.
Therefore $A$ is a discrete subgroup of $\mathbb{R}_+$ and so it is isomorphic with $\mathbb{Z}$. $\square$

\theorem{Let $J = [\eta,\zeta] \subseteq [0,1]$ be a closed interval with endpoints in $S$
and $z \in \PL_{S,G}(J)$, then:

\noindent (1) $C_{\PL_{S,G}(I)}(z)$ is isomorphic with a direct product of copies of $\mathbb{Z}$'s
and $\PL_2(J_i)$'s for some suitable intervals $J_i \subseteq I$.

\noindent(2) For every positive integer $n$ we can decide whether or not $\sqrt[n]{z}$ exists.
\label{thm:centralizers-PLSG}}

\noindent \emph{Proof}. The proofs of (1) and (2) follow from the proofs of Propositions
\ref{thm:compute-roots} and \ref{thm:centralizers-F}
by replacing every occurrence of $\partial_2$ with $\partial_S$ and by applying
the previous corollary to get the centralizers of elements in $\PL_{S,G}^0(J)$. Moreover,
to prove (2) we need to observe that, in order to start the procedure, we need to verify
whether or not $\sqrt[n]{z'(\eta^+)} \in S$. $\square$

\noindent The following is an immediate generalization of Proposition \ref{intersectionofcetralizers}:

\propositionname{Intersection of Centralizers}
{Let $J = [\eta,\zeta] \subseteq [0,1]$ be a closed interval with endpoints in $S$, let
$z_1, \ldots, z_k \in \PL_{S,G}(J)$ and define the subgroup
$C:=C_{\PL_{S,G}(I)}(z_1) \cap \ldots \cap C_{\PL_{S,G}(I)}(z_k)$.
If the interval $J$ is divided by the points in the union
$\partial_S D(z_1) \cup \dots \cup \partial_S D(z_k)$
into the intervals $J_i$ then
$$
C = C_{J_1} \cdot C_{J_2} \cdot \ldots \cdot C_{J_r},
$$
where $C_{J_i}:= \{f \in C \mid f(t)=t, \forall t \not \in J_i\} = C \cap \PL_{S,G}(J_i)$.
Moreover, each
$C_{J_i}$ is isomorphic to either $\mathbb{Z}$, or $\PL_{S,G}(J_i)$ or
the trivial group. $\square$ \label{thm:intersectionofcentralizers}}

\corollary{Let $J = [\eta,\zeta] \subseteq [0,1]$ be a closed interval with endpoints in $S$
and $y,z\in \PL_{S,G}^0(J)$. Then $C_{\PL_{S,G}(J)}(y,z)$ is either empty or countable. \label{thm:estimate}}

\noindent \emph{Proof.} Suppose that the set $C_{\PL_{S,G}(J)}(y,z)$ is not empty, then we have that
$C_{\PL_{S,G}(J)}(y,z)= g_0 \cdot C_{\PL_{S,G}(J)}(y)$ for a suitable $g_0 \in \PL_{S,G}(J)$. Thus
$\# C_{\PL_{S,G}(J)}(y,z) = \# C_{\PL_{S,G}(J)}(y) = \aleph_0$ which is countable by Theorem
\ref{thm:new-centralizer}. $\square$

\noindent In order to solve the conjugacy problem in $\PL_{S,G}(I)$, we need to check
whether or not there are candidate conjugators in a given interval of initial slopes.

\lemma{Let $J = [\eta, \zeta]$ be a closed interval with endpoints in $S$
and let $W = [w , 1] \cap G$ for some number $w \in \mathbb{R}$.
If $y,z \in \PL^0_{S,G}(J)$, then the set
\[
\{g'(\eta^+) \mid g \in C_{\PL_{S,G}(J)}(y,z) \} \cap W
\]
is contained in a finite set $V$ that can be constructed directly. \label{thm:possible-interval-initial-slopes-PLSG}}

\noindent \emph{Proof.} We will use the notation of Theorem \ref{thm:new-centralizer}.
Since the argument of this proof will be based on the Stair Algorithm, which works
in $\PL_+(J)$, we can restrict our attention on the interval between $\eta$ and the first
fixed point of $z$. Hence, we can assume $y,z \in \PL_+^<(J)$ without loss of generality.
We choose a positive integer $r$ following the proof of Proposition \ref{thm:stair-algorithm}:
that is, we choose the smallest integer $r$ such that
\[
\min\{z^{-r}(\alpha),y^{-r}(\eta + w(\alpha-\eta))\} > \beta.
\]
using the lowest possible initial number $w$.
Using the explicit conjugator formula for an initial slope $q \in W$ (see Corollary \ref{thm:explicit-conjugator}), we
know that the candidate conjugator has the shape $g_q:=y^{-r} g_{0,q} z^r$ on the interval
$[\eta,z^{-r}(\alpha)]$ for a suitable map $g_{0,q}$ that has initial slope $q \in W$. Our choice of $r$
guarantees that, for any $q \in W$, the map $g_q$ lies inside the final linearity
box at the point $z^{-r}(\alpha)$.

\noindent \emph{Claim:} Choose an integer $i$ such that $z^{-ir}(\beta)>z^{-r}(\alpha)$. Then
$z^r$ must have a breakpoint $p \in [z^{-ir}(\beta),z^{-(i+1)r}(\alpha)]$.

\noindent \emph{Proof of the Claim.} Let $a=\frac{d}{dt}z^{r}(t)\Big\vert_{t=\eta^+}<1$.
If $z^r$ were linear on $[z^{-ir}(\beta),z^{-(i+1)r}(\alpha)]$ then, by Step 3 of Theorem \ref{thm:centralizers-PLSG},
we would have that $z^r$ is linear on every interval $[z^{-kir}(\beta),z^{-k(i+1)r}(\alpha)]$
with slope $a$ for every positive integer $k \ge 2$. Arguing as in the conclusion
of Theorem \ref{thm:centralizers-PLSG}, this would imply that $\frac{d}{dt}z^{r}(t)\Big\vert_{t=\zeta^-}=a<1$
which is a contradiction. $\square$

\noindent By construction, the map $g_{0,q}$ can be built to be linear on the interval $[\eta,z^{-(i+1)r}(\alpha)]$.
We observe that $z^r$ has a breakpoint at $p$, hence $g_0 z^r$ must have a breakpoint at $p$.
Now, for the map $y^{-r} g_{0,q} z^r$ to be a candidate conjugator, it must be linear around the point $p$, thus
the breakpoints of $g_{0,q}z^r$ on the interval $[z^{-ir}(\beta),z^{-(i+1)r}(\alpha)]$ must be
canceled by the set $\{c_1,\ldots,c_v\}$ of all the breakpoints of $y^{-r}$ on $[\eta,\zeta]$, thus the image
of $p$ under $g_{0,q}z^r$ must go to a breakpoint of $y^{-r}$.
Since $g_{0,q}z^r(p)= q (z^r(p)-\eta)+\eta \in \{c_1,\ldots,c_v\}$, then there are only finitely
many choices for $q \in W$. $\square$

\remark{Since the finite set $V$ of Lemma \ref{thm:possible-interval-initial-slopes-PLSG}
can be computed directly, we can run the stair algorithm on all elements of $V$ as possible
initial slopes and thus find all possible conjugators with slopes in $[w,1] \cap G$.}

\section{Simultaneous Conjugacy Problem in $\PL_{S,G}(I)$ \label{sec:simultaneous}}

In this section we wrap up all the arguments of the Chapter to
solve the $k$-simultaneous conjugacy problem. We will first deal with the case $k=1$ and then
with the general case. Unlike the approach adopted in the first part of the paper for the case of $F$,
in this second part we have first
solved the conjugacy problem in the special case of $y=z$ before approaching the conjugacy problem for
two elements $y,z \in \PL_{S,G}(I)$.

\subsection{The Ordinary Conjugacy Problem for $\PL_{S,G}(I)$ \label{ssec:conjugacy-problem}}

\theorem{The conjugacy problem in $\PL_{S,G}(I)$ is solvable.}

\noindent \emph{Proof.} Let $y,z \in \PL_{S,G}(I)$, $y \ne z$.
We use Proposition \ref{thm:same-intersection} and suppose that
$\partial_S D(y) = \partial_S D(z) = \{ 0 = \alpha_0  < \alpha_1 < \ldots < \alpha_r < \alpha_{r+1} = 1\}$.
Now we restrict to an interval $J_i=[\alpha_i,\alpha_{i+1}]$. For simplicity, we still call $y|_{J_i}$ and
$z|_{J_i}$, $y$ and $z$. In order for $y$ and $z$ to be conjugate, we must have $y'(\alpha_i^+)=z'(\alpha_i^+)$
and $y'(\alpha_{i+1}^-)=z'(\alpha_{i+1}^-)$.
Up to taking inverses of $y$ and $z$, we can assume that $q= y'(\alpha_i^+)=z'(\alpha_i^+) < 1$.
Now observe that $g^{-1}yg=z$ if and only if $(y^v g)^{-1}y (y^v g)=z$ for every
$v \in \mathbb{Z}$. If $q^{\rho(g)}$ is the initial slope of $g$, then $q^{v + \rho(g)}$ is the initial
slope of $y^v g$. Thus, up to taking powers of $y$ we can assume
that the exponent of the initial slope of $g$ is in $[q,1]$. By Lemma \ref{thm:possible-interval-initial-slopes-PLSG},
the set of possible initial slopes inside $[q,1] \cap G$ is finite and can be directly
constructed, so we can apply the Stair Algorithm on each of them and verify if any of the obtained
maps is a conjugator. All the other conjugators
are found by taking the products $y^v g$ with $v \in \mathbb{Z}$. $\square$

\subsection{The $k$-Simultaneous Conjugacy Problem in $\PL_{S,G}(I)$\label{ssec:simultan-conjug}}

The algorithm used to solve the $k$-simultaneous problem in the case of the group $F$
can be extended in full generality, except for one of its steps.

\theorem{The $k$-simultaneous conjugacy problem in $\PL_{S,G}(I)$ is solvable \label{thm:simultaneous-conjug-PLSG}.}

\noindent \emph{Proof.} To prove the solvability of the $k$-simultaneous conjugacy problem we can mimic completely
the proof used for Thompson's group $F$. We need to replace every occurrence of $\partial_2$ with $\partial_S$
and speak of elements of $S$ instead of dyadic rational numbers.
The only part in which we need refine the argument is in the case of Subsection \ref{ssec:specialcase}
in which we reduce to solve the equation
\begin{equation} \label{eq:simult}
x^m=g_0 \widehat{y}^n
\end{equation}
where $x,y,g_0 \in \PL_{S,G}([\eta,\zeta])$ are given and we are looking for $m,n \in \mathbb{Z}$
satisfying the previous equation. We define $q = g_0'(\eta^+) \in \mathbb{R_+}$ and so
$x'(\eta^+)=q^\alpha, y'(\eta^+)=q^\beta,g_0'(\eta^+)=q$ for some $\alpha, \beta \in \mathbb{R}$.
Notice that in Subsection \ref{ssec:specialcase} we have
$\alpha,\beta,\gamma \in \mathbb{Z}$, while here not all of them are integers. We must then have
\begin{equation} \label{eq:simult2}
q^\alpha=x'(\eta^+)^m=(g_0 \widehat{y}^n)'(\eta^+)=q^{1 +\beta n}
\end{equation}
and therefore we need to solve the equation
\begin{equation} \label{eq:simult3}
\alpha m = 1 + \beta n
\end{equation}
for some $m,n \in \mathbb{Z}$. We observe that if equation (\ref{eq:simult3}) is solvable, then
$\alpha$ is rational if and only if $\beta$ is rational. Thus, if either $\alpha$ or $\beta$ is a rational
number it is immediate to check whether there is a solution to (\ref{eq:simult3}).
If $\alpha$ and $\beta$ are both irrational, then equation (\ref{eq:simult3})
becomes a $\mathbb{Q}$-linearity dependence relation and, if it is solvable,
then the dimension of the vector space
generated by $\alpha,\beta$ and $1$ over $\mathbb{Q}$ is exactly 2.
By Remark \ref{thm:rewrite-requirement} and
the last of the requirements in Remark \ref{thm:ring-requirements} we are able to detect if this last statement is true or
not. In case it is true, then there is a unique solution to (\ref{eq:simult3}) and it is given by
the coordinates of $1$ in the basis $\alpha$ and $\beta$, thus
it is now trivial to check if there is a integer solution or not. In case there is a solution to
equation (\ref{eq:simult3}), we do not need to find a bound for $m,n \in \mathbb{Z}$ as for the case of Thompson's group $F$,
because there is at most one solution.
The remaining part of the algorithm follows as before. $\square$

\section{Interesting Examples \label{sec:examples}}

Now that we have developed the general theory, we are going to see a few interesting examples where
the simultaneous conjugacy problem is solvable.
We will not dwell too much on the details here, sketching only why it is possible to
verify the requirements.

\example{$S=\mathbb{Q}$ and $G=\mathbb{Q}^*=\mathbb{Q} \cap (0,\infty)$.}

\noindent Since $\mathbb{Q}$
is a field, $S / \mathcal{I}=\{0\}$ so all the requirements of Remark \ref{thm:ring-requirements}
are satisfied. To solve the simultaneous conjugacy problem, we need to solve equation (\ref{eq:simult2}),
which becomes
\[
\frac{a_1^m}{b_1^m}=\frac{c a_2^n}{d b_2^n}
\]
where we can assume that all numerators are coprime with the denominators. By equating prime factors
in the equation to be solved, we get a system of equations of the type
$\alpha_i m = \gamma_i + \beta_i n$, for $\alpha_i, \beta_i, \gamma_i \in \mathbb{Z}$. All of them can be solved
in the same fashion as in Lemma \ref{thm:special-case} and we can reduce equation (\ref{eq:simult})
to the equation $X^k=G_0 Y^k$ and solve it as in Subsection \ref{ssec:specialcase}.

\example{$S=\mathbb{Z}\big[\frac{1}{n_1}, \ldots, \frac{1}{n_k}\big]$
and $G=\langle n_1, \ldots, n_k \rangle$ for $n_1,\ldots,n_k \in \mathbb{Z}$.\label{thm:example-rationals}}

\noindent We observe that $S=\mathbb{Z}\big[\frac{1}{n_1 \ldots n_k}\big]$ and
it can be shown that, if $r:=n_1 \ldots n_k$, then $S / \mathcal{I} \cong \mathbb{Z}/r\mathbb{Z}$
as rings and therefore the requirements of Remark \ref{thm:ring-requirements} are also satisfied. Equation
(\ref{eq:simult2}) can be treated as in the previous example.
For $k=1$, we recall that the groups $\PL_{S,G}(I)$ are known as \emph{generalized Thompson's groups}.

\example{$S=\mathbb{Z}\big[\frac{1}{n_1}, \ldots, \frac{1}{n_k}, \ldots \big]$
with $G=\langle \{n_i\}_{i \in \mathbb{N}} \rangle$ for a sequence $\{n_i\}_{i \in \mathbb{N}} \subseteq \mathbb{Z}$.}

\noindent This example is easily reducible to the previous one, since if we are given a finite set
$E$ of elements in $\PL_{S,G}(I)$ we can consider the set $\{n_{i_1}^{\alpha_{i_1}},\ldots,
n_{i_v}^{\alpha_{i_v}}\}$ of all slopes of elements of $E$. Then $E \subseteq \PL_{S',G'}(I)$
where $S':=\mathbb{Z}\big[\frac{1}{n_{i_1}}, \ldots, \frac{1}{n_{i_v}}\big]$ and
$G':=\langle n_{i_1}, \ldots, n_{i_v} \rangle$.

\example{$S$ finite algebraic extension over $\mathbb{Q}$
and $G=S^*:= S \cap (0,\infty)$ }

\noindent As with the first example, since $S$ is a finite algebraic extension it is not
difficult to verify that all the requirements of Remark \ref{thm:ring-requirements} are
satisfied.


\example{$S=\mathbb{R}$ and $G=\mathbb{R}_+$}

\noindent In order to verify
the requirements for this case, we need to discuss exactly what we
mean by real number and how we implement it in a machine.
In most cases, we work with numbers which are expressed as roots of polynomials in some subfields of $\mathbb{R}$ and we
are able to give a complete answer and the same is true for all the requirements of Remark
\ref{thm:ring-requirements}.

\chapter{Cryptanalysis of the Shpilrain-Ushakov protocol for Thompson's group}
\label{chapter5}

\section{Introduction}

Public Key Cryptography is involved in the exchange of information
between two parties $A$ and $B$, often labeled as ``Alice'' and ``Bob''.
Before they start sending data to the other party, they must agree on a way
to send it. The type of encryption they use is called ``public key'' because part of information they need to agree on (in
this context it is usually a group and some of its subgroups or elements) and the encryption scheme are in the public domain.
Alice and Bob each choose secretly some information, respectively $i_A$ and $i_B$. They both use their secret information to encrypt
some public data $w$, and send it to the other party. Alice receives the encrypted information $e(i_B,w)$ and
she encrypts it using her own information, obtaining some data $e(i_A,e(i_B,w))$. Similarly, Bob receives
$e(i_A,w)$ and encrypts it using $i_B$ to get $e(i_B,e(i_A,w))$.
The public protocol is usually chosen so that after this procedure they obtain the same information, that is
$e(i_A,e(i_B,w))=e(i_B,e(i_A,w))$.
This common information is now referred to as the \emph{shared secret key}
and it is now used to exchange messages between the two parties.
This ``commutativity'' of encryption comes from
generalizing the Diffie-Hellman cryptosystem based on the infinite cyclic group
(see \cite{cryptonotes} for details).

A third party $E$ (``Eve'') is listening and detecting anything the two parties are exchanging.
Thus Eve captures $e(i_A,w)$ and $e(i_B,w)$ and any message encrypted using the shared secret key. To break
the protocol Eve must try to extract $i_A$ and $i_B$ or equivalent elements.
This discussion
will be made precise in the next sections, by describing precisely the Shpilrain-Ushakov public key cryptography
protocol based on Thompson's group $F$ and how the information gets transmitted between the two parties.

The Chapter is organized as follows. In Section \ref{sec:protocol} and Section \ref{sec:platform}
we recall the protocol. In Section 4 we
recall the choice of parameters proposed in \cite{su}.
In section \ref{sec:recovering-keys} we give an efficient attack that always recovers
the secret key. In Sections \ref{sec:transitivity} and
\ref{sec:recovering-keys-alternatives} we show another type of attack.
In Section \ref{sec:restriction-centralizers} we make some comments on possible generalizations of this protocol.
The material of this Chapter is going to appear in the \emph{Journal of Cryptology} \cite{matucci1}.

\subsection*{History and related works.} The first attack on this protocol
was announced by Ruinskiy, Shamir and Tsaban in November 2005 at the Bochum Workshop
\emph{Algebraic Methods in Cryptography}, showing that the parameters given in \cite{su}
should be increased to have higher security of the system. Their attack was improved in other
announcements and was finalized in \cite{rst2} at the same time that the material of this Chapter was written.
Their attack describes a more general procedure which uses length functions.
We remark that the same authors have been developing new techniques involving
``subgroup distance functions'' and that they applied them to the same protocol for $F$ as a test case \cite{rst3}.
The approach of Ruinskiy, Shamir and Tsaban in their papers is heuristic,
and its success rates are good but not $100\%$. Our approach is deterministic, and provably succeeds
in all possible cases.

\section{The Protocol \label{sec:protocol}}

The protocol proposed in \cite{su} is based on the \emph{decomposition problem}:
given a group $G$, a subset $X \subseteq G$ and $w_1,w_2 \in G$, find $a,b \in X$
with $aw_1b=w_2$, given that such $a,b$ exist. Here is
the protocol in detail:

\subsection*{Public Data.} A finitely presented group $G$, an element $w \in G$ and two subgroups $A,B$ of $G$
such that $ab=ba$ for all $a \in A$, $b \in B$.

\subsection*{Private Keys and Communication.} Alice chooses $a_1 \in A$, $b_1 \in B$ and sends the element
$u_1=a_1 w b_1$ to Bob. Bob chooses $b_2 \in B$, $a_2 \in A$ and sends the element $u_2=b_2wa_2$ to Alice.
Alice then computes the element $K_A=a_1 u_2 b_1=a_1 b_2 w a_2 b_1$ and Bob computes the element
$K_B=b_2 u_1 a_2=b_2 a_1 w b_1 a_2$. Since $A$ and $B$ commute elementwise, $K=K_A=K_B$ becomes
Alice and Bob's shared secret key to send one bit. Alice and Bob need to generate and compute
a shared secret key for each bit they want to send.

To communicate bits, the two parties send elements $x \in G$.
If Alice wants to send a $1$, she uses the relations of $G$ to rewrite the word representing $K$,
``scrambling'' the way it appears, and sends it.
If she wants to send a $0$, she chooses a random element $x \in G$ and sends it; with overwhelming probability she will
pick an element different from $K$. Now Bob solves the word problem for
$x K^{-1}$, to identify whether he has received a $0$ or a $1$.
Hence, it is important for the word problem to be efficient in the group.

\subsection*{Eavesdropper's Data and Brute force Attack.} Eve has all the public data and the two elements $u_1$ and
$u_2$, observed during Alice and Bob's exchange.

\noindent \textbf{To break the protocol} To recover $a_1$, Eve needs to find a pair $(\overline{a},\overline{b})$ such that
\[
u_1 = \overline{a} w \overline{b}
\]
Thus Eve computes
\[
\overline{a} u_2 \overline{b} = \overline{a} b_2 w a_2 \overline{b} = b_2 \overline{a} w \overline{b} a_2 = b_2 u_1 a_2 = K.
\]
Eve can always use what is called a \emph{brute force attack}, that is, try all the elements of $A$ to get candidates for the
shared secret key to test on the exchanged message.
However, as the groups $A$ and $B$ are usually chosen to be infinite, this is a cumbersome and slow way
to look for a suitable pair $(\overline{a},\overline{b})$ and she has to look for something more efficient.

\section{The Subgroups $A_s,B_s$ \label{sec:platform}}

\begin{figure}[0.5\textwidth]
 \begin{center}
  \includegraphics[height=6cm]{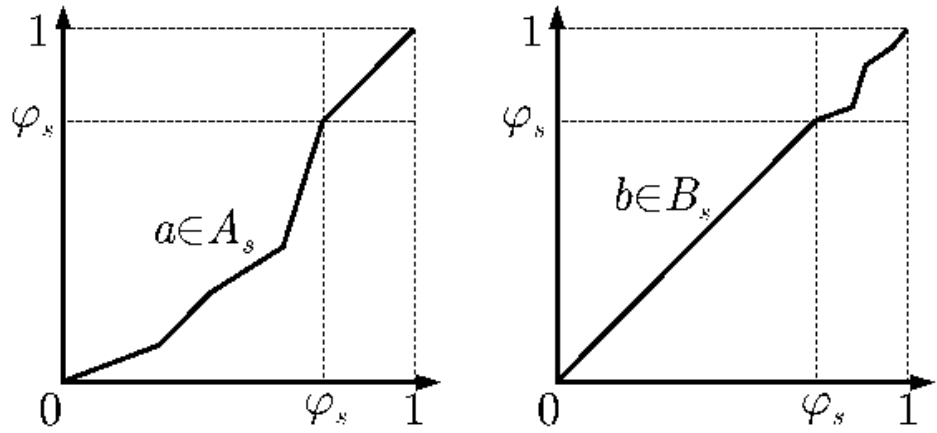}
 \end{center}
 \caption{An example of an element of $A_s$ and one of $B_s$.}
 \label{fig:commuting-elements}
\end{figure}

Now we apply the protocol described in the previous Section to the special case of $G=F$
Thompson's group with the standard generating set defined in Chapter \ref{chapter1}.
We introduce a notation which will be useful for the definition of the subgroups $A$ and $B$.
For every positive integer $k$ we call
\[
\varphi_k:=1-\frac{1}{2^{k+1}}.
\]
From the definition of $x_k$, we get
\[
x_k^{-1}\left(\left[\varphi_k,1\right]\right)=\left[\varphi_{k+1},1\right] \subseteq
\left[\frac{3}{4},1\right]
\]
implying that, for $t \in [\varphi_k,1]$, we have
\[
\frac{d}{dt}x_0x_k^{-1}(t)=x_0'(x_k^{-1}(t))(x_k^{-1})'(t)=2 \cdot \frac{1}{2}=1
\]
which means $x_0 x_k^{-1}$ is the identity in the interval $[\varphi_k,1]$.
For any $s \in \mathbb{N}$, Shpilrain and Ushakov define in \cite{su} the following sets
\[
S_{A_s}=\{x_0 x_1^{-1}, \ldots, x_0 x_s^{-1}\}
\]
\noindent and
\[
S_{B_s}=\{x_{s+1},\ldots,x_{2s}\}
\]
and then define the subgroups $A_s:=\langle S_{A_s} \rangle$ and $B_s:=\langle S_{B_s} \rangle$.
The previous argument immediately yields that all elements of $A_s$ commute with
all elements of $B_s$ (see figure \ref{fig:commuting-elements}), i.e.

\lemmaname{Shpilrain-Ushakov \cite{su}}{For every fixed $s \in \mathbb{N}$, $ab=ba$ for every
elements $a \in A_s$ and $b \in B_s$.}

\convention{For this Chapter only,
for every dyadic number $d \in [0,1]$ we denote by $\PL_2([0,d])$ the set of functions in $\PL_2(I)$
which are the identity on $[d,1]$. Moreover, if we are given a
piecewise linear map defined only on $[0,d]$ we will assume it is extended to $[0,1]$ by defining it
as the identity on $[d,1]$. Similar remarks apply to $\PL_2([d,1])$.}

\noindent Parts (i) and (iii) of the following Lemma are in \cite{su}, while part (ii) is a simple observation.

\lemma{(i) $A_s$ is the set of elements whose normal form is of the type
\[
x_{i_1}\ldots x_{i_m} x_{j_m}^{-1} \ldots x_{j_1}^{-1}
\]
where $i_k - k < s$ and $j_k - k < s$, for all $k=1,\ldots,m$.

\medskip
\noindent (ii) $B_s=\PL_2([\varphi_s,1])$.

\medskip
\noindent (iii) Let $a \in A_s$ and $b \in B_s$ be such that their normal forms are
\begin{eqnarray*}
a = x_{i_1}\ldots x_{i_m} x_{j_m}^{-1} \ldots x_{j_1}^{-1} \\
b = x_{c_1}\ldots x_{c_u} x_{d_v}^{-1} \ldots x_{d_1}^{-1}.
\end{eqnarray*}
Then the normal form of $ab$ is
\[
ab = x_{i_1}\ldots x_{i_m} x_{c_1+m}\ldots x_{c_u+m} x_{d_v+m}^{-1} \ldots x_{d_1+m}^{-1}
x_{j_m}^{-1} \ldots x_{j_1}^{-1}.
\]\label{thm:definition-A_s}}

\theoremname{Shpilrain-Ushakov \cite{su}}{In Thompson's group $F$, the normal form of a given word $w$
can be computed in time $O(|w| \log |w|)$,
where $|w|$ is the length of the normal form in the generators $x_0,x_1, x_2, \ldots$ \label{thm:normal-forms}}

\section{Suggested Parameters for the Encryption \label{sec:parameters}}

We now illustrate briefly the choice of parameters proposed in \cite{su}.
Alice and Bob select an integer $s \in [3,8]$ and an even integer
$M \in [256,320]$ uniformly and randomly.
Moreover, they also choose a random element $w \in \langle x_0,x_1, \ldots, x_{s+2} \rangle$ with $|w|=M$,
where $|w|$ is as in Theorem \ref{thm:normal-forms}.
The numbers $s,M$ and the element $w$ are now part of the public data.

To proceed with the protocol described in Section \ref{sec:protocol},
Alice chooses random elements $a_1 \in A_s,b_1 \in B_s$, with $|a_1|=|b_1|=M$, while Bob chooses
random elements $a_2 \in A_s,b_2 \in B_s$, with $|a_2|=|b_2|=M$.
Now they both compute the shared secret key:
\[
K=a_1 b_2 w a_2 b_1.
\]
Shpilrain and Ushakov remark that this choice of parameters gives a key space which
increases exponentially in $M$, i.e., $|A_s(M)| \ge \sqrt{2}^M$, thereby making it difficult for Eve to perform a brute force attack.

\section{Recovering the Shared Secret Key \label{sec:recovering-keys}}

We begin this section by providing the theoretical background
for the attack. We will use the piecewise-linear point of view to understand why the attack works
and then rephrase it combinatorially.
We will now describe how Eve, by knowing the elements $w,u_1,u_2$, can always
recover one of the two legitimate parties' private keys. She chooses whose key to crack, depending on whether the graph of $w$
is above or below the point $(\varphi_s,\varphi_s)$.

\subsection{Recovering Bob's Private Keys: $w(\varphi_s) \le \varphi_s$}

Since $w(t)\le \varphi_s$ for all $t \in [0,\varphi_s]$, we observe the following identity
\[
u_2(t)=b_2 w a_2(t)=w a_2(t), \, \, \forall t \in [0,\varphi_s].
\]
Therefore, Eve may apply $w^{-1}$ to the left of both sides of the previous equation to obtain
\[
w^{-1} u_2(t)=a_2(t), \, \, \forall t \in [0,\varphi_s]
\]
and so $w^{-1}u_2 \in A_s B_s$ and
\[
a_2(t)=\begin{cases}
w^{-1} u_2(t) & t \in [0,\varphi_s] \\
t & t \in [\varphi_s,1].
\end{cases}
\]
Now Eve has the elements $a_2$, $w$ and $u_2=b_2wa_2$ and she computes
\[
b_2=u_2 a_2^{-1}w^{-1}
\]
thereby detecting Bob's private keys and the shared secret key $K$.

\subsection{Recovering Alice's Private Key: $w(\varphi_s)>\varphi_s$}

Since $w^{-1}(t) < \varphi_s$ for all $t \in [0,\varphi_s]$, we have
\[
u_1^{-1}(t)=b_1^{-1} w^{-1} a_1^{-1}(t)= w^{-1} a_1^{-1}(t)  ,\, \, \forall t \in [0,\varphi_s].
\]
By applying the same technique as in the previous subsection Eve recovers $a_1^{-1}$ and obtains that $u_1w^{-1} \in A_s B_s$. Thus,
she is able to detect $a_1,b_1$ and the shared secret key $K$. Alternatively, Eve observes
\[
w^{-1}u_1(t)=w^{-1}a_1 w b_1(t)= b_1(t)  ,\, \, \forall t \in [\varphi_s,1]
\]
and so
\[
b_1(t)=\begin{cases}
t & t \in [0,\varphi_s] \\
w^{-1} u_1(t) & t \in [\varphi_s,1].
\end{cases}
\]

\subsection{Outline of the attack} We expand on the previous discussion to describe a combinatorial attack.
Assume that Eve has the elements $w,u_1,u_2$.

\begin{enumerate}
\item{Eve writes the normal forms of $z_1:=u_1 w^{-1}$ and $z_2:=w^{-1}u_2$.}
\item{By the previous discussion, either $z_1 \in A_s B_s$ or $z_2 \in A_s B_s$ (or both).
She can detect which one using Lemma \ref{thm:definition-A_s}(i) and selects this $z_i$.}
\item{She  computes the $A_s$-part $a_{z_i}$ of $z_i$.}
\item{If $i=1$, she computes $b_{z_1}:= w^{-1} a_{z_1}^{-1} u_1$. If $i=2$,
she computes $b_{z_2} := u_2 a_{z_2}^{-1}w^{-1}$.}
\item{Eve computes $K$ from $u_1,u_2,a_{z_{i}},b_{z_i}$.}
\end{enumerate}

The only point of this procedure which needs further explanation is (2). When we have
the normal forms of $z_1,z_2$, we know that one of them is in $A_s B_s$.
We write the normal form $z_i = x_{i_1}\ldots x_{i_e} x_{j_f}^{-1} \ldots x_{j_1}^{-1}$
and we look at the notation of Lemma \ref{thm:definition-A_s}(i): we need to find
the smallest index $r$ in $z_i$ such that either $i_{r+1}$ or $j_{r+1}$
does not satisfy the index condition in Lemma \ref{thm:definition-A_s}(i). To verify if $z_i \in A_s B_s$,
we need to check whether it has the form described in Lemma \ref{thm:definition-A_s}(iii): we remove the first
$r$ letters and the last $r$ letters of $z_i$ from the word and we lower all the indices of the remaining letters by $r$;
if what remains is a word whose indices are in $\{s+2, \ldots, 2s\}$, then we have an element of $B_s$, otherwise $z_i \not \in A_s B_s$.
If $z_i \in A_s B_s$, then $a_{z_i}$ will be
the product of the first $r$ elements of $z_i$ and the last $r$ ones.

\subsection{Complexity of the attack.} By Theorem \ref{thm:normal-forms}
we know that computing normal forms can be done in time $O(M \log M)$, where $M$
is the size of the inputs suggested in Section \ref{sec:parameters}. Part (2) of the attack can be executed in time $O(M)$,
by just reading the indices of the normal forms and finding when the relation of Lemma
\ref{thm:definition-A_s}(i) breaks down. Finally, the last steps are just
multiplications and then simplifications so they can again be performed in time
$O(M \log M)$. Therefore, Eve can recover the shared secret key in time $O(M \log M)$.

\remark{The previous discussion shows that there is no need to pass from words to piecewise-linear functions and back.
The attack can be performed entirely by using the combinatorial point of view which is used for encryption. The piecewise-linear
point of view is necessary only to prove that the combinatorial attack works. We also remark that the complexity of
the attack is independent of the parameter $s$.}

\section{Transitivity of $A_s$ and $B_s$ \label{sec:transitivity}}

The previous section showed how to recover the shared secret key of one of
the two involved parties, based on whether the graph of $w$ lies above or below
the point $(\varphi_s,\varphi_s)$. However, it is possible to find the shared
secret key even in the cases not studied in the previous section. More precisely,
it is possible to attack Alice's word in the case $w(\varphi_s) \le \varphi_s$
and Bob's word in the case $w(\varphi_s)>\varphi_s$. We need a better description of
the subgroups $A_s$. If $s=1$, we observe that $A_1=\langle x_0 x_1^{-1}\rangle$ is a cyclic group.
For larger values of $s$, $A_s$ becomes the full group of piecewise linear homeomorphism
on $[0,\varphi_s]$.

\theorem{$A_s=\PL_2([0,\varphi_s])$, for every $s \ge 2$. \label{thm:subgroups-are-F}}

\noindent In order to prove the Theorem we need the following two Lemmas.

\lemma{$x_0^{-1} \PL_2([0,\varphi_s])x_0=\PL_2([0,\varphi_{s+1}])$. Similarly,
$x_0^{-1}\PL_2([\varphi_s,1])x_0=\PL_2([\varphi_{s+1}])$. \label{thm:isomorphic-subgroups}}

\noindent \emph{Proof.} This result is a special case of Theorem \ref{thm:thompson-like},
but it is straightforward to verify it too. Observe that
$x_0^{-1}(\varphi_s)$ is fixed by $x_0^{-1}fx_0$ for every $f \in \PL_2([0,\varphi_s])$
and
\[
x_0^{-1}(\varphi_s)=\frac{1}{2}\left(1 - \frac{1}{2^{s+1}} - \frac{1}{2}\right) + \frac{3}{4} = \varphi_{s+1}
\]
therefore the result holds. The other result follows similarly. $\square$

\noindent The next corollary is also a special case of Theorem \ref{thm:thompson-like}.

\corollary{$\PL_2([0,\varphi_s]) \cong \PL_2([\varphi_s,1]) \cong F$, for every $s\ge0$. \label{thm:isomorphic-piecewise}}

\lemma{$A_2=\PL_2\left(\left[0,\frac{7}{8}\right]\right)$. \label{thm:first-isomorphism}}

\begin{figure}[0.5\textwidth]
 \begin{center}
  \includegraphics[height=6cm]{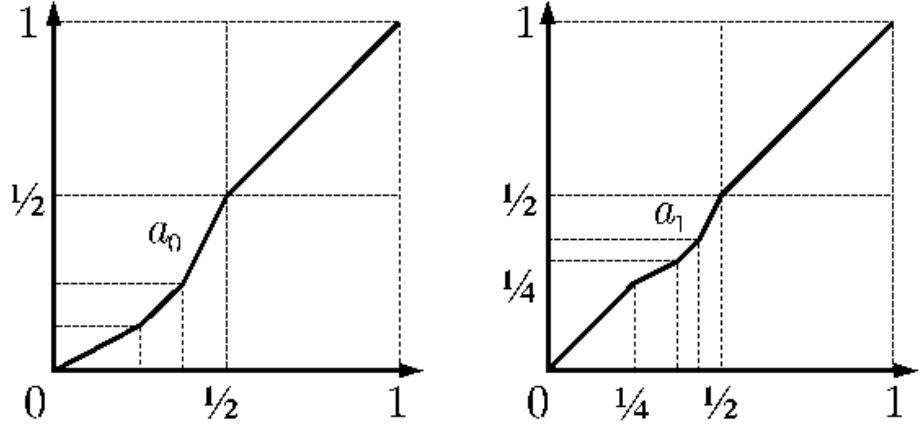}
 \end{center}
 \caption{The two standard generators for $\PL_2([0,\varphi_2])$.}
 \label{fig:generators-A_s}
\end{figure}

\noindent \emph{Proof.} Define $a=(\rho_2^*)^{-1}(x_0)$ and $b=(\rho_2^*)^{-1}(x_1)$ (see figure
\ref{fig:generators-A_s}). Then

\[
\begin{array}{cc}
a(t)=\begin{cases}
\frac{1}{2}t & t \in \left[0,\frac{1}{4}\right] \\
\left(t-\frac{1}{4}\right)+\frac{1}{8} & t \in \left[\frac{1}{4},\frac{3}{8}\right] \\
2\left(t - \frac{3}{8}\right) + \frac{1}{4} & t \in \left[\frac{3}{8},\frac{1}{2}\right] \\
t & t \in \left[\frac{1}{2},1\right]
\end{cases}

&

b(t)=\begin{cases}
t & t \in \left[0,\frac{1}{4}\right] \\
\frac{1}{2}\left(t- \frac{1}{4} \right) + \frac{1}{4} & t \in \left[\frac{1}{4},\frac{3}{8}\right] \\
\left(t-\frac{3}{8}\right)+\frac{5}{16} & t \in \left[\frac{3}{8},\frac{7}{16}\right] \\
2\left(t - \frac{7}{16}\right) + \frac{3}{8} & t \in \left[\frac{7}{16},\frac{1}{2}\right] \\
t & t \in \left[\frac{1}{2},1\right]
\end{cases}
\end{array}
\]

\noindent One sees that $a=x_0^2x_1^{-1}x_0^{-1}$
and that $b=x_0 x_1^2 x_2^{-1} x_1^{-1} x_0^{-1}$ (for example this can
be verified using tree diagrams for $F$, i.e. using Proposition 1.3.3 in \cite{belkthesis}
or by a direct computation). Since $\rho_2^*$ is an isomorphism,
$\PL_2([0,\frac{1}{2}])=\langle a,b\rangle$ and, by Lemma \ref{thm:isomorphic-subgroups},
$\PL_2([0,\frac{7}{8}])=\langle x_0^2 a x_0^{-2}, x_0^2 b x_0^{-2}\rangle$. By Lemma
\ref{thm:definition-A_s} we have
\begin{eqnarray*}
x_0^2 a x_0^{-2} = x_0^4 x_1^{-1}x_0^{-3} \in A_2 \\
x_0^2 b x_0^{-2} = x_0^3 x_1^2 x_2^{-1} x_1^{-1} x_0^{-3} \in A_2
\end{eqnarray*}
so that $\PL_2([0,\frac{7}{8}]) \subseteq A_2$. The other inclusion is obvious. $\square$

\medskip

\noindent \emph{Proof of Theorem \ref{thm:subgroups-are-F}}. By Lemma \ref{thm:first-isomorphism}
we have that $A_2=\PL_2([0,\varphi_2])$ and so, by applying Lemma \ref{thm:isomorphic-subgroups} and
the definition of $A_s$ we have

\[
\PL_2([0,\varphi_s])= x_0^{s-2} A_2 x_0^{2-s} \subseteq A_s \subseteq \PL_2([0,\varphi_s])
\]

\medskip

\noindent therefore implying that $A_s=\PL_2([0,\varphi_s])$. $\square$

\corollary{$A_s \cong B_s \cong F$, for every $s \ge 2$.}

\noindent The previous Theorem and Lemma \ref{thm:extension-partial-maps} in Chapter \ref{chapter4} yield the following corollaries:

\corollaryname{Transitivity of $A_s$}{For any $t_1,t_2 \in \mathbb{Z}\left[\frac{1}{2}\right] \cap [0,\varphi_s]$ we can construct an
$a \in A_s$ with $a(t_1)=t_2$. \label{thm:transitivity-A_s}}

\noindent \emph{Proof.} We appeal to the transitivity of $F$ on the dyadic numbers of $[0,1]$. We use
the map $\rho_s$ of Lemma \ref{thm:isomorphic-piecewise} and pull back $a_1,a_2$ to $[0,1]$ as
$\rho_s^{-1}(a_1),\rho_s^{-1}(a_2)$. We consider the two dyadic partitions $0<\rho_s^{-1}(a_1)<1$
and $0<\rho_s^{-1}(a_2)<1$ of the interval $[0,1]$ and apply Lemma \ref{thm:cfp-extension} to find
and $f \in F$ with $f(\rho_s^{-1}(a_1))=\rho_s^{-1}(a_2)$. Then the element $(\rho_s^*)^{-1}(f) \in A_s$
satisfies the requirements. $\square$

\corollaryname{Extendability of $A_s$}{Let $t_0 \in \mathbb{Z}\left[\frac{1}{2}\right] \cap [0,\varphi_s]$ and $\bar{a}(t)=a|_{[0,t_0]}$
for an element $a \in A_s$. Assume we know $\bar{a}$, but that we do not know $a$. Then we can construct
an $a_\sigma \in A_s$ such that $a_\sigma(t)=\bar{a}(t)$ for all $t \in [0,\varphi_s]$.
\label{thm:extension-A_s}}

\noindent \emph{Proof.} We observe that $\rho_s^*(\bar{a})$ is a piecewise linear map between
two intervals with dyadic extremes and contained in $[0,1]$. We apply Lemma \ref{thm:extension-partial-maps}
to extend $\rho_s^*(\bar{a})$ to a piecewise linear homeomorphism $f \in F$. Then we pull back $f$
to $A_s$ through $(\rho_s^*)^{-1}(f)$ which satisfies the requirements. $\square$

\remark{The analogues of the last two corollaries are true for the interval
$[\varphi_s,1]$ and $B_s$ too.}

\section{Using Transitivity to Attack the Shared Secret Key \label{sec:recovering-keys-alternatives}}

With the new description of $A_s$ and $B_s$ given in
section \ref{sec:transitivity}, it is now possible to attack the secret
keys in the cases left open from section \ref{sec:recovering-keys}.

\subsection{Attacking Alice's word for the case $w(\varphi_s) \le \varphi_s$}

We have
\begin{eqnarray*}
u_1(t)=a_1w(t), \forall t \in [0,\varphi_s],
\end{eqnarray*}
thus
\begin{eqnarray*}
a_1(t)=u_1w^{-1}(t), \forall t \in [0,w(\varphi_s)]
\end{eqnarray*}
and so $a_1$ is uniquely determined in $[0,w(\varphi_s)]$.
We apply corollary \ref{thm:extension-A_s} to find an element $a_\sigma \in A_s$ such
that $a_\sigma=a_1$ on the interval $[0,w(\varphi_s)]$. If we define
\[
b_\sigma:=w^{-1} a_\sigma^{-1} u_1
\]
then we have that
\begin{eqnarray*}
b_\sigma(t)=w^{-1}a_\sigma^{-1}a_1 w(t)= w^{-1} w(t)=t, \forall t \in [0,\varphi_s]
\end{eqnarray*}

\noindent Therefore $b_\sigma \in B_s$ and $a_\sigma w b_\sigma = u_1$
and so Eve can recover the shared secret key $K$ by using the pair $(a_\sigma,b_\sigma)$.

\remark{We observe
that any extension of $a_1|_{[0,w(\varphi_s)]}$ to an element $a_\sigma$ of $\PL_2([0,\varphi_s])$ will yield a suitable element
to attack Alice's key. Moreover, any element $a_1' \in A_s$
such that $a_1'wb_1'=u_1$, for some suitable $b_1' \in B_s$,
will be an extension of $a_1|_{[0,w(\varphi_s)]}$.}

\subsection{Attacking Bob's word for the case $w(\varphi_s) > \varphi_s$}

Eve considers $u_2^{-1}=a_2^{-1} w^{-1} b_2^{-1}$ and recovers a pair $(a_\sigma^{-1},b_\sigma^{-1})$
to get the shared secret key in the same fashion of the previous subsection.

\remark{Both the techniques of this section have been carried out using the transitivity of $A_s$ (Corollary \ref{thm:transitivity-A_s}).
They can also be solved by using the analogue of Corollary
\ref{thm:extension-A_s} for $B_s$ to get another pair $(a_\sigma,b_\sigma)$ which can be used to
retrieve the secret key.}

\section{Comments and Alternatives to the Protocol \label{sec:restriction-centralizers}}

This section analyzes possible alternatives and weaknesses of our methods. We observe that,
if instead of $\PL_2(I)$ we had used a larger group of piecewise linear homeomorphisms of the unit interval, the same
technique would have worked, as long as the commuting subgroups $A$ and $B$ had disjoint supports.
More generally, we can copy this idea if the given group $G$ acts on some space and we
have $A, B$ with disjoint support. We will now see some examples
of how this is possible.

\subsection{Choice of the subgroups $A$ and $B$}

We recall the following result from Chapter \ref{chapter4} (see Theorems \ref{thm:centralizers-F}
and \ref{intersectionofcetralizers}):

\theorem{Let $A=\langle a_1, \ldots, a_m \rangle \le F$ be a finitely generated subgroup. Then

\medskip
\noindent (i) There exists a dyadic partition of $[0,1]=I_1 \cup \ldots \cup I_n$ such that the centralizer
$C_F(A):=\{f \in F \, | \, af=fa, \forall a \in A\}$
is a product of subgroups $C_1, \ldots, C_n$, where $C_r \le \{f \in F \, | \, f(t)=t, \forall t \not \in I_r \}$.
Moreover, we have
\begin{itemize}
\item $C_r=\PL_2(I_r)$ if and only if of $a_i|_{I_r}=id$, for all $i=1, \ldots, r$.
\item $C_r \cong \mathbb{Z}$ if and only if $a_1|_{I_r}, \ldots, a_m|_{I_r}$ have a common root on $I_r$.
\item $C_r = 1$ if and only if there are $i \ne j$ such that $a_i|_{I_r}, a_j|_{I_r}$ have no common root on $I_r$.
\end{itemize}

\medskip
\noindent (ii) There exist two elements $g_1,g_2 \in F$ such that $C_F(A)=C_F(g_1) \cap C_F(g_2)$.}

Going back to the protocol introduced in Section \ref{sec:protocol} we observe that, after
we choose a finitely generated subgroup $A=\langle f_1,\ldots,f_m \rangle$,
we are very restricted in our choice of the subgroup $B$. Since $B \le C_F(A)$, we must make sure that the elements of
$B$, when restricted to $I_r$, are powers of common roots of the $a_i$'s, if at least one $a_i$ is non-trivial on $I_r$.
This gives a tight restriction on the subgroup $B$ whose support is essentially disjoint from that of $A$, except in the
intervals where they all are powers of a common root. An attack similar to that of Section \ref{sec:recovering-keys}
can thus be applied on each interval $I_r$: if their supports are disjoint on $I_r$, we can act as before, otherwise elements
of $A$ and $B$ are powers of a common root on $I_r$.

With more general commuting subgroups, the attack of Section \ref{sec:recovering-keys} does not immediately give either of the two keys.
However it is likely that a modification of the given algorithm can recover the shared
secret key for any choice of commuting subgroups $A$ and $B$.

\medskip

\subsection{Alternative Protocol and Attacks}

Ko-Lee et al.\ \cite{kolee} introduced a slightly different protocol based on the decomposition problem
(They worked with braid groups, but we will apply their protocol to Thompson's group).
In their protocol, Alice picks $a_1,a_2 \in A$ and sends $u_1 = a_1 w a_2$ to Bob, while Bob
chooses $b_1,b_2 \in B$ and sends $u_2=b_1 w b_2$ to Alice. We can still attempt to solve this new protocol, by
again dividing the problem into various cases. We assume that we use the same subgroups $A_s$ and $B_s$
and we work in the case $w(\varphi_s) \le \varphi_s$ to
show how to attack the private keys of Bob.
We apply the analogue for $B_s$ of Corollary \ref{thm:transitivity-A_s} and
find a $b_0$ such that $b_0^{-1}(w^{-1}(\varphi_s))=u_2^{-1}(\varphi_s)= b_2^{-1}w^{-1}(\varphi_s)$.
We define
\begin{eqnarray*}
b_1'=b_1 \\
b_2'= b_2 b_0^{-1} \\
u_2'=b_1'wb_2'
\end{eqnarray*}
so that $b_2'(w^{-1}(\varphi_s))=w^{-1}(\varphi_s)>\varphi_s$. Thus we have
\[
u_2'(t)=b_1'(t)wb_2'(t)=wb_2'(t), \forall t \in [0,w^{-1}(\varphi_s)]
\]
hence
\[
b_2'(t)=w^{-1}u_2'(t), \forall t \in [0,w^{-1}(\varphi_s)].
\]

\noindent Thus $b_2'$ is uniquely determined in $[0,w^{-1}(\varphi_s)]$.
We apply corollary \ref{thm:extension-A_s} for $B_s$ to find
a $b_{\sigma_2} \in B_s$ such that $b_{\sigma_2}=b_2'$ on $[0,w^{-1}(\varphi_s)]$
and we define
\[
b_{\sigma_1}:=u_2'b_{\sigma_2}^{-1}w^{-1}.
\]
Thus
\[
b_{\sigma_1}(t)=b_1' w b_2' b_{\sigma_2}^{-1}w^{-1}(t)=b_1'(t)=t, \forall t \in [0,\varphi_s]
\]
therefore $b_{\sigma_1} \in B_s$. Therefore the pair $(b_{\sigma_1},b_{\sigma_2})$ satisfies
$u_2'=b_{\sigma_1}wb_{\sigma_2}$ and so Eve can recover the shared secret key $K$. A similar argument
can be used to attack the element $a_1wa_2$, with the transitivity results for $A_s$.

\subsection{A comment on the Alternative Protocol}

The weakness in the protocol discussed in the previous
subsection arises from the fact that the chosen subgroups $A_s$ and $B_s$ are transitive on the intervals
on which they act nontrivially. This suggests that a possible way to
avoid such attacks is for $A$ and $B$ to be chosen to be not transitive on their support.

\remark{We observe that the attacks of section \ref{sec:recovering-keys-alternatives} and section
\ref{sec:restriction-centralizers} can be carried out in a fashion similar to that of Section \ref{sec:recovering-keys}, still producing
a solution in polynomial time.}

\chapter{Structure Theorems for Subgroups of Homeomorphisms Groups}
\label{chapter6}

Let $\Homeo(S^1)$ denote the full group of orientation-preserving homeomorphisms of the unit interval and $G$
be one of its subgroups. In this Chapter we recall the notion of \emph{rotation number} for an element of
$\Homeo(S^1)$. This number is invariant under conjugacy and describes the behavior of an element
under infinitely many iterations.
Loosely speaking, it describes how close to a rotation an element is, when we iterate it many times on a point of the circle.
It is a well known result that the rotation number map $rot:G \to \mathbb{R}/\mathbb{Z}$ is a group
homomorphism if the group $G$ is abelian. We give a direct proof that the same result is true if we assume that the group
$G$ has no non-abelian free subgroups. This was first deduced as a Corollary of a Theorem by Margulis in \cite{marg1}.
Our methods are independent and we recover Margulis's Theorem as a byproduct. We use our understanding of the rotation number
map to obtain a classification of subgroups of $\Homeo(S^1)$ and show how to build examples of such subgroups.

The Chapter is
organized as follows: Section \ref{sec:back-tools} recalls the necessary language and tools which will be used
in the Chapter; Section \ref{sec:rot-homo} shows that the rotation map is a homomorphism on subgroups as above;
Section \ref{sec:structure-embedding} explains the main structure theorem and
shows how to construct directly embeddings in $\Homeo(S^1)$ realizing the subgroups of the structure theorem;
Section \ref{sec:free-action} presents an analogue of Sacksteder's Theorem (see \cite{farbsha}) for fixed-point free subgroups,
showing that they must always be abelian;
Section \ref{sec:margulis}
uses the fact that the rotation map is a homomorphism to prove Margulis' Theorem on invariant measures on the unit circle.
The material of this Chapter represents
joint work with Collin Bleak and Martin Kassabov.

\section{Background and Tools \label{sec:back-tools}}

\begin{figure}[0.5\textwidth]
 \begin{center}
  \includegraphics[height=4.5cm]{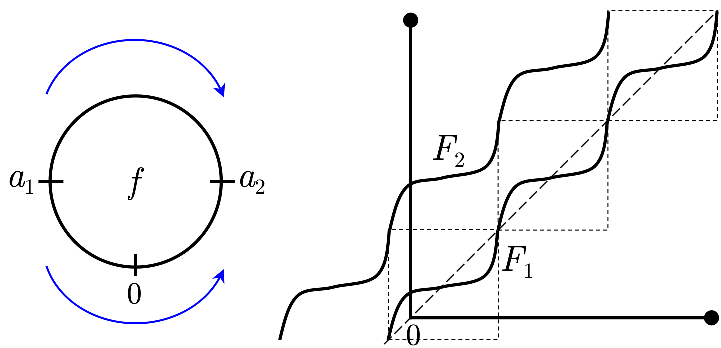}
 \end{center}
 \caption{Two lifts of a circle homeomorphism.}
 \label{fig:lift-3}
\end{figure}

In this section we collect some known results we will use throughout the Chapter. We begin by recalling the definition
of rotation number. Given $f \in \Homeo(S^1)$, let $F: \mathbb{R} \to \mathbb{R}$
represent one lift of $f$ (see figure \ref{fig:lift-3}) via the standard covering projection $\exp : \mathbb{R} \to S^1$,
where we think of $S^1$ as a subset of the complex
plane and use $\exp(t) = e^{2 \pi it}$.

\begin{figure}[0.5\textwidth]
 \begin{center}
  \includegraphics[height=4.5cm]{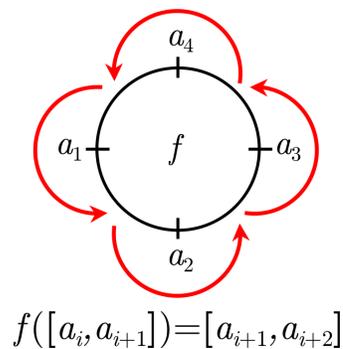}
 \end{center}
 \caption{A homeomorphism with rotation number $\frac{1}{4}$.}
 \label{fig:rot-number}
\end{figure}

\noindent We consider the limit
\begin{equation}\label{eq:rotation}
\lim_{n \to \infty} \frac{F^n(x)}{n}
\end{equation}
\noindent It is possible to prove that the previous limit exists and it is
is independent of the choice of $x$ used in the above calculation (see \cite{herman}).
Moreover, such a limit is independent of the choice of lift $\pmod{1}$.
\definition{We say that
\[
\lim_{n \to \infty} \frac{F^n(x)}{n} \pmod{1} := rot(f) \in  \mathbb{R}/\mathbb{Z}
\]
is the \textbf{rotation number} of $f$ (see figure \ref{fig:rot-number}.}

Since the rotation number is independent of the choice of the lift, we will work with
a preferred lift of elements and of functions.
For any element $x \in S^1$ we denote by $\wh{x}$ the lift of $x$ in $[0,1)$. If
$g \in \Homeo(S^1)$ and the fixed point set $Fix(g) = \emptyset$, we denote by $\wh{g}$
the lift to $\Homeo(\mathbb{R})$ such that $t < \wh{g}(t) < t+1$ for all $t \in \mathbb{R}$.
If $g \in \Homeo(S^1)$ and $Fix(g) \ne \emptyset$,
we denote by $\wh{g}$ the lift to $\Homeo(\mathbb{R})$ such that $Fix(\wh{g}) \ne \emptyset$.
We will use these definitions for lifts of elements and functions in Lemma \ref{thm:lift-control}(iv) and
throughout the proof of Theorem \ref{thm:rot-homo}.
If we use this lift to compute the limit defined in (\ref{eq:rotation}), the result is always in $[0,1)$.
A proof of the next three results can be found in \cite{herman} and \cite{mackay}.

\lemmaname{Properties of the Rotation Number}{Let $f,g \in \Homeo(S^1)$, $G\le \Homeo(S^1)$ and $k$ be a positive
integer. Then:

\noindent (i) $rot(f^g)=rot(f)$

\noindent (ii) $rot(f^k)=k \cdot rot(f)$

\noindent (iii) If $G$ is abelian then the map
\[
\begin{array}{cccc}

rot: & G & \longrightarrow & \mathbb{R}/\mathbb{Z} \\
     & f & \longmapsto     & rot(f)
\end{array}
\]
is a homomorphism

\noindent (iv) If $rot(g)=p/q \pmod{1} \in \mathbb{Q}/\mathbb{Z}$ and $s \in S^1$ is such that $g^q(s)=s$,
then $\wh{g}^q(\wh{s})=\wh{s}+p$. \label{thm:lift-control}}

\noindent Two of the most important results about the rotation number are stated below.

\theoremname{Poincar\'e}{Let $f \in \Homeo(S^1)$ be a homeomorphism. Then

\noindent (i) $f$ has a periodic orbit of length $q$ if and only if $rot(f)=p/q \pmod{1} \in \mathbb{Q}/\mathbb{Z}$ and
$p,q$ are coprime.

\noindent (ii) $f$ has a fixed point if and only if $rot(f) =0$.}

\theoremname{Denjoy}{Suppose  $f \in \Homeo(S^1)$ is piecewise-linear with finitely many breakpoints
or is a $C^1$ homeomorphism whose first derivative has bounded variation. If the rotation
number of $f$ is irrational, then $f$ is conjugate (by an element in
$\Homeo(S^1)$) to a rotation. Moreover, every orbit of $f$ is dense in $S^1$. \label{thm:denjoy}}

\noindent The following is a standard result proved by Fricke and Klein (independently)
which we will need in the proofs of section \ref{sec:rot-homo}. Our citation is to a more recent proof in English.

\begin{figure}[0.5\textwidth]
 \begin{center}
  \includegraphics[height=3cm]{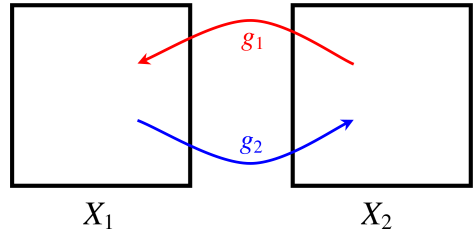}
 \end{center}
 \caption{A graphical description of the ping pong lemma.}
 \label{fig:ping-pong-lemma}
\end{figure}

\theoremname{Ping-Pong Lemma}{Let $G$ be a group of permutations on a set $X$, let $g_1, g_2$
be elements of $G$ of order at least three. If $X_1$ and $X_2$ are disjoint subsets of $X$ and for all integers
$n \ne 0$, $i \ne j$, $g_i^n(X_j) \subseteq X_i$ , then $g_1, g_2$ freely generate the free group $F_2$
on two generators (see figure \ref{fig:ping-pong-lemma} ). \label{thm:ping-pong}}
\noindent \emph{Proof.} See result 24 in section II.B of \cite{delaharpe1}. $\square$

\section{The Rotation Number Map is a Homomorphism \label{sec:rot-homo}}

Our main goal for this section is to prove the following result.

\theorem{Let $G \le \Homeo(S^1)$ with no non-abelian free subgroups. Then the $rot$ map is a homomorphism.
\label{thm:rot-homo}}

\begin{figure}[0.5\textwidth]
 \begin{center}
  \includegraphics[height=4.5cm]{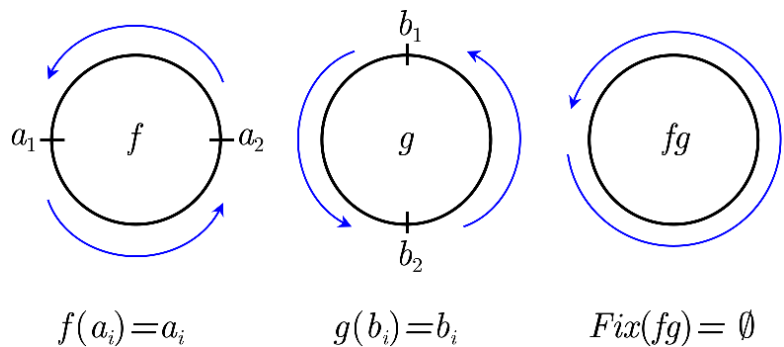}
 \end{center}
 \caption{The rotation map is not a homomorphism in general.}
 \label{fig:rot-number-zero-not-homomorphism-4}
\end{figure}

\noindent Before we begin with the proof of Theorem \ref{thm:rot-homo}, we want to give a short account of its history.
In his paper \cite{marg1}, Margulis proved a Theorem on the existence
of $G$-invariant measures on $S^1$ which yields Theorem \ref{thm:rot-homo} as a corollary.
Instead, we will give a direct proof of Theorem \ref{thm:rot-homo} and, in section \ref{sec:margulis},
we will derive the original Theorem of Margulis. We notice that the statement of Theorem \ref{thm:rot-homo}
does not hold in general: figure \ref{fig:rot-number-zero-not-homomorphism-4} shows two elements with rotation
number $0$ (hence they have fixed points) but whose product has no fixed points and must have non-zero rotation number.

Our proof divides naturally into several steps.

\lemma{Let $f,g \in \Homeo(S^1)$ such that $Fix(f) \ne \emptyset \ne Fix(g)$.
If the intersection $Fix(f) \cap Fix(g) = \emptyset$, then $\langle f,g \rangle$ contains a non-abelian free subgroup. Equivalently,
if $\langle f,g \rangle$ does not contain any non-abelian free subgroups, then $Fix(f) \cap Fix(g) \ne \emptyset$.
\label{thm:free-group}}

\begin{figure}[0.5\textwidth]
 \begin{center}
  \includegraphics[height=4.5cm]{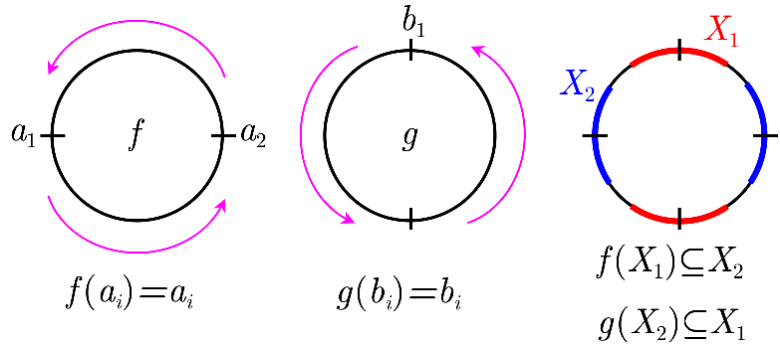}
 \end{center}
 \caption{Two elements generating a free subgroup.}
 \label{fig:ping-pong-4}
\end{figure}

\noindent \emph{Proof.} Let $S^1 \setminus Fix(f)=\bigcup I_\alpha$ and $S^1 \setminus Fix(g)=\bigcup J_\beta$,
for some suitable families of disjoint open intervals $\{I_\alpha\},\{J_\beta\}$. By construction,
$\partial I_\alpha \subseteq Fix(f)$ and $\partial J_\beta \subseteq Fix(g)$.
We assume $Fix(f) \cap Fix(g) = \emptyset$ so that
\[
S^1 \subseteq \left(\bigcup I_\alpha\right) \cup \left(\bigcup J_\beta\right).
\]
Since $S^1$ is compact, we can write $S^1=I_1 \cup \ldots \cup I_r \cup J_1 \cup \ldots \cup J_s$.
Define $I=I_1 \cup \ldots \cup I_r$ and $J=J_1 \cup \ldots \cup J_s$. Since each $x \in \partial J$
lies in the interior of $I$, then there is an open neighborhood $U_x$ of $x$ such that $U_x \subseteq I$.
Let $X_g = \bigcup_{x \in \partial J} U_x$. Similarly we build an open set $X_f$. If $x \in \partial J$,
then the sequence $\{f^n(x)\}_{n \in \mathbb{N}}$ accumulates at a point of $\partial I$ and so there is
an $n \in \mathbb{N}$ such that $f^n(U_x) \subseteq X_f$. By repeating this process for each $U_x$,
we can find a positive integer $n_0$ big enough so that $f^{n_0}(X_g) \cup f^{-n_0}(X_g) \subseteq X_f$.
We act similarly on $g$ and so we find an $N$ big enough so that for all $m \ge N$ we have
\[
f^m(X_g) \cup f^{-m}(X_g) \subseteq X_f, \; \; g^m(X_f) \cup g^{-m}(X_f) \subseteq X_g.
\]
If we define $g_1 = f^N,g_2=g^N,X_1=X_f,X_2=X_g$, we have satisfied
the hypothesis of Theorem \ref{thm:ping-pong} since both of the elements $g_1,g_2$ have infinite order.
Thus $\langle g_1,g_2 \rangle$ is a non-abelian free subgroup of $\langle f,g\rangle$
(see figure \ref{fig:ping-pong-4} to see an example of two elements generating a free group). $\square$

\corollary{Let $f,g \in \Homeo(S^1)$ such that $Fix(\wh{f}) \ne \emptyset \ne Fix(\wh{g})$.
If $Fix(\wh{f}) \cap Fix(\wh{g}) = \emptyset$, then $\langle f,g \rangle$ contains a non-abelian free subgroup.
\label{thm:free-group-2}}

\noindent If $G \le \Homeo(S^1)$ is a group, we define the set of homeomorphisms with fixed points
\[
G_0=\{g \in G \mid \exists s\in S^1, g(s)=s\} = \{ g \in G \mid rot(g)=0 \} \subseteq G.
\]

\corollary{Let $G \le \Homeo(S^1)$ with no non-abelian free subgroups.
The subset $G_0$ is a normal subgroup of $\Homeo(S^1)$.}

\noindent \emph{Proof.} Let $f,g \in G_0$ then, by Lemma \ref{thm:free-group}, they must have a common
fixed point, hence $fg^{-1} \in G_0$ and $G_0$ is a subgroup of $G$.
Moreover, if $f \in G,g \in G_0$
and $s \in Fix(g)$, we have that $f^{-1}(s) \in Fix(f^{-1}gf)$ and so that $f^{-1}gf \in G_0$
and therefore $G_0$ normal. $\square$

\noindent If $f$ has no fixed points then the support of $f$ is the whole circle $S^1$, otherwise the support can be broken
into open intervals on which $f$ is a one-bump function, that is $f(x) \ne x$ on them.
Given $f \in \Homeo(S^1)$, we define an \emph{orbital} as an open component of the support of $f$.

\noindent The following three Lemmas can be derived using techniques similar to those of \cite{bleakthesis}, however we
give a direct proof of them.

\begin{figure}[0.5\textwidth]
 \begin{center}
  \includegraphics[height=4cm]{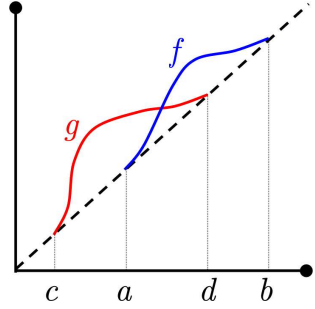}
 \end{center}
 \caption{Intersecting bumps.}
 \label{fig:build-large-weak}
\end{figure}

\lemma{Let $f,g \in \Homeo(S^1)$ and let $(a,b)$ be an orbital for $f$ and $(c,d)$ be an orbital for $g$
such that $c<a<d<b$ (see figure \ref{fig:build-large-weak}). For every $\varepsilon > 0$, there are two integers $M,N$ such that
$f^M g^N$ has an orbital containing $(c+\varepsilon,b-\varepsilon)$.
\label{thm:build-large-weak}}

\noindent \emph{Proof.} Without loss of generality we can assume that
$c + \varepsilon < a$ and $b-\varepsilon>d$. There is an integer $N \ne 0 $ such that
$s:=g^N(c+\varepsilon)>a$. Now choose a second integer $M \ne 0 $ such that $f^M(s) > b-\varepsilon$.
If $x \in (c+\varepsilon,b-\varepsilon)$, then $f^M g^N(x) > f^M g^N(c+\varepsilon)=f^M(s)>b-\varepsilon>x$,
hence $f^M g^N$ has an orbital containing $(c+\varepsilon,b-\varepsilon)$. $\square$

%
%
%

\lemma{Let $H\le \Homeo(S^1)$ and let $(a,b)$ be an interval
such that $Fix(H) \cap (a,b)=\emptyset$.
For every $\varepsilon >0$, there is an element $w \in H$ such that $w$ has an orbital
containing $(a+\varepsilon, b-\varepsilon)$. \label{thm:build-large}}

\noindent \emph{Proof.} We fix a point $c \in (a,b)$ and an $\varepsilon>0$.
We define the following notation: let $\mathcal{J}$ be the family of elements of $H$
that do not fix $c$. The set $\mathcal{J}$ is non-empty, since $c \not \in Fix(H)$.
For any element of $h \in \mathcal{J}$, let $O(h)$ denote the orbital of $h$
containing $c$ and let $l(h)$ and $r(h)$ respectively be the left and right endpoints of $O(h)$.

\noindent \emph{Claim 1.} The following supremum of all right endpoints of elements in $\mathcal{N}$
\[
R:=\sup_{h \in \mathcal{J}} \{ r(h) \}  \ge b.
\]
\noindent \emph{Proof.} Assume, by contradiction, that $R < b$.
Since $R \not \in Fix(H)$, there is an $f \in H$ with an orbital $(s,t)$
containing $R$ and we can assume that $f(x)>x$, for all $x \in (s,t)$. By definition of $R$, there is an element
$g \in \mathcal{J}$ with an orbital $O(g)$ such that its right endpoint $r(g)>s$. If $s \le  l(g)$, then
$f \in \mathcal{J}$ and its right endpoint is bigger than $R$, which is not possible. Hence we must have
$l(g)<s<r(g)<t$ and we can apply Lemma \ref{thm:build-large-weak} to find an element $f^Mg^N$ with an orbital
containing $(l(g)+\delta,t-\delta)$ and with $\delta$ chosen to be small enough that
$c,R \in (l(g)+\delta,t-\delta)$. This would imply that $f^Mg^N \in \mathcal{J}$ and
its right endpoint is bigger than $R$, which is a contradiction to the definition of supremum. $\square$

Thus we must have $R \ge b$ and so the family $\mathcal{J}$ has elements with ``large'' orbitals on the right.
We will now extend this procedure to make them ``large'' on the left. We define the following new subfamily of elements of
$\mathcal{J}$
\[
\mathcal{K}=\{h \in \mathcal{J} \; | \; r(h)> b-\varepsilon\}
\]
By Claim 1, the family $\mathcal{K}$ is non-empty.

\noindent \emph{Claim 2.} The following infimum of all right endpoints of elements in $\mathcal{K}$
\[
L:=\inf_{h \in \mathcal{K}} \{ l(h) \} \le a.
\]
\noindent \emph{Proof.} We repeat the idea of the previous Claim, by assuming that $L>a$ and then finding
two elements $f,g$ on which we can apply Lemma \ref{thm:build-large-weak}. $\square$

\noindent By Claim 2, we can choose an element $w \in \mathcal{K}$ that has an orbital with $l(w)<a+\varepsilon$. By
definition of $\mathcal{K}$, $w$ satisfies the thesis. $\square$

\lemma{Let $H \le \Homeo(S^1)$ and let $(a_1,b_1),\ldots,(a_r,b_r)$ be disjoint intervals
such that
\[
Fix(H) \cap \left( \bigcup_{i=1}^r(a_i,b_i) \right)=\emptyset.
\]
For every $\varepsilon >0$, there is an element $w \in H$ such that,
the element $w$ has a support containing $(a_i+\varepsilon, b_i-\varepsilon)$,
for each $i=1,\ldots,r$.}

\noindent \emph{Proof.} By induction on the number of intervals $r$. The case $r=1$ has been proven in Lemma
\ref{thm:build-large}. We assume the result holds for the $r-1$ intervals $(a_1,b_1),\ldots,(a_{r-1},b_{r-1})$.
Define the family of elements
\[
\mathcal{J}=\left\{h \in H \; | \; \mathrm{supp}(h) \supseteq \bigcup_{i=1}^{r-1}[a_i+\varepsilon,b_i-\varepsilon] \right\}
\]
By induction hypothesis, the family $\mathcal{J}$ is non-empty. We also observe that $\mathrm{supp}(h)$
is always an open set, hence it is a union of intervals that contains the closed set
$\bigcup_{i=1}^{r-1}[a_i+\varepsilon,b_i-\varepsilon]$ properly. We fix a point $c \in (a_r,b_r)$.
We now want to prove that $\mathcal{J}$ contains elements that do not fix $c$.

\begin{figure}[0.5\textwidth]
 \begin{center}
  \includegraphics[height=4cm]{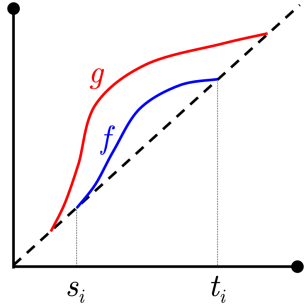}
 \end{center}
 \caption{Non-intersecting bumps.}
 \label{fig:build-large-strong}
\end{figure}

\noindent \emph{Claim.} The following subfamily of $\mathcal{J}$:
\[
\mathcal{K}=\left \{h \in \mathcal{J} \; \Big \vert \; h(c) \ne c \right \} \ne \emptyset.
\]
\noindent \emph{Proof.} Let $f \in \mathcal{J}$. If $f(c) \ne c$, we are done. Otherwise,
suppose that $f(c)=c$. Since $c \not \in Fix(H)$, there is a $g \in H$ such that $g(c) \ne c$.
For each $i=1,\ldots,r-1$ consider the interval $(s_i,t_i)$ of the support of $f$ containing $[a_i+\varepsilon,b_i-\varepsilon]$
properly. On each $(s_i,t_i)$ we have two cases: (i) an orbital of $g$ contains $s_i$ or $t_i$, so we
we can apply Lemma \ref{thm:build-large-weak}
to find integers $M_i,N_i$ such that $f^{M_i} g^{N_i}$ has an orbital containing $[a_i+\varepsilon,b_i-\varepsilon]$, or
(ii) an orbital of $g$ contains $(s_i,t_i)$ properly and again we can find powers
$M_i,N_i$ such that $f^{M_i} g^{N_i}$ has an orbital containing $[a_i+\varepsilon,b_i-\varepsilon]$
(see figure \ref{fig:build-large-strong}).
If we now take $M=\max\{M_1,\ldots,M_{r-1}\}$ and $N=\max\{N_1,\ldots,N_{r-1}\}$,
we have that $f^M g^M \in \mathcal{J}$ and, by construction, $f^M g^N \in \mathcal{J}$. $\square$

\noindent The proof now proceeds as in Lemma \ref{thm:build-large}. We first find elements with orbitals
whose right endpoint is near $b_r$ and then do the same on the left. As done in the previous Claim, this
procedure can be followed so that the support of the families of elements always contain properly
the union $\bigcup_{i=1}^{r-1}[a_i+\varepsilon,b_i-\varepsilon]$. $\square$

\lemmaname{Finite Intersection Property}{Let $G \le \Homeo(S^1)$ with no non-abelian free subgroups.
Then the family $\{Fix(g)| \, g \in G_0\}$ satisfies the finite intersection property,
i.e. for all $n$-tuples $g_1, \ldots,g_n \in G_0$, we have $Fix(g_1) \cap \ldots \cap Fix(g_n) \ne \emptyset$.}

\noindent \emph{Proof.} We use induction on $n$, with the case $n=2$ being true by Lemma \ref{thm:free-group}.
We assume that the result is true for any $(n-1)$-tuple of elements in $G_0$.
Let $g_1, \ldots g_n \in G_0$ and assume, by contradiction, that $Fix(g_1) \cap \ldots \cap Fix(g_n) = \emptyset$.
Let $H=\langle g_1, \ldots, g_{r-1} \rangle$ and notice that $Fix(H) \ne \emptyset$ by induction hypothesis.
Write $S^1 \setminus Fix(H)=\bigcup I_\alpha$ and $S^1 \setminus Fix(g_r)=\bigcup J_\beta$,
for some suitable families of open intervals $\{I_\alpha\},\{J_\beta\}$. By construction,
$\partial I_\alpha \subseteq Fix(H)$ and $\partial J_\beta \subseteq Fix(g_r)$.
Since $Fix(H) \cap Fix(g_r) = \emptyset$ we have
\[
S^1 \subseteq \left(\bigcup I_\alpha\right) \cup \left(\bigcup J_\beta\right).
\]
Since $S^1$ is compact, we can write $S^1=I_1 \cup \ldots \cup I_r \cup J_1 \cup \ldots \cup J_s$
and notice that $Fix(g_r) \subseteq \bigcup_{i=1}^r I_i$.
Since the intersection $Fix(H) \cap \left( \bigcup_{m=1}^s I_m \right) = \emptyset$. If $I_i=(a_i,b_i)$
we apply Lemma \ref{thm:build-large} to build an element $w \in H$ such that $w$ has an orbital containing
$(a_i+\varepsilon,b_i-\varepsilon)$, for every $i=1,\ldots,r-1$. We can choose $\varepsilon>0$
to be small enough so $Fix(g_r) \subseteq \bigcup_{i=1}^r (a_i+\varepsilon,b_i-\varepsilon) \subseteq \mathrm{supp}(w)$
thus implying that $Fix(w)\cap Fix(g_r)=\emptyset$.
We can again apply Lemma \ref{thm:free-group} to build a  non-abelian free group inside $\langle w,g_n\rangle$,
contradicting the assumption on $G$.
Thus, for every finite set $H \subset G_0$, we have
\[
\bigcap_{h \in H}Fix(h) \ne \emptyset
\]
which proves that the family $\{Fix(g) | \, g \in G_0\}$ has the finite intersection property. $\square$
%

\corollary{Let $G \le \Homeo(S^1)$ with no non-abelian free subgroups.
The subgroup $G_0$ admits a global fixed point, i.e. $Fix(G_0)\ne \emptyset$.}

\noindent \emph{Proof.} By the previous Lemma we have that the family $\{Fix(g)| \, g \in G_0\}$ has the finite intersection property. By compactness of the unit circle $S^1$ we have that:
\[
Fix(G_0) = \bigcap_{g \in G_0} Fix(g) \ne \emptyset. \; \;\square
\]

In order to prove Theorem \ref{thm:rot-homo} we observe that the element $(fg)^k$ can be rewritten $f^k g^k h_k$
for some suitable product of commutators $h_k \in [G,G]$. If we prove that every element $[G,G]$ has a global
fixed point $s$ we can compute the rotation number on $s$, so that $(fg)^k(s)=f^k g^k$.
We will prove that this is indeed the case.

\lemma{Let $G \le \Homeo(S^1)$ and let $f,g \in G$. Suppose one of the following two cases is true:

\noindent (i) $G$ has no non-abelian free subgroups and $rot(f)=rot(g) \in \mathbb{Q}/\mathbb{Z}$, or

\noindent (ii) $rot(f)=rot(g) \not \in \mathbb{Q} / \mathbb{Z}$.

\noindent Then $fg^{-1} \in G_0$. \label{thm:same-rot-fix-point}}

\noindent \emph{Proof.} (i) Assume $rot(f)=rot(g)=k/m \in \mathbb{Q} / \mathbb{Z}$ with $k,m$ positive integers
and that $G$ has no non-abelian free subgroups. Moreover, $f^m$ and $g^m$ have fixed points in $S^1$
and $\wh{f}^m(\wh{x})=\wh{x}+k$ and $\wh{g}^m(\wh{y})=\wh{y}+k$, for any $x \in Fix(f^m),y \in Fix(g^m)$ by Lemma
\ref{thm:lift-control}(iv).
Thus $f^m$ and $g^m$ must have a common fixed point $s \in S^1$ by Lemma \ref{thm:free-group}.
We argue, by contradiction, that $fg^{-1} \not \in G_0$, so that $\wh{f}>\wh{g}$ or $\wh{f}<\wh{g}$.
Suppose the former so $\wh{f}^m>\wh{g}^m$, but this is impossible as $\wh{f}^m(\wh{s})=\wh{s}+k=\wh{g}^m(\wh{s})$.

\begin{figure}[0.5\textwidth]
 \begin{center}
  \includegraphics[height=4cm]{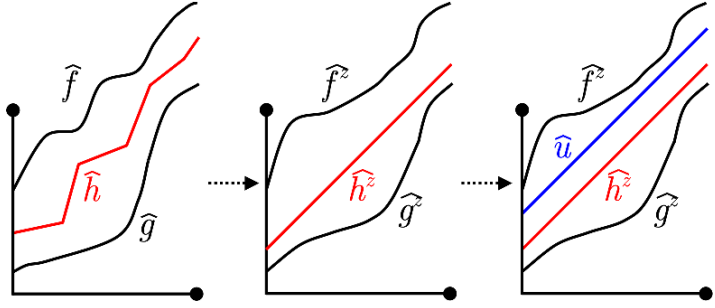}
 \end{center}
 \caption{Making room for straight lines between $\wh{f}$ and $\wh{g}$.}
 \label{fig:straighten-maps}
\end{figure}

(ii) Assume now that $rot(f)=rot(g) \not \in \mathbb{Q} / \mathbb{Z}$. Again, we argue by contradiction that
$fg^{-1} \not \in G_0$ and we suppose $\wh{f}>\wh{g}$. We observe that, for any map $h \in \Homeo(S^1)$
such that $\wh{f} \ge \wh{h} \ge \wh{g}$, we have $rot(f) \ge rot(h) \ge rot(g)=rot(f)$ and
so $rot(h) \not \in \mathbb{Q}/\mathbb{Z}$. By compactness of $S^1$ and the fact that $\wh{f}>\wh{g}$,
we can find an $h \in \PL_+(S^1)$
such that $\wh{f}>\wh{h}>\wh{g}$. We use Denjoy's theorem
on $h$ to find $z \in Homeo(S^1)$ such that $h^z$ ia a rotation. Therefore
$\wh{h^z}$ is a straight line $t \to t+rot(g)$. Now since
$\wh{f^z}>\wh{h^z}$, we can find a rotation $u \in Homeo(S^1)$ such that
its lift is a straight line $\wh{u}:t \to t+rot(u)$ with $rot(u)>rot(h)$ and $\wh{f^z}> \wh{u} > \wh{h^z}$.
(see figure \ref{fig:straighten-maps}). To conclude we observe that
\[
rot(f)=rot(f^z) \ge rot(u) > rot(h^z)=rot(h)=rot(f)
\]
yielding a contradiction. $\square$

\corollary{Let $G \le \Homeo(S^1)$ with no non-abelian free subgroups, then we have $[G,G] \le G_0$.
\label{thm:commutator-fixed}}


\corollary{Let $G \le \Homeo(S^1)$ with no non-abelian free subgroups. Suppose that
every element of $G \setminus \{id\}$ has no fixed points. Then $G$ is abelian. \label{thm:non-fixed-abelian-1}}


\noindent Corollary \ref{thm:non-fixed-abelian-1} is true in a greater generality. In fact, Theorem
\ref{thm:non-fixed-abelian-2} proves it without any
requirement on the subgroups of $G$.

\lemma{Let $G \le \Homeo(S^1)$ with no non-abelian free subgroups.
Let $f,g \in G$ and $s \in S^1$ be a fixed point of $[f,g]$.
Then $\wh{s}$ is a fixed point for $[U,V]$, for any $U$ lift of $f$ and $V$ lift of $g$ in $\Homeo(\mathbb{R})$.
\label{thm:lift-commutator-fixed}}

\noindent \emph{Proof.} If $T(x)=x+1$ then $U=T^m \wh{f}$ and $V=T^n \wh{g}$ for some suitable integers $m,n$.
Moreover, it is immediate to verify that
\[
[U,V](x) = [\wh{f},\wh{g}](x), \forall x \in \mathbb{R}.
\]
Thus, we need only to prove that $[\wh{f},\wh{g}](\wh{s})=\wh{s}$. We divide the proof into two cases.

\noindent \emph{Case 1: $Fix(\wh{f}) = \emptyset = Fix(\wh{g})$.} We have that $t < \wh{f}(t) < t+1$ and
$t < \wh{g}(t) < t+1$ for all $t \in \mathbb{R}$. Since $fg(s)=gf(s)$ we also have
$k:=\wh{f}\wh{g}(\wh{s})-\wh{g}\wh{f}(\wh{s}) \in \mathbb{Z}$.
If $\wh{f}\wh{g}(\wh{s}) \ge \wh{g}\wh{f}(\wh{s})$ then, since
$\wh{g}$ is increasing and $\wh{s}<\wh{f}(\wh{s})$, we have $\wh{g}\wh{f}(\wh{s})>\wh{g}(\wh{s})$
and so

\[
|k|= \wh{f}\wh{g}(\wh{s})-\wh{g}\wh{f}(\wh{s}) \le  \wh{f}\wh{g}(\wh{s}) - \wh{g}(\wh{s}) < 1
\]
implying that $k=0$. A similar argument holds if $\wh{f}\wh{g}(\wh{s}) \le \wh{g}\wh{f}(\wh{s})$.

\noindent \emph{Case 2: $Fix(\wh{f}) \ne \emptyset$ or $Fix(\wh{g}) \ne \emptyset$.}
We can assume that $Fix(\wh{f}) \ne \emptyset$. Then also $Fix(\wh{g}^{\,-1}\wh{f} \,\wh{g}) \ne \emptyset$.
By Corollary \ref{thm:free-group-2} we must have $Fix(\wh{f}) \cap Fix(\wh{g}^{\,-1}\wh{f} \,\wh{g}) \ne \emptyset$
and so this implies that $[\wh{f},\wh{g}]$ intersects the diagonal. Since $[f,g](s)=s$, then
$[\wh{f},\wh{g}](\wh{s})=\wh{s}$. $\square$

\noindent Now we are ready to prove the main theorem of this section.

\noindent \emph{Proof of Theorem \ref{thm:rot-homo}.} Let $f,g \in G$. We write the power
$(fg)^k=f^k g^k h_k$ where $h_k$ is a suitable product of commutators
used to shift the $f$'s and $g$'s to the left. Since $h_k \in [G,G] \le G_0$ for all positive integers $k$ then,
if $s \in S^1$ is a global fixed point for $G_0$, we have $h_k(s)=s$. Similarly, we observe that
$(\wh{f}\wh{g})^k=\wh{f}^k \, \wh{g}^k \, H_k$ where $H_k$ is a suitable product of commutators and
$H_k$ is a lift for $h_k$. By Lemma \ref{thm:lift-commutator-fixed} we must have that
$H_k(\wh{s})=\wh{s}$ for all positive integers $k$. Thus we observe that:
\[
(\wh{f}\wh{g})^n(\wh{s}) = \wh{f}^n \, \wh{g}^n \, H_n(\wh{s}) =
\wh{f}^n \, \wh{g}^n (\wh{s}).
\]
We now find upper and lower bounds for $\wh{f}^n \, \wh{g}^n (\wh{s})$. Observe that,
for any two real numbers $a,b$ we have that
\[
\wh{f}^n(a) + b - 1 < \wh{f}^n(a) + \lfloor b \rfloor \le  \wh{f}^n(a+b) <
\wh{f}^n(a)+ \lfloor b \rfloor +1 \le \wh{f}^n(a) + b +1
\]
where $\lfloor \cdot \rfloor$ denotes the floor function. By applying this inequality to
$\wh{f}^n \, \wh{g}^n (\wh{s}) = \wh{f}^n(\wh{s}+(\wh{g}^n(\wh{s})-\wh{s}))$ we get
\begin{eqnarray*}
\wh{f}^n(\wh{s})+\wh{g}^n(\wh{s})-\wh{s} - 1 \le
\wh{f}^n(\wh{s}+(\wh{g}^n(\wh{s})-\wh{s})) \le
\wh{f}^n(\wh{s})+\wh{g}^n(\wh{s})-\wh{s}+1.
\end{eqnarray*}
We divide the previous inequalities by $n$, and get
\[
\frac{\wh{f}^n(\wh{s})+\wh{g}^n(\wh{s})-\wh{s}-1}{n} \le \frac{(\wh{f} \, \wh{g})^n(\wh{s})}{n} \le \frac{\wh{f}^n(\wh{s})+\wh{g}^n(\wh{s})-\wh{s}+1}{n}.
\]
By taking the limit for $n \to \infty$ of the previous expression, we immediately
obtain $rot(fg)=rot(f)+rot(g)$. $\square$

\corollary{Let $G \le \Homeo(S^1)$ with no non-abelian free subgroups. Then $rot:G \to \mathbb{R}/\mathbb{Z}$
is a group homomorphism and

\noindent (i) $\ker(rot)=G_0$,

\noindent (ii) $G/G_0 \cong rot(G)$.

\noindent (iii) for all $f,g \in G$, $fg^{-1} \in G_0$ if and only if $rot(f)=rot(g)$. \label{thm:G_0-normal}}

\section{Applications: Margulis' Theorem \label{sec:margulis}}

In this section we show how the techniques developed in Section \ref{sec:rot-homo} yield two known results
for groups of homeomorphisms of the unit circle.

\theoremname{Margulis, \cite{marg1}}{Let $G \le \Homeo(S^1)$. Then at least one of the two following statements must be
true:

\noindent (i) $G$ has a non-abelian free subgroup, or

\noindent (ii) there is a $G$-invariant probability measure on $S^1$.}

\noindent \emph{Proof.} We assume that (i) does not hold. The proof must be divided into two cases.

\noindent \emph{Case 1: $G /G_0$ is finite.} Let $s \in Fix(G_0)$ and consider the finite orbit $s^G$.
Then for every subset $X \subseteq S^1$ assign
\[
\mu(X)=\frac{\# \, s^G \cap X}{\# \, s^G}.
\]
This obviously defines a probability measure on $S^1$.

\noindent \emph{Case 2: $G/G_0$ is infinite and therefore $rot(G)$ is dense in $\mathbb{R} / \mathbb{Z}$.}
Fix $s \in Fix(G_0)$ as an origin and write $S^1$ as $[0,1]$.
We regard $s^G$ as a subset of $[0,1]$ and define the map
$\varphi:s^G \to rot(G)$, given by $\varphi(s^g)=rot(g)$, for any $g \in G$. It is immediate that
$\varphi$ is well-defined and order-preserving. We take the ``closure'' of
this map, by defining
\[
\begin{array}{cccc}
\ov{\varphi}: &  [0,1] & \longrightarrow & [0,1] \\
              &    a   & \longmapsto     & \sup \{rot(g) \; | \; s^g \le a. \; g \in G \}.
\end{array}
\]
By construction, the map $\ov{\varphi}$ is order-preserving. Moreover, since the image of $\ov{\varphi}$
contains $rot(G)$, it is dense in $[0,1]$. Since $\ov{\varphi}$ is an order-preserving map whose image is dense in
$[0,1]$, then $\ov{\varphi}$ is a continuous map. This allows us to define the Lebesgue-Stieltjes
measure associated to $\ov{\varphi}$ on the Borel algebra of $S^1$ (see \cite{malliavin}),
that is, for every half-open interval $(a,b] \subseteq S^1$ we define
\[
\mu((a,b]):=\ov{\varphi}(b)-\ov{\varphi}(a).
\]
Since the map $rot$ is a homomorphism, it is straightforward to see
that the measure $\mu$ is $G$-invariant. For example,
consider an interval $(s^{g_1},s^{g_2}]$ such that $rot(g_1) \le rot(g_2) <1$ (that is, neither $s^{g_1}$ nor
$s^{g_2}$ wrap around the circle and pass $s$). If $g \in G$, we have that
\[
\mu(s^{gg_1},s^{gg_2}]=rot(gg_2)-rot(gg_1)=rot(g)+rot(g_2)-rot(g)-rot(g_1)=\mu(s^{g_1},s^{g_2}]
\]
The other cases are dealt similarly, by doing additions and subtractions in $\mathbb{R}/\mathbb{Z}$
The same can be verified for any other half-open interval. By definition of the measure, $\mu(S^1)=1$
and $\mu(\mathrm{p})=0$, for every point $p \in S^1$, and we are done. $\square$

The following result is strongly believed to be well known, but unfortunately we were unable to
find a reference for it.

\theorem{Suppose $G$ is a subgroup of $\Homeo(S^1)$ which contains no non-abelian free
subgroups, and there is an element $g \in G$ such that

\noindent (i) $rot(g) \not \in \mathbb{Q} / \mathbb{Z}$, and

\noindent (ii) $g$ is piecewise-linear with finitely many breakpoints, or $C^1$ with bounded variation in its first derivative,

\noindent then $G$ is topologically conjugate to a group of rotations.
In particular, $G$ is abelian.\label{thm:irrational-abelian}}

\noindent \emph{Proof.} By Denjoy's Theorem \ref{thm:denjoy} the orbits of $g$ are dense in $S^1$.
Suppose there is $id \ne h \in G_0$ and let $s \in Fix(G_0)$. Then $g(s) \in Fix(h)$ in fact
\[
h(g(s))=g g^{-1}hg(s)= g(h^g(s))=g(s)
\]
since $h^g \in G_0 \unlhd G$. Thus $h$ must fix the sequence of points
$\{g^k(s)\}$, which is dense in $S^1$ and so $h$ must fix the whole $S^1$,
giving a contradiction, since $h \ne id$. Thus $Fix(G_0) = S^1$ and so $G_0$ is trivial.
By Corollary \ref{thm:G_0-normal} we have $G \cong G/G_0 \cong rot(G) \le \mathbb{R} / \mathbb{Z}$.
By Denjoy's Theorem \ref{thm:denjoy}, there is a $z \in \Homeo(S^1)$ such that $g^z$ is a rotation. Thus
$G^z \le C_{\Homeo(S^1)}(g^z)=\{$all rotations$\}$. $\square$

\section{Structure and Embedding Theorems \label{sec:structure-embedding}}

We start the section with a result which classifies the structure of subgroups of $\Homeo(S^1)$
with no non-abelian free subgroups. We consider an orbit $s^G$ of a point $s$ of $Fix(G_0)$ under the
action of $G$ (so that any of the points of the closure $\ov{s^G} \subseteq Fix(G_0)$),
then we choose a fundamental domain $D$ for the action of $G$ on $S^1 \setminus \ov{O(x)}$.
Since $S^1 \setminus \ov{s^G}$ is open, the fundamental domain will be given by a union of intervals. By restricting
$G_0$ to this fundamental domain and we get a group $H_0$ which acts as a set of homeomorphisms of
a disjoint union of intervals. We will prove that $G$ is isomorphic to the wreath structure of $G/G_0$ over the group $H_0$ which acts on the fundamental domain.

\theorem{Let $G \le \Homeo(S^1)$ with no non-abelian free subgroups. Then:

\noindent (i) $G$ is abelian, or

\noindent (ii) $G \hookrightarrow H_0 \wr T$, the standard unrestricted wreath product,
where $T:= G/G_0$ is isomorphic to a countable subgroup
of $\mathbb{R}/\mathbb{Z}$ and $H_0 \le \prod \Homeo(I_i)$ has no non-abelian free subgroups.
\label{thm:nonfree-classification}}

\remark{If $G \le \PL_+(S^1)$ is non-abelian,
then $T \le \mathbb{Q} / \mathbb{Z}$ because of Denjoy's Theorem. This is also a consequence
of Theorem \ref{thm:irrational-abelian}.}

\noindent \emph{Proof.} (i) If $G_0=\{id\}$ then $G \cong G/G_0 \cong rot(G) \le \mathbb{R}/\mathbb{Z}$.
(ii) Suppose $G_0$ non-trivial, so that $Fix(G_0) \ne S^1$ and define $T=G/G_0$. Let $s \in Fix(G_0)$
and consider the open subset $S^1 \setminus \ov{s^T}$, where $s^T$ is the orbit of $s$ under the action of $T$.
The set $S^1 \setminus \ov{s^T}$ is a collection
of, at most countably many, disjoint open intervals. We can define a fundamental domain for the action
of $T$ on $S^1 \setminus \ov{s^T}$ as the union $D=\bigcup_{i \in \mathcal{I}} I_i$ of a
collection $\{I_i\}_{i \in \mathcal{I}}$  of at most countably many intervals $I_i$ such that
\begin{eqnarray*}
t_1(D) \cap t_2(D) = \emptyset, \, \, t_1 \ne t_2,  \\
\\
S^1 \setminus \ov{s^T} = \bigcup_{t \in T} t(D)
\end{eqnarray*}

\noindent \emph{Claim 1:} The fundamental domain $D$ exists.

\noindent \emph{Proof of Claim 1.} Let $T$ act on $S^1 \setminus \ov{s^T}$ and consider two intervals
$I_1,I_2$ to be equivalent if there is $t \in T$ such that $t(I_1)=I_2$. For each equivalence class $\mathcal{C}_i$,
we apply the Axiom of Choice to choose an interval $I_i$
representing the class. We define $D$ to be the union of these representatives. $\square$

\noindent Since $\ov{s^T} \subseteq Fix(G_0)$ we have
\[
S^1 \setminus \bigcup_{t \in T} t(D) \subseteq Fix(G_0).
\]

\noindent \emph{Claim 2:} Define $H_0:=G_0|_{D}$ restriction of $G_0$ to $D$. Then
there is an embedding $H_0^t=G_0|_{t^{-1}(D)} \hookrightarrow \prod_{i\in \mathcal{I}} \Homeo(t^{-1}(I_i))$.
In particular, $H_0 \hookrightarrow \prod_{i\in \mathcal{I}} \Homeo(I_i)$.

\noindent \emph{Proof of Claim 2.} This is immediate, once we observe that if
$h \in H_0$ and $t \in G\setminus G_0$, then $t^{-1}ht(t^{-1}(D))=t^{-1}(D)$,
since $G_0$ fixes the endpoints of the intervals $I_i$. $\square$

\noindent It is important to notice that $H_0$ is not necessarily contained in $G_0$, since $H_0$ has its
support in $D$, while an element of $G_0$ has support in $\bigcup_{t \in T} t(D)$.
From Claim 2 it is now obvious that there is an embedding
\[
\varphi: G_0 \hookrightarrow  \prod_{t \in T} H_0^t =\prod_{t \in T} \prod_{i \in \mathcal{I}} \Homeo(t^{-1}(I_i))
\]

\noindent We observe $\prod_{t \in T} H_0^t \cong \langle H_0^t \,| t \in T \rangle \le \Homeo(S^1)$
and we define $H := \langle H_0^t \mid t \in T \rangle$ and
\[
E:=\langle G, H_0^t \mid t \in T \rangle \le \Homeo(S^1).
\]
By definition of $E$ we get the following exact sequence
\[
1 \to H \overset{i}{\to} E \overset{\pi}{\to} E/H \to 1
\]
where $i$ is the inclusion map and $\pi$ is the natural projection $\pi:E \to E/H$.
Notice that $E/H \cong G/(G\cap H)$ and $G \cap H \le G_0$, by definition of $G_0$.
By the argument above $G_0 \le H$ and so $G \cap H = G_0$ thus implying that $E/H \cong G/G_0=T$,
so we can rewrite the sequence as
\[
1 \to H \overset{i}{\to} E \overset{\pi}{\to} T \to 1
\]
where $T$ acts on the base group by shifts.

\noindent \emph{Claim 3:} $E \cong E_0:=H_0 \wr T$

\noindent \emph{Proof of Claim 3.}
By a standard result in cohomology of groups (see Theorem 11.4.10 in \cite{rob}),
if we can prove that $H^2(T,Z(\prod H_0^t))=0$ (where $Z(\prod H_0^t)$ denotes the center of $\prod H_0^t$),
there can be only one possible extension of $\prod H_0^t$ by $T$.
We observe that $H_0 \wr T$ is one such extension, so it suffices to prove that $H^2(T,Z(\prod H_0^t))=0$ to show that
any other extension will be equivalent to the wreath product $H_0 \wr T$.
We use Shapiro's Lemma to compute this cohomology group (see Proposition 6.2 in \cite{brown1}). We have
\begin{eqnarray*}
H^2(T,Z(\prod H_0^t))=H^2(T,\prod Z(H_0)^t)= \\
=H^2(T,\mbox{Coind}^T_{\{id\}} Z(H_0)) = H^2(\{id\},Z(H_0))=0. \; \; \; \; \; \; \; \square
\end{eqnarray*}
\remark{We observe that the wreath product in the previous result is unrestricted,
because the elements of $\Homeo(S^1)$ can have infinitely many bumps and so the elements
of $G_0$ may be non-trivial on infinitely many intervals. Conversely, if we assume $G \le \PL_+(S^1)$,
this would imply that any element in $G_0$ is non-trivial only at finitely many intervals, and so that
$G_0$ can be embedded in the direct sum $\bigoplus$.
This argument explains why the wreath products in the following Theorem \ref{thm:first-embedding}
is unrestricted and the ones in Theorems \ref{thm:embedding-thompson} and \ref{thm:embedding-PL}
are restricted.}

\remark{We notice that, in the statement of the previous Theorem, we may replace the conclusion
``$H_0 \le \prod \Homeo(I_i)$'' with $H_0 \le \Homeo(I)$, because we can always build an embedding
$\prod \Homeo(I_i) \hookrightarrow H_0$.}

\noindent We now turn to prove existence results and show that subgroups with wreath product structure do exist
in $\Homeo(S^1)$ and in $\PL_+(S^1)$.

\theorem{For every $T \le \mathbb{R} / \mathbb{Z}$ countable and for every $H_0 \le \Homeo(I)$
there is an embedding $H_0 \wr T \hookrightarrow \Homeo(S^1)$,
where the wreath product $H_0 \wr T=\left(\prod H_0 \right) \rtimes T$ is unrestricted. \label{thm:first-embedding}}


\noindent \emph{Proof.} We divide the proof into two cases: $T$ infinite and $T$ finite.
If $T$ is infinite, we enumerate the elements of $T=\{t_1, \ldots, t_n, \ldots \}$ and we
build a sequence:
\[
\frac{1}{2}, \frac{1}{2^2}, \ldots, \frac{1}{2^n}, \ldots
\]
We identify $S^1$ with the interval $[0,1]$ to fix an origin and an orientation of the unit circle.
$T$ is countable subgroup of $\mathbb{R}/\mathbb{Z}$, so it is non-discrete and therefore it is dense in $S^1$.
Now define the following map:
\[
\begin{array}{cccc}
\varphi: & [0,1] =S^1 & \longrightarrow & [0,1] =S^1\\
         & x          & \longmapsto     & \sum_{t_i<x} \frac{1}{2^i}
\end{array}
\]
(where $t_i < x$ is written with respect to the order in $[0,1]$). It is immediate from
the definition to see that the map is order-preserving and it is injective, when restricted
to $T$. Observe now that
\begin{eqnarray*}
\varphi(t_1) = \sum_{t_i<t_1} \frac{1}{2^i} \\
\varphi(t_1 + \varepsilon) = \sum_{t_i< t_1 + \varepsilon} \frac{1}{2^i}
\end{eqnarray*}
If we let $\varepsilon \to 0$, we see
\begin{eqnarray*}
\alpha : = \varphi(t_1) \le \varphi(t_1 + \varepsilon) \underset{\varepsilon \to 0}{\longrightarrow}
\sum_{t_i \le t_1} \frac{1}{2^i} = \alpha + \frac{1}{2}
\end{eqnarray*}
Since $\varphi$ is order-preserving, we have $(\alpha, \alpha + \frac{1}{2}) \cap \varphi(T)=\emptyset$.
More generally, it can be seen that
\[
\bigcup_{i \in \mathbb{N}} \left(\varphi(t_i),\varphi(t_i)+\frac{1}{2^i}\right) \cap \ov{\varphi(T)} = \emptyset
\]
\noindent \emph{Claim.} The unit circle can be written as the disjoint union
\[
S^1=\bigcup_{i \in \mathbb{N}} \left(\varphi(t_i),\varphi(t_i)+\frac{1}{2^i}\right) \cup \ov{\varphi(T)}.
\]
\noindent \emph{Proof of Claim.} Let $X = \bigcup_{i \in \mathbb{N}} \left(\varphi(t_i),\varphi(t_i)+\frac{1}{2^i}\right)$
and let $x_0 \not \in X$. We want to prove that, for any $\varepsilon > 0$, there is a $t_\varepsilon \in T$
such that $x_0 - \varepsilon < \varphi(t_\varepsilon) < x_0$: thus if we take $\varepsilon_n = \frac{1}{n}$,
we can find a sequence $t_{\varepsilon_n} \in T$ such that $\varphi(t_{\varepsilon_n}) \to x_0$ and so $x_0 \in \ov{\varphi(T)}$.

Assume, by contradiction, that there is an $\varepsilon>0$ such that $\varphi(t) \not \in (x_0-\varepsilon,x_0)$,
for any $t \in T$. Observe that $(x_0-\varepsilon,x_0) \cap A=\emptyset$. If this were not true, there would be
a $t_i \in T$ such that $(x_0-\varepsilon,x_0) \cap \left(\varphi(t_i),\varphi(t_i)+\frac{1}{2^i}\right) \ne \emptyset$.
We have the following three cases:
\begin{itemize}
\item $\varphi(t_i) \in (x_0-\varepsilon,x_0)$. This is impossible, because of the definition of $(x_0-\varepsilon,x_0)$.
\item $\varphi(t_i) + \frac{1}{2^i} \in (x_0-\varepsilon,x_0)$. Let $\{t_{i_r}\} \subseteq T$ be a decreasing
sequence converging to $t_i^+$, then $\lim_{r \to \infty} \varphi(t_{i_r})=\varphi(t_i)+\frac{1}{2^i}$.
Thus there is an $r$ such that $\varphi(t_{i_r}) \in (x_0-\varepsilon,x_0)$, contradicting the assumption on
$(x_0-\varepsilon,x_0)$ so this case is not possible.
\item $(x_0-\varepsilon,x_0) \subseteq \left(\varphi(t_i),\varphi(t_i)+\frac{1}{2^i}\right)$. This is also
impossible, as it would imply that $x_0 \in \left(\varphi(t_i),\varphi(t_i)+\frac{1}{2^i}\right) \subseteq A$.
\end{itemize}
Thus $(x_0-\varepsilon,x_0) \cap A=\emptyset$ and so
\[
1=m([0,1]) \ge m((x_0-\varepsilon,x_0)) + m(A)= \varepsilon + 1 > 1
\]
where $m$ is the Lebesgue measure on $[0,1]$. Hence we have a contradiction and the Claim is proved. $\square$

\noindent We can visualize the set $C:=\ov{\varphi(T)}$ as a Cantor set.
If we regard $[0,1]$ as $S^1$, then 
the group $T$ acts on $[0,1]$ by rotations and so each $t \in T$ induces a map $t:C \to C$.
Now we extend this map to a map $t:S^1 \to S^1$ by sending an interval
$X_i:=\left[\varphi(t_i),\varphi(t_i)+\frac{1}{2^i}\right]
\subseteq S^1 \setminus C$ linearly
onto the interval
$t(X_i):=\left[\varphi(t_j),\varphi(t_j)+\frac{1}{2^j} \right]$, where $t_j=t+t_i$ according to the enumeration of $T$.
Thus we can identify $T$ as a subgroup of $\Homeo(S^1)$.

We stretch the interval $I$ into $\ov{X_1}$ and we can regard the group $H_0$
as a subgroup of $\{g \in \Homeo(S^1) \; | \; g(x)=x, \forall x \not \in X_1\} \cong \Homeo(X_1)$
(we still call $H_0$ this subgroup of $\Homeo(S^1)$).
We now consider the subgroup $\langle H_0^t \; | \; t \in T\rangle$
obtained by spreading $H_0$ on the circle through conjugation by elements of $T$.
Since supp$(H_0^t) \subseteq t(X_1)$ for any $t \in T$, the groups $H_0^t$ have
disjoint support hence
they commute elementwise pairwise and $\langle H_0^t \; | \; t \in T\rangle \cong
\prod_{t \in T} H_0^t$. Moreover, the conjugation action of $T$ on $\langle H_0^t \; | \; t \in T\rangle$ permutes the
subgroups $H_0^t$. If follows that
\[
\langle H_0^t,T \; | \; t \in T \rangle = H_0 \wr T \hookrightarrow \Homeo(S^1).
\]
In case $T=\{t_1,\ldots,t_k\}$ is finite, then it is a closed subset of $S^1$. We define $X_i:=(t_i,t_{i+1})$,
for $i=1,\ldots,k$, where $t_{k+1}:=t_1$.
We can copy the same procedure of the infinite case,
by noticing that $S^1=\bigcup_{i=1}^k X_i \cup T$ and embedding $H_0$ into subgroups of $\Homeo(S^1)$
isomorphic with $\Homeo(X_i)$. $\square$

\noindent We now follow the previous proof, but we need to be more careful in order to embed
Thompson's group $T=\PL_2(S^1)$ into $\PL_+(S^1)$.

\theorem{There is an embedding $\varphi:\mathbb{Q}/\mathbb{Z} \hookrightarrow \PL_2(S^1)$ such
that $rot(\varphi(x))=x$ for every $x \in \mathbb{Q}/\mathbb{Z}$ and there is an interval $I \subseteq S^1$
with dyadic endpoints such that $\varphi(x)I$ and $\varphi(y)I$ are disjoint, for all
$x,y \in \mathbb{Q}/\mathbb{Z}$ with $x \ne y$.}

\noindent \emph{Proof.} We consider the set of elements $\{x_n=1/n! \; | \; n \in \mathbb{N}\}$ of $\mathbb{Q}$ which are the primitive
$n!$-th roots of $1$ in $\mathbb{Q}$
with respect to the sum. Thus $n x_n = x_{n-1}$ for each $n$. We want to send each $x_n$ to a homeomorphism $X_n$
of $\PL_2(S^1)$ with $rot(X_n)=1/n!$ and such that $(X_n)^{n!}=id_{S^1}$ and so, since
$\langle x_n \; | \; n \in \mathbb{N}\rangle = \mathbb{Q}/\mathbb{Z}$, we will have an embedding
$\mathbb{Q}/\mathbb{Z} \hookrightarrow \PL_2(S^1)$.
For every positive integer $n$ we choose and fix a partition $P_n$ of the unit interval $[0,1]$ into $2n-1$
intervals whose length is a power of $2$. To set up a notation, we always assume to look at $S^1$ from
the origin of the axes: from this point of view right will mean clockwise and left will mean counterclockwise
and we will always read intervals clockwise.

If we have a partition of $S^1$ in $2m$ intervals, we define a ``shift by 2'' in $\PL_2(S^1)$ to be the homeomorphism
$X$ which permutes the intervals of the partition cyclically and such
that $rot(X)=1/m$ and $X^m=id$. In other words, the ``shift by 2'' sends linearly an interval $V$
to another interval $W$ which is 2 intervals to the right of $V$.

\begin{figure}[0.5\textwidth]
 \begin{center}
  \includegraphics[height=4cm]{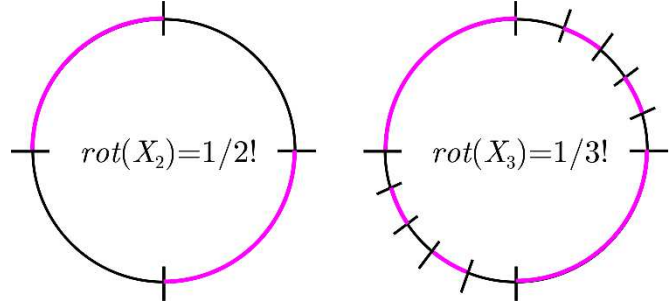}
 \end{center}
 \caption{How to build the map $X_3$ from $X_2$.}
 \label{fig:Q-mod-Z-3-x}
\end{figure}

We want to build a sequence of maps $\{X_n\}$ that acts on a partition of $S^1$ made by
$2(n!)$ intervals $J_{n,1},I_{n,1} \ldots, J_{n,n!},I_{n,n!}$, that are ordered so that
each one on the right of the previous one. The map $X_n$ acts as the ``shift by 2'' map on this partition.
We define $X_1=id$. To build $X_2$,
we cut $S^1$ in four intervals $I_{2,1},J_{2,1},I_{2,2},J_{2,2}$ of length $1/4$, each
one on the right of the previous one: $X_2$ is then defined to be the map
which shifts linearly all these intervals by 2, thus sending the $I$'s onto the $I$'s and the $J$'s onto the $J$'s.
$X_2$ is thus the rotation map by $\pi$. Assume now we have built $X_n$ and we want to build $X_{n+1}$.
We take the $2(n!)$ intervals of the partition associate to $X_n$ and we divide
each of the intervals $I_{n,i}$ according to the proportions given by
the partition $P_n$, and thus cutting each $I_{n,i}$ into $2n+1=2(n+1)-1$ intervals.
On the other hand, we leave all the $J_{n,i}$'s undivided. Now we have a partition of $S^1$ into
\[
n! + (2n+1) n! = 2[(n+1)!]
\]
intervals with dyadic endpoints. Starting $J_{n+1,i}:=J_{n,i}$ we relabel all the intervals of the new partition
by $I$'s and $J$'s. By shifting all the intervals by 2, we have defined
a new piecewise linear map $X_{n+1} \in \PL_2(S^1)$ (see figure \ref{fig:Q-mod-Z-3-x} to see the construction
of maps $X_2$ and $X_3$).

We need to verify that
$(X_{n+1})^{n+1}=X_n$. We observe that $Y_n:=(X_{n+1})^{n+1} \in \PL_2(S^1)$ shifts every interval linearly by $2n+2$.
By construction $Y_n$ sends $J_{n,i}$ linearly onto $J_{n,i+1}$, while it sends $J_{n,i}$
piecewise-linearly onto $J_{n,i+1}$. All the possible breakpoints of $Y_n$ on the interval $I_{n,i}$
occur at the points of the partition $P_n$, but it is a straightforward computation
to verify that the left and right slope coincide
at these points, thus giving that $Y_n$ sends $J_{n,i}$ linearly onto $J_{n,i+1}$.
To build the embedding $\varphi: \mathbb{Q}/\mathbb{Z} \to \PL_2(S^1)$ we
define $\varphi(x_n):=X_n$ and then extend it to a group homomorphism recalling that $\mathbb{Q}/\mathbb{Z}=\langle x_n\rangle$.
The map $\varphi$ must be injective since if $\varphi(x)=id$ then, by using the fact that $(X_{r+1})^{r+1}=X_r$
for any integer $r$, we can write $id=\varphi(x)=X_n^k$ for some suitable integers $n,k$, hence $k$ is a multiple of $n!$
and we can rewrite $x$ as $kx_n=(n!)x_n=0$.
Finally, we notice that for every $x,y \in \mathbb{Q}/\mathbb{Z}, x \ne y$
we have that $\varphi(x)(J_{2,1})$ and $\varphi(y)(J_{2,1})$ are disjoint.
In fact, if we define $V=\varphi(y)(J_{2,1})$, then the two intervals can be rewritten as $\varphi(x y^{-1})(V)$ and $V$,
and so, since $\varphi$ is an embedding and $xy^{-1} \ne 1$, they must be distinct. $\square$

\noindent As an immediate consequence of the previous theorem, we get the following result.

\theorem{For every $H \le \mathbb{Q} / \mathbb{Z}$ there is an embedding $F \wr H \hookrightarrow T$,
where $F$ and $T$ are the respective Thompson's groups and the wreath
product $F \wr H = \left( \bigoplus F \right) \rtimes H$ is restricted. \label{thm:embedding-thompson}}

\noindent \emph{Proof.}  We prove it for the full group $H = \mathbb{Q} / \mathbb{Z}$. We apply the previous
Theorem to build an embedding $\varphi:\mathbb{Q}/\mathbb{Z} \hookrightarrow \PL_2(S^1)$. Moreover,
by construction, the image $\varphi(\mathbb{Q}/\mathbb{Z})$ acts as permutations on the intervals $J_{n,i}$.
Hence, we recover that
\[
\PL_2(J_{2,1}) \wr \mathbb{Q}/\mathbb{Z} \hookrightarrow \PL_2(S^1). \; \; \; \square
\]


\theorem{For every $H \le \mathbb{Q} / \mathbb{Z}$
there is an embedding $\PL_+(I) \wr H \hookrightarrow \PL_+(S^1)$, where the
wreath product $\PL_+(I) \wr H= \left(\bigoplus \PL_+(I) \right) \rtimes H$ is restricted. \label{thm:embedding-PL}}

\noindent \emph{Proof.} The proof of this result is similar to the one of Theorem \ref{thm:embedding-thompson},
except that here we do not require the endpoints of the interval $I$ to be dyadic. $\square$

\remark{We remark that none of the embedding results require the groups to have no
non-abelian free subgroups. For Theorems \ref{thm:embedding-thompson} and \ref{thm:embedding-PL}
the absence of non-abelian free subgroups is guaranteed by the Brin-Squier Theorem \cite{brin1}.
However, we observe that in
Theorem \ref{thm:first-embedding} we can have non-abelian free subgroups and still build the embedding.}

%
%

\section{Fixed-Point Free Actions on the Circle \label{sec:free-action}}

Sacksteder's Theorem states that every fixed-point free action on the real line
must be abelian (see Theorem 2.3 in \cite{farbsha}).
The same result is true also for $G \le \Homeo(I)$ (see Lemma 4.4 by Plante-Thurston
in \cite{plathu}). Both Sacksteder's Theorem and Plante-Thurston's Lemma
are proved by observing that $G$ is an archimedean group and so they appeal to a Theorem of Holder
(for a proof, see \cite{farbsha}).
The following is an alternative proof of the well known version of
Sacksteder's Theorem for subgroups of the group $\Homeo(S^1)$.

\lemma{Let $G \le \Homeo(S^1)$. Suppose that every element of $G\setminus \{id\}$ has no fixed points.
If $f,g \in G$ are such that $rot(f)=rot(g)$, then $f=g$.
Moreover, for every $h,k \in G$, we have that $[h,k]=id$.}

\noindent \emph{Proof.} If $rot(f)=rot(g) \not \in \mathbb{Q}/\mathbb{Z}$ then it has already been proved in
Lemma \ref{thm:same-rot-fix-point}(ii). We suppose now that $rot(f)=rot(g) = k/m \in \mathbb{Q}/\mathbb{Z}$,
hence $f^m,g^m$ both have fixed points. By definition of $G$ this implies that $f^m=g^m=id_{S^1}$
and so $\wh{f}^m(t)=\wh{g}^m(t)+u$ for some integer $u$ and for all $t \in \mathbb{R}$.
We argue, by contradiction, that $fg^{-1}$ has no fixed points in $S^1$. Thus we can assume
$\wh{f}>\wh{g}$ on $\mathbb{R}$. If $u = 0$, then $\wh{g}^m(t) < \wh{f}^m(t)=\wh{g}^m(t)$,
which is a contradiction. If $u \ne 0$ then, for every positive integer $r$, we have
\begin{eqnarray*}
\wh{f}^{mr}(t) = \wh{f}^{m (r-1)}(\wh{f}^m(t))= \wh{f}^{m (r-1)}(\wh{g}^m(t)+u)= \\
=\wh{f}^{m (r-1)}(\wh{g}^m(t))+u=\wh{f}^{m (r-2)}(\wh{g}^m(t)+u)+u = \\
= \wh{f}^{m (r-2)}(\wh{g}^m(t))+ 2u = \ldots = \wh{g}^{mr}(t)+ur.
\end{eqnarray*}
Up to moving $ur$ on the other hand side of the equation, we can assume $u>0$ and divide by $mr$
\[
\frac{\wh{f}^{mr}(t)}{mr} = \frac{\wh{g}^{mr}(t)}{mr} + \frac{u}{m}.
\]
We send $r \to +\infty$ and we get
\[
rot(f)=rot(g)+ \frac{u}{m} > rot(g) = rot(f)
\]
yielding a contradiction. For the second part, we just note that $rot(h)=rot(h^k)$. $\square$

\corollary{Let $G \le \Homeo(S^1)$. Suppose that every element of $G\setminus \{id\}$
has no fixed points. Then $G$ is abelian. \label{thm:non-fixed-abelian-2}}

%
%

\chapter{Centralizers of subgroups of $\Homeo(S^1)$}
\label{chapter7}

In this Chapter we give a description of centralizers of elements in $\PL_+(S^1)$ and Thompson's group $T$.
This analysis is a first step toward a possible solution of the simultaneous conjugacy problem.
We recall that, if $G \le \Homeo(S^1)$, the subset of elements of $G$ with a fixed point does not
necessarily form a subgroup. In Chapter \ref{chapter6} we showed that this happens when $G$ satisfies some additional
requirements. In this Chapter we show that this is still the case when the group $G$ is the centralizer in
$\PL_+(S^1)$ of elements in $\PL_+(S^1)$ with rational rotation number.
In certain cases, we can embed the subgroup of elements with fixed points in $\PL_+(I)$ and use
the results of Chapter \ref{chapter4} to describe it and write $G$ as an extension of this subgroup.
We will identify $S^1$ with $\mathbb{R}/\mathbb{Z}$ to have a well defined origin $0$ on $S^1$.
The material of this Chapter represents
joint work with Collin Bleak and Martin Kassabov.

\section{Centralizers of torsion elements \label{sec:roots-identity}}

In this section we determine centralizers for torsion elements in $\Homeo(S^1)$, in
$\PL_+(S^1)$ and Thompson's group $T$.
In all of these cases, the strategy will be to conjugate the element to a rotation.

\convention{Let $G$ be either of the symbols $\Homeo$ or $\PL_+$ or $\PL_2$. To make the treatment
unified in this subsection, we will write $G(I),G(S^1),G(\mathbb{R})$ to refer to the groups
$\Homeo(I),\Homeo(S^1),\Homeo(\mathbb{R})$ (respectively $\PL_+(I),\PL_+(S^1),\PL_+(\mathbb{R})$)
and $\PL_2(I),\PL_2(S^1),\PL_2(\mathbb{R})$). \label{thm:G-convention}}

\lemma{Let $G$ be either of the symbols $\Homeo$ or $\PL_+$.
Every torsion element of $G(S^1)$ is conjugate to a rotation. \label{thm:torsion-conjugate-rotation}}

\noindent \emph{Proof.} Let $f\in G(S^1)$ be such that $rot(f)=m/k$ and $f^k=id \in G(S^1)$. Since
$(m,k)=1$ are coprime there exist two integers $\alpha,\beta$ such that $\alpha m + \beta k=1$ and so
\[
\frac{1}{k}=\alpha \frac{m}{k} + \beta.
\]
If we define $g=f^\alpha$, we have that $rot(g)=rot(f^\alpha)=\alpha \cdot
rot(f)=\alpha  \frac{m}{k} = \frac{1}{k} \pmod{1}$.
Since $rot(g)=\frac{1}{k}$, we have $\wh{g}^k(t)=t+1$ by Lemma \ref{thm:lift-control}(iv)
and the order of $g$ is $k$ so that
$\langle f \rangle = \langle g \rangle$. We want to find
$h \in G(\mathbb{R})$ such that $h(t+1)=h(t)+1$ and $\wh{g}h(t)=h(t+\frac{1}{k})$,
that is, a map $h$ such that $h^{-1} \wh{g} h$ is the translation by $\frac{1}{k}$.

Choose a real number $A \in [0,1)$ and choose an orientation-preserving
piecewise linear homeomorphism that sends the interval $[0,\frac{1}{k}]$
to the interval $[A,\wh{g}(A)]$, hence  $h(0)=A$ and $h(\frac{1}{k})=\wh{g}(A)$.
We extend $h$ to an element of $G(\mathbb{R})$ by defining
\[
h(t)=\wh{g}^r h \left(t-\frac{r}{k}\right)
\]
if $t \in \left[\frac{r}{k},\frac{r+1}{k}\right]$, for some integer $r$. By
construction, we have that $\wh{g}h(t)=h(t+\frac{1}{k})$ for all $t \in \mathbb{R}$.
Hence, for any $t \in \mathbb{R}$, we have
\[
h(t+1)=h(t + k\frac{1}{k})=\wh{g}^k h(t)=h(t)+1.
\]
If we define $v \in G(S^1)$ as $v(t):=h(t) \pmod{1}$ we have that $\wh{v}=h$ and $v$ is a conjugator between
$g$ and a pure rotation $T_{\frac{1}{k}}$ by $\frac{1}{k}$. $\square$

\theorem{Let $G$ be either of the symbols $\Homeo$ or $\PL_+$. Let $f \in G(S^1)$ be a torsion element.
Then there exist two subgroups $H,K \le C_{G(S^1)}(f)$
such that $H \cong G(I), K \cong \mathbb{R}/\mathbb{Z}$ and $H \cap K = 1$, $C_{G(S^1)}(f)=HK$
and neither $H$ nor $K$ is a normal subgroup.}

\noindent \emph{Proof.} By the previous result, there is an element $h \in G(S^1)$ such that $h^{-1}fh=T_\alpha$
rotation by $\alpha = \frac{m}{k}$. Up to taking a suitable power of $f$, we can assume
that $m=1$. Thus $C_{G(S^1)}(f)=C_{G(S^1)}(hT_\alpha h^{-1})=h C_{G(S^1)}(T_\alpha) h^{-1}$.
We want to find $C_{G(S^1)}(T_\alpha)$. We observe that $C_{G(S^1)}(T_\alpha)$ acts transitively on $S^1$,
since it contains the set $H$ of all rotations of $S^1$.
By the previous result,
we have that $K:=Stab_{C_{G(S^1)}(T_\alpha)}(0) \cong G(\left[0,\frac{1}{k}\right])$.
Elements of $K$ appear as follows: choose an element of $G(\left[0,\frac{1}{k}\right])$
and copy it on each interval $\left[\frac{r}{k},\frac{r+1}{k}\right]$, for $r=0,\ldots,k-1$.
For any $g \in C_{G(S^1)}(T_\alpha)$, there is a rotation $T_\beta$ such that $g T_\beta$ fixes $0$,
and so $g T_\beta \in K$. Hence $H \cdot K = C_{G(S^1)}(T_\alpha)$ and, by construction, we have
$H \cap K=1$. It is easily seen that neither $H$ nor $K$ is normal in $C_{G(S^1)}(T_\alpha)$. $\square$

\remark{The product in the previous Lemma is a special instance of
the Zappa-Szep product. An overview of this type of product can be found in \cite{brinzappa}.
We recall that a group $G$ is the \emph{Zappa-Szep product} of two subgroups $H,K$ if $H \cap K =1$, $HK=G$
but they are not necessarily normal in $G$.}

There is another way to classify the centralizers of the previous Theorem. We are now going to give a description
that will give information about centralizers in Thompson's group $T$ too. The idea is similar to the one of
Lemma \ref{thm:torsion-conjugate-rotation}: instead of conjugating the element to a rotation, we will
rescale the circle to get an isomorphic group where the torsion element is indeed a rotation.


\theorem{Let $G$ be either of the symbols $\Homeo$ or $\PL_+$ or $\PL_2$.
Let $g \in G(S^1)$ with $rot(g)=p/q \in \mathbb{Q}/\mathbb{Z}$, with $p,q$ coprime numbers,
and such that $g^q = id_{S^1}$. Then
$C_{G(S^1)}(g)$ is a central extension
\[
1 \to C_q \to C_{G(S^1)}(g) \to G(S^1) \to 1
\]
where $C_q$ is the cyclic group of order $q$. \label{thm:centralizer-torsion-thompson}}

\noindent \emph{Proof.} Since $p$ and $q$ are coprime numbers then, up
to taking a suitable power, we can assume that $g \in G(S^1)$ with $rot(g)=1/q$.
In the case of $G(S^1)=T$, since $0$ is dyadic then $g(0)$ is dyadic too.
We choose $D:=[0,g(0)]$ as a fundamental domain for the action of $g$,
since
\[
S^1=\bigcup_{i=0}^{q-1} g^i(D) \qquad \mbox{and} \qquad g|_{g^i(D)}=g^i (g|_D).
\]
We are now going to stretch the unit circle to a circle of ``length $q$'', by transforming
the fundamental domain into an interval of length 1 and then
reproducing $g$ in this new setting. If this stretching is done carefully,
using conjugation by a suitable map, the map $g$ becomes a rotation, which is then simpler to centralize.
We look for a homeomorphism $H:\mathbb{R} \to \mathbb{R}$ such that
\begin{itemize}
\item $H(g^k(0))=k$, for any integer $k$, and
\item $H(g(x))=t(H(x))=H(x)+1$, where $t(x)=x+1$.
\end{itemize}
To construct $H$, choose any piecewise-linear homeomorphism $H:[0,g(0)] \to [0,1]$
with finitely many breakpoints: it is immediate to find one such map
in the cases $G=\Homeo$ or $G=\PL_+$. For the case $G=\PL_2$ we apply theorem \ref{thm:thompson-like}
to find a piecewise-linear homeomorphism with the additional requirement of having dyadic rational
breakpoints and all slopes that are power of $2$.
Then we extend it to a map $H \in G(\mathbb{R})$ by defining
\[
H(x):=H(g^{-k}(x))+k \qquad \mbox{if} \; x \in [g^k(0),g^{k+1}(0)] \; \mbox{for some integer $k$}.
\]
It is immediate to see, from the definition of $H$, that $H(g(x))=t(H(x))$ for any real number $x$.
By passing to quotients in $S^1_1:=[0,1]/\{0,1\}$ and $S^1_q:=[0,q]/\{0,q\}$,
we get a map $h:S^1_1 \to S^1_q$, defined by $h(x \pmod{1}):= H(x) \pmod{q}$. Define a map
\[
\begin{array}{cccc}
\varphi:& G(S^1_1) & \longmapsto & G(S^1_q) \\
        &       f     & \longmapsto & hfh^{-1}.
\end{array}
\]
This map is clearly an isomorphism.
By construction we have that $\varphi(g)$ is the rotation map $r:S^1_q \to S^1_q$
defined by $r(x \pmod{q})=t(x) \pmod{q}= x+1 \pmod{q}$. We define two isomorphic copies of the group $G(S^1)$ by
putting $G_1:=C_{G(\mathbb{R})}(t)/\langle t \rangle \cong G(S^1_1)$
and $G_2:=C_{G(\mathbb{R})}(t^q)/\langle t^q \rangle \cong G(S^1_q)$.
Using the isomorphism $\varphi$ it follows $C_{G_1}(g) \cong C_{G_2}(r)$, so we can study the second
centralizer, as the rotation map $r$ is easier to deal with.

Inside the circle $S^1_q$ the map $r$ has rotation number $rot(r)=1/q$.
To compute the centralizer $C_{G_2}(r)$, we need to find all $v \in G_2$ that are induced
by some map $V \in G(\mathbb{R})$ satisfying $V(x+q)=V(x)+q$ (since $v$ is a map on the circle $S^1_q$)
as well as the equality $V(x+1)=V(x)+1$ (since $v$ centralizes the rotation $r$). In other words,
\[
C_{G_2}(r) \cong \{V \in G(\mathbb{R}) \; | \; V(x+1)=V(x)+1\}/\langle t^q \rangle
\]
By construction $\langle r \rangle$ is contained in the center of $C_{G_2}(r)$ and has order $q$. To conclude we
just observe that the quotient is
\[
\frac{C_{G_2}(r)}{\langle r \rangle} \cong \{V \in G(\mathbb{R}) \; | \; V(x+1)=V(x)\} \cong G(S^1). \; \; \; \; \square
\]

\remark{We observe that the extension of Theorem \ref{thm:centralizer-torsion-thompson} does not split. If the extension did split,
following the proof of the Theorem, we would be able to write $C_{G_2}(r)$ as the direct product $\langle r \rangle \times G(S^1)$
where the element
$(r,id)$ has no $q$-th root. In fact, for any element $(x,y) \in \langle r \rangle \times G(S^1)$, we have $(x,y)^q=(id,y^q)$.
However, we observe that the group $C_{G_2}(r)$ contains every rotation contained in $G_2$. Hence it is possible to build a suitable
rotation $v$ in $G_2$ with rotation number $\frac{1}{q^2}$ such that $v^q=r$, leading to a contradiction.}

%
%

\section{Non-torsion elements with rational $rot$ number}

We can use the procedure of Theorem \ref{thm:centralizer-torsion-thompson}
and adapt it to the case of non-torsion elements.
That is, we can stretch the fundamental domain of the action to become an interval of length 1 and then
centralizers will be determined by their behavior on the fundamental domain. We begin with
an elementary result that we will use in the rest of the sections of this chapter.

\lemma{Let $G \le Homeo(S^1)$ stabilizing a finite subset $X \subseteq S^1$. Then
\begin{enumerate}
\item If $g \in G$ fixes a point $x \in X$, then $g|_X=id_X$.
\item The restriction $G|_X$ is a cyclic group $C_k$, where $k$ divides the size $|X|$.
\item Let $G_0=\{g \in G \;|\; rot(g)=0\}$. Then $G_0$ is a normal subgroup of $G$ and $G/G_0 \cong rot(G) \cong C_k$.
\end{enumerate} \label{thm:non-fixed-cyclic}}

\noindent \emph{Proof}. We order the points of $X$ on the circle and label them as $\{1, 2, \ldots, n \}$.
Assume that $1<g(1)=r \le n$. We want to prove that $g$ shifts all the elements of $X$ by $r$ units in the same direction.
The map $g$ sends the interval $[1,2]$ into the interval $[r,g(2)]$. Since the map $g$ is order preserving, we must have
$g(2)>g(1)$ and $g(2)=r+1$. Otherwise, if $g(2)>r+1$, then we would have $1<g^{-1}(r+1) <2$ and $g^{-1}(r+1) \in X$.
Similarly we prove that $g(i)=r+i-1 \pmod{n}$. Hence $g=h^{r-1}$, where $h$ is the map $h(i)=i+1 \pmod{n}$.
This is true for any $g|_X \in G|_X$, therefore $G|_X \le \langle h \rangle$ and $G|_X$ is cyclic of order $k$, for
some integer $k$. If $v$ is the generator of $G|_X$, all of its orbits have size $k$ and so, by the class equation,
we have $|X|=km$, for some integer $m$. It also follows that if $g \in G$ fixes a point $x \in X$, then $g|_X=id_X$.
Moreover, it is now immediate to see that if $g(1)=r$, then $rot(g)=r/n$ and so $rot(G) \cong C_k$.

Let now $g \in G_0$ and we can assume $0 \in Fix(g)$. Write $g \in \Homeo([a,b])$ for some suitable
interval $[a,b]$ of length 1. Since $X$
is $g$-invariant and $g$ preserves the orientation, then $g(i)=i$. Otherwise, $g(i)>i$ for all $i$ and $g$ is a shift map.
But this would imply that $n<g(n)<b$ and $g(n) \in X$ and this is impossible, because $X$ has only $n$ elements.
Hence $G_0$ must be precisely the kernel of the action of $G$. We conclude that $G/G_0 \cong G|_X$ and we are done. $\square$

\note{For the remainder of this Chapter, if $H$ is a subgroup of $G(S^1)$, we denote by
$H_0$ the subset of elements of $H$ that have fixed points.}

%
%
%
%
%
%
%
%

\lemma{Let $G$ be either of the symbols 
$\PL_+$ or $\PL_2$.
Let $g \in G(S^1)$ with $rot(g)=p/q \in \mathbb{Q}/\mathbb{Z}$ and such that $g^q \ne id_{S^1}$.
Then $C_{G(S^1)}(g)_0$ is a subgroup and the group $C_{G(S^1)}(g)$ is an extension
\[
1 \to C_{G(S^1)}(g)_0 \to C_{G(S^1)}(g) \to C_k \to 1
\]
where $C_k$ is the cyclic group of order $k$. \label{thm:centralizers-subgroup-extensions}}

\noindent \emph{Proof.} Since the map $g^q$ has fixed points, it can be considered as an element of $\PL_+(J)$
for some interval $J$, hence the set $X:=\partial Fix(g^q)$ is finite. The conclusion follows via
Lemma \ref{thm:non-fixed-cyclic}, since $C_{G(S^1)}(g)$ stabilizes $X$. $\square$

\remark{The previous proof does not extend immediately to the case $\Homeo(S^1)$ since the set $\partial Fix(g^q)$
is not always finite.}

\theorem{Let $G$ be either of the symbols $\PL_+$ or $\PL_2$.
Let $g \in G(S^1)$ with $rot(g)=p/q \in \mathbb{Q}/\mathbb{Z}$ and such that $g^q \ne id_{S^1}$. Then
$C_{G(S^1)}(g)$ is an extension
\[
1 \to G(I)^r \times \mathbb{Z}^s \to C_{G(S^1)}(g) \to C_k \to 1
\]
where $C_k$ is the cyclic group of order $k$.}

\noindent \emph{Proof.} By Lemma \ref{thm:centralizers-subgroup-extensions}
it is sufficient to give a description of the subgroup $C_{G(S^1)}(g)_0$. We consider
again the action of $\langle g \rangle$ acts on the finite set $X:=\partial Fix(g^q)$.
We choose a fundamental domain for the action of $\langle g \rangle$ on $S^1$,
that is an interval $D$ such that $S^1 = \bigcup_{v=0}^{q-1} g^v(D)$. To
do so we select inequivalent elements of $X$ as endpoints of the intervals in such a way that they build a unique interval $D$.
By definition of $D$ we can write
\[
g^q|_{h(D)} \circ h =h \circ g^q|_D
\]
for every $h \in \langle g \rangle$, and so the structure of bumps of $g^q$ on $h(D)$ is the same as it has on $D$.
In particular, $\left( C_T(g)_0 \right)|_D$ and $\left( C_T(g)_0 \right)|_{h(D)}$ must be isomorphic, for every
element $h \in \langle g \rangle$. Moreover since every element $h \in \left( C_T(g)_0 \right)|_D$
centralizes $g^q |_D$ we
can use Theorem \ref{thm:centralizers-PLSG}(ii) to determine centralizers.
In the case $G=\PL_2$ we have two cases: (i) if $X \cap \mathbb{Z}\left[\frac{1}{2} \right] \ne \emptyset$,
then $g^q|_D$ can be seen as an element of $G(I)$ and, by Theorem \ref{thm:centralizers-PLSG}(ii),
we have $C_{G(S^1)}(g^q)_0 |_D \cong \mathbb{Z}^m \times F^n$.
(ii) If $X \cap \mathbb{Z}\left[\frac{1}{2} \right] = \emptyset$, then we can still see $g^q$ as an element of some
$G(J)$ for some interval $J$ with non-dyadic endpoints and apply Theorem \ref{thm:centralizers-PLSG}(ii)
because the Stair Algorithm is valid independently of the endpoints. Hence, since no fixed point of $g^q$ is dyadic
the centralizer is generated by a suitable root of $g^q$ and is thus infinite cyclic, i.e.
$C_{G(S^1)}(g^q)_0 |_D \cong \mathbb{Z}$. In the case $G=\PL_+$, we only have case (i) and we get again that
$C_{G(S^1)}(g^q)_0 |_D \cong \mathbb{Z}^m \times \PL_+(I)^n$. $\square$

%



\section{More results on Centralizers}

There are many more cases to be explored. We conclude with some results and a discussion
leading to the more difficult and general cases, i.e. the presence of elements with irrational
rotation number or the classification of centralizers of subgroups.

\lemma{Let $f \in \PL_+(S^1)$ be such that $rot(f) \not \in \mathbb{Q}/\mathbb{Z}$,
then $C_{\PL_+(S^1)}(f)$ can be embedded in $\mathbb{R}/\mathbb{Z}$.}

\noindent \emph{Proof.} By Denjoy's Theorem \ref{thm:denjoy}, there is an
$h \in \Homeo(S^1)$ such that $h^{-1} f h = T_\alpha$, the rotation
by $\alpha$. Then the map $\varphi(g):=hgh^{-1}$ sends $C_{\PL_+(S^1)}(f)$
injectively into $C_{\Homeo(S^1)}(T_\alpha)$. It is a well known fact that $C_{\Homeo(S^1)}(T_\alpha)
\cong  \mathbb{R} / \mathbb{Z}$. $\square$

\noindent The previous result is not a complete classification and it is left open for future work.
We describe now some possible directions toward a complete classification.
Suppose $H$ is a subgroup of $\PL_+(S^1)$. If $H$ has no non-abelian
free subgroups, we can use the fact that $H$ can be described as an extension by Theorem \ref{thm:nonfree-classification}
as a starting point for describing $C_{\PL_+(S^1)}(H)$.
However, if $H$ has non-abelian free subgroups, the structure of centralizers can partially
be described. Recall that if the set
\[
\{ h \in H \setminus \{id\} \; | \; Fix(h) \ne \emptyset \}
\]
is empty then $H$ must be abelian by Theorem \ref{thm:non-fixed-abelian-2}. This could
be again a starting point for describing $C_{\PL_+(S^1)}(H)$. If some element of $H$
has a fixed point, we have the following result.

\theorem{Let $H \le \PL_+(S^1)$. Assume that $H$ has a non-abelian free subgroup and contains an element $h_0$
such that $Fix(h_0) \ne \emptyset$.
Then $C_{\PL_+(S^1)}(H) \cong C_k$, for some finite cyclic group $C_k$.}

\noindent \emph{Proof.} Let $x_0 \in \partial Fix(h_0)$ and fix an orientation on the circle
so that we can define intervals.  The proof then divides naturally into several steps.

\noindent \emph{Claim 1.} $X:=\partial Fix(h_0)$ is a finite $C_{\PL_+(S^1)}(H)$-invariant subset of $S^1$.

\noindent \emph{Proof.} The first part follows from the fact that $H$ commutes with $h_0$. $\square$

\noindent \emph{Claim 2.} $H$ has no global fixed point.

\noindent \emph{Proof.} Assume, by contradiction, that $y \in Fix(H)$,
then we can fix $y$ as an origin and write $H \le \PL_+(I)$. This is impossible because $\PL_+(I)$
has no-non abelian free subgroups by the Brin-Squier Theorem \cite{brin1}. $\square$

\noindent \emph{Claim 3.} The set $C_{\PL_+(S^1)}(H)_)$ is a normal subgroup and, for
any $\gamma \in C_{\PL_+(S^1)}(H)_0$, we have $\gamma|_X=id_X$.

\noindent This follows from Lemma \ref{thm:non-fixed-cyclic}(iii), because $C_{\PL_+(S^1)}(H)_0$ is the kernel
of the action of $C_{\PL_+(S^1)}(H)$ on $X$. $\square$

\noindent \emph{Claim 4.} The subgroup $C_{\PL_+(S^1)}(H)_0$ is trivial.

\noindent \emph{Proof.} Let now $C$ be the largest connected set containing $x_0$ on which $\gamma$ is the identity.
Assume by contradiction that $C \ne S^1$, so that $C$ is a suitable interval $(a,b)$. By Claim 2 $b \not \in Fix(H)$,
so there is an $h_1 \in H$ with an orbital $(x_-,x_+)$ containing $b$,
that is $x_-,x_+ \in Fix(h_1)$ and $(x_-,x_+) \subseteq S^1 \setminus Fix(h_1)$.
Consider the restrictions $\ov{\gamma}:=\gamma|_{(x_-,x_+)}$ and $\ov{h_1}:=h_1|_{(x_-,x_+)}$.
Since $[\ov{\gamma},\ov{h_1}]=id|_{(x_-,x_+)}$ and $\ov{h_1}$ is a one-bump
function, Theorem \ref{thm:centralizers-PLSG}(i) implies that $\ov{\gamma}$ is a power of some root of $\ov{h_1}$ and,
in particular, $\ov{\gamma}$ is a one-bump function on $(x_-,x_+)$. The map $\ov{\gamma}$
fixes the midpoint of $(x_-,b)$, therefore $\ov{\gamma}=id|_{(x_-,x_+)}$, hence $\gamma$ is the identity on $(a,x_+)$,
contradicting the maximality of $C$. Therefore $C=S^1$ and, more generally, $K=\{id\}$. $\square$

\noindent The conclusion now follows via Lemma \ref{thm:non-fixed-cyclic}. $\square$

\chapter{A Growth Formula for Thompson's Group $F$}
\label{chapter8}

In this final Chapter we provide an algorithm to compute the size of the balls
in Thompson's group $F$ with respect to the standard $2$-element generating set. We briefly
review the definition of ``growth of a group''.
Let $G$ be a finitely generated group with a fixed generating set $S$. For
$n \in \mathbb{N}$, let $B_n$ denote the ball of radius $n$ in the Cayley graph
of $G$. The \emph{growth function} for $G$ is defined as

\[
\begin{array}{cccc}
\gamma: & \mathbb{N} & \longrightarrow & \mathbb{N} \\
        &     n      & \longmapsto     & |B_n|
\end{array}
\]
Moreover, the \emph{growth rate} of $G$ is defined to be the limit
\[
\lim_{n \to \infty} \sqrt[n]{ | B_n |}.
\]
Guba \cite{guba1} and Burillo
\cite{burillo1} give estimates for lower bounds of the growth rate.
We recall that Thompson's group $F$ has the following infinite presentation
\[
\langle x_0, x_1, x_2, \ldots | x_k^{-1} x_n x_k = x_{n+1}, \forall k<n \rangle
\]
Since $x_k = x_0^{1-k}x_1 x_0^{k-1}$ for $k \ge 2$ the group $F$ is generated
by the elements $x_0$ and $x_1$.
We will use \emph{forest diagrams} introduced by Belk and Brown in \cite{bebr}
to give a procedure to compute the size of the $n$-balls with respect to $x_0$ and $x_1$, which may lead
to information on the growth rate, or at least provide better bounds than those already
known. We remark that it is an ongoing open question to
determine the \emph{growth series} of Thompson's group $F$ (the generating function of the
sequence $\{|B_n|\}_{n \in \mathbb{N}}$), or at least to detect whether or not it is rational or algebraic.

\section{Forest Diagrams}

\begin{figure}
 \begin{center}
  \includegraphics{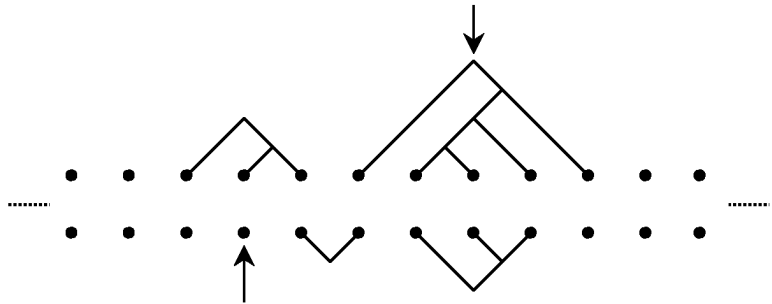} 
 \end{center}
 \caption{A forest diagram for an element of $F$.}
 \label{fig:intro-forestdiagram}
\end{figure}

The idea of this Chapter is to use forest diagrams introduced by Belk and Brown in \cite{bebr}
to find a recursion formula for a partition of the ball of radius $n$.
We will use their length formula for forest diagrams
to calculate distances in the Cayley graph of $F$.
Let $\Gamma$ denote the \emph{Cayley graph} of $F$. This graph has a vertex for each element of
$F$ and an edge from $f$ to $xf$ for every $x \in \{x_0, x_1\}$. The \emph{distance}
between two points in the Cayley graph is the length of a minimal path between them.
The \emph{norm} $\ell(v)$ of a vertex $v \in \Gamma$
is the distance from $v$ to the identity vertex of $\Gamma$.
Each vertex of $\Gamma$ can be represented by a \emph{forest diagram} as shown in figure 1.
Such a diagram consists of a pair of bounded, bi-infinite binary forests (\emph{the top forest}
and the \emph{bottom forest}) together with an order-preserving bijection of their leaves.
Let us be a bit more precise about these definitions. A \emph{bi-infinite binary forest} is
a sequence $(\ldots, T_{-1}, T_0, T_1, \ldots)$ of finite binary trees. We can represent such a forest
as a line of trees, together with a pointer at the tree $T_0$ (as in figure \ref{fig:intro-forestdiagram}).

\begin{figure}
 \begin{center}
  \includegraphics{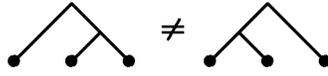} 
 \end{center}
 \caption{Different trees}
 \label{fig:intro-different-trees}
\end{figure}

A forest is
\emph{bounded} if only finitely many of its trees are nontrivial. Note that our binary trees
are planar, i.e., the trees showed in figure \ref{fig:intro-different-trees} are different.

\begin{figure}
 \begin{center}
  \includegraphics{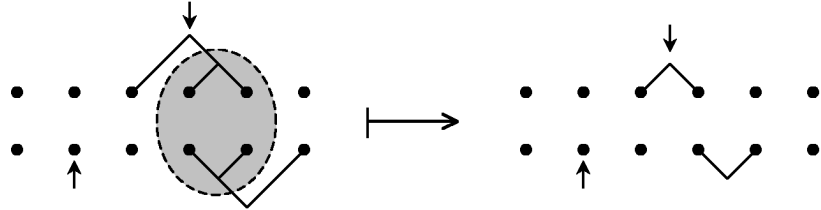}  
 \end{center}
 \caption{Reductions in a forest diagram.}
 \label{fig:intro-reductions}
\end{figure}

In particular, any binary tree comes with a linear ordering on its leaves; and this in
turn induces a linear ordering on the leaves of a bi-infinite binary forest.
A   \emph{caret} is a pair of edges in a forest that join two vertices to a common parent. We
call a caret \emph{grounded} if it joins two leaves. A \emph{reduction} of a forest diagram consists
of removing an opposing pair of grounded carets (see figure \ref{fig:intro-reductions}).

The inverse of a reduction is called an \emph{expansion}. Two forest diagrams are \emph{equivalent}
if one can be transformed into the other by a sequence of reductions and expansions.
A forest diagram is \emph{reduced} if it does not have any opposing pairs of grounded
carets. It turns out that every forest diagram is equivalent to a unique reduced forest
diagram.
\propositionname{\cite{bebr}, Section 4}{There is a one-to-one correspondence between
vertices of $\Gamma$ and equivalence classes of forest diagrams. Therefore, every element
of $F$ can be represented uniquely by a reduced forest diagram.}

\remark{We will frequently identify vertices of $\Gamma$, elements of $F$, and reduced
forest diagrams. For example, if $f \in F$, we might talk about the ``top forest of $f$''.
We hope this will not cause any confusion.}

\begin{figure}
 \begin{center}
  \includegraphics{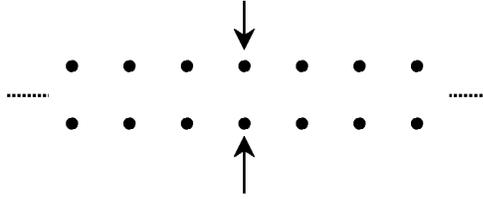}  
 \end{center}
 \caption{The trivial forest diagram.}
 \label{fig:intro-trivial-tree}
\end{figure}

\noindent The forest diagram for the identity is showed in figure \ref{fig:intro-trivial-tree}.

\noindent Given a forest diagram for a vertex $v \in \Gamma$, it is easy to find forest diagrams for the
neighbors of $v$:
\propositionname{\cite{bebr}, Section 4}{Let $\mf{f}$ be a reduced forest diagram
representing the vertex $v \in \Gamma$. Then:

\noindent 1. A forest diagram for $x_0 v$ can be obtained by moving the top pointer of $\mf{f}$ one tree
to the right.

\noindent 2. A forest diagram for $x_1 v$ can be obtained by ``dropping a caret at the current
position''. That is, the forest diagram for $x_1 v$ can be obtained by attaching a
caret to the roots of the top trees in $\mf{f}$ indexed by $0$ and $1$. Afterward, the top
pointer points to the root of the new, combined tree.}

\noindent The bottom forest remains unchanged in either case. Note that the given forest
diagram for $x_1 v$ will need to be reduced if the new caret opposes a grounded caret
from the bottom tree. In this case, left-multiplication by $x_1$ effectively deletes a
grounded caret from the bottom tree.

\corollary{Again, let $\mf{f}$ be a reduced forest diagram for a vertex $v \in \Gamma$. Then:

\noindent 1. A forest diagram for $x^{-1}_0 v$ can be obtained by moving the top pointer of $\mf{f}$ one
tree to the left.

\noindent 2. A forest diagram for $x^{-1}_1 v$ can be obtained by deleting the top caret of the current
tree. The top pointer ends at the resulting left-child tree. If the current tree is
trivial, one must first perform an expansion. In this case, left-multiplication by $x^{-1}_1$
effectively creates a new grounded caret in the bottom tree.}

\begin{figure}
 \begin{center}
  \includegraphics{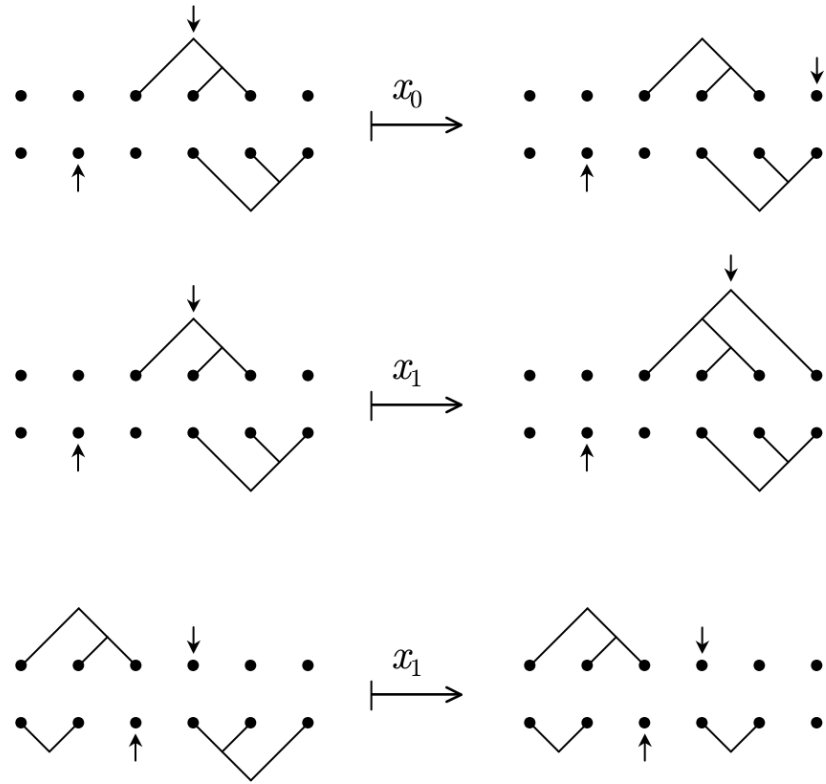}  
 \end{center}
 \caption{Some sample edges from the Cayley graph of $F$.}
 \label{fig:intro-cayley graph}
\end{figure}

\section{The Length Formula}

Since the action of $x_0$ and $x_1$ is relatively simple, it comes as no surprise that one
can find the length $\ell(f)$ of an element $f \in F$ directly from a forest diagram. Our
treatment of the length formula is based on \cite{bebr}.
We begin with some terminology. A \emph{space} is the region between two leaves in a
forest. A space is \emph{interior} if it lies between two leaves from the same tree, and \emph{exterior}
if it lies between two trees. Note that every exterior space in a forest is either to the
left or the right of the pointer.
Given a forest diagram for an element $f \in F$, we label the spaces between the
leaves of each forest as follows. Label a space:

\textbf{L} (for left) if it is exterior and to the left of the corresponding pointer,

\textbf{N} (for necessary) if it is not of type \textbf{L} and if the leaf to the right of the space is a
left leaf in its caret,

\textbf{I} (for interior) if it is interior and not of type \textbf{N}, or

\begin{figure}
 \begin{center}
  \includegraphics{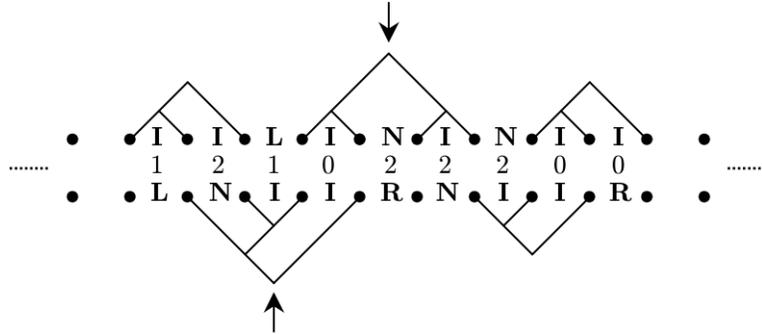}  
 \end{center}
 \caption{The length of this element is 22}
 \label{fig:intro-length}
\end{figure}

\textbf{R} (for right) if it is exterior, to the right of the corresponding pointer, and not of type \textbf{N}.
See figure \ref{fig:intro-length} for an example.
The spaces of a forest diagram come in pairs: one from the top forest and one
from the bottom forest. The \emph{support} of a forest diagram is the minimum interval
that contains both pointers and all nontrivial trees. We only label space pairs in the
support of a forest diagram.
The \emph{weight} of a space pair is determined by the table \ref{fig:length-formula-table}.

\begin{table}
\begin{center}
\begin{tabular}{|c|cccc|}
\hline
         & \bf{N} & \bf{I} & \bf{R} & \bf{L} \\
\hline

\bf{N}   &   2    &   2    &   2    &   1    \\

\bf{I}   &   2    &   0    &   0    &   1    \\

\bf{R}   &   2    &   0    &   2    &   1    \\

\bf{L}   &   1    &   1    &   1    &   2    \\

\hline
\end{tabular}
\end{center}
\caption{The table of weight of the spaces }
\label{fig:length-formula-table}
\end{table}

We can now state the length formula for elements of $F$. We follow the exposition
in Section 5 of \cite{bebr}, which is a simplification of Fordham's original formula \cite{ford1}.
Viewing Thompson's group as a diagram group, Guba \cite{guba1} has
recently obtained a different version of the length formula.

\theoremname{Length Formula \cite{bebr}}{Let $f \in F$, and let $\mf{f}$ be its reduced forest diagram.
Let $\ell_1(f)$ be the total number of carets of $\mf{f}$, and let $\ell_0(f)$ be the sum of the
weights of all space pairs in the support of $\mf{f}$. Then $f$ has length
\[
\ell(f) = \ell_1(f) + \ell_0(f)
\]}

\noindent An element of $F$ is called \emph{positive} if it lies in the submonoid generated by
$\{x_0, x_1, x_2, \ldots\}$. Using the length formula, it is possible to estimate the growth of the elements of the positive monoid
with respect to the $\{x_0,x_1\}$ generating set:

\theoremname{Belk-Brown \cite{bebr}; Burillo \cite{burillo1}}{Let $p_n$ denote the number of positive elements of length $n$, and let:
\begin{equation*}
p(x) = \sum_{n=0}^\infty p_n x^n
\end{equation*}
Then:
\begin{equation*}
p(x) = \frac{1 - x^2}{1 - 2x - x^2 + x^3}
\end{equation*}
In particular, $p_n$ satisfies the recurrence relation:
\begin{equation*}
p_n = 2 p_{n-1} + p_{n-2} - p_{n-3}
\end{equation*}
for large $n$.}

\noindent Using a proof similar to that given by Belk-Brown in \cite{bebr} we obtain a recurrence
formula for a partition of the $n$-sphere of $F$, thus getting estimates for the growth rate.

\section{Partitioning the $n$-sphere}

\begin{figure}
 \begin{center}
  \includegraphics{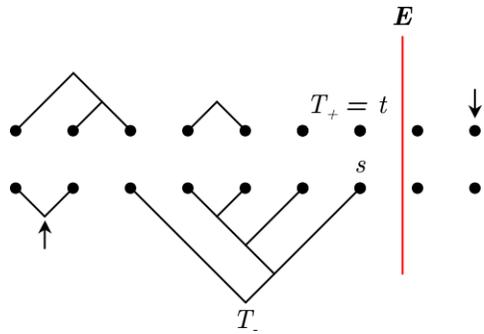}  
 \end{center}
 \caption{The trees $T_+,T_-$ and the line $E=E(f)$.}
 \label{fig:1-critical-line}
\end{figure}

We wish to cut the $n$-sphere so that there is a recurrence formula between the sizes of the slices.
In order to do so, we need to
establish some notations to define how to cut the ball of radius $n$. Throughout this section,
$\mf{f}$ will be a reduced forest diagram for an element $f \ne x_0^k$ for all
$k \in \mathbb{Z}$: this assumption will assure that there exists a nontrivial tree in the diagram $\mf{f}$.
By looking at the diagram $\mf{f}$ we define $\left(\begin{smallmatrix} t \\ s \end{smallmatrix}\right)$ to
be the rightmost pair of corresponding leaves of $\mf{f}$ such that at least one of the two leaves
belongs to a non-trivial tree. Let $T_+$ be the tree of the top forest with $t$ as its rightmost leaf
and let $T_-$ be the tree of the bottom forest with $s$ as its rightmost leaf.
We define the \emph{critical space}
$E=E(f)$ be the space to the right of the pair $\left(\begin{smallmatrix} t \\ s \end{smallmatrix}\right)$
and we denote it by a vertical line passing through it (see figure \ref{fig:1-critical-line}).

\subsection{Weights of the indicated trees}

\begin{figure}
 \begin{center}
  \includegraphics{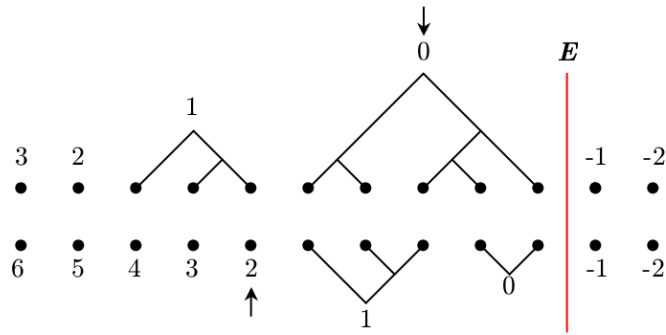} 
 \end{center}
 \caption{In this case $u(f)=0$ and $w(f)=5$.}
 \label{fig:2-weight}
\end{figure}

We order both the forests from right to left by placing an integer on each tree such that
the two trees $T_+,T_-$ correspond to the zeroes. Then we define (see figure \ref{fig:2-weight}):

\begin{eqnarray*}
u(f)= \; \textrm{index of the tree corresponding to the top pointer} \\
w(f)= \; \textrm{index of the tree corresponding to the bottom pointer}
\end{eqnarray*}

\subsection{Length of arcs of right edges}

\begin{figure}
 \begin{center}
  \includegraphics{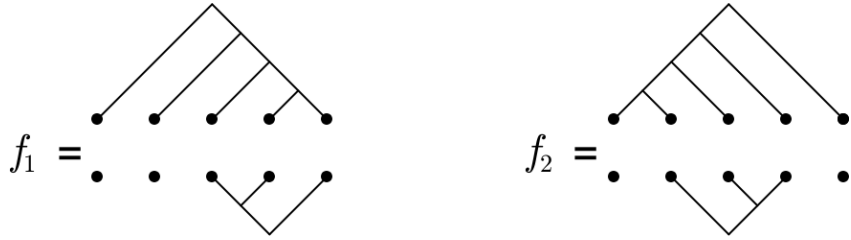} 
 \end{center}
 \caption{Here $b(f_1)=4,c(f_1)=1, b(f_2)=1,c(f_2)=0$.}
 \label{fig:3-arcs}
\end{figure}

Given a binary tree $T$, whether it is oriented upward or downward,
we define its \emph{right arm} to be the longest
path in $T$ starting from the root and made only of right edges.
Then define:
\begin{eqnarray*}
b(f)= \; \textrm{number of edges of the right arm of} \; T_+ \\
c(f)= \; \textrm{number of edges of the right arm of} \; T_-
\end{eqnarray*}
Notice that we always have $\max\{b(f),c(f)\}>0$ (see figure \ref{fig:3-arcs}).

\subsection{Slices of the $n$-sphere} For $(i,j) \in \mathbb{N} \times \mathbb{N}
\setminus \{(0,0)\}$, $p,q \in \mathbb{Z}$ and $n \in \mathbb{N}$ we define the subsets
\[
Z_{i,j,p,q,n}=\left\{f \in F\setminus\ \langle x_0 \rangle \Bigg\vert
\begin{array}{cc}b(f)=i, & c(f)=j, \\
u(f)=p, & w(f)=q, \\
\ell(f)=n & {}
\end{array}
\right\}.
\]
For a fixed
$n$, this family of sets forms a partition of the $n$-sphere.
%
%
%
%
%
%
%
%
%
\noindent It is immediate that the inverses are related by
\[
(Z_{i,j,p,q,n})^{-1}=Z_{j,i,q,p,n}
\]
because inverting a forest diagram means turning it upside down. Thus we can always assume that
$q\le p$, that is ``the bottom pointer is always on the right of the top pointer''.
The interesting cases happen when we find $\max\{i,j,|p|,|q|\} \le n$, in fact:

\medskip

\noindent (a) If $i>n$, $T_+$ has more than $n$ carets so by the length formula $\ell(f)>n$ and
$Z_{i,j,p,q,n}=\emptyset$. Similarly the same is true for $j>n$.

\medskip

\noindent (b) If $p>n$, there are more than $n$ spaces of type $({\bf X,Y})$ with ${\bf X} ={\bf N}$ or
${\bf Y} ={\bf R}$.
Moreover, a space of type $({\bf I,R})$, must have a caret in the top forest.
Each of the spaces in the support has weight $\ge 1$ and so $\ell(f)>n$.

\medskip

\noindent (c) If $q<-n$, there are more than $n$ spaces of type $({\bf L,Y })$, for any ${\bf Y}$ and so their weight
is at least $1$ and so $\ell(f)>n$ again.

\section{A recurrence formula for the slices in 5 variables}

We define a map which shortens the length
of elements and we will show how to keep track of this reduction. This will provide the desired relation.

\subsection{The Shortening Map $\lambda$.} 
\begin{figure}
 \begin{center}
  \includegraphics{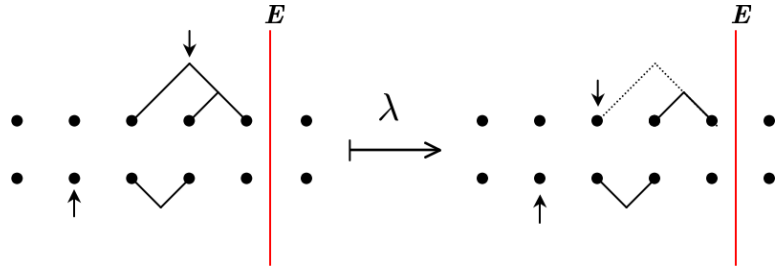} 
 \end{center}
 \caption{The action of $\lambda$ when the top pointer is on $T_+$.}
 \label{fig:4-maps}
\end{figure}

On each slice $Z_{i,j,p,q,n}$
we define a map $\lambda$, provided that $p,q$ are positive integers, $i>0$ for
the map $\lambda$:
\[
\begin{array}{cccc}
\lambda:& Z_{i,j,p,q,n} & \longrightarrow & F \\
        &       f       & \longmapsto     &   x_{p+1}^{-1} f,\\
\end{array}
\]
where $x_{r+1}=x_0^{-r}x_1 x_0^r$. By construction, we see that $\lambda(f)$
is already reduced as a forest diagram and $\lambda$ is an injective map.
The map $\lambda$ removes the top caret of the tree $T_+(f)$
and adjusts the top pointer in the following way:
\begin{itemize}
\item{if the top pointer of $\mf{f}$ is not on $T_+(f)$, then we do not move it}
\item{if the top pointer of $\mf{f}$ is on $T_+(f)$, then the top pointer of $\lambda(f)$ falls on the left
of the two subtrees of $T_+(f)$ (see figure \ref{fig:4-maps}).}
\end{itemize}

\noindent In other words, if $f \in Z_{i,j,p,q,n}$
then the $p$ is the number of the tree where we must add a caret to
get back $f$ from $\lambda(f)$.

\subsection{Length reduction of $\lambda$}
\begin{figure}
 \begin{center}
  \includegraphics{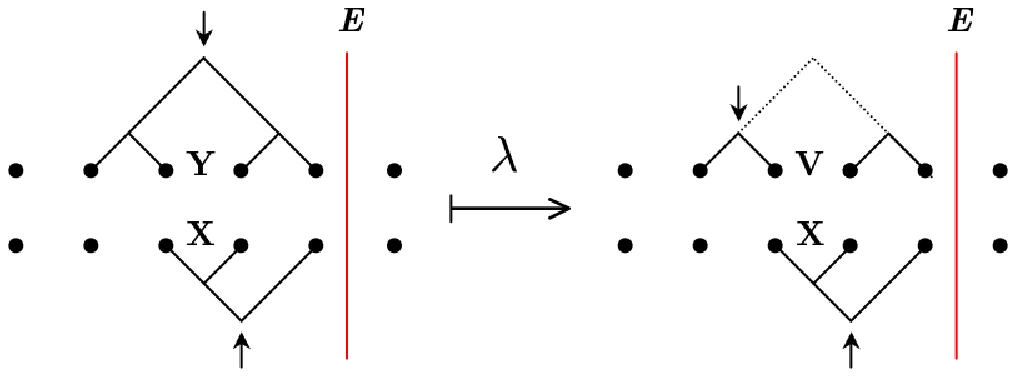} 
 \end{center}
 \caption{The various possibilities for $({\bf Y,X})$ and $({\bf V,X})$}
 \label{fig:5-lengths}
\end{figure}

Let $({\bf Y,X})$ the type of the space under the root of $T_+(f)$ and
$( {\bf V,X })$ the type of the space under the root of $T_+(\lambda(f))$.
By definition of $\lambda$, the top pointer of $\lambda(f)$ is always on the left of the space $( {\bf V,X })$.
We look at all the possibilities for the weights of the spaces $( {\bf Y,X})$ and $({\bf V,X})$
(see figures \ref{fig:5-lengths} and table \ref{fig:lambda-length-reduction}).

\begin{figure}
 \begin{center}
  \includegraphics{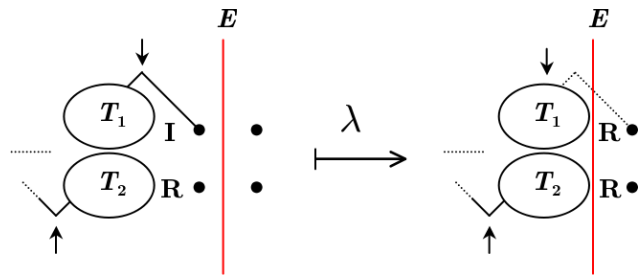} 
 \end{center}
 \caption{In each case $\ell(\lambda(f)) \le \ell(f) - 1$.}
 \label{fig:6-mapseffect}
\end{figure}

\begin{table}
\begin{center}
\begin{tabular}{|r|c|r|c|}
\hline
$({\bf Y,X})$  & weight & $({\bf V,X})$ & weight\\
\hline
$({\bf N,N})$ & 2 & $({\bf N,N})$ & 2 \\
$({\bf I,N})$ & 2 & $({\bf R,N})$ & 2 \\
$({\bf N,I})$ & 2 & $({\bf N,I})$ & 2 \\
$({\bf I,I})$ & 0 & $({\bf R,I})$ & 0 \\
$({\bf N,R})$ & 2 & $({\bf N,R})$ & 2 \\
$({\bf I,R})$ & 0 & $({\bf R,R})$ & 2 \\
$({\bf N,L})$ & 1 & $({\bf N,L})$ & 1 \\
$({\bf I,L})$ & 1 & $({\bf R,L})$ & 1 \\
\hline
\end{tabular}
\end{center}
\caption{How $\lambda$ reduces the length of elements}
\label{fig:lambda-length-reduction}
\end{table}
\noindent We notice that, if the type of space $({\bf Y,X})$ is different from $({\bf I,R})$,
then $\ell(\lambda(f)) = \ell(f) -1$, because the map $\lambda$ removes only a caret and it does not
move the space $E=E(f)$. In the case $({\bf Y,X})$ is of type $({\bf I,R})$, then the space $({\bf V,X})$
is now of type $({\bf R,R})$ and it is out of the support of $\lambda(f)$ so we do not
count it for the evaluation of length $\lambda(f)$, therefore
$\ell(\lambda(f)) \le \ell(f) - 1$ (see figure \ref{fig:6-mapseffect}).

\subsection{The case $p < 0$ \label{ssec:p<0}}

Define
\[
\begin{array}{cccc}
\theta_1:& Z_{i,j,p,q,n} & \longrightarrow & Z_{i,j,p+1,q+1,n-2} \\
         &       f       & \longmapsto     & x_0^{-1}f x_0
\end{array}
\]
\begin{figure}
 \begin{center}
  \includegraphics{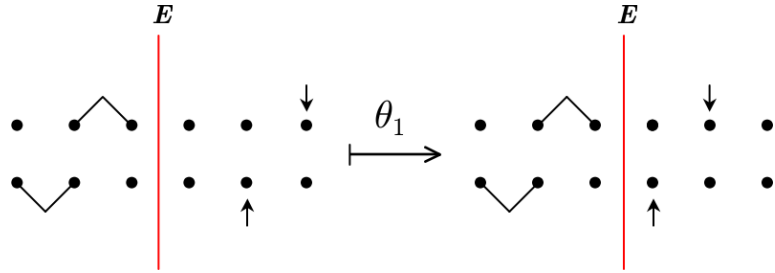} 
 \end{center}
 \caption{The map $\theta_1$.}
 \label{fig:7-case-q-less}
\end{figure}
\noindent so that $\theta_1$ moves both pointers by one space to the left. By construction $\theta_1$ is bijective
and $\ell (\theta_1(f))=\ell(f)-2$. In fact, either a space of type $({\bf L,L})=2$ is lost, or a space of type
$({\bf R,L})=1$ is lost and one of type $({\bf L,L})=2$ becomes an $({\bf R,L })=1$. Therefore:
\[
|Z_{i,j,p,q,n}|=|Z_{i,j,p+1,q+1,n-2}|
\]

\subsection{The case $q<0\le p$ \label{ssec:q<0<p}}

Define
\[
\begin{array}{cccc}
\theta_2:& Z_{i,j,p,q,n} & \longrightarrow & Z_{i,j,p,q+1,n-1} \\
         &       f       & \longmapsto     & f x_0
\end{array}
\]
\begin{figure}
 \begin{center}
  \includegraphics{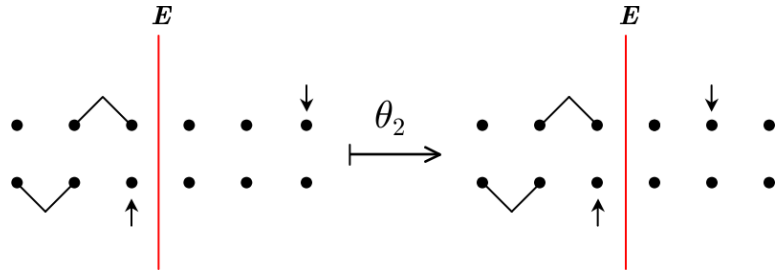} 
 \end{center}
 \caption{The map $\theta_2$.}
 \label{fig:8-case-q-more}
\end{figure}
so that $\theta_2$ moves the bottom pointer by one space to the left. This map is bijective and
$\ell(\theta_2(f))=\ell(f)-1$, since only a space of type $({\bf R,L})=1$ is lost. Thus:
\[
|Z_{i,j,p,q,n}|=|Z_{i,j,p,q+1,n-1}|
\]

\subsection{The case $q \ge 0$ and $i+j \ge 2$ \label{ssec:q>0-i+j>1}}

By definition
of the shortening map $\lambda$ it is easy to verify the following three equalities:
\[
\begin{array}{lc}
\lambda(Z_{i,j,p,q,n})= Z_{i-1,j,p+1,q,n-1}     & \forall i\ge1, j \ge 1 \\
 &  \\
\lambda(\left(Z_{i,j,p,q,n}\right)^{-1})= \left(Z_{i,j-1,p,q+1,n-1}\right)^{-1}     & \forall i \ge 1, j \ge 1 \\
 &  \\
\lambda(Z_{i,0,p,q,n})= Z_{i-1,0,p+1,q,n-1} & \forall i \ge 2 \\
 &  \\
\end{array}
\]

\subsection{The case $q=j=0$ and $i=1$ \label{ssec:qj=zero-i=1}}

\begin{figure}
 \begin{center}
  \includegraphics{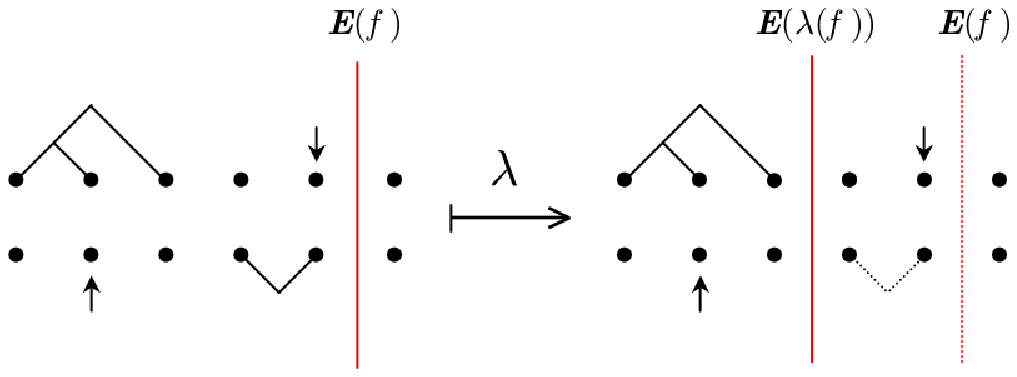} 
 \end{center}
 \caption{supp$(f)=$ supp$(\lambda(f))$.}
 \label{fig:9-case-p-zero}
\end{figure}

Using the map $\lambda$ it can be seen that:
\[
\lambda(Z_{1,0,p,0,n})= \left(\bigcup_{\substack{c,d=0,
\\ \max\{c,d\}>0}}^{n} \bigcup_{r = -n}^{-1} Z_{c,d,r+p+1,r,n-1} \right)
\cup \{ x_0^{1-n}\}
\]
In fact, when we apply $\lambda$ to an element of $f \in Z_{1,0,p,0,n}$ then $\lambda(f)$ can have
no $E$-line, and so $\lambda(f)=x_0^{1-n}$ or it can still have an $E$-line, which is moved to the left
by a suitable number of spaces. In this second case, what happens is that by applying the map $\lambda$
we remove a top caret on the right, but since the bottom pointer of $f$ was on the rightmost leaf inside
supp$(f)$ then supp$(\lambda(f))$ still contains the same number of spaces, that is supp$(f)=$ supp$(\lambda(f))$
(see figure \ref{fig:9-case-p-zero}).

\subsection{The case $q>0=j$ and $i=1$ \label{ssec:q>j=zero-i=1}}

The final case to observe by using the shortening maps is this:

\begin{eqnarray*}
\lambda(Z_{1,0,p,q,n})= \left(\bigcup_{\substack{c,d=0, \\ \max\{c,d\}>0}}^{n}
\bigcup_{r=0}^{q-1} Z_{c,d,r+p+1-q,r,n-2(r-(q-1))-1} \right) \cup \\
\cup \left(\bigcup_{\substack{c,d=0, \\ \max\{c,d\}>0}}^{n}
\bigcup_{r=-n}^{-1} Z_{c,d,r+p+1-q,,r,n-2(q-1)-1} \right) \cup \{ x_0^{p+1-q}\}
\end{eqnarray*}

\begin{figure}
 \begin{center}
  \includegraphics{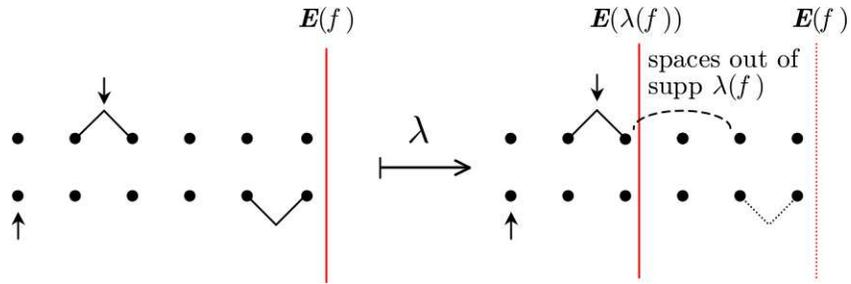} 
 \end{center}
 \caption{supp$(\lambda(f))$ may be reduced.}
 \label{fig:10-case-p-more}
\end{figure}

In fact, when we apply $\lambda$ to an element of $f \in Z_{1,0,p,q,n}$ then $\lambda(f)$ can have
no $E$-line, and so $\lambda(f)=x_0^{p+1-q}$ or it can still have an $E$-line, which is moved to the left
by a suitable number of spaces. In this second case, what happens is that by applying the map $\lambda$
we remove a bottom caret on the right, but we have to consider the fact that the support of $\lambda(f)$ is
now reduced with respect from that of $f$. There might be empty spaces in $f$, between the rightmost caret
and the first tree we find immediately to the left of it, and so we need to keep track of this in the union on the left.
We need to consider all the possibilities for the position of the new $E$-line (see figure
\ref{fig:10-case-p-more}).

\subsection{Commuting parameters \label{ssec:commute}} Define the following map
\[
\begin{array}{cccc}
\varphi:& Z_{1,0,p,q+1,n} & \longrightarrow & Z_{1,0,q,p+1,n} \\
         &       f       & \longmapsto     & x_{q+1} f^{-1} x_{p+1}
         \end{array}
\]
It is immediate that this is a bijection.

\section{Reducing the recurrence to 3 variables}

In this section we put together all the information about the slices and lower the number of parameters
from 5 to 3. We will use the following notation: if we write a sum $\sum_a^b(\ldots)$
where $b < a$ then this sum symbol will denote zero. For example $\sum_{i=3}^2 i^2=0$. Now we define,
for positive $p,q,n \in \mathbb{Z}$:

\[
z(p,q,n)=|Z_{1,0,p,q+1,n}|.
\]
Moreover we define $z(p,q,n)=0$ for all negative $p,q,n \in \mathbb{Z}$. It is immediate from the definition and
the bijection of Subsection \ref{ssec:commute} that
\[
z(p,q,n)=z(q,p,n)
\]
If we let $v=\min\{p,q\}$ and $s=\max\{p,q\}$, then we can decompose the $n$-ball
in slices of the types described in the previous section. We start by applying
the formula of Subsection \ref{ssec:qj=zero-i=1} or Subsection \ref{ssec:q>j=zero-i=1}.
Then we apply the map $\lambda$ to the remaining pieces and use the formulas of
Subsections \ref{ssec:p<0}, \ref{ssec:q<0<p} and \ref{ssec:q>0-i+j>1} to remove carets and move the pointers in order to obtain
a slice of the type $Z_{1,0,a,b+1,c}$. It is a straightforward computation to see that
\begin{eqnarray*}
z(p,q,n)= \sum_{i,j=0}^n \sum_{r=0}^{v-1} z(p-r+i-1,q-r+j-1,n-2r-i-j)-\\
+ \sum_{r=0}^{v-1} z(p-r-1,q-r-1,n-2r) + \\
+ \sum_{\substack{i,j=0 \\ i \ge1 }}^n \sum_{r=v}^{s-1} z(i-1,s-r+j-1,n-v-r-i-j) +\\
+ \sum_{j=1}^n \sum_{r=v}^{s-1} z(0,s-r+j, n- v-r-j-1) + \\
+ \sum_{i.j=1}^{n} \sum_{r=s}^{s+n} z(i-1, j-1, n-v+s-2r -i -j) + \\
+ \sum_{i=1}^n \sum_{r=s}^{s+n} z(i-1,0, n-v+s-2r-i-j-1) + \\
+ \sum_{j=1}^n \sum_{r=s}^{s+n} z(0,j-1,n-v+s-2r-i-j-1)
\end{eqnarray*}

\section{Open Question: the Growth Series of $F$}

Let $G$ be a finitely generated group with a fixed generating set $S$ and let $\gamma$ denote
the associated growth function. We recall that the limit
\[
\gamma:=\lim_{n \to \infty} \sqrt[n]{\gamma(n)}
\]
is called \emph{growth rate} of $G$ with respect to $S$. We say that $G$ has \emph{exponential growth}
if this limit is positive. It can be shown that Thompson's group $F$ has exponential
growth with respect to the generating set $\{x_0,x_1\}$ (see \cite{cfp}). Although, the precise
growth rate is not known, some estimates have been given by Burillo \cite{burillo1},
proving that $\gamma$ cannot be less than the largest root of the equation
$x^3 - x^2 - 2x + 1 = 0$, that is $\gamma \ge 2.2469796 \ldots$, and later were
improved by Guba \cite{guba1}, showing that $\gamma \ge \frac{3+\sqrt{5}}{2}$.
The \emph{growth series} of $G$ with respect to $S$ is the generating
function
\[
\Gamma(t)= \sum_{n=0}^\infty \gamma(n) t^n
\]

\question{Is the function $\Gamma(t)$ rational? Is it algebraic?}

\noindent It has been suggested to use the recurrence formula derived in this chapter
to build a language such that the language growth function
coincides with the growth function of $F$ with respect to $\{x_0,x_1\}$ or with the recurrence formula
that we have derived. It is known that if a language
is regular, the growth series is rational and hence this would be an interesting
direction to try. It would be possible to say something even in the case that the language were
proved to be context-free or indexed (see \cite{wordprocessing} for the definitions
and the standard results about languages and growth functions).

\appendix
\chapter{Omitted Proofs}

\section{Chapter \ref{chapter2} Appendix: Positive Cochains}

\theorem{Let $G$ be a directed graph, and let $c\in
H^1\left(G,\mathbb{Z}\right)$.  Suppose that:
\begin{equation*}
c\left(\ell\right)\geq 0
\end{equation*}
for every directed cycle $\ell$ in $G$.  Then $c$ can be represented by a
cochain that takes a non-negative value on each directed edge. \label{thm:appendix}}

We shall prove this statement using a version of the Farkas lemma.
Call a vector $v\in\mathbb{R}^n$ \emph{non-negative} if each of its entries is
non-negative.

\lemmaname{Farkas}{Let $S$ be a subspace of $\mathbb{R}^n$, and let $a\in\mathbb{R}^n$. \
Then either:

\begin{enumerate}
\item The affine subspace $a+S$ contains a non-negative vector, or
\item There exists a non-negative $v\in S^\perp$ such that $\langle v,a\rangle<0$. $\square$
\end{enumerate}}

\begin{figure}[0.5\textwidth]
 \centering
  \includegraphics[height=4cm]{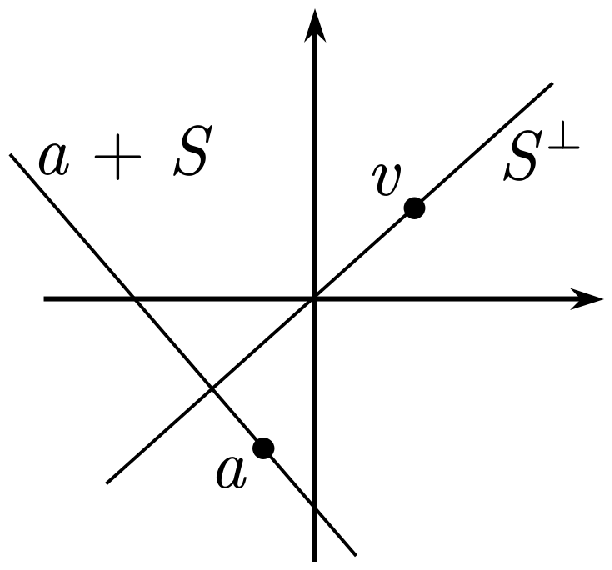}
  \caption{The Farkas Lemma}
  \label{fig:farkas-lemma}
\end{figure}

\noindent Figure \ref{fig:farkas-lemma} illustrates this fact.

Because $a+S$ does not intersect
the first quadrant, $S^\perp$ contains a vector $v$ in the first quadrant
with $\langle v,a\rangle<0$. \ See \cite{zieg1} for more information
on the Farkas lemma, including alternate versions and a simple proof.

\noindent \emph{Proof of Theorem \ref{thm:appendix}:} Let $E$ be the set of directed edges in $G$, and
let $V$ be the set of vertices.  We will begin by producing a non-negative
cocycle in $\mathbb{R}^E$ that represents~$c$.

The set of all cochains representing $c$ is the affine subspace
\begin{equation*}
\alpha+\text{im}\left(\delta\right)\subset\mathbb{R}^E
\end{equation*}
where $\alpha\in\mathbb{R}^E$ is any cocycle representing $c$ and
$\delta\colon\mathbb{R}^V\rightarrow\mathbb{R}^E$ is the coboundary map. \
The orthogonal complement to im$\left(\delta\right)$ is the space of cycles:
\begin{equation*}
\text{im}\left(\delta\right)^\perp=\ker\left(\partial\right)
\end{equation*}
where the boundary map $\partial\colon\mathbb{R}^E\rightarrow\mathbb{R}^V$ is
the adjoint to $\delta$. \ By hypothesis,
$\langle\alpha,\ell\rangle=c\left(\ell\right)\geq 0$ for every directed cycle
$\ell$, and therefore $\langle\alpha,\sigma\rangle\geq 0$ for every positive
cycle $\sigma\in\ker\left(\partial\right)$. \ From the Farkas lemma, we
conclude that the affine subspace $\alpha+$im$\left(\delta\right)$ contains a
non-negative vector $\beta$.

So far, we have proved the existence of a non-negative real cochain $\beta$
representing~$c$. \ We wish to modify $\beta$ to have integer entries. \
Consider the image cochain
$\pi\left(\beta\right)\in\left(\mathbb{R}/\mathbb{Z}\right)^E$. \ Since
$\langle\beta,\ell\rangle=c\left(\ell\right)\in\mathbb{Z}$ for any cycle
$\ell$ with integer coefficients, the image $\pi\left(\beta\right)$ evaluates
to $0\in\mathbb{R}/\mathbb{Z}$ on any cycle, and is therefore a coboundary. \
Choose a function $f\colon V\rightarrow\mathbb{R}/\mathbb{Z}$ so that $\delta
f=\pi\left(\beta\right)$, and let $\overline{f}\colon V\rightarrow[0,1)$ be
the lift of $f$. \ Then the difference $\beta-\delta\overline{f}$ must have
integer entries. \ Since $\beta$ is non-negative and
$\left|\left(\delta\overline{f}\right)\left(e\right)\right|<1$ for any
directed edge $e$, the entries of $\beta-\delta\overline{f}$ must be
non-negative integers, and so $\beta-\delta\overline{f}$ is the desired
representative for $c$. $\square$

\section{Chapter \ref{chapter4} Appendix: Some Computations}

Here is the proof of Proposition \ref{thm:equiv-cond-trans}:

\proposition{Let $J \subseteq [0,1]$ be a closed interval with endpoints in $S$
and let $u,v \in J \cap S$. Then $\pi(u)=\pi(v)$
if and only if there is a $g \in \PL_{S,G}(J)$ such that
$g(u)=v$. \label{thm:equiv-cond-trans-appendix}}

\begin{figure}[0.5\textwidth]
 \begin{center}
  \includegraphics[height=6.5cm]{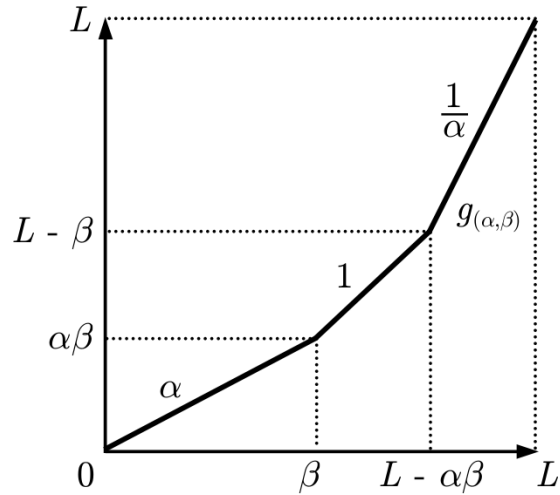}
 \end{center}
 \caption{The basic function to get transitivity.}
 \label{fig:transitive}
\end{figure}

\noindent \emph{Proof.} The sufficient condition is implied by Lemma \ref{thm:nec-cond-trans}.
Suppose now that $J=[\eta,\zeta]$ and let $L=\zeta-\eta$.
We recenter the axis at $(\eta,\eta)$, so that interval $J$ is now $[0,L]$.
For $\alpha \in G, \beta \in J \cap S$ such that $\alpha \beta < L-\beta$
define (see figure \ref{fig:transitive})
\[
g_{\alpha,\beta}(t):=\begin{cases}
\alpha t & t \in [0,\beta] \\
t- (1-\alpha)\beta & t \in [\beta, L-\alpha\beta] \\
\frac{1}{\alpha}(t-L) + L  & t \in [L-\alpha\beta,L]
\end{cases}.
\]

Using the maps $g_{(\alpha,\beta)}$ or $g_{(\alpha,\beta)}^{-1}$ we can send any number $\beta \le t \le L-\alpha\beta$
to $t - (1-\alpha)\beta$ and any number $\alpha \beta \le t \le L -\beta$ to $t + (1-\alpha)\beta$.
We define a relation on $J \cap S$ by saying that $t_1 \sim t_2$,
if either $t_2 = g_{(\alpha,\beta)}(t_1)$ for some $\alpha \in G, \beta \in J \cap S$ such that $\beta \le t \le L-\alpha\beta$
or $t_2 = g_{(\alpha,\beta)}^{-1}(t_1)$ for some $\alpha \in G, \beta \in J \cap S$ such that $\alpha \beta \le t \le L -\beta$.
Then we take the transitive closure of this relation, to get an equivalence relation.
Now, since $\pi(u)=\pi(v)$ then we have that $v-u \in \mathcal{I}$ and so
\[
v - u =  (1-\alpha_1)\beta_1+ \ldots + (1-\alpha_k)\beta_k
\]
for some $\alpha_i \in G,\beta_i \in J \cap S$. We want to rewrite $v-u$ as a sum of terms with $\beta_i$'s
small enough so that we can use the defined equivalence relation. We will define a suitable sequence of numbers
$m_i$ and $\beta_{i,j}$ with $1 \le j \le m_i$, for each $i=1,\ldots, k$. Take $\beta_1$ and choose a number
$\beta_{i,1} \in J \cap S$ small enough such that $g_{(\alpha_i,\beta_{i,1})}$ can be defined. Then choose
inductively a number $\beta_{i,j} \in J \cap S$ small enough such that it satisfies all the following three properties
\begin{itemize}
\item{$g_{(\alpha_i,\beta_{i,j})}$ can be defined}
\item{the number $\beta_{i,j+1}^0:=\beta_i-\beta_{i,1}-\ldots-\beta_{i,j}>0$ is strictly positive}
\item{the number
\[
u+(1-\alpha_1)\sum_{s=1}^{m_1}\beta_{1,s} + \ldots + (1-\alpha_i)\sum_{s=1}^{j-1}\beta_{i,s}
\]
lies in the interval $[\beta_{i,j},L-\alpha_i\beta_{i,j}]$.}
\end{itemize}

\noindent We stop when we find an index $m_i$ such that the number $\beta_{i,m_i}^0$ has the property
that $g_{(\alpha_i,\beta_{i,m_i}^0)}$ can be defined and
\[
u+(1-\alpha_1)\sum_{s=1}^{m_1}\beta_{1,s} + \ldots + (1-\alpha_i)\sum_{s=1}^{m_i-1}\beta_{i,s}
\]
lies in the interval $[\beta_{i,m_i}^0,L-\alpha_i \beta_{i,m_i}^0]$
and so we define $\beta_{i,m_i}:=\beta_{i,m_i}^0$.
We iterate this argument for each $i=1,\ldots,k$ and thus we can rewrite
\[
v-u = (1-\alpha_1)\sum_{j=1}^{m_1}\beta_{1,j} + \ldots + (1-\alpha_k)\sum_{j=1}^{m_k}\beta_{k,j}
\]
and so
\begin{eqnarray*}
u \sim u + (1-\alpha_1)\beta_{1,1} \sim \\
u + (1-\alpha_1)(\beta_{1,1}+\beta_{1,2}) \sim \ldots \sim \\
u + (1-\alpha_1)\sum_{j=1}^{m_1}\beta_{1,j} + \sim \ldots \sim \\
u + (1-\alpha_1)\sum_{j=1}^{m_1}\beta_{1,j} + \ldots + (1-\alpha_k)\sum_{j=1}^{m_k}\beta_{k,j}= v
\end{eqnarray*}
implying that there exists an element $g \in \PL_{S,G}(J)$ such that $g(u)=v$. $\square$

\bibliography{go2}

\end{document}